\documentclass[11pt, oneside, article]{memoir}

\settrims{0pt}{0pt}
\setlrmarginsandblock{1.2in}{1.2in}{*}
\setulmarginsandblock{.98in}{.98in}{*}
\setheadfoot{\onelineskip}{2\onelineskip}
\setheaderspaces{*}{1.5\onelineskip}{*}
\checkandfixthelayout

\setsecheadstyle{\normalfont\large\bfseries}
\setsubsecheadstyle{\normalfont\large\bfseries}
\setsubsubsecheadstyle{\normalfont\normalsize\bfseries}

\usepackage{newpxtext}

\usepackage{amssymb}

\usepackage[varg,bigdelims]{newpxmath}

\usepackage{amsmath}

\usepackage{amsthm}
\usepackage{thmtools}
\usepackage{mathtools}
\usepackage{aliascnt}
\usepackage{array}
\usepackage{dutchcal}
\usepackage{array}
\usepackage{capt-of}
\usepackage{accents}
\usepackage{chngcntr}
\usepackage[multiple]{footmisc}
\usepackage{tikz}
\usepackage{quiver}
\usepackage{xcolor}
\definecolor{darkgreen}{rgb}{0,0.35,0}
\usepackage{enumitem}
\usepackage[colorlinks,citecolor=darkgreen,linkcolor=darkgreen,urlcolor=blue]{hyperref}
\usepackage[nameinlink,capitalize]{cleveref}

\usetikzlibrary{arrows.meta}
\tikzcdset{arrow style=tikz, diagrams={>={Straight Barb[angle'=90,length=1mm]}}}
\tikzset{>={Straight Barb[angle'=90,length=1mm]}}

\AddToHook{cmd/appendix/before}{%
  \crefalias{chapter}{appendix}%
  \crefalias{section}{appendix}%
  \crefalias{subsection}{appendix}
}

\setlist[itemize,enumerate]{
  topsep=3pt,
  itemsep=0pt,
  parsep=0pt,
  partopsep=0pt 
}

\DeclareSymbolFont{stmry}{U}{stmry}{m}{n}
\DeclareMathSymbol\fatsemi\mathop{stmry}{"23}

\DeclareFontFamily{U}{dmjhira}{}
\DeclareFontShape{U}{dmjhira}{m}{n}{ <-> dmjhira }{}
\DeclareRobustCommand{\yoneda}{\text{\usefont{U}{dmjhira}{m}{n}\symbol{"48}}}

\newcommand{\sslash}{/\!\!/}

\DeclareMathAlphabet{\mathbbold}{U}{bbold}{m}{n}

\DeclareFontFamily{U}{mathx}{}
\DeclareFontShape{U}{mathx}{m}{n}{<-> mathx10}{}
\DeclareSymbolFont{mathx}{U}{mathx}{m}{n}
\DeclareMathAccent{\widehat}{0}{mathx}{"70}
\DeclareMathAccent{\widecheck}{0}{mathx}{"71}

\theoremstyle{plain}
\newtheorem{theorem}{Theorem}[section]
\newtheorem{corollary}[theorem]{Corollary}
\newtheorem{proposition}[theorem]{Proposition}
\newtheorem{lemma}[theorem]{Lemma}

\theoremstyle{definition}
\newaliascnt{def}{theorem}
\newtheorem{definition}[def]{Definition}
\aliascntresetthe{def}

\theoremstyle{remark}
\newaliascnt{rem}{theorem}
\newtheorem{myremark}[rem]{Remark}
\aliascntresetthe{rem}
\newtheorem{myexample}[rem]{Example}
\newtheorem{mywarning}[rem]{Warning}
\newtheorem{myconvention}[rem]{Convention}
\newtheorem{mycounterexample}[rem]{Counterexample}

\def\theHsection{\theHchapter.\arabic{section}}
\def\theHtheorem{\theHsection.\arabic{theorem}}
\def\theHcorollary{\theHsection.\arabic{theorem}}
\def\theHproposition{\theHsection.\arabic{theorem}}
\def\theHlemma{\theHsection.\arabic{theorem}}
\def\theHconjecture{\theHsection.\arabic{theorem}}
\def\theHdefinition{\theHsection.\arabic{theorem}}
\def\theHmyremark{\theHsection.\arabic{theorem}}
\def\theHmyexample{\theHsection.\arabic{theorem}}
\def\theHmyconvention{\theHsection.\arabic{theorem}}
\def\theHmywarning{\theHsection.\arabic{theorem}}
\def\theHmycounterexample{\theHsection.\arabic{theorem}}
\def\theHdef{\theHsection.\arabic{theorem}}
\def\theHrem{\theHsection.\arabic{theorem}}

\makeatletter

\let\@origchapter\chapter
\let\@origappendix\appendix

\renewcommand{\thechapter}{\arabic{chapter}}

\newif\if@inappendixmode
\@inappendixmodefalse

\renewcommand{\appendix}{%
  \@origappendix
  \@inappendixmodetrue
  \counterwithout{theorem}{section}%
  \@addtoreset{theorem}{chapter}%
  \gdef\thetheorem{\thechapter.\arabic{theorem}}%
  \global\let\thedef\thetheorem
  \global\let\therem\thetheorem
  \global\let\thedefinition\thetheorem
  \global\let\thelemma\thetheorem
  \global\let\thecorollary\thetheorem
  \global\let\theproposition\thetheorem
  \global\let\theconjecture\thetheorem
  \global\let\themyremark\thetheorem
  \global\let\themyexample\thetheorem
  \global\let\themyconvention\thetheorem
  \global\let\themywarning\thetheorem
  \global\let\themycounterexample\thetheorem
  \gdef\theHtheorem{\theHchapter.\arabic{theorem}}%
  \gdef\theHcorollary{\theHchapter.\arabic{theorem}}%
  \gdef\theHproposition{\theHchapter.\arabic{theorem}}%
  \gdef\theHlemma{\theHchapter.\arabic{theorem}}%
  \gdef\theHconjecture{\theHchapter.\arabic{theorem}}%
  \gdef\theHdefinition{\theHchapter.\arabic{theorem}}%
  \gdef\theHmyremark{\theHchapter.\arabic{theorem}}%
  \gdef\theHmyexample{\theHchapter.\arabic{theorem}}%
  \gdef\theHmyconvention{\theHchapter.\arabic{theorem}}%
  \gdef\theHmywarning{\theHchapter.\arabic{theorem}}%
  \gdef\theHmycounterexample{\theHchapter.\arabic{theorem}}%
  \gdef\theHdef{\theHchapter.\arabic{theorem}}%
  \gdef\theHrem{\theHchapter.\arabic{theorem}}%
}

\newcommand{\@syncthmdisplay}{%
  \global\let\thedef\thetheorem
  \global\let\therem\thetheorem
  \global\let\thedefinition\thetheorem
  \global\let\thelemma\thetheorem
  \global\let\thecorollary\thetheorem
  \global\let\theproposition\thetheorem
  \global\let\theconjecture\thetheorem
  \global\let\themyremark\thetheorem
  \global\let\themyexample\thetheorem
  \global\let\themyconvention\thetheorem
  \global\let\themywarning\thetheorem
  \global\let\themycounterexample\thetheorem
}
\newcommand{\@syncthmH}[1]{%
  \gdef\theHtheorem{#1}\gdef\theHcorollary{#1}\gdef\theHproposition{#1}%
  \gdef\theHlemma{#1}\gdef\theHconjecture{#1}\gdef\theHdefinition{#1}%
  \gdef\theHmyremark{#1}\gdef\theHmyexample{#1}\gdef\theHmyconvention{#1}%
  \gdef\theHmywarning{#1}\gdef\theHmycounterexample{#1}%
  \gdef\theHdef{#1}\gdef\theHrem{#1}%
}
\newcommand{\appendixsectionnumbering}{%
  \counterwithin*{theorem}{section}
  \gdef\thetheorem{\thechapter.\arabic{section}.\arabic{theorem}}%
  \@syncthmdisplay
  \@syncthmH{\theHchapter.\arabic{section}.\arabic{theorem}}%
}

\renewcommand{\chapter}{%
  \if@inappendixmode%
    \renewcommand{\thechapter}{\Alph{chapter}}%
    \renewcommand{\theHchapter}{\Alph{chapter}}%
    \counterwithout*{theorem}{section}
    \gdef\thetheorem{\thechapter.\arabic{theorem}}%
    \@syncthmdisplay
    \@syncthmH{\theHchapter.\arabic{theorem}}%
  \else
    \renewcommand{\thechapter}{\arabic{chapter}}%
    \renewcommand{\theHchapter}{\arabic{chapter}}%
  \fi
  \@origchapter
}

\crefname{chapter}{Section}{Sections}
\newcounter{bgchapter}
\crefname{bgchapter}{Background}{Background}
\newcounter{mypart}
\crefname{mypart}{Part}{Parts}

\newcommand{\mypart}[2]{%
  \clearpage%
  \renewcommand{\themypart}{#1}%
  \refstepcounter{mypart}%
  \addcontentsline{toc}{part}{Part #1: #2}%
  \part*{Part #1: #2}%
  \label{part:#1}%
  \renewcommand{\thebgchapter}{#1}%
  \refstepcounter{bgchapter}%
  \phantomsection%
  \addcontentsline{toc}{chapter}{Background}%
  \chapter*{Background}%
  \renewcommand{\theHchapter}{bg#1}%
  \renewcommand{\thechapter}{#1}
  \label{bg:#1}%
  \setcounter{theorem}{0}%
  \setcounter{section}{0}%
}

\makeatother

\newenvironment{example}
{%
  \pushQED{\qed}\myexample}
{\popQED\endmyexample}

\newenvironment{remark}
{%
  \pushQED{\qed}\myremark}
{\popQED\endmyremark}

\newenvironment{convention}
{%
  \pushQED{\qed}\myconvention}
{\popQED\endmyconvention}

\newenvironment{warning}
{%
  \pushQED{\qed}\mywarning}
{\popQED\endmywarning}

\newenvironment{counterexample}
{%
  \pushQED{\qed}\mycounterexample}
{\popQED\endmycounterexample}

\DeclareMathOperator{\id}{id}
\DeclareMathOperator{\im}{im}

\DeclareMathOperator{\Sh}{Sh}
\DeclareMathOperator{\Psh}{Psh}

\DeclareMathOperator{\El}{El}

\DeclareMathOperator{\Lan}{Lan}
\DeclareMathOperator{\colimw}{colim}
\DeclareMathOperator{\ophelper}{op}
\DeclareMathOperator{\cohelper}{co}

\DeclareMathOperator*{\colim}{colim}
\newcommand{\flim}{\mathop{\widetilde{\mathrm{lim}}}}

\DeclareMathOperator\yo{\!\yoneda\!}

\DeclareMathOperator\groth{\int\!}

\DeclarePairedDelimiter{\gen}{\langle}{\rangle}
\DeclarePairedDelimiter{\abs}{\lvert}{\rvert}
\DeclarePairedDelimiter{\norm}{\lVert}{\rVert}

\newcommand{\dash}{{\mathord{\text{--}}}}
\newcommand{\then}{\mathbin{\fatsemi}}
\newcommand{\op}{^{\ophelper}}
\newcommand{\co}{^{\cohelper}}

\newcommand{\ot}{\leftarrow}
\newcommand{\xto}[1]{\xrightarrow{#1}}
\newcommand{\xot}[1]{\xleftarrow{#1}}
\newcommand{\tn}[1]{\textnormal{#1}}

\newcommand{\imp}{\Rightarrow}

\newcommand{\bang}{\mathrm{!}}

\newcommand{\Pow}[1]{2^{#1}}

\newcommand{\namedcat}[1]{\mathbf{#1}}
\newcommand{\tcat}[1]{\mathbf{#1}}
\newcommand{\cat}[1]{\mathbf{#1}}
\newcommand{\spa}[1]{\mathbf{#1}}
\newcommand{\pre}[1]{\cat{#1}}
\newcommand{\ion}[1]{\spa{#1}}
\newcommand{\Halo}{\namedcat{Halo}}
\newcommand{\LHalo}[1]{\Halo_{#1}}
\newcommand{\Set}{\namedcat{Set}}
\newcommand{\SET}{\namedcat{SET}}
\newcommand{\Cat}{\namedcat{Cat}}
\newcommand{\CAT}{\namedcat{CAT}}
\newcommand{\Span}{\namedcat{Span}}

\newcommand{\Prof}{\namedcat{Prof}}

\renewcommand{\Top}{\namedcat{Top}}
\newcommand{\FinSet}{\namedcat{FinSet}}
\newcommand{\MultiGph}{\namedcat{MultiGph}}

\newcommand{\Part}{\namedcat{Part}}
\newcommand{\Surj}{\namedcat{Surj}}

\newcommand{\Unlink}{\namedcat{Unlink}}
\newcommand{\Shear}{\namedcat{Shear}}

\newcommand{\Int}{\namedcat{Int}}

\newcommand{\Cmd}{\mathbb{C}\namedcat{md}}
\newcommand{\CmdSet}{\Cmd}

\newcommand{\CmdId}{\Cmd_=}
\newcommand{\CmdIdSet}{\CmdId}

\newcommand{\Com}{\namedcat{Cmd}^\sharp}
\newcommand{\ComSet}{\Com}
\newcommand{\Con}{\namedcat{Cmd}^\natural}

\newcommand{\ConSet}{\Con}

\newcommand{\ConO}{\Con}
\newcommand{\ConOSet}{\ConO}

\newcommand{\MndInS}[1]{\namedcat{EM}^{\mathrm{S}}}

\newcommand{\SSpan}{\mathbb{S}\namedcat{pan}}

\newcommand{\MMod}{\mathbb{M}\namedcat{od}}
\newcommand{\Fam}{\mathbf{Fam}}

\newcommand{\CoalgOn}[2]{{#2}^{#1}}
\newcommand{\CoalgSet}[1]{\CoalgOn{#1}{\Set}}
\newcommand{\CoalgS}[1]{\CoalgOn{#1}{\bS}}

\newcommand{\KlOn}[2]{{#2}_{#1}}

\newcommand{\KlS}[1]{\KlOn{#1}{\bS}}

\newcommand{\Comon}[1]{\mathbf{Comon}(#1)}

\newcommand{\clift}[1]{\overline{#1}}

\newcommand{\pts}[1]{{\abs{#1}}}
\newcommand{\obs}[1]{\pts{#1}}
\newcommand{\fpts}[1]{{\norm{#1}}}
\newcommand{\qpts}[1]{\pts{#1}}
\newcommand{\one}[1]{\pts{#1}}
\newcommand{\lift}[1]{\widecheck{#1}}
\newcommand{\ptlift}[1]{\overline{#1}}
\newcommand{\extend}[1]{\widehat{#1}}

\newcommand{\nerve}[1]{(\dash)^{#1}}
\newcommand{\bsheaf}[1]{\nerve{\lift{#1}}}

\newcommand{\halo}[1]{\accentset{\circ}{#1}}
\newcommand{\lhalo}[1]{\accentset{\overline{\circ}}{#1}}

\makeatletter
\newcommand{\doubleoverline}[1]{\mathpalette\dbl@overline{#1}}
\newcommand{\dbl@overline}[2]{\sbox0{$\m@th#1#2$}\dimen0=\fontdimen8\textfont3\vbox{\hrule height\dimen0\kern.75pt\hrule height\dimen0\kern2pt\box0}}
\makeatother

\newcommand{\fptlift}[1]{\doubleoverline{#1}}

\newcommand{\carr}{\Sigma}
\newcommand{\rcarr}{\Delta}

\newcommand{\leindex}[1]{\Sigma_{#1}}
\newcommand{\reindex}[1]{\Delta_{#1}}

\newcommand{\none}{\mathbf{0}}
\newcommand{\point}{\mathbf{1}}
\newcommand{\arro}{\mathbf{2}}
\newcommand{\compo}{\mathbf{3}}

\newcommand{\sier}{\arro}

\newcommand{\slicel}{\Sigma_\bang}
\newcommand{\slicer}{\Delta_\bang}
\newcommand{\sliced}{\Sigma}

\newcommand{\car}[1]{\carr^{#1}}
\newcommand{\rcar}[1]{\rcarr^{#1}}

\newcommand{\carA}{\car{\cA}}

\newcommand{\carB}{\car{\cB}}

\newcommand{\carC}{\car{\cC}}
\newcommand{\rcarC}{\rcar{\cC}}
\newcommand{\carD}{\car{\cD}}

\newcommand{\carP}{\car{\cP}}
\newcommand{\rcarP}{\rcar{\cP}}

\newcommand{\ph}[1]{#1}

\newcommand{\smbullet}{\scalebox{0.6}{$\bullet$}}

\newcommand{\com}[2]{{#1}^{#2}}
\newcommand{\con}[2]{{#1}^{#2}}
\newcommand{\comr}[2]{R}
\newcommand{\connat}[1]{\vec{#1}}

\newcommand{\conion}[1]{{#1}^{-1}}
\newcommand{\contop}[1]{{#1}^{-1}}

\newcommand{\carion}[1]{\car{#1}}

\newcommand{\fun}[2]{{#2}^{#1}}

\newcommand{\topbasis}[1]{\pP_{#1}}
\newcommand{\catbasis}[1]{\pP_{#1}}

\newcommand{\Ps}[1]{\Psh(#1)}

\newcommand{\Coco}[1]{\mathcal{P}(#1)}

\newcommand{\Ho}[1]{#1}

\newcommand{\Ionad}{\namedcat{Ion}}
\newcommand{\IonadO}{\namedcat{Ion}}

\newcommand{\End}[1]{\fun{#1}{#1}}

\newcommand{\SetSet}{\End{\Set}}

\newcommand{\AccCom}{\Com_{\mathrm{small}}}

\newcommand{\AccCon}{\Con_{\mathrm{small}}}
\newcommand{\AccCmd}{\Cmd_{\mathrm{small}}}
\newcommand{\AccCmdId}{(\CmdId)_{\mathrm{small}}}
\newcommand{\AccIonad}{\Ionad_{\mathrm{small}}}
\newcommand{\AccComSet}{\AccCom}

\newcommand{\AccConSet}{\AccCon}
\newcommand{\AccConO}{\AccCon}
\newcommand{\AccConOSet}{\AccConO}
\newcommand{\AccCmdSet}{\AccCmd}
\newcommand{\AccCmdIdSet}{\AccCmdId}

\newcommand{\setof}[1]{\{#1\}}

\newcommand{\cmnd}[1]{{#1}}

\newcommand{\eA}{A}
\newcommand{\eE}{E}
\newcommand{\eS}{S}

\newcommand{\eX}{X}
\newcommand{\xA}{A}
\newcommand{\xB}{B}
\newcommand{\xS}{S}
\newcommand{\xT}{T}
\newcommand{\xU}{U}
\newcommand{\xV}{V}
\newcommand{\xW}{W}
\newcommand{\xX}{X}
\newcommand{\xY}{Y}
\newcommand{\xZ}{Z}

\newcommand{\hK}{k}
\newcommand{\hH}{h}
\newcommand{\cA}{\cmnd{A}}
\newcommand{\cB}{\cmnd{B}}
\newcommand{\cC}{\cmnd{C}}
\newcommand{\cD}{\cmnd{D}}
\newcommand{\cE}{\cmnd{E}}

\newcommand{\bA}{\cat{A}}
\newcommand{\bB}{\cat{B}}
\newcommand{\bC}{\cat{C}}
\newcommand{\bD}{\cat{D}}
\newcommand{\bE}{\cat{E}}
\newcommand{\bG}{\cat{G}}

\newcommand{\bI}{\cat{I}}
\newcommand{\bJ}{\cat{J}}
\newcommand{\bL}{\cat{L}}
\newcommand{\bM}{\cat{M}}
\newcommand{\bS}{\cat{S}}

\newcommand{\pP}{\ph{P}}
\newcommand{\cP}{\gen{\pP}}
\newcommand{\lP}{\lift{\pP}}

\newcommand{\tG}{\spa{G}}
\newcommand{\tE}{\spa{E}}
\newcommand{\tX}{\spa{X}}
\newcommand{\tY}{\spa{Y}}
\newcommand{\iX}{\ion{X}}
\newcommand{\iY}{\ion{Y}}
\newcommand{\cX}{\toptocmd{\tX}}
\newcommand{\cY}{\toptocmd{\tY}}
\newcommand{\cbC}{\cattocmd{\bC}}
\newcommand{\cbD}{\cattocmd{\bD}}

\newcommand{\sS}{\namedcat{S}}
\newcommand{\ASet}{\CoalgSet{\cA}}

\newcommand{\AS}{\CoalgS{\cA}}
\newcommand{\BSet}{\CoalgSet{\cB}}

\newcommand{\BS}{\CoalgS{\cB}}
\newcommand{\CSet}{\CoalgSet{\cC}}

\newcommand{\CS}{\CoalgS{\cC}}
\newcommand{\DSet}{\CoalgSet{\cD}}

\newcommand{\DS}{\CoalgS{\cD}}
\newcommand{\PSet}{\CoalgSet{\cP}}

\newcommand{\PS}{\CoalgS{\cP}}

\newcommand{\KlCS}{\KlS{\cC}}

\newcommand{\SetA}{\Set/\pts{\cA}}

\newcommand{\SetC}{\Set/\pts{\cC}}

\newcommand{\SetP}{\Set/\pts{\pP}}
\newcommand{\SetpA}{\Set^\pts{\cA}}
\newcommand{\SetpB}{\Set^\pts{\cB}}
\newcommand{\SetpC}{\Set^\pts{\cC}}
\newcommand{\SetpD}{\Set^\pts{\cD}}
\newcommand{\SetpP}{\Set^\pts{\pP}}
\newcommand{\SetpX}{\Set^\pts{\iX}}
\newcommand{\SetpY}{\Set^\pts{\iY}}
\newcommand{\SA}{\bS/\pts{\cA}}
\newcommand{\SB}{\bS/\pts{\cB}}
\newcommand{\SC}{\bS/\pts{\cC}}
\newcommand{\SD}{\bS/\pts{\cD}}

\newcommand{\qSC}{\bS/\qpts{\cC}}

\newcommand{\idSet}{\id_\Set}

\newcommand{\La}[2]{\Lan_{#1}#2}

\newcommand{\act}{\mathbin{\cdot}}

\newcommand{\cattocmd}[1]{{{#1}^\bullet}}
\newcommand{\cattocmdp}[1]{{{#1}^\circ}}
\newcommand{\toptocmd}[1]{{{#1}^\circ}}
\newcommand{\cmdtotop}[1]{{\tX_{#1}}}
\newcommand{\topreflect}[1]{(\cmdtotop{#1})^\circ}

\newcommand{\pointcmd}{\point}
\newcommand{\arrocmd}{\arro}

\newcommand{\opens}[1]{\mathcal{O}(#1)}
\newcommand{\neighbors}[1]{\mathcal{N}_{#1}}
\newcommand{\ionopens}[1]{\Sh(#1)}

\newcommand{\coto}{\to}
\newcommand{\conto}{\to}

\newcommand{\bb}{\mathbb{B}}
\newcommand{\nn}{\mathbb{N}}
\newcommand{\zz}{\mathbb{Z}}
\newcommand{\qq}{\mathbb{Q}}
\newcommand{\rr}{\mathbb{R}}

\newcommand{\final}{terminal-object-preserving}
\newcommand{\Final}{Terminal-object-preserving}
\newcommand{\cartesian}{cartesian}
\newcommand{\slice}{slice}

\newcommand{\hh}[2][]{#1 \tn{#2} #1}
\newcommand{\qqand}{\hh[\qquad]{and}}
\newcommand{\qand}{\hh[\quad]{and}}
\renewcommand{\iff}[1][\;\;]{#1\Leftrightarrow#1}

\newcommand{\qor}{\hh[\quad]{or}}

\definecolor{arrowcolor}{RGB}{0,0,255}
\definecolor{pointcolor}{RGB}{0,0,100}

\tikzset{ob/.style={font=\small}}
\tikzset{a/.style={->,shorten <=-1,shorten >=-1}}
\tikzset{ab/.style={<-,shorten <=-1,shorten >=-1}}
\tikzset{long/.style={shorten <=-3,shorten >=-3}}
\tikzset{arr/.style={font=\tiny}}
\tikzset{cell/.style={font=\footnotesize}}
\tikzset{eq/.style={line width=2.4pt,postaction={draw=white, line width=1.6pt, shorten <=-0.05pt, shorten >=-0.05pt}}}
\tikzset{ar/.style={above right=-2}}
\tikzset{al/.style={above left=-2}}
\tikzset{br/.style={below right=-2}}
\tikzset{bl/.style={below left=-2}}
\tikzset{arcl/.style={above right=-3}}
\tikzset{alcl/.style={above left=-3}}
\tikzset{brcl/.style={below right=-3}}
\tikzset{blcl/.style={below left=-3}}
\tikzset{acl/.style={above=-2}}
\tikzset{bcl/.style={below=-2}}

\newcommand{\eqquad}{\hspace{3pt}}
\newcommand{\equad}{\hspace{-2pt}}
\newcommand{\eqq}{\eqquad=\eqquad}
\newcommand{\eq}{\equad=\equad}

\makeatletter
\newcommand{\thickhline}{%
  \noalign {\ifnum 0=`}\fi \hrule height 1pt
  \futurelet \reserved@a \@xhline
}
\newcolumntype{?}{!{\vrule width 1pt}}
\makeatother

\newcommand{\thanksAFOSR}[1]{This material is based upon work supported by the Air Force Office of Scientific Research under award number #1}

\newcommand{\anoteNI}[1]{}

\begin{document}
\title{COMONADS AS SPACES}
\author{Aaron David Fairbanks, Kevin Carlson, and David I.\ Spivak}
\date{}

\maketitle

\vspace{-3.5em}
\begin{abstract}
  Comonads on $\Set$ generalize both categories and topological
  spaces. Expanding upon Garner's work on ionads, we develop aspects
  of the theory of topological spaces for arbitrary comonads on
  arbitrary categories. Our approach is centered around density
  comonads, which provide an abstraction of subbases. We study
  subbases as well as bases in terms of density comonads, and we study
  continuous maps of comonads in terms of functors between coalgebra
  categories, with definitions that recover the usual notions for
  topological spaces.

  Whereas Ahman and Uustalu characterized categories as precisely the
  polynomial comonads on $\Set$, we characterize topological spaces as
  precisely the density comonads of diagrams of subsets of a set,
  which are familiar as topological subbases. We show that every
  comonad on $\Set$ has an underlying topological space, and that this
  construction is a reflection with respect to continuous maps;
  similarly, every comonad on $\Set$ has an underlying small category,
  and this construction is a coreflection. We also show that the
  category of all comonads on $\Set$ with continuous maps is complete,
  and that its full subcategory of accessible comonads is cocomplete.
  Continuous maps and ordinary comonad morphisms form a double
  category, which, in the case of the polynomial comonads on $\Set$,
  recovers the double category of functors and retrofunctors of Clarke
  and Di Meglio. We find topological intuition for these concepts in
  terms of ``halos'', an abstraction of infinitesimal neighborhoods of
  points, defined as formal limits of neighborhood systems.

  We include a long appendix of counterexamples, many applicable to
  general (co)monad theory rather than the particular concerns of this
  text.
\end{abstract}

\clearpage

\hypertarget{toc}{}
\tableofcontents*

\clearpage

\addtocontents{toc}{\vspace{-.5em}}

\newcommand{\chaptertocspace}{}

\addtocontents{toc}{\protect\setcounter{tocdepth}{0}}
\phantomsection
\addcontentsline{toc}{chapter}{Introduction}

\chapter*{Introduction}

Monads provide an abstract notion of \emph{algebra}: groups, rings,
and vector spaces are all examples of algebraic structures
corresponding to monads on the category of sets and functions
$\Set$. In this paper we explore the dual idea that comonads provide
an abstract notion of \emph{space}.

\section*{Topological spaces}

The first justification of this is that topological spaces are
identified with certain comonads on $\Set$. For any topological space
$\tX$, we have the following adjunction.
\[
  \begin{tikzcd}[column sep=50pt]
    \Sh(\tX)\ar[r, shift left=6pt,
    "\parbox{50pt}{\centering\tiny sum of\\ stalks}"]\ar["\bot"{anchor=center},r, draw=none]&
    \Set\ar[l, shift left=6pt, "\parbox{50pt}{\centering\tiny sheaf of\\ functions}"]
  \end{tikzcd}
\]
This adjunction is comonadic: the coalgebras of the induced comonad on
$\Set$ are sheaves on $\tX$. More explicitly, the comonad is given
by
\[\xA \mapsto \sum_{x\in \pts{\tX}} \colim_{U \ni x} \xA^U\] where
$\pts{\tX}$ is the set of points of $\tX$, and $U$ ranges over
neighborhoods of $x$, so that the colimit is the formula for the stalk
at $x$ of the sheaf of arbitrary functions $U\to A$. Moreover, we can
recover the topological space $\tX$ from this comonad as the terminal
object of the category of coalgebras $\Sh(\tX)$, with set of points
given by its underlying set and topology given by subobjects in
$\Sh(\tX)$.

The comonads corresponding to topological spaces also admit a more
abstract characterization. Starting with a small category $\bB$ and a
functor $\pP\colon\bB\to\Set$, we obtain the so-called \emph{density
  comonad} $\cP$, the left Kan extension of $\pP$ along itself.
\[
  \qquad
  \begin{tikzcd}[column sep=15pt, row sep=15pt]
    {\bB} &{}& {\Set} \\
    & {\Set}
    \arrow["\pP", from=1-1, to=1-3]
    \arrow["\pP"', from=1-1, to=2-2]
    \arrow["\cP \,\coloneqq\, \La{\pP}{\pP}"', dashed, from=2-2, to=1-3]
    \arrow["\Downarrow"{description, pos=0.3}, draw=none, from=1-2, to=2-2]
  \end{tikzcd}
\]
We will show that $\cP$ corresponds to a topological space if and only
if $\pP$ constitutes a diagram of subsets of a set\footnote{More
  precisely, to say that $\pP$ constitutes a diagram of subsets of a
  set means each colimit inclusion is injective.}  --- namely the
colimit of $\pP$ --- and in this case the topological space is the one
generated by that diagram of sets, considered as a subbasis
(\cref{thm:tfaetop}).

Thus the concept of topological space falls out of comonads on $\Set$,
without any mention of unions or finite intersections. Moreover,
several central definitions from topology have clean expression in
this form, and they generalize to arbitrary comonads.

\begin{center}
  \renewcommand{\arraystretch}{1.5}
  \begin{tabular}{?>{\centering\arraybackslash}p{11em}?>{\centering\arraybackslash}p{11em}?}%
    \thickhline
    Topological spaces & Comonads on $\Set$ \\
    \thickhline
    Subbasis & $\pP \colon \bB \to \Set$ \\
    \hline
    Space & Density comonad $\cP$\\
    \hline
    Points & $\cP(1)$ \\
    \hline
    Sheaf & Coalgebra \\
    \hline
    Basis & $\bB\to\PSet$ dense \\
    \hline
    Continuous map & $\underset{\phantom{1}}{\overset{\phantom{1}}{\begin{tikzcd}[column sep=15pt, row sep=25,ampersand replacement=\&]
          \CSet \&\& \DSet \\
          \& \Set
          \arrow["\text{\tiny preserves $1$}"', dashed, from=1-3, to=1-1]
          \arrow[""{name=0, anchor=center, inner sep=0}, from=1-1, to=2-2]
          \arrow[""{name=1, anchor=center, inner sep=0}, from=1-3, to=2-2]
          \arrow["\Rightarrow", "\text{\tiny cartesian}"',outer sep=-1pt, shift right=5pt, draw=none, from=1, to=0]
        \end{tikzcd}}}$ \\
    \thickhline
  \end{tabular}\vspace{.1in}
  \captionof{table}{Basic vocabulary for comonads as spaces}\label{fig:table}
\end{center}

For instance, it is not too difficult to show that a subbasis
$\pP \colon \bB \to \Set$ of a topological space $\tX$ is a basis if
and only if $\Sh(\tX)$ canonically embeds as a full subcategory of
$\Ps{\bB}$.  More abstractly, this is to say that the inclusion
$\bB \to \Sh(\tX)$ is \emph{dense}.  We therefore define a
\emph{basis} of a comonad on a general category $\bS$ to be a functor
$\pP \colon \bB \to \bS$ admitting a density comonad
$\cP\colon\bS\to\bS$ such that the canonical lift\footnote{The lift
  $\bB \to \PS$ exists by abstract properties of density comonads,
  which we cover in \cref{sec:density}.} $\bB \to \PS$, where $\PS$
denotes the category of coalgebras, is dense.

It is important to note that the standard definition of morphism of
comonads, when applied to comonads corresponding to topological spaces,
does \emph{not} yield continuous maps.
Whereas a comonad map from $\cC$ to $\cD$ induces a functor from
$\cC$-coalgebras to $\cD$-coalgebras, a continuous map of spaces from
$\tX$ to $\tY$ induces an \emph{inverse image} functor from $\Sh(\tY)$
to $\Sh(\tX)$.\footnote{A continuous map also induces a \emph{direct
    image} functor $\Sh(\tX) \to \Sh(\tY)$, but this does not commute
  with taking underlying sets (that is, sums of stalks), and therefore
  does not correspond to a comonad map either.}  Such a functor comes
from a continuous map just when it is given on underlying sets by
pullback along a function $\one{f} \colon \pts{\tX} \to \pts{\tY}$.

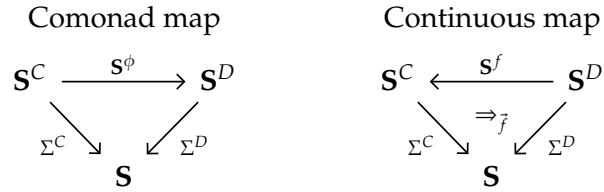
\[
  \overset{\raise.5em\hbox{\normalsize Comonad map}}{
    \begin{tikzcd}[column sep=15pt,ampersand replacement=\&]
      \CS \&\& \DS \\
      \& \bS
      \arrow["{\com{\bS}{\phi}}", from=1-1, to=1-3]
      \arrow["\carC"', from=1-1, to=2-2]
      \arrow["\carD", from=1-3, to=2-2]
    \end{tikzcd}
  }
  \qquad\qquad
  \overset{\raise.5em\hbox{\normalsize Continuous map}}{
    \begin{tikzcd}[column sep=15pt,ampersand replacement=\&]
      \CS \&\& \DS \\
      \& \bS
      \arrow["{\con{\bS}{f}}"', from=1-3, to=1-1]
      \arrow["\carC"', ""{name=0, anchor=center, inner sep=0}, from=1-1, to=2-2]
      \arrow["\carD", ""{name=1, anchor=center, inner sep=0}, from=1-3, to=2-2]
      \arrow["\Rightarrow_{\connat{f}}", shift right=10pt, draw=none, from=1, to=0]
    \end{tikzcd}
  }
\]

Indeed, ordinary comonad maps $\phi$ correspond to commutative
triangles shown above left, where $\carC$ and $\carD$ denote the
comonadic functors and $\com{\bS}{\phi}$ is arbitrary; continuous maps
$f$ correspond to triangles shown above right,\footnote{We identify
  such continuous maps up to natural isomorphism commuting with the
  cartesian triangles.} where $\con{\bS}{f}$ preserves the terminal
object and $\connat{f}$ is a cartesian natural
transformation.\footnote{Here we are assuming $\DS$ has a terminal
  object for simplicity. However, this assumption can be dropped, as
  discussed in \cref{rem:noterminal}.} The map
$\one{f} \colon \pts{\cC} \to \pts{\cD}$ is given by the component of
$\connat{f}$ at the terminal object.


The above definitions make sense for arbitrary comonads on arbitrary
categories. But for concreteness, our examples of focus will be
comonads on $\Set$.

\section*{Ionads}

Not all comonads on $\Set$ come from topological spaces: just as
monads provide a rather general notion of algebra, comonads provide a
rather general notion of space. For instance, more general than
topological spaces --- but still not as general as arbitrary comonads
on $\Set$ --- are the \emph{ionads} of~\cite{garner:ionads}. Ionads
are ``toposes with points'': they generalize topological spaces in the
same way that toposes generalize locales.

An ionad consists of a set $\xX$ and a finite-limit-preserving comonad
on $\Set^\xX$. This provides a category-theoretic analogue of the
order-theoretic definition of a topological space as an \emph{interior
  operator}, or finite-meet-preserving comonad on the poset of subsets
$2^\xX$ of a set $\xX$. Moreover, an ionad is indeed equivalent to a
topos equipped with a set of \emph{enough points}, a.k.a.\ a
\emph{separating set of points}~\cite[C2.2.12]{johnstone:elephant}, as
was noted in~\cite{garner:ionads} and is straightforward to show
(\cref{lem:ionadtopos}).


It was observed by Simon Henry and the second author\footnote{See
  \url{https://mathoverflow.net/questions/457580/examples-of-non-polynomial-comonads-on-set}.}
that ionads are also equivalent to pullback-preserving comonads on
$\Set$ itself (\cref{cor:ionads}).  Moreover, every small category
$\bC$ corresponds to an ionad, whose points are the objects of $\bC$,
and whose topos of coalgebras is $\fun{\bC}{\Set}$. Translating from
ionads to comonads on $\Set$ then recovers the result
of~\cite{ahman-uustalu} that categories are identified with certain
comonads on $\Set$, namely those carried by polynomial
functors.\footnote{Polynomial functors can be equivalently described
  as wide-pullback-preserving functors or as
  connected-limit-preserving functors; similarly, pullback-preserving
  functors are equivalent to finite-connected-limit-preserving
  functors~\cite[Lemma A1.2.9]{johnstone:elephant}.} We thus have the
following diagram of inclusions of objects of various kinds relevant
to the paper.
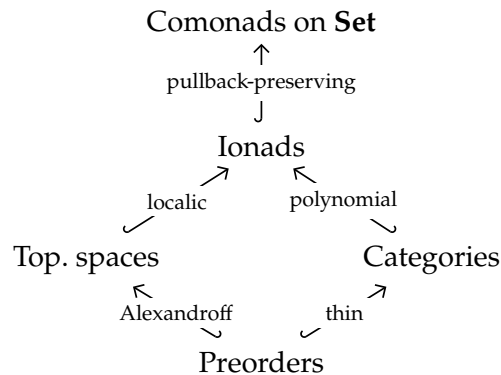
\[
  \begin{tikzcd}[column sep=-15pt, row sep=25pt]
    &\text{Comonads on $\Set$}\\[5pt]
    &\text{Ionads}\ar[hook, u, "\text{pullback-preserving}" description]\\
    \text{Top.\ spaces}\ar[hook, ur, "\text{localic}" description]&&\text{Categories}\ar[hook',ul, "\text{polynomial}" description]\\[-5pt]
    &\text{Preorders}\ar[hook', ul, "\text{Alexandroff}" description]\ar[hook,ur, "\text{thin}" description]
  \end{tikzcd}
\]

When we view a category as a generalized space, we view the objects as
the points, and we view the arrows out of an object as the
``infinitesimal neighborhood'' of the point (a notion we make precise
in \cref{sec:halos}). Each arrow is then a witness of a point being
near another. Intuitively, an ionad generalizes the notion of
topological space by allowing a point to be near another point in
multiple ways.

Our definition of continuous map of comonads is a generalization of
Garner's definition of \emph{continuous map of ionads}, which
corresponds to a geometric morphism of toposes respecting the given
sets of points. Garner also defines 2-cells between continuous maps of
ionads, called \emph{specializations}, which in particular recover
specialization preorders of topological spaces. As observed by Garner,
continuous maps between ionads corresponding to categories are
precisely functors, and specializations between these are precisely
natural transformations.

Just like continuous maps, specializations readily generalize beyond
ionads to arbitrary comonads. The 2-category of comonads on $\Set$,
continuous maps, and specializations thus includes Garner's 2-category
of ionads, which in turn includes both the thin 2-category $\Top$
(whose 2-cells are pointwise specializations) and the 2-category
$\Cat$. Even more satisfying, the 2-category extends to a double
category of comonads, comonad maps, and continuous maps, which
recovers (\cref{thm:doublecat}) the double category of categories,
retrofunctors (a.k.a.\
cofunctors~\cite{higgins-mackenzie,aguiar,dimeglio:lenses,ahman-uustalu:updates,niu-spivak}),
and functors from~\cite{clarke-dimeglio}.

\section*{Universal coalgebra}

Not all comonads of interest on $\Set$ come from ionads (toposes $+$
points). Still, the theory and intuition of comonads as spaces can be
put to work. This is easiest to illustrate with an example.

Whereas the category of \emph{directed multigraphs} is the presheaf
topos $\fun{\rightrightarrows}{\Set}$ where $\rightrightarrows$ is the
walking parallel pair of arrows, the category $\MultiGph$ of
\emph{undirected multigraphs} is not a topos. By an undirected
multigraph, we mean a set $E$ of edges, a set $V$ of vertices, and a
function $E \to \mathbb{P}_{1,2}(V)$ from $E$ to the set of subsets of
$V$ with one or two elements, i.e.\ assigning to each edge an
unordered pair of (not necessarily distinct) vertices; a morphism
consists of a function between vertex sets and a function between edge
sets respecting assigned vertices.  The comonadic functor
$\fun{\rightrightarrows}{\Set} \to \Set$, with induced polynomial
comonad corresponding to the category $\rightrightarrows$, sends a
directed multigraph to the disjoint union $E + V$ of edges and
vertices (its ``elements''). The analogous functor
$\MultiGph \to \Set$, likewise sending an undirected multigraph to the
disjoint union $E + V$ of edges and vertices, is comonadic as well,
inducing a non-pullback-preserving comonad on $\Set$, which we will
call $\cC$.



Running through some of the basic vocabulary from \cref{fig:table},
this comonad $\cC$ represents a ``space'' whose ``sheaves'' are
undirected multigraphs. It has two ``points'', the edge and the
vertex, just like the category $\rightrightarrows$. It has a ``basis''
$\pP \colon \bB \to \Set$ where $\bB$ is the full subcategory of
$\MultiGph$ consisting of the graphs
\[
\begin{tikzpicture}[baseline={([yshift=-3.5pt]current bounding box.center)}]
  \node [inner sep=0.25] (s) at (0,0) {\tiny$\bullet$};
  \node [inner sep=0.25] (t) at (.45,0) {\tiny$\bullet$};
  \draw [line width=0.25mm] (s) -- (t);
\end{tikzpicture}
\;\; \text{(edge)}
\qquad\qqand\qquad
\begin{tikzpicture}[baseline={([yshift=-3.5pt]current bounding box.center)}]
  \node [inner sep=0] at (0,0) {\tiny$\bullet$};
\end{tikzpicture}
\;\;\text{(vertex)}
\]
and the functor $\pP$ sends each to the disjoint union of its edges
and vertices. (There is a similar basis for the comonad corresponding
to $\rightrightarrows$.)
We have a bijective-on-points ``continuous map'' from $\cC$ to the
comonad corresponding to $\rightrightarrows$, whose inverse image
functor $\fun{\rightrightarrows}{\Set} \to \MultiGph$ sends a directed
multigraph to its underlying undirected multigraph. Intuitively, this
continuous map is a refinement that remembers the distinction between
the source and the target of an edge.

\begin{center}
  \vspace*{-55pt}
  \[
    \begin{tikzpicture}[scale=.75]
      \path (0,3) -- (0,-3);
      \fill [rounded corners=24, fill=black, opacity=.25] (-1.875,-1.125) rectangle (1.875,1.125);
      \fill [rounded corners=18, fill=black, opacity=.25] (-.125,-.875) rectangle (1.625,.875);
      \coordinate (e) at (-.75,0);
      \coordinate (v) at (.75,0);
      \node [pointcolor] at (-.8,-.55) {$E$};
      \node [pointcolor] at (.8,-.55) {$V$};
      \begin{scope}[yscale=.75]
        \draw [-stealth, line width=2,shorten <=3,shorten >=6,arrowcolor] (e) .. controls +(.5,.75) and +(-.5,.5) .. node [pos=.45,above] {\footnotesize$s$} (v);
        \draw [-stealth, line width=2,shorten <=3,shorten >=6,arrowcolor] (e) .. controls +(.5,-.75) and +(-.5,-.5) .. node [pos=.45,below] {\footnotesize$t$} (v);
      \end{scope}
      \draw[-stealth, line width=2,shorten <=3,shorten >=6,arrowcolor] (e) .. controls +(-.5,1) and +(-1.125,0) .. node [yshift=7,xshift=9] {\footnotesize$\id_E$} (e);
      \draw[-stealth, line width=2,shorten <=3,shorten >=6,arrowcolor] (v) .. controls +(.5,1) and +(1.125,0) .. node [yshift=6,xshift=-12] {\footnotesize$\id_V$} (v);
      \node [circle, fill=pointcolor, inner sep=3] at (e) {};
      \node [circle, fill=pointcolor, inner sep=3] at (v) {};
      \node at (0,1.625) {$\rightrightarrows$};
      \node at (0,-1.625) {$\xX^3 + \xX$};
    \end{tikzpicture}
    \qquad\qquad\qquad
    \begin{tikzpicture}[scale=.75]
      \path (0,3) -- (0,-3);
      \fill [rounded corners=24, fill=black, opacity=.25] (-1.875,-1.125) rectangle (1.875,1.125);
      \fill [rounded corners=18, fill=black, opacity=.25] (-.125,-.875) rectangle (1.625,.875);
      \coordinate (e) at (-.75,0);
      \coordinate (v) at (.75,0);
      \node [pointcolor] at (-.8,-.55) {$E$};
      \node [pointcolor] at (.8,-.55) {$V$};
      \begin{scope}[yscale=.75]
        \draw [-stealth, line width=2,shorten <=3,shorten >=6,arrowcolor] (e) .. controls +(.5,.75) and +(-.5,.5) .. node [pos=.45,above] {\footnotesize$?$} (v);
        \draw [-stealth, line width=2,shorten <=3,shorten >=6,arrowcolor] (e) .. controls +(.5,-.75) and +(-.5,-.5) .. node [pos=.45,below] {\footnotesize$?$} (v);
      \end{scope}
      \draw[<->,line width=.75] (-.075,-.3) -- (-.075,.3);
      \draw[-stealth, line width=2,shorten <=3,shorten >=6,arrowcolor] (e) .. controls +(-.5,1) and +(-1.125,0) .. node [yshift=7,xshift=9] {\footnotesize$\id_E$} (e);
      \draw[-stealth, line width=2,shorten <=3,shorten >=6,arrowcolor] (v) .. controls +(.5,1) and +(1.125,0) .. node [yshift=6,xshift=-12] {\footnotesize$\id_V$} (v);
      \node [circle, fill=pointcolor, inner sep=3] at (e) {};
      \node [circle, fill=pointcolor, inner sep=3] at (v) {};
      \node at (0,1.625) {$\cC$};
      \node at (0,-1.625) {$\xX \times \binom{\xX+1}{2} + \xX$};
    \end{tikzpicture}
    \vspace{-45pt}
  \]
  \nopagebreak 
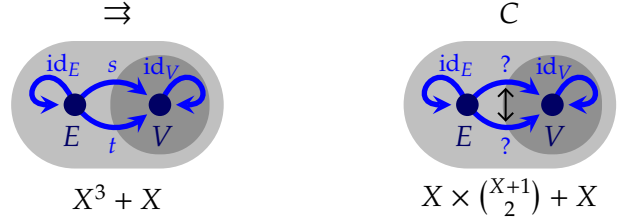
\captionof{figure}{Illustrations of the above comonads
    as spaces, with formulas displayed. Both have two points, $E$
    (edge) and $V$ (vertex). Points correspond to the formula's
    coproduct summands, equivalently the formula's value at $1$. It is
    also possible to extract from the formula of a comonad on $\Set$
    the ``infinitesimal neighborhoods'' of points, as we will explain
    in more detail in \cref{sec:halos}. Whereas the ``infinitesimal
    neighborhood'' of the edge point in the former comonad consists of
    itself and an ordered pair of vertices $s$ and $t$ (represented by
    the functor $\xX^3 \cong \xX \times \xX^2$), the ``infinitesimal
    neighborhood'' of the edge point in the latter comonad consists of
    itself and an unordered pair of vertices (represented by the
    functor $\xX \times \binom{\xX+1}{2}$).}\label{fig:multigraphcomonads}
\end{center}

Although the above comonad $\cC$ does not preserve pullbacks (so does
not correspond to an ionad), it comes quite close: it preserves
\emph{weak pullbacks}. This is a standard assumption in
\emph{universal
  coalgebra}~\cite{rutten,johnstone-power-tsujishita-watanabe-worrell,jacobs,adamek-velebil,gumm-schroder,gumm:elements},
in which coalgebras are used to model transition systems, dynamical
systems, automata, and various infinite data types. Though we do not
cover this angle in detail in this paper, whereas monads correspond to
\emph{varieties} of algebras specified by equations, comonads
correspond to \emph{covarieties} of coalgebras specified by
coequations~\cite{adamek-porst}.

A potential application of the theory we develop is to approach
coalgebra using tools from topology. As a start in this direction, we
characterize weak-pullback-preserving comonads on $\Set$ as the
density comonads of diagrams for which every cospan in the category of
elements extends to a square (\cref{cor:wpb}) --- a straightforward
generalization of the usual concrete definition of basis --- and we
also show that such diagrams are indeed bases in our more abstract
categorical sense.\footnote{This condition that every cospan in the
  category of elements extends to a square, a.k.a.\
  \emph{co-confluence}, is more specific than our general abstract
  definition of basis. All comonads on all categories admit some basis
  (possibly large, for example the comonadic functor itself), but not
  all comonads on $\Set$ admit a basis satisfying the more specific
  co-confluence condition.} (For instance, the above basis of the
comonad $\cC$ satisfies this condition.) Thus a standard assumption
from coalgebra is justified naturally from a topological perspective.


Not only are pullback-preserving comonads on $\Set$ identified with
finite-limit-preserving comonads on categories $\Set^\xX$ (ionads),
but more generally than this, arbitrary comonads on $\Set$ are
identified with comonads on categories $\Set^\xX$ that preserve the
terminal object and finite intersections
(\cref{cor:setslices}). Notice these are two reasonable
generalizations of topological space as finite-meet-preserving
interior operator: preservation of finite limits, versus preservation
of the terminal object and finite intersections. From this
perspective, both ionads and comonads on $\Set$ serve as reasonable
category-theoretic analogues of topological spaces, the latter more
general.

\section*{Outline and contributions}

\cref{bg:I} is a condensed survey of comonads.  \cref{sec:comonads}
reviews basic definitions. \cref{sec:onset} discusses what is known
about comonads on $\Set$. \cref{sec:density} introduces density
comonads.

\cref{sec:spaces} is about the relationship between comonads on $\Set$
and topological spaces. In \cref{sec:toptocomonad}, we show that an
arbitrary comonad on $\Set$ has an underlying topological space; this
may be viewed as a corollary of the result
from~\cite{gumm-schroder:bounded} that every coalgebra of an
endofunctor on $\Set$ has a canonical topology. In
\cref{sec:comonadtotop}, we characterize topological spaces as the
density comonads of diagrams of subsets.

\cref{sec:bases} is about generalized bases. In \cref{sec:basessub},
we introduce the definition of basis for a comonad on an arbitrary
category. In \cref{sec:subbases}, we study the relationship between
generalized bases and subbases. In \cref{sec:coconfluence}, we relate
the usual topological definition of basis to bases of
weak-pullback-preserving comonads on $\Set$.

\cref{bg:II} gives prerequisites around the topic of
ionads. \cref{sec:comonadfunctors} introduces colax comonad functors and
2-cells between them. \cref{sec:ionads} is a review of
ionads. \cref{sec:slices} establishes the equivalence between
pullback-preserving comonads on $\Set$ and ionads, as well as the
equivalence between arbitrary comonads on $\Set$ and comonads on
categories $\Set^X$ that preserve the terminal object and finite
intersections.

\cref{sec:set} is about generalized continuous maps between comonads
on $\Set$. In \cref{sec:continuousset} we give definitions and
generalize some results from~\cite{garner:ionads}. In
\cref{sec:doubleset} we define two double categories $\CmdIdSet$ and
$\CmdSet$ of comonads on $\Set$ (the former a thin sub-double-category
of the latter), which when restricted to just the comonads
corresponding to categories (i.e.\ polynomial comonads) recover the
two known double categories of retrofunctors (a.k.a.\ cofunctors) and
functors from~\cite{clarke} and~\cite{clarke-dimeglio}. In
\cref{sec:limits}, we show that the category of comonads on $\Set$ and
continuous maps is complete, and its full subcategory of small
(a.k.a.\ accessible) comonads is both complete and cocomplete.

In \cref{sec:halos} we introduce ``infinitesimal neighborhoods'', or
\emph{halos}, of points of comonads on $\Set$.  Various concepts
encountered thus far are explained and justified through this
perspective: subbases, continuous maps, comonad maps, specializations,
and identifications all obtain intuitive spatial interpretations in
terms of infinitesimal neighborhoods.

\cref{sec:maps} generalizes continuous maps, as well as the double
categories $\CmdIdSet$ and $\CmdSet$, to comonads on categories
besides $\Set$.
\cref{app:comprehensive} reviews the comprehensive factorization
system, which can be used to generalize results of \cref{sec:maps}
further to categories without terminal objects. \cref{app:universal}
establishes a universal property of density comonads on $\Set$ with
respect to continuous maps. \cref{sec:comonadicity} recalls
comonadicity results.  \cref{app:counterexamples} collects
counterexamples, which show that various results of the paper cannot
be strengthened.

We also refer the reader to~\cite{fairbanks:monads}, originally
written as an appendix to this paper, for elaboration on the formal
category theory underpinning \cref{part:II}. Specifically, one has two
double categories of comonads in an abstract 2-category. In the case
of comonads in $\CAT$, these include the above $\CmdId(\bS)$ and
$\Cmd(\bS)$ as sub-double-categories.

We hope to investigate bicomodules between comonads in future work.

\paragraph{Acknowledgements}
This work was primarily developed by the first author, who plans to
incorporate it in his forthcoming thesis. We would like to thank
Nathanael Arkor, Clark Barwick, Bryce Clarke, Theo Johnson-Freyd,
Simon Henry, Hayato Nasu, Bob Par\'e, Dorette Pronk, Sridhar Ramesh,
Alyssa Renata, and Peter Selinger for useful conversations.
\thanksAFOSR{FA9550-23-1-0376}.

\vspace{35pt}

\begin{center}
  \renewcommand{\arraystretch}{1.5}%
  \begin{tabular}{|>{\centering\arraybackslash}p{5em}|>{\centering\arraybackslash}p{30em}|}
    \hline
    $\Set$ & Category of sets and functions \\
    \hline
    $\Cat$ & Category (or 2-category) of categories and functors \\
    \hline
    $\Top$ & Category (or 2-category) of topological spaces and continuous maps \\
    \hline
    $\Ionad$ & Category (or 2-category) of ionads and continuous maps \\
    \hline
    $\point$ & Terminal category or terminal topological space \\
    \hline
    $\arro$ & Walking arrow category or Sierpinski topological space \\
    \hline
    $\CS$ & Category of coalgebras of comonad $\cC$ on category $\bS$ \\
    \hline
    $\KlCS$ & Category of cofree coalgebras of comonad $\cC$ on category $\bS$ \\
    \hline
    $\fun{\bC}{\bS}$ & Category of functors from category $\bC$ to category $\bS$ \\
    \hline
    $\Ps{\bC}$ & Category of presheaves on category $\bC$ \\
    \hline
    $\Sh(\tX)$ & Category of sheaves on topological space (or ionad) $\tX$ \\
    \hline
    $\opens{\tX}$ & Poset of open subsets of topological space $\tX$ \\
    \hline
    $\bS/\xX$ & Slice category of category $\bS$ over object $\xX$ \\
    \hline
    $(F / G)$ & Comma category of functors $F\colon \bC \to \bS$ and $G \colon \bD \to \bS$ \\
    \hline
    $\El(\pP)$ & Category of elements of functor $\pP \colon \bB \to \Set$ \\
    \hline
    $\Com(\bS)$ & Category (or 2-category) of comonads on $\bS$ and comonad maps \\
    \hline
    $\Con(\bS)$ & Category (or 2-category) of comonads on $\bS$ and continuous maps \\
    \hline
    $\Cmd(\bS)$ & Double category of comonads on $\bS$ with specialization cells \\
    \hline
    $\CmdId(\bS)$ & Double category of comonads on $\bS$ with identification cells \\
    \hline
    $\CmdSet$ & $\Cmd(\Set)$ (and similarly for $\ComSet$, $\ConSet$, and $\CmdIdSet$) \\
    \hline
  \end{tabular}
  \captionof{table}{Notation used throughout the paper}
\end{center}

\setcounter{table}{1}
\begin{center}
  \renewcommand{\arraystretch}{1.5}%
  \begin{tabular}{|>{\centering\arraybackslash}p{4em}|>{\centering\arraybackslash}p{30em}|}
    \hline
    $\cP$ & Density comonad of functor $\pP$ \\
    \hline
    $\pts{\pP}$ & Colimit of functor $\pP$ \\
    \hline
    $\pts{\bC}$ & Objects of category $\bC$ \\
    \hline
    $\pts{\tX}$ & Points of topological space (or ionad) $\tX$ \\
    \hline
    $\one{\phi}$ & Map $\pts{\cC} \to \pts{\cD}$ underlying comonad map $\phi \colon \cC \to \cD$ \\
    \hline
    $\one{f}$ & Map $\pts{\cC} \to \pts{\cD}$ underlying continuous map $f \colon \cC \to \cD$ \\
    \hline
    $\bang$ & Map from initial object or into terminal object \\
    \hline
    $\clift{\cC}$ & Comonad on $\bS/\pts{\cC}$ induced by comonad $\cC$ on $\bS$\\
    \hline
    $\ptlift{\pP}$ & Functor $\bB \to \bS/\pts{\pP}$ induced by functor $\pP \colon \bB \to \bS$\\
    \hline
    $\lP$ & Functor $\bB \to \PS$ induced by functor $\pP \colon \bB \to \bS$\\
    \hline
    $\extend{\pP}$ & Functor $\Ps{\bB} \to \bS$ induced by functor $\pP \colon \bB \to \bS$\\
    \hline
    $\nerve{\pP}$ & Functor $\bS \to \Ps{\bB}$ induced by functor $\pP \colon \bB \to \bS$\\
    \hline
    $\topbasis{\tX}$ & Inclusion $\opens{\tX} \hookrightarrow \Set$ for topological space $\tX$\\
    \hline
    $\cmdtotop{\cC}$ & Topological space underlying comonad $\cC$ on $\Set$ \\
    \hline
    $\cbC$ & Comonad on $\Set$ corresponding to category $\bC$ \\
    \hline
    $\cattocmdp{\bC}$ & Comonad on $\Set$ corresponding to category $\bC\op$ \\
    \hline
    $\cX$ & Comonad on $\Set$ corresponding to topological space $\tX$ \\
    \hline
    $\carC$ & Comonadic forgetful functor $\CS \to \bS$ \\
    \hline
    $\rcarC$ & Right adjoint  $\bS \to \CS$ to $\carC$ \\
    \hline
    $\slicel$ & Slice category projection $\bS/\xX \to \bS$ \\
    \hline
    $\slicer$ & Right adjoint $ \bS \to \bS/\xX$ to $\slicel$ \\
    \hline
    $\leindex{\one{f}}$ & Functor $\bS/\xX \to \bS/\xY$ given by composing with map $\one{f} \colon \xX \to \xY$ \\
    \hline
    $\reindex{\one{f}}$ & Right adjoint $\bS/\xY \to \bS/\xX$ to $\leindex{\one{f}}$ \\
    \hline
    $\Pi_{\one{f}}$ & Right adjoint $\bS/\xX \to \bS/\xY$ to $\reindex{\one{f}}$ \\
    \hline
    $\com{\bS}{\phi}$ & Functor $\CS \to \DS$ underlying comonad map $\phi \colon \cC \to \cD$ \\
    \hline
    $\con{\bS}{f}$ & Functor $\DS \to \CS$ underlying continuous map $f \colon \cC \to \cD$ \\
    \hline
    $\connat{f}$ & Cartesian natural transformation underlying continuous map $f$ \\
    \hline
    $\halo{x}$ & Halo of point $x \in \pts{\cC}$ of comonad $\cC$ on $\Set$ \\
    \hline
  \end{tabular}
  \captionof{table}{Notation used throughout the paper (continued)}
\end{center}

\vspace{1em}

\addtocontents{toc}{\protect\setcounter{tocdepth}{1}}

\mypart{I}{Comonads}

\section{Comonads}\label[bgchapter]{sec:comonads}

\tikzset{every picture/.style={baseline={([yshift=-3.5pt]current bounding box.center)},scale=.125}}

\begin{definition}\label{def:comonad}
  A \emph{comonad} on a category $\bS$ is an endofunctor
  $\cC \colon \bS \to \bS$, equipped with a natural transformation
  $\delta \colon \cC \Rightarrow \cC \circ \cC$, called the
  \emph{comultiplication}, and a natural transformation
  $\varepsilon \colon \cC \Rightarrow \id_{\bS}$, called the
  \emph{counit}, such that the following associativity and unit laws
  hold:
  \[
    \begin{tikzpicture}[xscale=.75,yscale=-1.125]
      \node [ob] (s) at (-14,0) {$\bS$};
      \node [ob] (mm) at (-9,7) {$\bS$};
      \node [ob] (m) at (2,7) {$\bS$};
      \node [ob] (t) at (7,0) {$\bS$};
      \draw [a] (s) -- node[arr,bl] {$\cC$} (mm);
      \draw [a] (mm) -- node[arr,below] {$\cC$} (m);
      \draw [a] (s) -- node[arr,ar,pos=.3] {$\cC$} (m);
      \draw [a] (m) -- node[arr,br] {$\cC$} (t);
      \draw [a] (s) -- node[arr,above] {$\cC$} (t);
      \node [cell] at (0,2.5) {$\delta$};
      \node [cell] at (-6.75,4.75) {$\delta$};
    \end{tikzpicture}
    \eq
    \begin{tikzpicture}[xscale=.75,yscale=-1.125]
      \node [ob] (s) at (-7,0) {$\bS$};
      \node [ob] (m) at (-2,7) {$\bS$};
      \node [ob] (mm) at (9,7) {$\bS$};
      \node [ob] (t) at (14,0) {$\bS$};
      \draw [a] (m) -- node[arr,below] {$\cC$} (mm);
      \draw [a] (mm) -- node[arr,br] {$\cC$} (t);
      \draw [a] (s) -- node[arr,bl] {$\cC$} (m);
      \draw [a] (m) -- node[arr,al,pos=.7] {$\cC$} (t);
      \draw [a] (s) -- node[arr,above] {$\cC$} (t);
      \node [cell] at (0,2.5) {$\delta$};
      \node [cell] at (6.75,4.75) {$\delta$};
    \end{tikzpicture}
    \qquad\qquad
    \begin{tikzpicture}[xscale=1.125,yscale=-1.125]
      \path (0,11) -- (0,-4);
      \node [ob] (s) at (-7,0) {$\bS$};
      \node [ob] (m) at (3,7) {$\bS$};
      \node [ob] (t) at (7,0) {$\bS$};
      \draw [a] (s) -- node[arr,ar,pos=.3] {$\cC$} (m);
      \draw [a] (m) -- node[arr,br] {$\cC$} (t);
      \draw [a] (s) -- node[arr,above] {$\cC$} (t);
      \draw [eq] (s) .. controls +(1,6.5) and +(-6.5,1) .. (m);
      \node [cell] at (1.5,2.5) {$\delta$};
      \node [cell] at (-3,4.75) {$\varepsilon$};
    \end{tikzpicture}
    \eq
    \begin{tikzpicture}[xscale=1.125,yscale=-1.125]
      \path (0,11) -- (0,-4);
      \node [ob] (s) at (-6,0) {$\bS$};
      \node [ob] (t) at (6,0) {$\bS$};
      \draw [a] (s) .. controls +(3,6) and +(-3,6) .. node[arr,below] {$\cC$} (t);
      \draw [a] (s) -- node[arr,above] {$\cC$} (t);
      \node [cell] at (0,2.5) {$=$};
    \end{tikzpicture}
    \eq
    \begin{tikzpicture}[xscale=1.125,yscale=-1.125]
      \path (0,11) -- (0,-4);
      \node [ob] (s) at (-7,0) {$\bS$};
      \node [ob] (m) at (-3,7) {$\bS$};
      \node [ob] (t) at (7,0) {$\bS$};
      \draw [a] (s) -- node[arr,bl] {$\cC$} (m);
      \draw [a] (m) -- node[arr,al,pos=.7] {$\cC$} (t);
      \draw [a] (s) -- node[arr,above] {$\cC$} (t);
      \draw [eq] (m) .. controls +(6.5,1) and +(-1,6.5) .. (t);
      \node [cell] at (-1.5,2.5) {$\delta$};
      \node [cell] at (3,4.75) {$\varepsilon$};
    \end{tikzpicture}
  \]
  We refer to comonads by their carrying endofunctors $\cC$.
  
  A \emph{comonad map} between comonads $\cC$ and $\cD$ on $\bS$ is a
  natural transformation $\phi \colon \cC \Rightarrow \cD$ such that
  the following homomorphism laws hold:
  \vspace{-5pt}
  \[
    \begin{tikzpicture}[xscale=.75,yscale=-1]
      \node [ob] (s) at (-10,0) {$\bS$};
      \node [ob] (m) at (0,7) {$\bS$};
      \node [ob] (t) at (10,0) {$\bS$};
      \draw [a] (s) -- node[arr,arcl,pos=.3] {$\cC$} (m);
      \draw [a] (m) -- node[arr,alcl,pos=.7] {$\cC$} (t);
      \draw [a] (s) -- node[arr,above] {$\cC$} (t);
      \draw [a] (s) .. controls +(-1,6) and +(-6,2) .. node[arr,bl] {$\cD$} (m);
      \draw [a] (m) .. controls +(6,2) and +(1,6) .. node[arr,br] {$\cD$} (t);
      \node [cell] at (0,2.5) {$\delta$};
      \node [cell] at (-6.5,5) {$\phi$};
      \node [cell] at (6.5,5) {$\phi$};
      \path [a] (s) .. controls +(6,-6) and +(-6,-6) .. node[arr,above] {$\phantom{\cC}$} (t);
    \end{tikzpicture}
    \eqq
    \begin{tikzpicture}[xscale=.75,yscale=-1]
      \node [ob] (s) at (-10,0) {$\bS$};
      \node [ob] (m) at (0,7) {$\bS$};
      \node [ob] (t) at (10,0) {$\bS$};
      \draw [a] (s) -- node[arr,bl] {$\cD$} (m);
      \draw [a] (m) -- node[arr,br] {$\cD$} (t);
      \draw [a] (s) -- node[arr,above=-2,pos=.2] {$\cD$} (t);
      \draw [a] (s) .. controls +(6,-6) and +(-6,-6) .. node[arr,above] {$\cC$} (t);
      \node [cell] at (0,2.5) {$\delta$};
      \node [cell] at (0,-2.5) {$\phi$};
      \path (s) .. controls +(-1,6) and +(-6,2) .. node[arr,bl] {$\phantom{\cD}$} (m);
      \path (m) .. controls +(6,2) and +(1,6) .. node[arr,br] {$\phantom{\cD}$} (t);
    \end{tikzpicture}
    \qquad\qquad
    \begin{tikzpicture}[xscale=.75,yscale=-1]
      \node [ob] (s) at (-10,0) {$\bS$};
      \node [ob] (t) at (10,0) {$\bS$};
      \draw [a] (s) -- node[arr,above] {$\cC$} (t);
      \draw [eq] (s) .. controls +(6,6) and +(-6,6) .. (t);
      \node [cell] at (0,2.3) {$\varepsilon$};
      \path (s) .. controls +(6,-6) and +(-6,-6) .. node[arr,above] {$\phantom{\cC}$} (t);
      \node [ob] (m) at (0,7) {$\phantom{\bS}$};
      \path (s) .. controls +(-1,6) and +(-6,2) .. node[arr,bl] {$\phantom{\cD}$} (m);
    \end{tikzpicture}
    \eqq
    \begin{tikzpicture}[xscale=.75,yscale=-1]
      \node [ob] (s) at (-10,0) {$\bS$};
      \node [ob] (t) at (10,0) {$\bS$};
      \draw [a] (s) -- node[arr,above=-2,pos=.2] {$\cD$} (t);
      \draw [eq] (s) .. controls +(6,6) and +(-6,6) .. (t);
      \draw [a] (s) .. controls +(6,-6) and +(-6,-6) .. node[arr,above] {$\cC$} (t);
      \node [cell] at (0,2.3) {$\varepsilon$};
      \node [cell] at (0,-2.5) {$\phi$};
      \node [ob] (m) at (0,7) {$\phantom{\bS}$};
      \path (s) .. controls +(-1,6) and +(-6,2) .. node[arr,bl] {$\phantom{\cD}$} (m);
    \end{tikzpicture}
  \]
\end{definition}

\tikzset{every picture/.style={}}

We denote the category of comonads on $\bS$ and comonad maps by
$\Com(\bS)$.

\begin{remark}
  Every adjunction
  \[
    \begin{tikzcd}[column sep=30pt]
      \bB\ar[r, shift left=5pt, "L"]\ar["\bot"{anchor=center},r, draw=none]&
      \bS\ar[l, shift left=5pt, "R"]
    \end{tikzcd}
  \]
  induces a comonad $\cC \coloneqq L \circ R$ on $\bS$, with the comonad
  counit given by the adjunction counit
  $\varepsilon \colon L \circ R \Rightarrow \id_{\bS}$ and the comonad
  comultiplication obtained from the adjunction unit
  $\eta \colon \id_{\bB} \Rightarrow R \circ L$ as
  $\delta \coloneqq L \circ \eta \circ R \colon L \circ R \Rightarrow L \circ R \circ L
  \circ R$.
\end{remark}

\begin{example}[Topological spaces]\label{ex:top}
  As mentioned in the introduction, topological spaces induce comonads
  on $\Set$. For every topological space $\tX$ we have an adjunction
  \[
    \begin{tikzcd}[column sep=50pt]
      \Sh(\tX)\ar[r, shift left=6pt,
      "\parbox{50pt}{\centering\tiny sum of\\ stalks}"]\ar["\bot"{anchor=center},r, draw=none]&
      \Set\ar[l, shift left=6pt, "\parbox{50pt}{\centering\tiny sheaf of\\ functions}"]
    \end{tikzcd}
  \]
  
  The left adjoint $\Sh(\tX) \to \Set$ sends a sheaf $F$ to the sum of
  its stalks $\sum_{x\in \pts{\tX}} \colim_{U \ni x} F(U)$ (where $U$
  ranges over neighborhoods\footnote{It does not matter whether we use
    open neighborhoods or arbitrary neighborhoods in this formula; in
    either case the colimit is the same.} of $x$). This is the set of
  points in the \'etal\'e space corresponding to the sheaf. The right
  adjoint $\Set \to \Sh(\tX)$ sends a set $\xA$ to the sheaf
  $U \mapsto \xA^U$ of $\xA$-valued functions. Points in the \'etal\'e
  space corresponding to this sheaf are then precisely germs of
  $\xA$-valued functions.

  We denote the induced comonad on $\Set$ by $\cX$, given
  by the formula
  \[\cX(\xA) = \sum_{x\in \pts{\tX}} \colim_{U \ni x} \xA^U\]
  Thus $\cX(\xA)$ is the set of all germs of functions from
  $\tX$ into $\xA$. We will sometimes refer to $\cX$ as the
  \emph{germs comonad} of $\tX$.

  The counit $\varepsilon_\xA \colon \cX(\xA) \to \xA$ sends
  the germ $[f]_x$ of a function $f \colon U \to \xA$ about a point
  $x$ in $U$ to the element $f(x)$ in $\xA$. The comultiplication
  $\delta_\xA \colon \cX(\xA) \to
  \cX(\cX(\xA))$ sends the germ $[f]_x$ to the
  germ $[[f]_{(\dash)}]_x$ of the function
  $[f]_{(\dash)} \colon U \to \cX(A)$ that sends each point
  $y$ in $U$ to the germ $[f]_y$.

  Comonad maps do \emph{not} correspond to continuous functions of
  topological spaces. We will see later in \cref{cor:topcom} that a
  comonad map from $\cX$ to $\cY$ instead
  corresponds to an \emph{\'etale map from a refinement of $\tX$ to
    $\tY$}.
\end{example}

\begin{example}[Categories]\label{ex:cat}
  As also mentioned in the introduction, small categories induce comonads on
  $\Set$. For every category $\bC$ we have an adjunction
  \[
    \begin{tikzcd}[column sep=60pt]
      \fun{\bC}{\Set}\ar[r, shift left=6pt,
      "\obs{\El(\dash)}"]\ar["\bot"{anchor=center},r, draw=none]&
      \Set\ar[l, shift left=6pt, "c \,\mapsto\, (\dash)^{\obs{c / \bC}}"]
    \end{tikzcd}
  \]

  The left adjoint $\fun{\bC}{\Set} \to \Set$ sends
  $F \colon \bC \to \Set$ to the sum of elements
  $\sum_{c\in \obs{\bC}} F(c)$. This is the set of objects of the
  category of elements $\obs{\El(F)}$. The right adjoint
  $\Set \to \fun{\bC}{\Set}$ sends a set $\xA$ to the functor $\bC \to \Set$
  defined by $c \mapsto \xA^{\obs{c / \bC}}$, where $\obs{c / \bC}$ denotes
  the set of all arrows out of $c$. If we regard the set of arrows out
  of an object as its ``infinitesimal neighborhood'', then we may view
  the elements of this functor as ``germs'' of $\xA$-valued functions.

  We denote the induced comonad on $\Set$ by $\cbC$, given
  by the formula
  \[\cbC(\xA) = \sum_{c\in \obs{\bC}} \xA^{\obs{c / \bC}}.\]

  The counit $\varepsilon_\xA \colon \cbC(\xA) \to \xA$ returns
  values at identities, sending a function
  $f \colon \obs{c / \bC} \to \xA$ to the element $f(\id_{c})$ in
  $\xA$. The comultiplication
  $\delta_\xA \colon \cbC(\xA) \to
  \cbC(\cbC(\xA))$ returns values at compositions,
  sending $f$ to the function $\obs{c / \bC} \to \cbC(\xA)$
  that takes an arrow $p \colon c \to d$ to the function that takes an
  arrow $q \colon d \to e$ and returns $f(q \circ p)$ in $\xA$.


  Comonad maps do \emph{not} correspond to functors of categories. We
  will see later in \cref{thm:doublecat} (as already observed in
  \cite{ahman-uustalu:updates}) that a comonad map from
  $\cbC$ to $\cbD$ instead corresponds to a
  \emph{retrofunctor} (a.k.a.\
  \emph{cofunctor}~\cite{higgins-mackenzie,aguiar,dimeglio:lenses,ahman-uustalu:updates,niu-spivak})
  from $\bC$ to $\bD$. A retrofunctor consists of a map
  $\one{\phi} \colon \obs{\bC} \to \obs{\bD}$ between sets of objects
  and maps \emph{backwards} from the arrows out of each
  $\one{\phi}(c)$ in $\bD$ to arrows out of $c$ in $\bC$, satisfying
  unit and associativity conditions; it can also be equivalently
  described as a span of functors
  $\bC \xot{\lambda} \bB \xto{\rho} \bD$ where $\lambda$ is bijective
  on objects and $\rho$ is a discrete opfibration.
\end{example}

\tikzset{every picture/.style={baseline={([yshift=-3.5pt]current bounding box.center)},scale=.125}}

\begin{definition}
  A \emph{coalgebra} of a comonad $\cC$ on $\bS$ is an object $\eE$ of
  $\bS$ and a map $\hH \colon \eE \to \cC(\eE)$ such that the following
  associativity and unit laws hold:
  \[
    \begin{tikzpicture}[xscale=.75,yscale=-1.125]
      \node [ob] (s) at (-14,0) {$\point$};
      \node [ob] (mm) at (-9,7) {$\bS$};
      \node [ob] (m) at (2,7) {$\bS$};
      \node [ob] (t) at (7,0) {$\bS$};
      \draw [a] (s) -- node[arr,bl] {$\eE$} (mm);
      \draw [a] (mm) -- node[arr,below] {$\cC$} (m);
      \draw [a] (s) -- node[arr,ar,pos=.3] {$\eE$} (m);
      \draw [a] (m) -- node[arr,br] {$\cC$} (t);
      \draw [a] (s) -- node[arr,above] {$\eE$} (t);
      \node [cell] at (0,2.5) {$\hH$};
      \node [cell] at (-7,5) {$\hH$};
    \end{tikzpicture}
    \eqq
    \begin{tikzpicture}[xscale=.75,yscale=-1.125]
      \node [ob] (s) at (-7,0) {$\point$};
      \node [ob] (m) at (-2,7) {$\bS$};
      \node [ob] (mm) at (9,7) {$\bS$};
      \node [ob] (t) at (14,0) {$\bS$};
      \draw [a] (m) -- node[arr,below] {$\cC$} (mm);
      \draw [a] (mm) -- node[arr,br] {$\cC$} (t);
      \draw [a] (s) -- node[arr,bl] {$\eE$} (m);
      \draw [a] (m) -- node[arr,al,pos=.7] {$\cC$} (t);
      \draw [a] (s) -- node[arr,above] {$\eE$} (t);
      \node [cell] at (0,2.5) {$\hH$};
      \node [cell] at (7,5) {$\delta$};
    \end{tikzpicture}
    \qquad\qquad
    \begin{tikzpicture}[xscale=.625,yscale=-1.125]
      \path (0,11) -- (0,-4);
      \node [ob] (s) at (-14,0) {$\point$};
      \node [ob] (t) at (7,0) {$\bS$};
      \draw [a] (s) .. controls +(4,6) and +(-4,6) .. node[arr,below] {$\eE$} (t);
      \draw [a] (s) -- node[arr,above] {$\eE$} (t);
      \node [cell] at (-3.5,2.5) {$=$};
    \end{tikzpicture}
    \eqq
    \begin{tikzpicture}[xscale=.75,yscale=-1.125]
      \path (0,11) -- (0,-4);
      \node [ob] (s) at (-7,0) {$\point$};
      \node [ob] (m) at (-2,7) {$\bS$};
      \node [ob] (t) at (14,0) {$\bS$};
      \draw [a] (s) -- node[arr,bl] {$\eE$} (m);
      \draw [a] (m) -- node[arr,al,pos=.7] {$\cC$} (t);
      \draw [a] (s) -- node[arr,above] {$\eE$} (t);
      \draw [eq] (m) .. controls +(8,1.5) and +(-1.5,8) .. (t);
      \node [cell] at (0,2.5) {$\hH$};
      \node [cell] at (7.25,5.25) {$\varepsilon$};
    \end{tikzpicture}
  \]
  That is, $\cC(\hH) \circ \hH = \delta_\eE \circ \hH$ and
  $\varepsilon_\eE \circ \hH = \id_\eE$. We refer to coalgebras by
  their carrying objects $\eE$.
  
  A \emph{coalgebra map} from a coalgebra $\eE$ to a coalgebra $\eE'$ is a
  map $\hK \colon \eE \to \eE'$ in $\bS$ such that the following homomorphism law
  holds: \vspace{-5pt}
  \[
    \begin{tikzpicture}[xscale=.75,yscale=-1]
      \node [ob] (s) at (-10,0) {$\point$};
      \node [ob] (m) at (0,7) {$\bS$};
      \node [ob] (t) at (10,0) {$\bS$};
      \draw [a] (s) -- node[arr,arcl,pos=.2] {$\eE$} (m);
      \draw [a] (m) -- node[arr,br] {$\cC$} (t);
      \draw [a] (s) -- node[arr,above] {$\eE$} (t);
      \draw [a] (s) .. controls +(-1,6) and +(-6,2) .. node[arr,bl] {$\eE'$} (m);
      \node [cell] at (0,2.5) {$\hH$};
      \node [cell] at (-6.5,5) {$\hK$};
      \path (s) .. controls +(6,-6) and +(-6,-6) .. node[arr,above] {$\phantom{\eE}$} (t);
    \end{tikzpicture}
    \eqq
    \begin{tikzpicture}[xscale=.75,yscale=-1]
      \node [ob] (s) at (-10,0) {$\point$};
      \node [ob] (m) at (0,7) {$\bS$};
      \node [ob] (t) at (10,0) {$\bS$};
      \draw [a] (s) -- node[arr,bl] {$\eE'$} (m);
      \draw [a] (m) -- node[arr,br] {$\cC$} (t);
      \draw [a] (s) -- node[arr,above=-2,pos=.2] {$\eE'$} (t);
      \draw [a] (s) .. controls +(6,-6) and +(-6,-6) .. node[arr,above] {$\eE$} (t);
      \node [cell] at (0,2.5) {$\hH'$};
      \node [cell] at (0,-2.5) {$\hK$};
      \path (s) .. controls +(-1,6) and +(-6,2) .. node[arr,bl] {$\phantom{B}$} (m);
    \end{tikzpicture}
  \]
  That is, $\cC(\hK) \circ \hH = \hH' \circ \hK$.
\end{definition}

\tikzset{every picture/.style={}}

We denote the category of coalgebras of a comonad $\cC$ on $\bS$ by
$\CS$, also known as the \emph{Eilenberg-Moore category} of $\cC$.
Extracting the carrier set induces a faithful functor
$\carC\colon\CS\to\bS$, and this has a canonical right adjoint
$\rcarC\colon \bS \to \CS$ sending each object $\xA$ to
$\delta \colon \cC(\xA) \to \cC(\cC(\xA))$, the \emph{cofree coalgebra
  on $\xA$}. We denote by $\KlCS \hookrightarrow \CS$ the full
subcategory of cofree coalgebras, also known as the \emph{Kleisli
  category} of $\cC$. Either adjunction induces the comonad $\cC$.

A functor is called \emph{comonadic} if it is essentially of the form
$\carC \colon \CS \to \bS$, the forgetful functor from a category of
coalgebras. In \cref{sec:comonadicity} we review some standard results
on comonadicity, which are used throughout the paper.

\begin{example}[Sheaves and presheaves]\label{ex:sheaves}\label{ex:copresheaves}
  The left adjoint sum of stalks functor $\Sh(\tX) \to \Set$ from
  \cref{ex:top}, which induces the comonad $\cX$
  corresponding to the topological space $\tX$, is comonadic. That is
  to say, $\cX$-coalgebras are equivalent to sheaves over
  $\tX$, or equivalently \'etal\'e spaces over $\tX$.

  To prove this, it suffices to observe that the sum of stalks functor
  factors as
  \[\Sh(\tX) \xto{i^*} \Set^{\pts{\tX}} \xto{\slicel} \Set\] where the first factor
  $i^*$ sends an \'etal\'e space $E \to \tX$ to its
  $\pts{\tX}$-indexed set of fibers (itself an \'etal\'e space over
  the discrete space $\pts{\tX}$, as obtained by pullback along the
  inclusion of the set of points $i \colon \pts{\tX} \to \tX$), and
  the second factor $\slicel$ sums up all the values over all
  points. In topos-theoretic terminology, $i^*$ arises as the
  \emph{inverse image} of the \emph{geometric surjection of
    toposes}~\cite[A4.2]{johnstone:elephant} given by the inclusion of
  the points $i \colon \pts{\tX} \to \tX$. It is well-known and
  readily verified that both factors are conservative,
  pullback-preserving, and left adjoint, and these conditions are
  sufficient to deduce comonadicity from the \emph{crude comonadicity
    theorem} (\cref{prop:crude}).

  Explicitly, a coalgebra $\eE \to \cX(\eE)$ of the comonad
  $\cX$ on $\Set$ corresponding to the topological space
  $\tX$ assigns to each element $e$ in $\eE$ the germ of an
  $\eE$-valued function about some point $x$ in $\tX$. In other words,
  each element $e$ in $\eE$ is assigned its own ``infinitesimal
  neighborhood'' of elements from $\eE$, locally resembling the
  neighborhood of a point $x$ in $\tX$. The unit law ensures $e$ in
  $\eE$ is itself assigned to the center of its own neighborhood, and
  the associativity law then does the rest to ensure the neighborhoods
  cohere together as an \'etal\'e space.

  Likewise, the left adjoint elements functor
  $\fun{\bC}{\Set} \to \Set$ from \cref{ex:cat}, which induces the
  comonad $\cbC$ corresponding to the category $\bC$, is
  comonadic. That is to say, $\cbC$-coalgebras are
  equivalent to functors $\bC \to \Set$ (a.k.a.\ \emph{$\bC$-sets} or
  \emph{$\bC$-actions} or \emph{$\bC$-copresheaves}).

  As before, we observe that the elements functor factors as
  \[\fun{\bC}{\Set} \xto{\reindex{i}} \Set^{\obs{\bC}} \xto{\slicel} \Set\]
  where the first factor $\reindex{i}$ simply reindexes along the
  inclusion of the set of objects $i \colon \obs{\bC} \to \bC$, and
  the second factor $\slicel$ sums up all the values over all
  objects. (Again, $\reindex{i}$ is the inverse image of a geometric
  surjection of toposes; in \cref{bg:II} we take a closer look at the
  general case of this construction, which applies to any topos
  equipped with a surjection from a set of points, a.k.a.\
  \emph{ionad}.)  Both functors are conservative,
  pullback-preserving,\footnote{In fact they are
    \emph{wide-pullback-preserving} functors between presheaf
    categories, which are equivalently the
    \emph{connected-limit-preserving} or \emph{parametric right
      adjoint} functors, also sometimes called \emph{polynomial}
    functors~\cite{spivak-garner-fairbanks}.} and left adjoint, hence
  the composite is comonadic by
  the crude comonadicity theorem (\cref{prop:crude}).

  Explicitly, a coalgebra $\eE \to \cbC(\eE)$ of the comonad
  $\cbC$ on $\Set$ from \cref{ex:cat} corresponding to the
  category $\bC$ assigns to each element $e$ in $\eE$ a function
  $\obs{c / \bC} \to \eE$ for some object $c$ in $\bC$. In other
  words, each element $e$ in $\eE$ is assigned an object $c$ in $\bC$
  and actions on $e$ by all arrows out of $c$. The unit law ensures
  $e$ in $\eE$ is sent to itself by the assigned identity in $\bC$,
  and the associativity law then does the rest to ensure the actions
  by arrows of $\bC$ on elements of $\eE$ are associative (and in
  particular compatible along assigned objects of $\bC$).
\end{example}

For any comonad map $\phi\colon\cC\to\cD$, there is an
associated functor
\[
  \begin{tikzcd}[column sep=15pt, row sep=15pt]
    \CS\ar[rr, dashed, "\com{\bS}{\phi}"]\ar[dr, "\carC"']&&\DS\ar[dl, "\carD"]\\
    &\bS
  \end{tikzcd}
\]
sending each coalgebra $\hH:\eE\to\cC(\eE)$ to its composite with
$\phi_\eE \colon \cC(\eE) \to \cD(\eE)$.


\tikzset{every picture/.style={baseline={([yshift=-3.5pt]current bounding box.center)},scale=.125}}
\begin{lemma}\label{lem:componentwisemonic}
  For any comonad map $\phi$, the functor $\com{\bS}{\phi}$ is
  faithful. If $\phi$ is componentwise monic, then $\com{\bS}{\phi}$ is
  also full.
\end{lemma}
\begin{proof}
  Since $\carC$ is faithful, so is $\com{\bS}{\phi}$ by
  cancellation. For the second statement, consider an arbitrary
  coalgebra map
  $\hK \colon \com{\bS}{\phi}(\eE) \to \com{\bS}{\phi}(\eE')$ and cancel
  the monic component of $\phi$ as follows:
  \[
    \begin{tikzpicture}[xscale=.75,yscale=-1]
      \node [ob] (s) at (-10,0) {$\point$};
      \node [ob] (m) at (0,7) {$\bS$};
      \node [ob] (t) at (10,0) {$\bS$};
      \draw [a] (s) -- node[arr,arcl,pos=.2] {$\eE$} (m);
      \draw [a] (m) -- node[arr,alcl,pos=.7] {$\cC$} (t);
      \draw [a] (s) -- node[arr,above] {$\eE$} (t);
      \draw [a] (s) .. controls +(-1,6) and +(-6,2) .. node[arr,bl] {$\eE'$} (m);
      \draw [a] (m) .. controls +(6,2) and +(1,6) .. node[arr,br] {$\cD$} (t);
      \node [cell] at (0,2.5) {$\hH$};
      \node [cell] at (-6.5,5) {$\hK$};
      \node [cell] at (6.5,5) {$\phi$};
      \path (s) .. controls +(6,-6) and +(-6,-6) .. node[arr,above] {$\phantom{\eE}$} (t);
    \end{tikzpicture}
    \eqq
    \begin{tikzpicture}[xscale=.75,yscale=-1]
      \node [ob] (s) at (-10,0) {$\point$};
      \node [ob] (m) at (0,7) {$\bS$};
      \node [ob] (t) at (10,0) {$\bS$};
      \draw [a] (s) -- node[arr,bl] {$\eE'$} (m);
      \draw [a] (m) -- node[arr,alcl,pos=.7] {$\cC$} (t);
      \draw [a] (s) -- node[arr,above=-2,pos=.2] {$\eE'$} (t);
      \draw [a] (s) .. controls +(6,-6) and +(-6,-6) .. node[arr,above] {$\eE$} (t);
      \draw [a] (m) .. controls +(6,2) and +(1,6) .. node[arr,br] {$\cD$} (t);
      \node [cell] at (0,2.5) {$\hH$};
      \node [cell] at (0,-2.5) {$\hK$};
      \node [cell] at (6.5,5) {$\phi$};
    \end{tikzpicture}
    \qquad\implies\qquad
    \begin{tikzpicture}[xscale=.75,yscale=-1]
      \node [ob] (s) at (-10,0) {$\point$};
      \node [ob] (m) at (0,7) {$\bS$};
      \node [ob] (t) at (10,0) {$\bS$};
      \draw [a] (s) -- node[arr,arcl,pos=.2] {$\eE$} (m);
      \draw [a] (m) -- node[arr,br] {$\cC$} (t);
      \draw [a] (s) -- node[arr,above] {$\eE$} (t);
      \draw [a] (s) .. controls +(-1,6) and +(-6,2) .. node[arr,bl] {$\eE'$} (m);
      \node [cell] at (0,2.5) {$\hH$};
      \node [cell] at (-6.5,5) {$\hK$};
      \path (s) .. controls +(6,-6) and +(-6,-6) .. node[arr,above] {$\phantom{\eE}$} (t);
    \end{tikzpicture}
    \eqq
    \begin{tikzpicture}[xscale=.75,yscale=-1]
      \node [ob] (s) at (-10,0) {$\point$};
      \node [ob] (m) at (0,7) {$\bS$};
      \node [ob] (t) at (10,0) {$\bS$};
      \draw [a] (s) -- node[arr,bl] {$\eE'$} (m);
      \draw [a] (m) -- node[arr,br] {$\cC$} (t);
      \draw [a] (s) -- node[arr,above=-2,pos=.2] {$\eE'$} (t);
      \draw [a] (s) .. controls +(6,-6) and +(-6,-6) .. node[arr,above] {$\eE$} (t);
      \node [cell] at (0,2.5) {$\hH$};
      \node [cell] at (0,-2.5) {$\hK$};
      \path (s) .. controls +(-1,6) and +(-6,2) .. node[arr,bl] {$\phantom{\eE'}$} (m);
    \end{tikzpicture}
  \]
\end{proof}
\tikzset{every picture/.style={}}

The comonad map $\phi$ can conversely be recovered from the functor
$\com{\bS}{\phi}$:

\begin{proposition}[\cite{street:monads}]\label{prop:respect}
  Comonad maps $\phi \colon \cC \to \cD$ are in one-to-one
  correspondence with functors
  $\com{\bS}{\phi} \colon \CS \to \DS$ making the triangle
  \[
    \begin{tikzcd}[column sep=15pt, row sep=15pt]
      \CS\ar[rr, "\com{\bS}{\phi}"]\ar[dr, "\carC"']&&\DS\ar[dl, "\carD"]\\
      &\bS
    \end{tikzcd}
  \]
  commute.\qed
\end{proposition}

The following standard fact will also be useful.



\begin{lemma}[{\cite[Theorem VI.3.1]{maclane}}]\label{lem:comparison}
  For any adjunction
  \[
    \begin{tikzcd}[column sep=30pt]
      \bB\ar[r, shift left=5pt, "L"]\ar["\bot"{anchor=center},r, draw=none]&
      \bS\ar[l, shift left=5pt, "R"]
    \end{tikzcd}
  \]
  inducing $\cC$, there is a unique \emph{comparison functor}
  $K \colon \bB \to \CS$ such that the following triangles
  commute:
  \[
    \begin{tikzcd}[column sep=15pt, row sep=15pt]
      {\bB} && \CS \\
      & \bS
      \arrow["K", from=1-1, to=1-3]
      \arrow["L"', from=1-1, to=2-2]
      \arrow["\carC", from=1-3, to=2-2]
    \end{tikzcd}
    \qquad\qquad
    \begin{tikzcd}[column sep=15pt, row sep=15pt]
      {\bB} && \CS \\
      & \bS
      \arrow["K", from=1-1, to=1-3]
      \arrow["R", from=2-2, to=1-1]
      \arrow["\rcarC"', from=2-2, to=1-3]
    \end{tikzcd}
    \tag*{\qed}
  \]
\end{lemma}

\section{Comonads on $\Set$}\label[bgchapter]{sec:onset}

It is not obvious that the definition of topological space, featuring
arbitrary unions but only finite intersections, should fall out of
the definition of comonad on $\Set$. The relevance of finite
intersections traces to the following fact:

\begin{proposition}[\cite{trnkova}]\label{prop:trnkova}
  Every endofunctor on $\Set$ preserves nonempty pairwise
  intersections, i.e.\ pullback squares consisting of monomorphisms
  where the pullback is nonempty.
\end{proposition}
\begin{proof}
  We will show the stronger fact that every such pullback is absolute,
  as observed in \cite{manes:classes}. More specifically, a
  \emph{split pullback} is a commutative square
  \[\begin{tikzcd}[column sep=20]
      \xU & {\xS_1} \\
      {\xS_2} & \xT
      \arrow["{m_1}", from=1-1, to=1-2]
      \arrow["{m_2}"', from=1-1, to=2-1]
      \arrow["{f_1}", from=1-2, to=2-2]
      \arrow["{f_2}"', from=2-1, to=2-2]
    \end{tikzcd}
  \]
  such that there exist retractions $e_1 \circ m_1 = \id_\xU$,
  $r_1 \circ f_1 = \id_{\xS_1}$, and $r_2 \circ f_2 = \id_{\xS_2}$. It is
  straightforward to verify that split pullbacks are
  pullbacks~\cite[Lemma 1.8(4)]{manes:classes}. Since the definition
  is equational, split pullbacks are preserved by all functors.


  Now suppose given a pullback diagram as above in $\Set$. If $f_1$
  and $f_2$ are monomorphisms and $\xU$ is nonempty, then in particular
  $m_1$ admits a retraction $e_1$ and $f_1$ admits a retraction
  $r_1$. For this pullback to be a split pullback, we require a
  retraction $r_2 \colon \xT \to \xS_2$ of $f_2$ satisfying
  $r_2 \circ f_1 = m_2 \circ e_1$.

  Observe that the object $\xT$ decomposes as
  $\xU + (\xS_1 \setminus \xU) + (\xS_2 \setminus \xU) + (\xT \setminus (\xS_1 +_\xU
  \xS_2) )$, as depicted by the following Venn diagram:
  
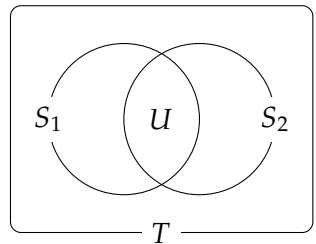
\[
    \begin{tikzpicture}
      \draw (-.5,0) circle (1);
      \draw (.5,0) circle (1);
      \draw[rounded corners] (-2,-1.5) rectangle (2,1.5);
      \node at (0,0) {$\xU$};
      \node [fill=white] at (-1.5,0) {$\xS_1$};
      \node [fill=white] at (1.5,0) {$\xS_2$};
      \node [fill=white] at (0,-1.5) {$\xT$};
    \end{tikzpicture}
  \]
  We define $r_2$ piecewise, given by $m_2 \circ e_1$ as desired on
  $\xS_1$, given by identity as desired on $\xS_2$, and chosen arbitrarily
  elsewhere (e.g.\ given by $m_2 \circ e_1 \circ r_1$).\footnote{More
    generally by the same argument, any pullback square consisting of split
    monomorphisms is a split pullback in any Boolean topos. One can
    also show that in any category, a pullback square consisting of monomorphisms
    between \emph{injective} objects is a split pullback if the
    comparison map from the corresponding pushout square is monic (as
    is the case when pairwise unions of subobjects are calculated as
    pushouts under intersections --- such as in any \emph{coherent
      category}~\cite[Proposition A1.4.3]{johnstone:elephant}).}
\end{proof}

\begin{corollary}\label{prop:intersections}
  Every comonad on $\Set$ preserves finite intersections, i.e.\
  pullback squares consisting of monomorphisms.\footnote{One can also
    show that an endofunctor on $\Set$ preserves finite intersections
    if and only if it preserves the equalizer of the two distinct maps
    $1 \rightrightarrows 2$.}
\end{corollary}
\begin{proof}
  In light of Trnkova's result \cref{prop:trnkova}, it remains to show
  that comonads on $\Set$ preserve pairwise intersections that are
  empty.

  In general, if $\cC$ is an endofunctor on $\Set$ (or any category
  with an object $\varnothing$ that has only isomorphisms into it)
  admitting a natural transformation to identity, then $\cC$ preserves
  $\varnothing$ and (co)limits for which the (co)limit object is
  $\varnothing$. Indeed, given a diagram $D \colon J \to \bS$ with
  (co)limit $\varnothing$, we have a natural transformation
  $\varepsilon \circ D \colon \cC \circ D \Rightarrow D$, inducing a map
  from the (co)limit of $\cC \circ D$ to
  $\varnothing \cong \cC(\varnothing)$.
\end{proof}

It is a standard fact (which we recall in \cref{lem:create}) that any
comonadic functor $\carC \colon \CS \to \bS$ creates all colimits, as
well as all limits that are preserved by $\cC$ and $\cC \circ
\cC$. Hence \cref{prop:intersections} implies that comonadic functors
into $\Set$ respect all the structure relevant to general topology:
they create colimits --- which are analogous to unions --- and finite
intersections.

A useful tool in the general study of comonads is the \emph{crude
  comonadicity theorem} (\cref{prop:crude}), already appealed to in
\cref{ex:sheaves}, which identifies \emph{preservation of coreflexive
  equalizers} as a sufficient condition for a conservative left
adjoint functor from a category with coreflexive equalizers to be
comonadic. While this may seem a rather technical condition, we will
see that under mild assumptions, it is simply equivalent to
preservation of finite intersections.

\begin{definition}\label{def:crep}
  A \emph{coreflexive pair} is a pair of parallel arrows
  $f, g \colon \xS \to \xT$ with a mutual retraction $r$. That is,
  $r \circ f = \id_\xS = r \circ g$.
  \[\begin{tikzcd}
      \xS & \xT
      \arrow["g"', shift right=3, from=1-1, to=1-2]
      \arrow["f", shift left=3, from=1-1, to=1-2]
      \arrow["r"{description}, from=1-2, to=1-1]
    \end{tikzcd}\]
  A \emph{coreflexive equalizer} is an equalizer of a coreflexive
  pair.
  
  For short, we will often call a comonad that preserves coreflexive
  equalizers \emph{crude}.
\end{definition}

The following tells us that coreflexive equalizer diagrams may
equivalently be viewed as certain pullback squares of monomorphisms,
and so indeed they are a special case of finite intersections.

\begin{lemma}\label{lem:coreflpull}
  Suppose $f, g \colon \xS \to \xT$ are a coreflexive pair. Then
  $m \colon \xU \to \xS$ is an equalizer of $f$ and $g$ if and only if
  \[\begin{tikzcd}[column sep=20pt]
      \xU & \xS \\
      \xS & \xT
      \arrow["m", from=1-1, to=1-2]
      \arrow["m"', from=1-1, to=2-1]
      \arrow["\lrcorner"{anchor=center, pos=0.125}, draw=none, from=1-1, to=2-2]
      \arrow["f", from=1-2, to=2-2]
      \arrow["g"', from=2-1, to=2-2]
    \end{tikzcd}\]
  is a pullback square. 
\end{lemma}
\begin{proof}
  The equalizer and the pullback satisfy the same universal property,
  since there no difference between the two types of cones: if
  $f \circ m_1 = g \circ m_2$, then composing with the mutual
  retraction of $f$ and $g$ yields $m_1 = m_2$.
\end{proof}

Therefore \cref{prop:intersections} recovers the following result from
\cite{linton-pare}, due to which the theory of comonads on $\Set$
works particularly smoothly.

\begin{corollary}[{\cite{linton-pare}}]\label{lem:corefl}
  Every comonad on $\Set$ is crude, i.e.\ preserves coreflexive
  equalizers.\footnote{In~\cite{linton-pare}, Linton and Par{\'e}
    characterize precisely the toposes on which every comonad is
    crude: those in which every non-initial object is injective. All
    numbered results about $\Set$ in this section hold more generally
    replacing $\Set$ with any topos satisfying that property. In
    particular, every well-pointed topos satisfies that property.}\qed
  %
\end{corollary}

We now recall the standard concept of regular monomorphism, which will
be relevant in our study of comonads on $\Set$ as spaces: in
\cref{sec:toptocomonad} we will define open subsets of comonads on
$\Set$ as regular subobjects of the terminal object.

\begin{definition}\label{def:coregular}
  A \emph{regular monomorphism} is an arrow that arises as an
  equalizer. Assuming its cokernel pair (i.e.\ pushout along itself)
  exists, an arrow is a regular monomorphism if and only if it is the
  equalizer of its cokernel pair.

  A \emph{coregular category} is a category with finite colimits in
  which every arrow factors (necessarily uniquely up to isomorphism)
  as an epimorphism followed by a regular monomorphism and such that
  these factorizations are stable under pushout.
\end{definition}

For example, $\Set$, or any topos, is coregular, and here all
monomorphisms are regular. Regular monomorphisms are always stable
under pullback, assuming the category has pullbacks.

We say that a functor \emph{preserves (regular) finite intersections}
if it preserves (regular) monomorphisms and their pairwise pullbacks.

\begin{lemma}\label{lem:functorpreserve}
  A functor out of a category with pushouts preserves regular finite
  intersections if and only if it preserves coreflexive equalizers.
\end{lemma}
\begin{proof}
  Split monomorphisms (in particular coreflexive pairs) are always
  regular monomorphisms, so the left-to-right direction follows from
  \cref{lem:coreflpull}. For the right-to-left direction, suppose
  given a pullback square consisting of regular monomorphisms:
  \[\begin{tikzcd}
      \xU & {\xS_2} \\
      {\xS_1} & \xT
      \arrow[hook, from=1-1, to=1-2]
      \arrow[hook, from=1-1, to=2-1]
      \arrow["\lrcorner"{anchor=center, pos=0.125}, draw=none, from=1-1, to=2-2]
      \arrow[hook, from=1-2, to=2-2]
      \arrow[hook, from=2-1, to=2-2]
    \end{tikzcd}\]
  We may extend any such diagram to the following larger diagram by
  forming pushouts:
  \[\begin{tikzcd}
      \xU & {\xS_2} & \xT \\
      {\xS_1} & \xT & \cdot \\
      \xT & \cdot & \cdot
      \arrow[hook, from=1-1, to=1-2]
      \arrow[hook, from=1-1, to=2-1]
      \arrow["\lrcorner"{anchor=center, pos=0.125}, draw=none, from=1-1, to=2-2]
      \arrow[hook, from=1-2, to=1-3]
      \arrow[hook, from=1-2, to=2-2]
      \arrow[hook, from=1-3, to=2-3]
      \arrow[hook, from=2-1, to=2-2]
      \arrow[hook, from=2-1, to=3-1]
      \arrow[hook, from=2-2, to=2-3]
      \arrow[hook, from=2-2, to=3-2]
      \arrow["\lrcorner"{anchor=center, pos=0.125, rotate=180}, draw=none, from=2-3, to=1-2]
      \arrow[hook, from=2-3, to=3-3]
      \arrow[hook, from=3-1, to=3-2]
      \arrow["\lrcorner"{anchor=center, pos=0.125, rotate=180}, draw=none, from=3-2, to=2-1]
      \arrow[hook, from=3-2, to=3-3]
      \arrow["\lrcorner"{anchor=center, pos=0.125, rotate=180}, draw=none, from=3-3, to=2-2]
    \end{tikzcd}\] Equivalently, this yields the colimit of the
  zig-zag
  $\xT \hookleftarrow \xS_1 \hookrightarrow \xT \hookleftarrow \xS_2
  \hookrightarrow \xT$; there is a canonical map from the colimit into
  $\xT$ that is a mutual retraction of the three inclusions from
  $\xT$. All of these pushout squares are pullbacks: each regular
  monomorphism is the pullback (equivalently, coreflexive equalizer)
  of its cokernel pair, and, as is readily verified, a pushout of span
  of split monomorphisms is a split pullback (as defined in the proof
  of \cref{prop:trnkova}). Thus by pullback pasting, the outer square
  is a pullback, corresponding to a coreflexive equalizer by
  \cref{lem:coreflpull}. The split pullback is preserved since it is
  absolute, so any functor preserving coreflexive equalizers then
  preserves the original pullback by pullback cancellation.
  %
\end{proof}

Thus \cref{prop:intersections} can conversely be deduced from Linton and
Par\'e's result \cref{lem:corefl}.
We also have that every category of coalgebras of a comonad on $\Set$
is coregular:

\begin{lemma}[{\cite[Regular Categories, Theorem 2.2]{barr-grillet-osdol}}]\label{lem:coregular}\label{lem:factor}
  Let $\cC$ be a comonad on a coregular category $\bS$. Then $\cC$
  preserves regular monomorphisms if and only if
  $\carC \colon \CS \to \bS$ does so.

  Moreover, in this case $\CS$ is coregular, $\carC$ preserves and
  reflects the factorization system, and factorings in $\bS$ of arrows
  in $\CS$ through regular monomorphisms in $\CS$ lift to $\CS$.
\end{lemma}
\begin{proof}
  The statement of~\cite[Regular Categories, Theorem
  2.2]{barr-grillet-osdol} includes all but the last claim about
  factorings, but this too is shown in the proof given there.
\end{proof}

\begin{corollary}\label{lem:regularmono}
  If $\cC$ is a comonad on $\Set$, then $\CSet$ is coregular. The
  regular monomorphisms are the coalgebra maps carried by injections.
\end{corollary}
\begin{proof}
  By \cref{lem:corefl}, every comonad on $\Set$ preserves all
  regular monomorphisms (i.e.\ all monomorphisms), so the conditions
  of \cref{lem:coregular} are satisfied.
\end{proof}

We now expand on the question of comonads preserving weak pullbacks, a
commonly studied class in coalgebra as mentioned in the
introduction. We first recall the definition of weak limits.

\begin{definition}
  A cone to a diagram $D \colon \bJ \to \bC$ is called a
  \emph{weak limit} of $D$ if all other cones to $D$ factor through it
  (but not necessarily uniquely).
\end{definition}

We note a couple of standard facts about weak limits. If a diagram
$D \colon \bJ \to \bC$ has a limit $\lim D$, then a cone to $D$ is a
weak limit if and only if its induced comparison map into $\lim D$ is
a split epimorphism.  Also, if $\bC$ has limits of shape $\bJ$, then a
functor $F \colon \bC \to \bD$ sends weak limits of shape $\bJ$ in
$\bC$ to weak limits in $\bD$ (``\emph{preserves weak limits} of shape
$\bJ$'') if and only if it sends limits of shape $\bJ$ in $\bC$ to
weak limits of shape $\bJ$ in $\bD$ (``\emph{weakly preserves limits}
of shape $\bJ$'').

We also note that preservation of weak pullbacks implies the following
useful condition.

\begin{definition}[\cite{manes},\cite{pare:taut}]
  A functor is \emph{taut} if it preserves preimages, i.e.\ pullbacks
  of monomorphisms along arbitrary arrows.
\end{definition}

\begin{lemma}[{\cite[Lemma 1.2]{johnstone-power-tsujishita-watanabe-worrell}}]\label{lem:weaktaut}
  Every weak-pullback-preserving functor is taut.
\end{lemma}
\begin{proof}
  First we show every weak-pullback-preserving functor $F$ preserves
  monomorphisms. Suppose $m \colon \xU \hookrightarrow \xX$ is a
  monomorphism, as shown below left.
  \[
    \begin{tikzcd}
      \xU & \xU \\
      \xU & \xX
      \arrow[equals, from=1-1, to=1-2]
      \arrow[equals, from=1-1, to=2-1]
      \arrow["\lrcorner"{anchor=center, pos=0.125}, draw=none, from=1-1, to=2-2]
      \arrow["m", hook, from=1-2, to=2-2]
      \arrow["m"', hook, from=2-1, to=2-2]
    \end{tikzcd}
    \qquad\mapsto\qquad
    \begin{tikzcd}[column sep=18pt]
      {F(\xU)} & {F(\xU)} \\
      {F(\xU)} & {F(\xX)}
      \arrow[equals, from=1-1, to=1-2]
      \arrow[equals, from=1-1, to=2-1]
      \arrow["{F(m)}", from=1-2, to=2-2]
      \arrow["{F(m)}"', from=2-1, to=2-2]
    \end{tikzcd}\]
  Then by assumption the cone shown above right is
  a weak pullback, i.e.\ all cones factor through it. But also, since
  its legs are identities, cones factor through it
  uniquely. Thus it is a pullback, i.e.\ $F(m)$ is a monomorphism.

  Now suppose given a pullback square as shown below left.
  \[
    \begin{tikzcd}
      {\xU'} & {\xX'} \\
      \xU & \xX
      \arrow["{m'}", hook, from=1-1, to=1-2]
      \arrow["f"', from=1-1, to=2-1]
      \arrow["\lrcorner"{anchor=center, pos=0.125}, draw=none, from=1-1, to=2-2]
      \arrow["g", from=1-2, to=2-2]
      \arrow["m"', hook, from=2-1, to=2-2]
    \end{tikzcd}
    \qquad\mapsto\qquad
    \begin{tikzcd}[column sep=18pt]
      {F(\xU')} & {F(\xX')} \\
      {F(\xU)} & {F(\xX)}
      \arrow["{F(m')}", hook, from=1-1, to=1-2]
      \arrow["{F(f)}"', from=1-1, to=2-1]
      \arrow["{F(g)}", from=1-2, to=2-2]
      \arrow["{F(m)}"', hook, from=2-1, to=2-2]
    \end{tikzcd}
  \]
  Then by assumption the cone shown above right is a weak pullback,
  i.e.\ all cones factor through it. But also, since one leg is a
  monomorphism, cones factor through it uniquely. Thus it is a pullback.
\end{proof}

\begin{lemma}[{\cite{johnstone-power-tsujishita-watanabe-worrell}}]\label{lem:weakpreserve}
  Let $\cC$ be a comonad on a category with pullbacks $\bS$. If $\cC$
  preserves weak pullbacks, then so does the comonadic functor
  $\carC \colon \CS \to \bS$.\footnote{This is stated in~\cite[Remark
    2.9]{johnstone-power-tsujishita-watanabe-worrell}, but the argument
    supplied there apparently does not work. We instead combine the
    arguments of~\cite[Lemma
    2.8]{johnstone-power-tsujishita-watanabe-worrell}
    and~\cite[Proposition
    3.4]{johnstone-power-tsujishita-watanabe-worrell}.}
\end{lemma}
\begin{proof}
  Let $A \to X \ot B$ be a cospan of $\cC$-coalgebras, and let $P$ be
  their pullback in $\bS$. Since by assumption $\cC$ preserves weak
  pullbacks, we have that $\cC(P)$ is a weak pullback of the cospan
  $\cC(A) \to \cC(X) \ot \cC(B)$. The coalgebra structure maps
  ($A \to \cC(A)$, $X \to \cC(X)$, $B \to \cC(B)$) form a morphism of
  cospans, which therefore induces a map $P \to \cC(P)$ making the
  following commute.
  \[\begin{tikzcd}[row sep=15, column sep=10]
      P && B \\
      & {\cC(P)} & {\cC(B)} \\
      A & {\cC(A)}
      \arrow[from=1-1, to=1-3]
      \arrow[from=1-1, to=2-2]
      \arrow[from=1-1, to=3-1]
      \arrow[from=1-3, to=2-3]
      \arrow[from=2-2, to=2-3]
      \arrow[from=2-2, to=3-2]
      \arrow[from=3-1, to=3-2]
    \end{tikzcd}\]
  Note that by \cref{lem:weaktaut} and \cref{lem:create}, the functor
  $\carC \colon \CS \to \bS$ creates pullbacks of monomorphisms, and
  so in particular it creates the following pullbacks:
  \[
    \begin{tikzcd}[row sep=15, column sep=15]
      Z & \cdot & B \\
      \cdot & {\cC(P)} & {\cC(B)} \\
      A & {\cC(A)}
      \arrow[hook, from=1-1, to=1-2]
      \arrow[hook, from=1-1, to=2-1]
      \arrow["\lrcorner"{anchor=center, pos=0.125}, draw=none, from=1-1, to=2-2]
      \arrow[from=1-2, to=1-3]
      \arrow[hook, from=1-2, to=2-2]
      \arrow["\lrcorner"{anchor=center, pos=0.125}, draw=none, from=1-2, to=2-3]
      \arrow[hook, from=1-3, to=2-3]
      \arrow[hook, from=2-1, to=2-2]
      \arrow[from=2-1, to=3-1]
      \arrow["\lrcorner"{anchor=center, pos=0.125}, draw=none, from=2-1, to=3-2]
      \arrow[from=2-2, to=2-3]
      \arrow[from=2-2, to=3-2]
      \arrow[hook, from=3-1, to=3-2]
    \end{tikzcd}
  \]
  The above $Z$ is equivalently the limit of the zig-zag
  $A \to \cC(A) \ot \cC(P) \to \cC(B) \ot B$.\footnote{This $Z$
    is in fact the product of $A$ and $B$ in $\CS$, as noted in
    \cite[Proposition
    3.4]{johnstone-power-tsujishita-watanabe-worrell}.} Thus we get a
  map $P \to Z$ in $\bS$. Note also that the following commutes:
  \[
    \begin{tikzcd}[row sep=15, column sep=15]
      Z & B && \\
      A & {\cC(P)} & {\cC(B)} \\
      & {\cC(A)} & P & B \\
      && A & X
      \arrow[from=1-1, to=1-2]
      \arrow[from=1-1, to=2-1]
      \arrow[from=1-1, to=2-2]
      \arrow[from=1-2, to=2-3]
      \arrow[curve={height=-25pt}, equals, from=1-2, to=3-4]
      \arrow[from=2-1, to=3-2]
      \arrow[curve={height=25pt}, equals, from=2-1, to=4-3]
      \arrow[from=2-2, to=2-3]
      \arrow[from=2-2, to=3-2]
      \arrow["\varepsilon", from=2-2, to=3-3]
      \arrow["\varepsilon"', from=2-3, to=3-4]
      \arrow["\varepsilon", from=3-2, to=4-3]
      \arrow[from=3-3, to=3-4]
      \arrow[from=3-3, to=4-3]
      \arrow[from=3-4, to=4-4]
      \arrow[from=4-3, to=4-4]
    \end{tikzcd}
  \]
  Therefore given any weak pullback $P'$ of $A \to X \ot B$ in $\CS$,
  we obtain a map $Z \to P'$ in $\CS$, and hence a composite map
  $P \to Z \to P'$ in $\bS$. Since $P$ is a pullback, we also have a
  comparison map $P' \to P$ in $\bS$, and moreover the roundtrip
  $P \to P' \to P$ is identity. Thus the comparison map $P' \to P$ is
  a split epimorphism, so $P'$ is a weak pullback in $\bS$ as desired.
\end{proof}

\begin{corollary}\label{cor:preservemono}
  Let $\cC$ be a comonad on a category with pullbacks $\bS$. If $\cC$
  preserves weak pullbacks, then the comonadic functor
  $\carC \colon \CS \to \bS$ preserves monomorphisms.
\end{corollary}
\begin{proof}
  Follows from \cref{lem:weakpreserve} and \cref{lem:weaktaut}.
\end{proof}

The assumption that a comonad $\cC$ preserves weak pullbacks ensures
especially pleasant properties of the category of
$\cC$-coalgebras. (The definition of \emph{regular category} is dual
to the definition of coregular category \cref{def:coregular}.)

\begin{proposition}[{\cite{johnstone-power-tsujishita-watanabe-worrell}}]\label{prop:regular}
  If $\cC$ is a weak-pullback-preserving comonad on a regular category
  $\bS$, then $\CS$ is regular.
\end{proposition}
\begin{proof}
  Follows from~\cite[Lemma 3.9]{johnstone-power-tsujishita-watanabe-worrell} and
  \cref{cor:preservemono}.
\end{proof}

The following provides a natural source of categories with a subobject
classifier that are not necessarily toposes (for instance, the
category $\MultiGph$ of undirected multigraphs from the introduction,
which we will encounter again later in \cref{ex:undgraphs} and
\cref{ex:symcat}).

\begin{proposition}[{\cite[Proposition 3.6]{johnstone-power-tsujishita-watanabe-worrell}}]\label{prop:classifier}
  Let $\cC$ be a comonad on a category $\bS$ with finite limits and a
  subobject classifier. If $\cC$ is taut and
  $\carC \colon \CS \to \bS$ preserves monomorphisms, then $\CS$ has a
  subobject classifier.

  In particular, if $\cC$ preserves weak pullbacks, then $\CS$ has a
  subobject classifier.\qed
\end{proposition}

\begin{remark}
  There is a more general analogue of this result for
  regular-subobject classifiers, replacing monomorphisms by regular
  monomorphisms. In particular, if $\cC$ is a taut comonad on $\Set$,
  then $\CSet$ has a regular-subobject classifier.
\end{remark}

We conclude with some general facts about categories of coalgebras for
comonads on $\Set$. Comonadic functors create
colimits (\cref{lem:create}), so all categories $\CSet$ are cocomplete.  Although
calculation of limits is less straightforward, categories of coalgebras
for comonads on $\Set$ are also complete.

\begin{lemma}[{\cite{linton}}]\label{lem:linton}
  If $\cC$ is a comonad on a category $\bS$ with (finite) products such
  that its category of coalgebras $\CS$ has coreflexive
  equalizers, then $\CS$ is (finitely) complete.\qed
\end{lemma}

In particular if $\cC$ is a crude comonad on a (finitely) complete
category $\bS$, then $\CS$ is (finitely) complete.

\begin{corollary}\label{lem:complete}
  If $\cC$ is a comonad on $\Set$, then its category of coalgebras
  $\CSet$ is complete.
\end{corollary}
\begin{proof}
  Follows from \cref{lem:corefl} and \cref{lem:linton}.
\end{proof}

The following lemma will also be useful.

\begin{lemma}
  In the situation of \cref{lem:factor}, the comonadic functor
  $\carC \colon \CS \to \bS$ reflects limits for which at least one of
  the projections is regular monic.
\end{lemma}
\begin{proof}
  By faithfulness, it suffices to show that the comparison maps in
  $\bS$ induced by cones in $\CS$ are coalgebra maps. This follows from
  \cref{lem:factor}.
\end{proof}

\begin{corollary}\label{lem:reflect}
  If $\cC$ is a comonad on $\Set$, then $\carC$ reflects equalizers,
  as well as pullbacks of monomorphisms along arbitrary maps.\qed
\end{corollary}

A functor into $\Set$ is called \emph{small} if it is a small colimit
of representable functors. A comonad on $\Set$ (or more generally any
functor into $\Set$ from an accessible category) is small if and only
if it is accessible~\cite[Proposition 2.4.2(i)]{makkai-pare}. It
follows from the theory of accessible
categories~\cite{makkai-pare,adamek-rosicky} that the category of
coalgebras associated to any accessible comonad on $\Set$ is locally
presentable (a.k.a.\ a category of models of a theory of limits):

\begin{lemma}[{\cite[Proposition A.2 (4)]{henry}}]\label{lem:accessibility}\label{lem:presentability}
  The category of (co)algebras of an accessible (co)monad on an
  accessible category is accessible.

  Thus the category of (co)algebras of an accessible (co)monad on a
  locally presentable category is locally presentable.\qed
\end{lemma}

In particular, this implies that the category of coalgebras $\CSet$ of
any small comonad $\cC$ on $\Set$ is total (meaning the Yoneda
embedding $\CSet \hookrightarrow \Ps{\CSet}$ has a left adjoint) and has a small
dense subcategory.

Many of the pleasant properties of categories of coalgebras of
comonads on $\Set$ may also be derived as consequences of the
following powerful result:

\begin{proposition}[{\cite{linton-pare}}]\label{lem:opmonadic}
  If $\cC$ is a comonad on $\Set$, then $(\CSet)\op$ is monadic over
  $\Set$.
\end{proposition}
\begin{proof}
  We have monadic functors
  \[\begin{tikzcd}
      {(\CSet)\op} & {\Set\op} & \Set
      \arrow["(\carC)\op", from=1-1, to=1-2]
      \arrow["{\Pow{(\dash)}}", from=1-2, to=1-3]
    \end{tikzcd}\] and their composite is monadic since the composite
  of a reflexive-coequalizer-preserving monadic functor with a monadic
  functor is monadic (\cref{lem:crudecompose}).
\end{proof}

In particular, $\CSet$ is cototal, coexact, and has a
regular-injective coseparator~\cite{vitale} given by the cofree
coalgebra $\cC(2)$.

It is called ``the fundamental theorem of topos theory'' that every
slice category of a topos is a topos. Similarly, every slice category
of a category of coalgebras of a comonad on $\Set$ is a category of
coalgebras of a comonad on $\Set$:

\begin{proposition}\label{prop:fundamental}
  Let $\cC$ be a comonad on $\Set$. For any coalgebra $A$, the
  composite forgetful functor
  \[\CSet / \eE \xto{\slicel} \CSet \xto{\carC} \Set\]
  where $\slicel$ denotes the slice category projection, is
  strictly\footnote{A left adjoint functor $\bB \to \bS$ inducing the
    comonad $\cC$ is called comonadic (resp. \emph{strictly}
    comonadic) when the comparison functor $K \colon \bB \to \CS$ from
    \cref{lem:comparison} is an equivalence (resp.\
    \emph{isomorphism}).} comonadic.
\end{proposition}
\begin{proof}
  The functor $\slicel \colon \CSet / \eE \to \CSet$ is crudely and
  strictly comonadic, with right adjoint defined via cartesian
  products in $\CSet$ which exist by \cref{lem:complete}. Thus the
  composite $\CSet \to \Set$ is also strictly comonadic by
  \cref{lem:crudecompose}.
\end{proof}

\begin{definition}\label{def:comonadofelements}\label{def:etale}
  We refer to the comonad on $\Set$ induced by $\CSet / \eE \to \CSet$ as in
  \cref{prop:fundamental} as the \emph{comonad of elements} of $\eE$,
  denoted $\El(\eE)$. This defines a functor
  $\El \colon \CSet \to \ComSet/\cC$.

  By \cref{prop:respect} the slice projection $\CSet / \eE \to \CSet$
  corresponds to a comonad map $\El(\eE) \to \CSet$. Note that functor
  of coalgebra categories $\DSet \to \CSet$ is a slice category
  projection (up to composing an isomorphism on the domain) if and
  only if it is a discrete fibration (\cref{def:df}), since $\DSet$
  always has a terminal object.

  We refer to such comonad maps $\phi \colon \bD \to \bC$ for which
  $\com{\Set}{\phi} \colon \DSet \to \CSet$ is a discrete fibration as
  \emph{\'etale}.
\end{definition}

\begin{example}\label{ex:dopf}
  Recall from \cref{ex:sheaves} that coalgebras of the comonad $\cX$
  corresponding to the topological space $\tX$ are equivalently
  \'etal\'e spaces (or sheaves) over $\tX$. It is a standard fact that
  the slice category over any particular \'etal\'e space $\eE$ is the
  category of \'etal\'e spaces over $\eE$, so the comonad $\El(\eE)$
  is $\toptocmd{\eE}$, the comonad corresponding to the topological
  space $\eE$. Thus the \'etale comonad maps between comonads
  corresponding to topological spaces are the same as \'etale maps
  between topological spaces.

  Recall also from \cref{ex:copresheaves} that the category of
  coalgebras of the comonad $\cbC$ corresponding to the
  category $\bC$ is equivalent to $\fun{\bC}{\Set}$. It is a standard
  fact that given a functor $\eE \colon \bC \to \Set$, the slice
  category $\fun{\bC}{\Set} / \eE$ is equivalent to
  $\fun{\El(\eE)}{\Set}$ where $\El(\eE)$ denotes the category of
  elements. This justifies the notation $\El(\eE)$ for the induced
  comonad. Recall also that a functor $\eE \colon \bC \to \Set$ may
  equivalently be described as a discrete opfibration
  $\El(\eE) \to \bC$. Thus the \'etale comonad maps between comonads
  corresponding to categories are the same as discrete opfibrations
  between categories.
\end{example}

We will give some more examples of comonads on $\Set$ after we have
introduced some tools to generate them as one would spaces.

\section{Density comonads}\label[bgchapter]{sec:density}

As discussed in the introduction, comonads on $\Set$ generalize
topological spaces, and functors $\pP\colon\bB\to\Set$ take the role
of generalized subbases. The density comonad $\cP$ of $\pP$ then takes
the role of the space generated by the subbasis. This section collects
definitions and results around density comonads.

\begin{definition}
  The \emph{left Kan extension}\footnote{The definitions of left Kan
    extension and density comonad also make sense in an abstract
    2-category.} of $\pP \colon \bB \to \bS$ along $F \colon \bB \to \bC$
  is, if it exists, the universal functor
  $\La{F}{\pP} \colon \bC\to \bS$ equipped with a natural
  transformation
  \[
    \begin{tikzcd}[column sep=15pt, row sep=15pt]
      {\bB} &{}& {\bS} \\
      & {\bC}
      \arrow["\pP", from=1-1, to=1-3]
      \arrow["F"', from=1-1, to=2-2]
      \arrow["\La{F}{\pP}"', dashed, from=2-2, to=1-3]
      \arrow["\Downarrow"{description, pos=0.3}, draw=none, from=1-2, to=2-2]
    \end{tikzcd}
  \]
  That is, $\La{F}{\pP}$ is defined by the natural isomorphism
  $\Ho{\fun{\bC}{\bS}}(\La{F}{\pP}, \dash) \cong
  \Ho{\fun{\bB}{\bS}}(\pP, \dash \circ F)$.
\end{definition}

\begin{definition}
  The \emph{density comonad} of a functor $\pP\colon \bB \to \bS$ is, if it exists, the
  left Kan extension of $\pP$ along itself
  $\cP \coloneqq \La{\pP}{\pP}$.
  \[
    \begin{tikzcd}[column sep=15pt, row sep=15pt]
      {\bB} &{}& {\bS} \\
      & {\bS}
      \arrow["\pP", from=1-1, to=1-3]
      \arrow["\pP"', from=1-1, to=2-2]
      \arrow["\cP"', dashed, from=2-2, to=1-3]
      \arrow["\Downarrow"{description, pos=0.3}, draw=none, from=1-2, to=2-2]
    \end{tikzcd}
  \]
\end{definition}

Density comonads are indeed comonads: one obtains the counit and
comultiplication from the universal property of the Kan extension.
If $\bB$ is small and $\bS$ is a cocomplete category such as $\Set$,
then $\cP$ necessarily exists, calculated as follows.

\begin{proposition}[{\cite[Proposition A.7]{diliberti}}]
  Given a diagram
  \[\begin{tikzcd}[column sep=small, row sep=small]
      \bB && \bS \\
      & {\bC}
      \arrow["\pP", from=1-1, to=1-3]
      \arrow["F"', from=1-1, to=2-2]
    \end{tikzcd}\] with $\bB$ small and $\bS$ cocomplete (and locally
  small), the left Kan extension $\La{F}{\pP} \colon \bC \to \bS$ of
  $\pP$ along $F$ is given by the composite
  \[\begin{tikzcd}[column sep=small, row sep=small]
      {\Ps{\bB}} && \bS \\
      & {\bC}
      \arrow["{\extend{\pP}}", from=1-1, to=1-3]
      \arrow["{\nerve{F}}", from=2-2, to=1-1]
    \end{tikzcd}
    \qquad
  \]
  where $\nerve{F}$ is the \emph{nerve functor}
  $\xA \mapsto \Ho{\bC}(F(\dash), \xA)$ and $\extend{\pP}$ is the
  colimit-preserving extension of $\pP$ to the free cocompletion $\Ps{\bB}$ of its domain
  $\bB$.\qed
\end{proposition}

The functor $\extend{\pP} \colon \Ps{\bB} \to \bS$ sends each
$\bB$-presheaf $W$ to the colimit of the diagram
\[\El(W) \to \bB \xto{\pP} \bS\] (i.e.\ the \emph{$W$-weighted
  colimit} of $\pP$, which may be calculated as the coend
$\int^{\xU \in \bB} W(\xU) \cdot \pP(\xU)$). Therefore the left Kan
extension $\La{F}{\pP} = \extend{\pP} \circ \nerve{F}$ sends an object
$\xA$ in $\bS$ to the colimit of the diagram
\[\El(\xA^F) \cong (F / \xA) \to \bB \xto{\pP} \bS.\]

Note also that $\extend{\pP} \colon \Ps{\bB} \to \bS$ is the left
adjoint of the nerve functor $\nerve{\pP} \colon \bS \to
\Ps{\bB}$. Thus in particular we have:

\begin{corollary}\label{prop:densityadj}
  Let $\pP \colon \bB \to \bS$ with $\bB$ small and $\bS$ cocomplete
  (and locally small). Then the density comonad $\cP$ is the comonad
  associated to the adjunction:
  \[
    \begin{tikzcd}[column sep=30pt]
      \cat{\Ps{\bB}}\ar[r, shift left=5pt, "\extend{\pP}"]\ar["\bot"{anchor=center},r, draw=none]&
      \cat{\bS}\ar[l, shift left=5pt, "\nerve{\pP}"]
    \end{tikzcd}
    \tag*{\qed}
  \]
\end{corollary}

In the more general situation that $\bS$ is not necessarily cocomplete
or $\bB$ is not necessarily small, but nevertheless the specific
relevant colimits exist, then the left Kan extension exists, given by
the same formula. Such left Kan extensions --- those given by the
expected colimit formula --- are called \emph{pointwise}.  We will see
shortly that all Kan extensions (even those involving large
categories) into $\Set$ are pointwise.

\begin{definition}\label{def:pointwise}
  The \emph{pointwise left Kan extension} of $\pP \colon \bB \to \bS$
  along $F \colon \bB \to \bC$, if it exists, is the functor
  $\bC \to \bS$ given by the colimit formula
  \[\xA \mapsto \colim_{F(\xU) \to \xA} \pP(\xU)\] defined just in case all
  such colimits exist.

  Pointwise left Kan extensions are in particular left Kan extensions;
  the defining natural transformation
  $\pP \Rightarrow \La{F}{\pP} \circ F$ is given by the colimit
  inclusion maps \[\pP(\xU) \to \colim_{F(\xV) \to F(\xU)} \pP(\xV)\]
  indexed by identities $\id_{F(\xU)}$. A left Kan extension is
  pointwise if and only if it is preserved by all representable
  functors $\Ho{\bS}(\dash, \xA)$. See~\cite[Theorem X.5.3]{maclane}
  for details.

  We likewise call the density comonad $\cP$ of $P$ \emph{pointwise}
  if the left Kan extension $\La{\pP}{\pP}$ is pointwise.
  In this case the counit and comultiplication
  \[\varepsilon_\xA \colon \colim_{\pP(\xU) \to \xA} \pP(\xU) \to \xA
    \qqand \delta_\xA \colon \colim_{\pP(\xU) \to \xA} \pP(\xU) \to
    \colim_{\pP(\xV) \to \left(\colim_{\pP(\xW) \to \xA}
        \pP(\xW)\right)} \pP(\xV)\] are respectively induced by the
  canonical cocone into $\xA$ and the map on colimit indices sending
  $f \colon \pP(\xU) \to \xA$ to the colimit inclusion
  $\pP(\xU) \to \left(\colim_{\pP(\xU) \to \xA} \pP(\xU)\right)$
  indexed by $f$.
\end{definition}


\begin{proposition}\label{lem:setpointwise}
  Let $\bS$ be a (locally small) category with small powers, e.g.\
  $\Set$. Any left Kan extension
  \[
    \begin{tikzcd}[column sep=15pt, row sep=15pt]
      {\bB} &{}& {\bS} \\
      & {\bC}
      \arrow["\pP", from=1-1, to=1-3]
      \arrow["F"', from=1-1, to=2-2]
      \arrow["\La{F}{\pP}"', dashed, from=2-2, to=1-3]
      \arrow["\Downarrow"{description, pos=0.3}, draw=none, from=1-2, to=2-2]
    \end{tikzcd}
  \]
  that exists is pointwise.
\end{proposition}
\begin{proof}
  A left Kan extension is pointwise if and only if it is preserved by
  all representable functors $\bS(\dash, \xA)$. It is well-known that
  left adjoints preserve left Kan extensions (see also the more
  general \cref{lem:kancomposeleft} below), and the functor
  $\bS(\dash, \xA) \colon \bS \to \Set\op$ is a left adjoint if and
  only if all powers of $\xA$ exist.
\end{proof}

In the case $\bS$ is $\Set$, we have additional ways of expressing
$\La{F}{\pP}$:

\begin{proposition}\label{prop:setformula}
  If the 
  left Kan extension of $\pP \colon \bB \to \Set$ along
  $F \colon \bB \to \bC$ exists, then its underlying functor is given
  by
  \begin{enumerate}[label=(\roman*)]
  \item the coend $\int^{\xU\in\bB}\pP(\xU) \times (\dash)^{F(\xU)}$, or
    equivalently\label{item:pcoend}
  \item the $\pP$-weighted colimit of
    $\bB\op\xto{F\op}\bC\op\xto{\yo}\fun{\bC}{\Set}$, or
    equivalently\label{item:pweighted}
  \item the colimit of the diagram
    $(\El(\pP))\op\to\bB\op\xto{F\op}\bC\op\xto{\yo}\fun{\bC}{\Set}$\label{item:pcolimit}
  \end{enumerate}
  where $\yo \colon \bC\op \hookrightarrow \fun{\bC}{\Set}$ denotes the Yoneda
  embedding.
\end{proposition}
\begin{proof}
  Given a covariant functor $\pP \colon \bB \to \Set$ and a
  contravariant functor $W \colon \bB\op \to \Set$, the $\pP$-weighted
  colimit of $W$ coincides with the $W$-weighted colimit of
  $\pP$. This is the tensor product of functors $\pP \otimes W$
  (a.k.a.\ profunctor composition), given by the coend
  $\int^{\xU \in \bB} W(\xU) \times \pP(\xU)$, where $W$ and $\pP$ evidently
  play symmetrical roles.

  By definition, the pointwise left Kan extension $\La{F}{\pP}(\xA)$ is
  the $\xA^F$-weighted colimit of $\pP$, whereas the above formulas are
  various ways of expressing the $\pP$-weighted colimit of $\xA^F$.
\end{proof}

\begin{remark}\label{rem:densitygerms}
  In analogy with \cref{ex:top}, for a general comonad $\cC$ on
  $\Set$, we still refer to the elements of $\cC(\xA)$ as $\xA$-valued
  \emph{germs}, and we still employ the notation $[f]_x$ for a
  germ. Indeed, the above coend formula for the density comonad $\cP$
  of $\pP \colon \bB \to \Set$ tells us that a germ
  $[f]_x \in \cP(\xA)$ comprises an element $x \in \pP(U)$ of an
  ``open set'' $\xU$ and a function $f \colon \pP(\xU) \to \xA$,
  modulo the equivalence relation generated by ``restricting'', i.e.\
  $[f]_x = [g]_y$ whenever there exists an arrow
  $i \colon \xU \to \xV$ in $\bB$ such that $\pP(i)(x) = y$ and
  $f = g\circ \pP(i)$. Thus $\cP(\xA)$ is reasonably interpreted as
  the set of ``germs'' of $\xA$-valued functions on some abstract
  ``space'' $\cP.$

  Moreover, just as in \cref{ex:top}, the counit
  $\varepsilon_\xA \colon \cP(\xA) \to \xA$ is always given by the
  formula $[f]_x \mapsto f(x)$, and the comultiplication
  $\delta_\xA \colon \cP(\xA) \to \cP(\cP(\xA))$ is always given by
  the formula $[f]_x \mapsto [[f]_{(\dash)}]_x$, where $[f]_{(\dash)}$
  denotes the $\cP(\xA)$-valued function on the same domain as $f$
  defined by $y \mapsto [f]_y$. This is simply a translation of the
  colimit formulas from \cref{def:pointwise} into the germ
  notation.

  Given $g \colon \xA \to \xB$ in $\Set$, the map
  $\cP(g) \colon \cP(\xA) \to \cP(\xB)$ is given by
  $[f]_x \mapsto [g\circ f]_x$. A coalgebra
  $\hH \colon \eE \to \cP(\eE)$ then consists of a germ
  $\hH(e) = [f_e]_{x_e}$ for each $e$ in $\eE$ satisfying
  $f_e(x_e) = e$ (the unit law) and
  $[\hH \circ f_e]_{x_e} = [[f_e]_{(\dash)}]_{x_e}$ (the associativity
  law). Thus each element $e$ in $\eE$ is assigned an ``infinitesimal
  neighborhood'' of elements in $\eE$, where the unit law ensures $e$
  is at the center of its own neighborhood, and the associativity law
  ensures the map from points to neighborhoods $\hH$ has a germ at
  each point equal to its neighborhood viewed as a germ of germs. A
  coalgebra map $\hK \colon \eE \to \eE'$ is such that
  $[\hK \circ f_e]_{x_e} = [f_{\hK(e)}]_{x_{\hK(e)}}$, i.e.\
  it respects assigned infinitesimal neighborhoods.
\end{remark}

Let us apply these tools to obtain some examples of comonads on $\Set$
besides topological spaces or categories.

\begin{example}[Undirected multigraphs]\label{ex:undgraphs}
  We construct the comonad from the introduction whose coalgebras are
  undirected multigraphs. Consider the inclusion $\pP$ of the
  subcategory of $\Set$ generated by the following arrows, where
  $\sigma$ is the involution interchanging $l$ and $r$:
  \[\begin{tikzcd}
      1 & {\setof{l,c,r}}
      \arrow["r"', shift right, from=1-1, to=1-2]
      \arrow["l", shift left, from=1-1, to=1-2]
      \arrow["\sigma", from=1-2, to=1-2, loop, in=55, out=125, distance=10mm]
    \end{tikzcd}\]
  The formula for the density comonad 
  \cref{prop:setformula} used as in \cref{rem:densitygerms} tells us that a germ
  $[f]_x \in \cP(\xA)$ is either 
  \begin{itemize}
  \item the germ of a function $f \colon \{l,c,r\} \to \xA$ (an ordered triple)
    about one of the elements of the domain, which we denote by
    $[a_l, a_c, a_r]_l$, $[a_l, a_c, a_r]_c$, or $[a_l, a_c, a_r]_r$, or
  \item the germ of a function $f \colon 1 \to \xA$ (an element) about
    the unique element of $1$, which we denote by $[a]_\ast$
  \end{itemize}
  subject to the equivalence relation generated by:
  \begin{align*}
    [a_l]_\ast&=[a_l,a_c,a_r]_l&&&[a_l,a_c,a_r]_l&=[a_r,a_c,a_l]_r\\
    [a_r]_\ast&=[a_l,a_c,a_r]_r&&&[a_l,a_c,a_r]_c&=[a_r,a_c,a_l]_c
  \end{align*}
  Note that no singleton germs will be identified with each other, so
  every germ reduces to either a germ of the form $[a]_\ast$ or a germ
  of the form $[a_r,a_c,a_l]_c = [a_l,a_c,a_r]_c$, picking out an
  element $a_c$ and an \emph{unordered pair} $\setof{a_l, a_r}$. In
  effect one is left with $\xA$ many germs about $\ast$ together with
  $\xA\left(\xA + {\xA\choose 2}\right)=\xA{\xA+1\choose 2}$ many
  germs about $c$. Thus the formula for $\cP(\xA)$ is given
  cardinality-wise by $\xA{\xA+1\choose 2}+\xA$.

  A coalgebra $\hH\colon \eE\to \cP(\eE)$ then assigns to each
  $e\in \eE$ either the germ $[e]_\ast$ or
  $[e_l,e,e_r]_c = [e_r,e,e_l]_c$ (using the unit law to pick out
  values that must be $e$) such that in the second case,
  $h(e_l)=[e_l]_\ast$ and $h(e_r)=[e_r]_\ast$, using the associativity
  law. Thus $\eE$ decomposes as $\eE_\ast+\eE_c$, and we have a
  function from $\eE_c$ to unordered pairs in $\eE_\ast$ given by
  $e\mapsto \setof{e_l,e_r}$. This is precisely an undirected multigraph,
  as desired.

  The comonad $\cP$ preserves weak pullbacks, and hence by
  \cref{prop:regular} and \cref{prop:classifier}, $\MultiGph$ is a
  regular category with a subobject classifier. (Later in
  \cref{ex:symcat} will be able to deduce that this comonad preserves
  weak pullbacks from general theory, but this is also straightforward
  to check.)  Explicitly, the subobject classifier in $\MultiGph$ is
  given by forgetting the direction of the edges of the subobject
  classifier in $\Set^{\rightrightarrows}$:
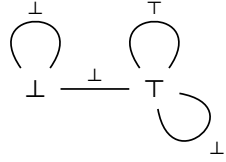
  \[\begin{tikzcd}
      \bot & \top
      \arrow["\bot", no head, from=1-1, to=1-1, loop, in=55, out=125, distance=10mm]
      \arrow["\bot", no head, from=1-1, to=1-2]
      \arrow["\top", no head, from=1-2, to=1-2, loop, in=55, out=125, distance=10mm]
      \arrow["\bot"', no head, from=1-2, to=1-2, loop, in=350, out=280, distance=10mm]
    \end{tikzcd}\]
  Despite this, as mentioned in the introduction, $\MultiGph$ is not a
  topos (see \cref{cex:undgraphs}).
\end{example}

Whereas a subbasis of a topological space is a diagram of subset
inclusions, the following is the complete opposite situation.

\begin{example}[Partitions]\label{ex:part}
  Let $\Surj$ be the category of sets and surjections.  Although
  $\Surj$ is large, the density comonad $\cP$ of the inclusion
  $\pP \colon \Surj \hookrightarrow \Set$ exists.

  The formula for the density comonad \cref{prop:setformula} used as
  in \cref{rem:densitygerms} tells us that a germ $[f]_x \in \cP(\xA)$
  consists of an arbitrary set $\xX$ with a distinguished element
  $x \in \xX$ and a function $f \colon \xX \to \xA$, modulo the
  equivalence relation generated by $[f]_x = [f\circ s]_{x'}$ whenever
  $s \colon \xX' \to \xX$ is a surjection sending $x'$ to $x$. Observe
  that the image $\im f = \xS \subseteq \xA$ of $f$ and the point
  $f(x) = a \in \xA$ assigned $x$ are invariant under this equivalence
  relation. Each equivalence class then has a canonical representative
  $[i]_a$ where $i \colon \xS \hookrightarrow \xA$ is a subset of
  $\xA$ and $a \in \xS$, since any germ $[f]_x$ with the same image
  $\im f = i$ and assigned point $f(x) = a$ is seen to satisfy
  $[f]_x = [i]_a$ by factoring the function $f$ as a surjection
  followed by an injection $f = i \circ s$.

  Thus the elements of $\cP(\xA)$ are precisely \emph{pointed subsets}
  \[\cP(\xA) = \setof{(\xS, a) \mid \xS \subseteq \xA, a \in \xS} = \sum_{\xS \subseteq \xA} \xS.\]

  A coalgebra $\hH \colon \eE \to \cP(\eE)$ then assigns each element
  $e$ a pointed subset $\xS$ of $\eE$ where the point is $e$ (the unit
  law) and such that all elements of $\xS$ are also assigned $\xS$ (the
  associativity law).
  In other words, a coalgebra is a \emph{partitioned set}, and a
  coalgebra map is seen to be a function such that the image of each
  part is a part. For future convenience we denote this category by
  $\Part \coloneqq \PSet$.  This comonad $\cP$ also evidently
  preserves weak pullbacks, so like $\MultiGph$ from
  \cref{ex:undgraphs} above, $\Part$ is regular with a subobject
  classifier.
  %
  %
  %
  %
  %
\end{example}

The following is one of the simplest diagrams in $\Set$ with a
nontrivial density comonad, yet the result is wholly unfamiliar.

\begin{example}\label{ex:shear}
  Let $\arro$ be the walking arrow, and let
  $\pP \colon \arro \to \Set$ be the diagram singling out the unique
  function of sets $2 \to 1$. The formula for the density comonad
  \cref{prop:setformula} used as in \cref{rem:densitygerms} tells us
  that a germ $[f]_x \in \cP(\xA)$ is either
  \begin{itemize}
  \item the germ of a function $f \colon 2 \to \xA$ (an ordered pair)
    about the first or second element of $2$, which we denote by
    $[a_1, a_2]_1$ or $[a_1, a_2]_2$, or
  \item the germ of a function $f \colon 1 \to \xA$ (an element) about
    the unique element of $1$, which we denote by $[a]_\ast$
  \end{itemize}
  modulo the equivalence relation $[a, a]_1 = [a]_\ast = [a,
  a]_2$. Each equivalence class then has a canonical representative:
  either $[a_1, a_2]_i$ with $a_1 \neq a_2$ or $[a]_\ast$.

  Thus the formula for $\cP(\xA)$ is given cardinality-wise by
  \[\pts{\cP(\xA)} = \pts{2 \xA^2 - \xA}\]
  since we have identified the two instances of an ordered pair
  whenever they contain the same repeated element. More specifically,
  the endofunctor $\cP$ is a pushout of endofunctors
  $(\dash)^2 \xot{\bang} (\dash) \xto{\bang} (\dash)^2$.

  Intuitively, we think of $\cP$ as a ``space'' whose diagram of open
  sets is $2 \to 1$. Since the ordering is by ``inclusion'', here the
  two-element set takes the role of a ``smaller'' open set $U$, and
  the one-element set takes the role of a ``larger'' open set $V$. We
  like to imagine the point dividing in two as we zoom in, revealing
  more detail on closer inspection.
  
  
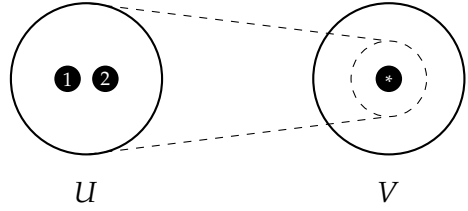
\[
    \begin{tikzpicture}
      \draw[thick] (-2,0) circle (1);
      \draw[thick] (2,0) circle (1);
      \fill (-2.25,0) circle (.175);
      \fill (-1.75,0) circle (.175);
      \fill (2,0) circle (.175);
      \draw[dashed,dash phase=-1.15pt] (2,0) circle (.5);
      \draw[dashed,shorten <=2.9pt] (2,.5) -- (-2,1);
      \draw[dashed,shorten <=2.9pt] (2,-.5) -- (-2,-1);
      \node at (-2,-1.5) {$U$};
      \node at (2,-1.5) {$V$};
      \node at (-2,1.5) {\phantom{$U$}};
      \node at (2,1.5) {\phantom{$V$}};
      \node[white] at (-2.25,0) {\tiny$1$};
      \node[white] at (-1.75,0) {\tiny$2$};
      \node[white] at (2,-.015) {\tiny$\ast$};
    \end{tikzpicture}
  \]

  A coalgebra $\hH \colon \eE \to \cP(\eE)$ assigns each element $e_0$ an
  $\eE$-valued germ of the form $[e_0, e_2]_1$, $[e_1, e_0]_2$, or
  $[e]_\ast$ (as guaranteed by the unit law) and such that in the
  first case $e_2$ is always assigned $[e_0, e_2]_2$ and in the second
  case $e_1$ is always assigned $[e_1, e_0]_1$ (as guaranteed by the
  associativity law). In other words, a coalgebra is a set partitioned
  into ordered pairs and singletons. Writing $\eE_1$ for the set of
  singletons and $\eE_2$ for the set of ordered pairs, a coalgebra map
  $\eE \to \eE'$ then amounts to a function
  $\eE_1 + \eE_2 \to \eE'_1 + \eE'_2$ sending $\eE_1$ entirely within
  $\eE'_1$. At this point we see one can entirely replace the coalgebra category 
  of $\cP$ with an isomorphic category whose objects 
  are arbitrary pairs $(E_1,E_2)$ of sets.

  For future convenience, we denote this category, whose
  objects are pairs of sets and whose arrows are functions between
  their coproducts satisfying this shear condition
  \[
    \begin{tikzcd}
      {\eE_1} & {\eE'_1} \\
      {\eE_2} & {\eE'_2}
      \arrow[from=1-1, to=1-2]
      \arrow["{+}"{description}, draw=none, from=1-1, to=2-1]
      \arrow["{+}"{description}, draw=none, from=1-2, to=2-2]
      \arrow[from=2-1, to=1-2]
      \arrow[from=2-1, to=2-2]
    \end{tikzcd}
  \]
  by $\Shear \cong \PSet$. This comonad does not preserve weak
  pullbacks, and proves useful for counterexamples; for instance the
  comonadic forgetful functor into $\Set$ does not preserve
  monomorphisms (\cref{cex:shear}).
\end{example}

Now we recall the category-theoretic concept of density, after which
density comonads were named.

\begin{definition}\label{def:dense}
  A functor $\pP \colon \bB \to \bS$ is \emph{dense} if the nerve
  functor $\nerve{\pP} \colon \bS \to \Ps{\bB}$, defined by
  $\xA \mapsto \Ho{\bS}(\pP, \xA)$, is fully faithful.
\end{definition}

The name density comonad is justified by the following fact:

\begin{proposition}[{\cite[Proposition X.6.2]{maclane}}]\label{prop:densecolimit}
  A functor $\pP \colon \bB \to \bS$ is dense if and only if its
  density comonad exists, is pointwise, and is identity, i.e.\ every
  object $\xA$ in $\bS$ is the canonical colimit of shape $(\pP / \xA)$.\qed
\end{proposition}
%

The following lemma will help us to calculate density comonads in practice.

\begin{lemma}\label{lem:kancomposeleft}
  Suppose given a diagram
  \[\begin{tikzcd}[row sep=3pt]
      & {\bC_1} & {\bS_1} \\
      {\bB}\ar[ru, "F_1"]\ar[rd, "F_2"']\\
      & {\bC_2} & {\bS_2}
      \arrow[""{name=0, anchor=center, inner sep=0}, "{L_1}", shift left=2, from=1-2, to=1-3]
      \arrow[""{name=1, anchor=center, inner sep=0}, "{R_1}", shift left=2, from=1-3, to=1-2]
      \arrow[""{name=2, anchor=center, inner sep=0}, "{L_2}", shift left=2, from=3-2, to=3-3]
      \arrow[""{name=3, anchor=center, inner sep=0}, "{R_2}", shift left=2, from=3-3, to=3-2]
      \arrow["\bot"{anchor=center},r, draw=none, from=1-2, to=1-3]
      \arrow["\bot"{anchor=center},r, draw=none, from=3-2, to=3-3]
    \end{tikzcd}\] where $L_1$ is left adjoint to $R_1$ and $L_2$ is
  left adjoint to $R_2$. Then
  \[\La{L_2 \circ F_2}{(L_1 \circ F_1)} \cong L_1 \circ \La{F_2}{F_1}
    \circ R_2\] if either side exists.\footnote{This holds in an
    abstract 2-category.}
\end{lemma}
\begin{proof}
  We have natural bijections between 2-cells of the following forms:
  \[
    \begin{tikzcd}[row sep=small, column sep=small]
      {\bB} & {\bC_1} & {\bS_1} \\
      {\bC_2} \\
      {\bS_2}
      \arrow["{F_1}", from=1-1, to=1-2]
      \arrow["{F_2}"', from=1-1, to=2-1]
      \arrow["{L_1}", from=1-2, to=1-3]
      \arrow["{L_2}"', from=2-1, to=3-1]
      \arrow[""{name=0, anchor=center, inner sep=0}, dashed, from=3-1, to=1-3]
      \arrow["\Downarrow"{description, pos=0.55}, draw=none, from=1-1, to=0]
    \end{tikzcd}
    \;\;
    \begin{tikzcd}[row sep=small, column sep=small]
      {\bB} && {\bC_1} \\
      {\bC_2} & {\bS_1} \\
      {\bS_2}
      \arrow["{F_1}", from=1-1, to=1-3]
      \arrow["{F_2}"', from=1-1, to=2-1]
      \arrow["\Downarrow"{description, pos=0.75}, draw=none, from=1-1, to=2-2]
      \arrow["{L_2}"', from=2-1, to=3-1]
      \arrow["{R_1}"', from=2-2, to=1-3]
      \arrow[dashed, from=3-1, to=2-2]
    \end{tikzcd}
    \;\;
    \begin{tikzcd}[row sep=small, column sep=small]
      {\bC_2} && {\bC_1} \\
      & {\bS_1} \\
      {\bS_2}
      \arrow["{\La{F_2}{F_1}}", from=1-1, to=1-3]
      \arrow["\Downarrow"{description, pos=0.75}, draw=none, from=1-1, to=2-2]
      \arrow["{L_2}"', from=1-1, to=3-1]
      \arrow["{R_1}"', from=2-2, to=1-3]
      \arrow[dashed, from=3-1, to=2-2]
    \end{tikzcd}
    \;\;
    \begin{tikzcd}
      {\bC_2} & {\bC_1} \\
      {\bS_2} & {\bS_1}
      \arrow["{\La{F_2}{F_1}}", from=1-1, to=1-2]
      \arrow["\Downarrow"{description}, draw=none, from=1-1, to=2-2]
      \arrow["{L_1}", from=1-2, to=2-2]
      \arrow["{R_2}", from=2-1, to=1-1]
      \arrow[dashed, from=2-1, to=2-2]
    \end{tikzcd}
    \qedhere
  \]
\end{proof}

\begin{remark}\label{rem:composepointwise}
  In the situation of the above lemma, if $\La{F_2}{F_1}$ is
  pointwise, then $\La{L_2 \circ F_2}{(L_1 \circ F_1)}$ is pointwise.
  This follows from the fact that a left Kan extension is pointwise if and
  only if it is preserved by the functor of bicategories
  $\Cat \to \Prof\co$ from categories to profunctors sending
  $F\colon \bC \to \bD$ to $\Ho{\bD}(F(\dash), \dash)$, and that
  functors of bicategories preserve adjoints.
\end{remark}

\begin{corollary}\label{lem:composeleft}
  Suppose that $F\colon\bB\to\bC$ has a density comonad
  $\gen{F}$ and that $L\colon\bC\to\bS$ has a right adjoint
  $R$. Then we have $\gen{L\circ F}=L\circ\gen{F}\circ R$.\qed
\end{corollary}

\begin{corollary}\label{cor:easybases}
  If $F \colon \bB \to \bC$ is dense and $L \colon \bC \to \bS$ has a
  right adjoint $R$, then the induced comonad $L \circ R$ is the
  pointwise density comonad $\gen{L \circ F}$ of the composite
  $\bB \xto{F} \bC \xto{L} \bS$.
\end{corollary}
\begin{proof}
  Since $F$ is dense, its density comonad on $\bC$ is pointwise and
  identity $\gen{F} \cong \id_{\bC}$ by \cref{prop:densecolimit}. By
  \cref{lem:composeleft}, the density comonad $\gen{L \circ F}$ is
  $L \circ \gen{F} \circ R \cong L \circ R$. This is moreover
  pointwise by \cref{rem:composepointwise}.
\end{proof}

\begin{example}[Canonical basis of topological space]\label{ex:topbasis}
  Let $\opens{\tX}$ denote the poset of open subsets of a topological space
  $\tX$, ordered by inclusion. The canonical embedding
  $\opens{\tX} \hookrightarrow \Sh(\tX)$ sending an open set to the
  corresponding representable sheaf is dense. Indeed, by definition,
  this is simply to say that sheaves embed canonically and fully
  faithfully in presheaves over $\opens{\tX}$.

  Now recall from \cref{ex:sheaves} that the sum of stalks functor
  $\Sh(\tX) \to \Set$, inducing the comonad $\cX$ on
  $\Set$, is comonadic. The composite
  \[\opens{\tX} \hookrightarrow \Sh(\tX) \to \Set\]
  is simply the forgetful functor
  $\topbasis{\tX} \colon \opens{\tX} \hookrightarrow \Set$, which
  sends an open subset of $\tX$ to that set of points itself.  Thus by
  \cref{cor:easybases}, the germs comonad $\cX$ is the
  density comonad $\gen{\topbasis{\tX}}$ of $\topbasis{\tX}$, as
  expected. Indeed, this agrees with the explicit description of the
  density comonad in terms of germs from \cref{rem:densitygerms}.
\end{example}

\begin{example}[Canonical ``basis'' of category]\label{ex:catbasis}
  Let $\bC$ be a (locally small) category. The Yoneda embedding
  $\yo \colon \bC\op \hookrightarrow \fun{\bC}{\Set} \hookrightarrow
  \Set$ sending an object to the corresponding representable functor
  is dense. Indeed, this is the rather tautologous statement that
  $\fun{\bC}{\Set}$ canonically embeds in $\fun{\bC}{\Set}$.

  Now recall from \cref{ex:copresheaves} that (assuming small $\bC$)
  the elements functor
  $\pts{\El(\dash)} \colon \fun{\bC}{\Set} \to \Set$, inducing the
  comonad $\cbC$ on $\Set$, is comonadic. The composite
  \[\bC\op \xhookrightarrow{\yo} \fun{\bC}{\Set} \xto{\pts{\El(\dash)}} \Set\]
  is the functor $\catbasis{\bC} \colon \bC\op \to \Set$ given by the
  formula
  \[
    \catbasis{\bC}(c) = \pts{\El(\yo_c)} = \pts{c/\bC} = \sum_{d \in \bC} \bC(c, d).
  \]
  Now by \cref{cor:easybases}, the comonad $\cbC$ may be
  constructed as the density comonad $\gen{\catbasis{\bC}}$ of this
  functor $\catbasis{\bC}$.
\end{example}

The following tells us in particular that every comonad admits a ``big
basis''\footnote{The word ``basis'' will be defined in general in \cref{def:basis}}, where $\bB$ is the entire category of coalgebras and the
diagram is simply the comonadic functor itself.

\begin{corollary}\label{cor:densityofleft}
  If $L \colon \bB \to \bS$ has a right adjoint $R$, then its induced
  comonad $L \circ R$ is its pointwise density comonad $\gen{L}$.
\end{corollary}
\begin{proof}
  Follows from \cref{cor:easybases}, noting that the identity functor
  is always dense.
  %
\end{proof}

Density comonads also satisfy another universal property: the comonadic functor associated to the density comonad of
$\pP$ is the free comonadic functor extending $\pP$.

\begin{proposition}[{\cite[Section II.1]{dubuc}}]\label{prop:density}
  Let $\pP \colon \bB \to \bS$ have density comonad $\cP$.
  Then we have a bijection between commutative triangles of the
  following forms, natural in comonads $\cC$ on
  $\bS$:\footnote{This holds in an abstract 2-category (assuming the
    relevant Eilenberg-Moore objects exist).}
  \[\begin{tikzcd}[row sep=10pt,column sep=15pt]
      {\PS} && {\CS} \\
      & {\bS}
      \arrow[dashed, from=1-1, to=1-3]
      \arrow["\carP"', from=1-1, to=2-2]
      \arrow["\carC", from=1-3, to=2-2]
    \end{tikzcd}
    \qquad
    \leftrightarrows
    \qquad
    \begin{tikzcd}[row sep=10pt,column sep=15pt]
      {\bB} && {\CS} \\
      & {\bS}
      \arrow[dashed, from=1-1, to=1-3]
      \arrow["\pP"', from=1-1, to=2-2]
      \arrow["\carC", from=1-3, to=2-2]
    \end{tikzcd}
    \qedhere
  \]
\end{proposition}

Note that by \cref{prop:respect}, the former triangles are
equivalently comonad maps from $\cP$ to $\cC$. \cref{prop:density}
allows us to replace the functor $\carC$ in the commutative triangle
defining a comonad map $\cC \to \cD$ by any functor whose density
comonad is $\cC$.

\begin{example}[Indiscrete comonads]\label{ex:indiscrete}
  Let $\xX \colon \point \to \bS$ be an object of a category $\bS$. By
  \cref{prop:density}, the density comonad $\gen{\xX}$, if it exists,
  satisfies the universal property that comonad maps
  $\gen{\xX} \to \cC$ are in bijection with $\cC$-coalgebra structures
  on $\xX$. This density comonad exists and is pointwise assuming
  $\bS$ has copowers, defined via the formula
  $\gen{\xX}(\xA) = \bS(\xX, \xA) \cdot \xX$.

  In the case $\bS\coloneqq\Set$, this comonad
  $\gen{\xX}(\xA) = \xX \times \xA^\xX$ corresponds to the indiscrete
  topological space with points $\pts{\xX}$, or equally well to the
  indiscrete category with objects $\obs{\xX}$.
\end{example}

\begin{definition}\label{def:lift}
  Let $\pP \colon \bB \to \bS$ have density comonad $\cP$. We denote
  by $\lP \colon \bB \to \bS$ the lift of $\pP$ that corresponds under
  \cref{prop:density} to the identity on $\PS$:
  \[\begin{tikzcd}[row sep=10pt,column sep=15pt]
      {\PS} && {\PS} \\
      & {\bS}
      \arrow[equals, from=1-1, to=1-3]
      \arrow["\carP"', from=1-1, to=2-2]
      \arrow["\carP", from=1-3, to=2-2]
    \end{tikzcd}
    \qquad
    \leftrightarrows
    \qquad
    \begin{tikzcd}[row sep=10pt,column sep=15pt]
      {\bB} && {\PS} \\
      & {\bS}
      \arrow["\lP", dashed, from=1-1, to=1-3]
      \arrow["\pP"', from=1-1, to=2-2]
      \arrow["\carP", from=1-3, to=2-2]
    \end{tikzcd}
  \]
\end{definition}

For any $\pP \colon \bB \to \bS$ with a density comonad, $\bB$ is thus
canonically interpreted as a diagram of coalgebras via $\lP$.

\begin{remark}\label{rem:lift}\label{rem:pointwiseformula}\label{rem:liftiscomparison}
  More explicitly, for any object $\xU$ of $\bB$, the coalgebra
  $\lP(\xU)$ is the map $\pP(\xU)\to \cP(\pP(\xU))$ provided by the
  defining natural transformation of
  the left Kan extension $\cP$.
  \[
    \begin{tikzcd}[column sep=15pt, row sep=15pt]
      {\bB} &{}& {\bS} \\
      & {\bS}
      \arrow["\pP", from=1-1, to=1-3]
      \arrow["\pP"', from=1-1, to=2-2]
      \arrow["\cP"', dashed, from=2-2, to=1-3]
      \arrow["\Downarrow"{description, pos=0.3}, draw=none, from=1-2, to=2-2]
    \end{tikzcd}
  \]
  In the case that $\cP$ is pointwise, this coalgebra structure
  $\pP(\xU)\to \cP(\pP(\xU))$ is the colimit inclusion map
  \[\pP(\xU) \to \colim_{\pP(\xV) \to \pP(\xU)} \pP(\xV)\]
  indexed by the identity
  $\id_{\pP(\xU)} \colon \pP(\xU) \to \pP(\xU)$.

  The functor $\lP$ will appear frequently in our study of comonads as
  spaces. If $\topbasis{\tX} \colon \opens{\tX} \to \Set$ is the
  canonical basis from \cref{ex:topbasis} of the comonad $\cX$
  corresponding to the topological space $\tX$, then the lift
  $\lift{\topbasis{\tX}} \colon \opens{\tX} \to \CoalgSet{\cX} \simeq
  \Sh(\tX)$ is the canonical functor
  $\opens{\tX} \hookrightarrow \Sh(\tX)$ viewing open sets of a space
  $\tX$ as representable (equivalently, subterminal) sheaves over
  $\tX$. More generally for a functor $\pP \colon \bB \to \Set$ with
  density comonad $\cP$, when we view elements of $\cP(\xA)$ as
  $\xA$-valued germs as in \cref{rem:densitygerms}, the coalgebra
  structure on $\cP(\xA)$ is given by sending each element $x$ in
  $\cP(\xA)$ to the identity germ $[\id_{\cP(\xA)}]_x$ about $x$.

  Moreover, in the situation of \cref{cor:easybases}, where $\pP$ is
  given by the composite of a dense functor $F \colon \bB \to \bC$ and
  a left adjoint functor $L \colon \bC \to \bS$, the functor $\lP$ is
  simply the composite of the dense functor $F$ and the comparison
  functor $K \colon \bC \to \PS$ from \cref{lem:comparison}. Indeed,
  the comparison functor $K$ is defined using the coalgebra structures
  $L(\xX) \to L(R(L(\xX)))$ given by the natural transformation
  $L \circ \eta$, which is also the defining natural transformation
  $L \rightarrow L \circ R \circ L$ for the density comonad
  $\gen{L} \cong L \circ R$; and running through the proof of
  \cref{lem:kancomposeleft} establishes that we obtain the defining
  natural transformation
  $L \circ F \rightarrow L \circ R \circ L \circ F$ for the density
  comonad $\gen{L \circ F} \cong L \circ R$ simply by pre-composing
  the dense functor $F$.

  In particular if $\pP$ is the composite of a dense functor $F$ and a
  comonadic functor $L$, the functor $\lP$ is simply the dense functor
  $F$. This applies to both \cref{ex:topbasis} (canonical bases of
  topological spaces) and \cref{ex:catbasis} (canonical bases of
  categories), where $\lP$ is the dense inclusion of representable
  sheaves or functors.
  
  We note also that in the situation of \cref{prop:densityadj}, where
  the density comonad $\cP$ is exhibited as the comonad induced by the
  left adjoint extension $\extend{\pP} \colon \Ps{\bB} \to \bS$ of
  $\pP$ to its free cocompletion, the induced lift
  $\lift{\extend{\pP}} \colon \Ps{\bB} \to \PS$, which is the same as
  the comparison functor $K$, is also the same as the extension of
  $\lP$ to the free cocompletion of its domain $\bB$, since the
  comonadic functor $\carP \colon \PS \to \bS$ creates colimits.
  %
\end{remark}

The following tells us that we may calculate the ``points''
$\pts{\cP} = \cP(1)$ of a density comonad $\cP$ as the colimit
$\pts{\pP}$ of $\pP$. This is as expected for topological spaces: the
colimit of $\topbasis{\tX} \colon \opens{\tX} \to \Set$ is simply the
set of points $\pts{\tX}$.  
\begin{remark}\label{rem:pts}
  Suppose $\pP \colon \bB \to \bS$ has a pointwise density comonad
  $\cP$. If $\bS$ has a terminal object, then $\pts{\cP} = \cP(1)$ is
  given by the colimit $\pts{\pP}$ of $\pP$:
  \[
    \pts{\pP} \coloneqq \colim_{\xU\in\bB}\pP(\xU)\cong\cP(1) = \pts{\cP}.
  \]
  This is because the density comonad $\cP$ applied to $1$ is the
  canonical colimit of shape $(\pP/1) \cong \bB$, i.e.\ the colimit of
  $\pP$.
\end{remark}

In fact, the same holds without assuming pointwiseness or that $\bS$ has a terminal object:

\begin{lemma}\label{lem:pts}
  Let $\pP \colon \bB \to \bS$ have density comonad $\cP$. Then $\PS$
  has a terminal object if and only if $\pP$ has a colimit
  $\pts{\pP}$.
\end{lemma}
\begin{proof}
  By~\cite[Proposition 7.6]{fairbanks:monads} (which we also recall in
  \cref{app:universal} as \cref{lem:densityformaltfae}), given another
  comonad $\cC$ on a category $\bS'$ and commutative squares
  \[
    \begin{tikzcd}[column sep=20]
      {\PS} & \CoalgOn{\cC}{\bS'} \\
      {\bS} & {\bS'}
      \arrow["G_1", from=1-1, to=1-2]
      \arrow["\carP"', from=1-1, to=2-1]
      \arrow["\carC", from=1-2, to=2-2]
      \arrow["F_1"', from=2-1, to=2-2]
    \end{tikzcd}
    \qquad\qqand\qquad
    \begin{tikzcd}[column sep=20]
      {\PS} & \CoalgOn{\cC}{\bS'} \\
      {\bS} & {\bS'}
      \arrow["G_2", from=1-1, to=1-2]
      \arrow["\carP"', from=1-1, to=2-1]
      \arrow["\carC", from=1-2, to=2-2]
      \arrow["F_2"', from=2-1, to=2-2]
    \end{tikzcd}
  \]
  where $F_1$ preserves the left Kan extension $\cP$, the canonical
  map from natural transformations $G_1 \Rightarrow G_2$ to natural
  transformations $G_1 \circ \lP \Rightarrow G_2 \circ \lP$, given by
  pre-composing $\lP \colon \bB \to \PS$, is bijective.

  Taking $F_1$ and $G_1$ to be identities $\id_{\bS}$ and $\id_{\PS}$,
  and taking $F_2$ and $G_2$ to be a constant functors, we obtain a
  natural bijection from cocones under the diagram $\id_{\PS}$ to
  cocones under the diagram $\lP$. Thus $\id_{\PS}$ has a colimit if
  and only if $\lP$ has a colimit; and in turn $\lP$ has a colimit if
  and only if $\pP$ has a colimit since comonadic functors preserve
  and create colimits. The result now follows from the fact that a
  colimit of an identity functor is the same as a terminal object.
\end{proof}

\begin{warning}
  In the special case of a topological space $\tX$, the elements in
  each open set $\xU$ are identified with a subset of the points
  $\pts{\tX}$ of $\tX$. That is, if we denote the forgetful functor
  from the poset of open subsets by
  $\topbasis{\tX} \colon \opens{\tX} \to \Set$, then each map
  $\topbasis{\tX}(\xU) \to \pts{\tX}$ is injective.

  This need not be the case for general comonads on $\Set$. As in
  \cref{rem:pts}, given an arbitrary diagram $\pP \colon \bB \to \Set$
  with density comonad $\cP$, the ``points'' of $\cP$ are the colimit
  $\pts{\pP}$ of $\pP$. But clearly the maps $\pP(\xU) \to \pts{\pP}$
  forming the colimiting cocone of a diagram $\pP$ are generally not
  injective. This is no real surprise, since we are precisely
  generalizing topological spaces by replacing the diagrams of open
  subset inclusions with arbitrary diagrams of sets.

  Therefore the concept of an element in some $\pP(\xU)$ is distinct
  from the concept of a point in $\pts{\pP}$. One might call the
  former a ``local point'' (or ``fine point'') and the latter a
  ``global point'' (or ``coarse point''), though having pointed out
  this subtlety here we will not dwell on it further in later
  sections.

  Of particular importance, when as in \cref{rem:densitygerms} we
  describe elements of $\cP(\xA)$ as \emph{germs} $[f]_x$ of functions
  $\pP(\xU) \to \xA$ about some $x \in \pP(\xU)$, here $x$ is
  \emph{not} simply a point in $\pts{\pP}$. (For example, the density
  comonad $\cP$ from \cref{ex:part} whose coalgebras are partitioned
  sets $\Part$ has but one point $\pts{\pP} \cong 1$, but it is not
  sufficient to consider germs of functions $1 \to \xA$, as there are
  more germs $\cP(\xA)$ than simply elements of $\xA$.) In contrast,
  for a topological space $\tX$, a germ $[f]_x$ of a function
  $\xU \to \xA$ from a neighborhood of $x$ is always equivalent to a
  germ $[f']_x$ of some extended function
  $f' \colon \pts{\tX} \to \xA$ by restricting along the inclusion
  $\xU \hookrightarrow \pts{\tX}$, and so in this case one can assume
  germs of functions are centered about points in $\pts{\tX}$.
\end{warning}

\chaptertocspace
\chapter{Spaces}\label[section]{sec:spaces}

In this section we characterize the comonads associated to topological spaces as the density
comonads of diagrams of subsets, i.e.\ topological
subbases. We construct for each comonad $\cC$ on $\Set$ an underlying
topological space $\cmdtotop{\cC}$ (\cref{def:underlyingspace}) along
with a canonical comonad map to $\cC$ from the comonad whose
coalgebras are sheaves on $\cmdtotop{\cC}$ (\cref{prop:sheavesff}). In
\cref{thm:tfaetop} we characterize when this canonical comonad map is
an isomorphism and show that this is the only case in which the
category of $\cC$-coalgebras is equivalent to the category of sheaves
on a topological space.

\section{Topological spaces from comonads}\label{sec:toptocomonad}

First we show that every coalgebra of every comonad on $\Set$ is
endowed with a canonical topology. Intuitively, this generalizes the
\'etal\'e space construction for sheaves. The underlying topological
space of the comonad $\cC$ is then given by the topology of the
terminal coalgebra $\pts{\cC}$.

The open subsets of the \'etal\'e space associated to a sheaf are precisely its
subsheaves, and the open subsets of a general $\cC$-coalgebra will be
similarly defined.

\begin{definition}\label{def:open}
  Let $\cC$ be a comonad on $\Set$. An \emph{open subset} of a $\cC$-coalgebra is a
  subobject given by a regular monomorphism --- equivalently by
  \cref{lem:regularmono}, carried by a monomorphism in $\Set$.
\end{definition}

Note that open subsets are indeed merely certain subsets, since by
\cref{lem:factor}, any subset of the carrier of a coalgebra admits at
most one coalgebra structure preserved by the subset inclusion.

\begin{lemma}
  Let $\cC$ be a comonad on $\Set$. The comonadic functor
  $\carC \colon \SetC \to \Set$ creates pullback squares consisting of
  injections in $\Set$.
\end{lemma}
\begin{proof}
  In \cref{prop:intersections} we saw that all comonads on $\Set$
  preserve finite intersections. The result follows because a
  comonadic functor $\carC$ creates limits preserved by $\cC$ and
  $\cC \circ \cC$ (\cref{lem:create}).
\end{proof}


\begin{corollary}\label{lem:intersection}
  The open subsets of any coalgebra of a comonad on $\Set$ are closed
  under element-wise finite intersections.\qed
\end{corollary}

\begin{lemma}\label{lem:coconfluentinjections}
  Let $\pP \colon \bB \to \Set$ be the diagram of subset inclusions
  constituting a basis of a topological space $\tX$. Then the colimit
  $\pts{\pP}$ is the set of points $\pts{\tX}$, and the inclusions
  into the colimit $\pP(\xU) \to \pts{\pP}$ are the inclusions of basic
  open subsets.
\end{lemma}
\begin{proof}
  The condition that $\pP$ constitutes a topological basis ensures,
  for each point $x$ of $\pts{\tX}$, that first of all there is a
  basic $\xU$ with $\pP(\xU)$ containing a representative of $x$, and
  that second of all if $x_1 \in \pP(\xU_1)$ and $x_2 \in \pP(\xU_2)$ both
  represent $x$, then there is a basic $\xU' \subseteq \xU_1 \cap U_2$ and
  $x' \in \pP(\xU')$ also representing $x$, and thus $x_1$ and $x_2$ are
  identified in the colimit.
\end{proof}

\begin{corollary}\label{cor:subcoalgebrabasis}
  Let $\eE\to\cC(\eE)$ be a coalgebra of a comonad $\cC$ on $\Set$.
  Suppose given a collection of open sets (regular subobjects)
  constituting a topological basis on $\eE$. Then the colimit of the
  diagram in $\CSet$ of inclusions between them is their element-wise
  union.
\end{corollary}
\begin{proof}
  Note that by \cref{lem:factor}, the diagram of inclusions between
  open subsets in $\CSet$ agrees with the diagram of subset inclusions
  in $\Set$. The result follows from \cref{lem:coconfluentinjections},
  since the comonadic functor creates colimits.
\end{proof}

\begin{corollary}\label{lem:union}
  The open subsets of any coalgebra of a comonad on $\Set$ are closed
  under arbitrary unions.
\end{corollary}
\begin{proof}
  Apply \cref{cor:subcoalgebrabasis} to the diagram of pairwise
  intersections of open subsets given by \cref{lem:intersection}.

  Alternatively, we may form the coproduct (disjoint union) of the
  open subsets, created by the comonadic functor, and take its regular
  subobject image, which exists and is the image in $\Set$ by
  \cref{lem:regularmono}.
\end{proof}

\begin{proposition}\label{prop:topology}
  If $\cC$ is a comonad on $\Set$ and $\eE\to\cC(\eE)$ is a
  $\cC$-coalgebra, then the open subsets form a topology on $\eE$.
\end{proposition}
\begin{proof}
  They are closed under finite intersections and
  unions (\cref{lem:intersection,lem:union}).
\end{proof}


\begin{remark}
  The above \cref{prop:topology} is not really new: an analogous
  result was shown in~\cite{gumm-schroder:bounded} for coalgebras of
  endofunctors on $\Set$ rather than coalgebras for comonads. The
  result for comonads (\cref{prop:topology}) follows without much work
  from the result for endofunctors: using the fact that all comonads
  on $\Set$ preserve monomorphisms, one can show that every
  endofunctor coalgebra admitting an injective coalgebra map into a
  comonad coalgebra is itself a comonad coalgebra. Therefore the two
  notions of open subset, defined via comonad coalgebra maps carried
  by injections versus endofunctor coalgebra maps carried by
  injections, coincide for comonad coalgebras.
\end{remark}

\begin{definition}\label{def:underlyingspace}
  The \emph{underlying topological space} $\cmdtotop{\cC}$ of a
  comonad $\cC$ on $\Set$ is $\pts{\cC}$, the terminal coalgebra,
  equipped with the topology from \cref{prop:topology}.

\end{definition}

\begin{remark}
  We may also view the construction of the underlying space of a
  comonad on $\Set$ through another lens, similar to \cite[Example
  4.8]{garner:ionads}. A topology on the set $\xX$ may equivalently be
  characterized as an \emph{interior operator}: a
  finite-meet-preserving comonad on the poset of subsets $2^\xX$. As
  we will show later in \cref{cor:setslices}, a comonad on $\Set$
  sending the terminal object to $\xX$ may equivalently be viewed as a
  comonad on $\Set^\xX$ preserving the terminal object and finite
  intersections. This restricts to a comonad on the full subcategory
  of subterminal objects in $\Set^\xX$, which is (up to equivalence)
  the poset $2^\xX$, and preservation of the terminal object and
  finite intersections here is the same as preservation of finite
  meets.

  This agrees with the above description of the underlying space,
  since an open subset with respect to the interior operator is a
  coalgebra of the comonad on $2^\xX$, which amounts to a coalgebra of
  the comonad on $\Set^\xX$ carried by a subterminal object, which in
  turn amounts to an open subset of the terminal coalgebra in the
  sense of \cref{def:open}.
\end{remark}

In passing, we note that the underlying topological space (as in
\cref{prop:topology}) of a coalgebra $\eE$ agrees with the underlying
topological space (as in \cref{def:underlyingspace}) of its comonad of
elements $\El(\eE)$ from \cref{def:comonadofelements}:

\begin{proposition}
  Let $\cC$ be a comonad on $\Set$ and let $\eE$ be a
  coalgebra. The underlying topological space of the comonad
  $\El(\eE)$ is the underlying topological space of the coalgebra
  $\eE$.
\end{proposition}
\begin{proof}
  The set of points of $\El(\eE)$ is given by the image of the
  terminal object under the comonadic forgetful functor
  $\CoalgSet{\El(\eE)} \to \CSet \to \Set$, i.e.\ the underlying
  set of $\eE$. The topology of $\eE$ is the same as the
  underlying topology of $\El(\eE)$, since they are both given by
  coalgebra maps into $\eE$ carried by injections.
\end{proof}

In general, maps of coalgebras are not necessarily continuous, i.e.\
the preimage of an open subset under a map of coalgebras need not be
open. The comonads $\cC$ on $\Set$ for which this does hold, yielding
a forgetful functor $\CSet \to \Top$, are precisely the taut
comonads. In particular, by \cref{lem:weaktaut}, every comonad
preserving weak pullbacks satisfies this condition.

\begin{proposition}\label{prop:tautpreimages}
  A comonad $\cC$ on $\Set$ is taut (i.e.\ preserves preimages) if and
  only if all coalgebra maps induce continuous maps of underlying
  spaces.
\end{proposition}
\begin{proof}
  The left-to-right direction is because the comonadic functor
  $\carC \colon \CSet \to \Set$ creates limits preserved by $\cC$ and
  $\cC \circ \cC$ (\cref{lem:create}). Conversely, assume preimages of
  open subsets are open. Then any pullback in $\Set$ of an open subset
  along a coalgebra map lifts to a square of coalgebra maps by
  \cref{lem:factor} which is a pullback square in $\CSet$ by
  \cref{lem:reflect}.
\end{proof}

Recall that a \emph{wide pullback} is a limit indexed by a category
obtained by freely adjoining a terminal object to a discrete
category. We say that a comonad on $\Set$ \emph{preserves arbitrary
  intersections} if it preserves wide pullbacks consisting of
monomorphisms.

\begin{proposition}\label{prop:alextop}
  A comonad $\cC$ on $\Set$ preserves arbitrary intersections if and
  only if open subsets of coalgebras are closed under arbitrary
  intersections, i.e.\ all coalgebras have Alexandroff underlying
  topological spaces.
\end{proposition}
\begin{proof}
  The comonadic functor $\carC \colon \CSet \to \Set$ creates limits
  preserved by $\cC$ and $\cC \circ \cC$ (\cref{lem:create}).
\end{proof}

\section{Comonads from topological subbases}\label{sec:comonadtotop}

We are now equipped to show that density comonads of functors into
$\Set$ generalize topological spaces generated by subbases.

\begin{proposition}\label{prop:subbases}
  Suppose $\pP \colon \bB \to \Set$ is a diagram with density comonad
  $\cP$ and such that all inclusion maps into the colimit
  $\pP(U) \to \pts{\pP}$ are injective. Then $\cP\cong\cX$, where
  $\tX$ is the topological space with points $\pts{\pP}$ and topology
  generated by the subbasis $\pP$.
\end{proposition}
\begin{proof}
  As described in \cref{ex:topbasis}, $\cX$ is the density comonad
  $\gen{\topbasis{\tX}}$ of
  $\topbasis{\tX} \colon \opens{\tX} \to \Set$. Recall the universal
  property of the density comonad from \cref{prop:density}:
  commutative triangles of the following forms are in natural
  bijection.
  \[
    \begin{tikzcd}[row sep=10pt,column sep=15pt]
      \CoalgSet{\gen{\topbasis{\tX}}} && \CSet \\
      & \Set
      \arrow["{\car{\gen{\topbasis{\tX}}}}"', from=1-1, to=2-2]
      \arrow[dashed, from=1-1, to=1-3]
      \arrow["{\carC}", from=1-3, to=2-2]
    \end{tikzcd}
    \qquad
    \leftrightarrows
    \qquad
    \begin{tikzcd}[row sep=10pt,column sep=15pt]
      \opens{\tX} && \CSet \\
      & \Set
      \arrow["{\topbasis{\tX}}"', hook, from=1-1, to=2-2]
      \arrow["F", dashed, hook, from=1-1, to=1-3]
      \arrow["{\carC}", from=1-3, to=2-2]
    \end{tikzcd}
  \]
  We must show that $\cX$ is also the density comonad of the subbasis
  $\pP \colon \bB \to \Set$. That is, we must show that $\cX$
  satisfies the same universal property with respect to $\pP$, so that
  commutative triangles as above are also in correspondence with
  commutative triangles
  \[
    \begin{tikzcd}[row sep=10pt,column sep=15pt]
      \bB && \CSet \\
      & \Set
      \arrow["{\pP}"', ""{name=0, anchor=center, inner sep=0}, from=1-1, to=2-2]
      \arrow["F'", dashed, from=1-1, to=1-3]
      \arrow["\carC", ""{name=1, anchor=center, inner sep=0}, from=1-3, to=2-2]
    \end{tikzcd}
  \]
  Suppose given such a functor $F' \colon \bB \to \CSet$. By
  assumption $\pP$ has colimit whose inclusion maps are
  injective. Since $\carC$ creates colimits, $F'$ is a diagram of open
  subsets of the colimit of $F'$, carried by $\pts{\pP}$, in
  $\CSet$. By \cref{prop:topology}, open subsets are closed under
  finite intersections and unions. Therefore this uniquely determines
  a map $F \colon \opens{\tX} \to \CSet$, as desired.
\end{proof}

\begin{remark}\label{rem:subbasiswarning}
  Let $\bB \subseteq \opens{\tX}$ be a subbasis of a topological space
  $\tX$ and let $\pP \colon \bB \to \Set$ denote the forgetful
  functor. Assuming the poset of subbasic neighborhoods of each point
  is connected, so that $\pts{\pP} \cong \pts{\tX}$, then we have
  $\cP \cong \cX$ by \cref{prop:subbases}. If we take a subbasis of $\tX$
  that does not satisfy this connectedness condition, the density
  comonad of the forgetful functor will instead correspond to the
  generated topology on $\pts{\pP}$, not $\tX$. One can ensure
  $\pts{\tX} \cong \pts{\pP}$ by including the whole space, i.e.\ the
  terminal object of $\opens{\tX}$, in the subbasis.
\end{remark}

\begin{remark}\label{rem:subbasisgerms}
  \cref{prop:subbases} also has a more concrete interpretation in terms
  of germs. The claim $\cP(\xA)\cong\cX(\xA)$ is that the set of all
  germs of $\xA$-valued functions on $\tX$ is equally well the colimit
  \[\sum_{x\in \pts{\tX}} \colim_{U \ni x} \xA^U\]
  where $U$ ranges over a subbasis of $\tX$ for which the poset of
  subbasic neighborhoods of each point is connected, as discussed in
  \cref{rem:subbasiswarning}.

  Indeed, we can alternatively give a direct proof that the germs of
  functions $\xX \to \xA$ with respect to a connected diagram of
  subsets satisfying the finite intersection property
  $\mathcal{B} \subseteq \Pow{\xX}$ of $\xX$ are the same as germs
  with respect to the generated proper filter $\mathcal{F}$ (the
  closure under finite intersections and supersets):
  \[\colim_{U \in \mathcal{B}} \xA^U \cong \colim_{U \in \mathcal{F}} \xA^U\]
  First we may assume without loss of generality that $\mathcal{B}$
  includes the maximal subset $X$, since by assumption of
  connectedness $\mathcal{B}$ is initial (dual to final,
  \cref{def:final}) in the poset given by adjoining $X$ to
  $\mathcal{B}$ and therefore the contravariantly-indexed colimits are
  the same. Now for any pair $U, V \in \mathcal{B}$, we may extend any
  map $U \cap V \to \xA$ to all of $\xX$, and all such extensions
  become equal in the colimit by alternately restricting to $U$ and
  $V$ and extending back again, establishing the claim.\footnote{If
    $\xA$ is $\varnothing$ then there is no such extension to all of
    $X$. Hence if we drop assumption that the members of $\mathcal{B}$
    pairwise intersect, so that an improper filter is generated, the
    colimit may then incorrectly evaluate to $\varnothing$ in this
    exceptional case.} Hence in particular if $\tX$ is a topological
  space and $x$ is a point, then germs of functions
  $\pts{\tX} \to \xA$ at $x$ can be calculated with respect to any
  connected neighborhood subbasis of $x$.\footnote{This claim that
    germs can be calculated with respect to connected neighborhood
    subbases is distinct from the more standard and straightforward
    one assuming a neighborhood \emph{basis}: a neighborhood basis is
    by definition an initial subpreorder of the complete system of
    neighborhoods, so any colimit indexed contravariantly by the
    complete system of neighborhoods may be calculated using a
    neighborhood basis in is place.}
\end{remark}

By \cref{def:underlyingspace}, the poset of opens of the
underlying topological space $\cmdtotop{\cC}$ of a comonad $\cC$ on
$\Set$ embeds in $\CSet$ as the regular subobjects (open subsets)
of the terminal coalgebra $\pts{\cC}$. In fact, more is true: the
entire category of sheaves $\Sh(\cmdtotop{\cC})$ embeds fully
faithfully in $\CSet$.

\begin{proposition}\label{prop:sheavesff}
  The canonical comonad map $\phi \colon \topreflect{\cC} \coto \cC$,
  obtained via \cref{prop:density} by the evident commutative triangle
  \[\begin{tikzcd}[sep=small]
      {\opens{\cmdtotop{\cC}}} && \CSet \\
      & \Set
      \arrow[hook, from=1-1, to=1-3, "i"]
      \arrow["\topbasis{\cmdtotop{\cC}}"', from=1-1, to=2-2]
      \arrow["\carC", from=1-3, to=2-2]
    \end{tikzcd}\] is componentwise monic.

  Therefore, the induced functor
  $\com{\Set}{\phi} \colon \Sh(\cmdtotop{\cC}) \to \CSet$ is fully
  faithful.
\end{proposition}
\begin{proof}
  The second statement follows from the first by
  \cref{lem:componentwisemonic}. For the first statement, note that
  the natural transformation carrying the comonad map $\phi$ is
  obtained by whiskering the counit of the adjunction induced by
  $\com{\Set}{\phi}$ with $\rcarC \colon \Set \to \CSet$ and
  $\carC \colon \CSet \to \Set$.
  \[
    \begin{tikzcd}[column sep=30pt]
      \Sh(\cmdtotop{\cC})\ar[r, shift left=6pt, "\com{\Set}{\phi}"]\ar["\bot"{anchor=center},r, draw=none]&
      \CSet\ar[l, shift left=6pt, "\comr{\Set}{\phi}"]\ar[r, shift left=6pt, "\carC"]\ar["\bot"{anchor=center},r, draw=none]& \Set\ar[l, shift left=6pt, "\rcarC"]
    \end{tikzcd}\] Hence it is enough to show that this counit is
  componentwise injective.
  Indeed, the right adjoint is given by
  \[\comr{\Set}{\phi}(\eE) \cong \colim_{U \to \comr{\Set}{\phi}(\eE)} U \cong \colim_{i(U) \to \eE} U\] where $\xU$ ranges over $\opens{\cmdtotop{\cC}}$,
  since $\opens{\cmdtotop{\cC}}$ is dense in
  $\Sh(\cmdtotop{\cC})$. Thus the induced comonad on $\CSet$ sends
  $\eE$ to
  \[\com{\Set}{\phi}\left(\colim_{i(U) \to \eE} U\right) \cong \colim_{i(U) \to \eE} i(U).\]
  The counit of the adjunction is the induced map into $\eE$ from
  this colimit, which is injective by \cref{cor:subcoalgebrabasis}.
\end{proof}

This full subcategory $\Sh(\cmdtotop{\cC})$ may be recognized within
$\CSet$ in the following various ways:

\begin{proposition}\label{lem:bangetale}
  Let $\cC$ be a comonad on $\Set$. The following are equivalent for a
  coalgebra $\eE$ in $\CSet$.
  \begin{enumerate}[label=(\roman*)]
  \item $\eE$ is the canonical colimit in $\CSet$ of shape
    $(i / \eE)$, where
    $i \colon \opens{\cmdtotop{\cC}} \hookrightarrow \CSet$ is the
    inclusion of open subsets (i.e.\ regular subobjects) of the terminal
    coalgebra.\label{item:cancolimit}
  \item $\eE$ is covered by a jointly epimorphic family of maps in
    $\CSet$ from open subsets of the terminal
    coalgebra.\label{item:iscolimit}
  \item $\eE$ is in the essential image of the fully faithful functor
    $\Sh(\cmdtotop{\cC}) \to \CSet$ (from
    \cref{prop:sheavesff}).\label{item:issheaf}
  \item The unique coalgebra map $\bang\colon \eE\to\pts{\cC}$ to
    the terminal coalgebra constitutes an \'etale map of the
    associated topological spaces (from
    \cref{prop:topology}).\label{item:isetale}
  \end{enumerate}
\end{proposition}
\begin{proof}
  \labelcref{item:cancolimit} $\imp$ \labelcref{item:iscolimit} is
  immediate. For \labelcref{item:iscolimit} $\imp$
  \labelcref{item:issheaf}, note that the fully faithful inclusion
  $\Sh(\cmdtotop{\cC}) \to \CSet$ from \cref{prop:sheavesff} is a
  coreflective subcategory by \cref{lem:adjointlift} and thus is
  closed under colimits. Moreover, every coalgebra admitting a jointly
  epimorphic family of coalgebra maps from regular subterminal
  coalgebras, which note are necessarily regular monic, is also a
  colimit of subterminal coalgebras by \cref{lem:intersection} and
  \cref{lem:coconfluentinjections}.
  
  \labelcref{item:issheaf} $\imp$ \labelcref{item:isetale} is
  immediate because maps between sheaves are \'etale.

  For \labelcref{item:isetale} $\imp$ \labelcref{item:cancolimit},
  suppose $\bang \colon \eE \to \pts{\cC}$ is \'etale (a local
  homeomeorphism). Then every point of the domain has an open
  neighborhood mapped homeomorphically to an open subset of
  $\pts{\cC}$ via $\bang$. Thus there is a jointly epimorphic family
  $\xU_i \colon \eE_i \hookrightarrow \eE$ consisting of open
  subsets such that each $\bang \circ U_i$ is an open subset of
  $\pts{\cC}$. Now by \cref{cor:subcoalgebrabasis}, $\eE$ is the
  canonical colimit of the diagram of inclusions between all $\xU_i$.
\end{proof}

We may thus characterize topological comonads in the following various
ways:

\begin{theorem}\label{thm:tfaetop}
  Let $\cC$ be a comonad on $\Set$. The following are equivalent.
  \begin{enumerate}[label=(\roman*)]
  \item\label{item:topological} $\cC$ is the germs comonad $\cX$ of a
    topological space $\tX$.
  \item\label{item:owntopological} $\cC$ is the germs comonad
    $\toptocmd{(\cmdtotop{\cC})}$ of its underlying topological space
    $\cmdtotop{\cC}$, i.e.\ it is the density comonad of the diagram
    $\topbasis{\cmdtotop{\cC}} \colon \opens{\cmdtotop{\cC}} \hookrightarrow
    \Set$.
  \item\label{item:subsets} $\cC$ is the density comonad of a diagram
    of subsets of a set --- that is, a diagram
    $\pP \colon \bB \to \Set$ where all maps in the colimiting cocone
    to $\pts{\pP}$ are injective.
  \item\label{item:subsetsterm} $\cC$ is the density comonad of a
    diagram $\pP \colon \bB \to \Set$ consisting of injections where
    $\bB$ has a terminal object.
  \item\label{item:densesubterms} The inclusion
    $\opens{\cmdtotop{\cC}} \hookrightarrow \CSet$ of open subsets
    (i.e.\ regular subobjects) of the terminal coalgebra is dense.
  \item\label{item:subterms} Every object of $\CSet$ is covered by
    open subsets of the terminal coalgebra.
  \item\label{item:alletale} All maps in $\CSet$ induce \'etale maps of
    the coalgebras' underlying topological spaces.
  \item\label{item:localic} $\CSet$ is a localic topos, i.e.\ the
    category of sheaves on a locale.
  \end{enumerate}
\end{theorem}
\begin{proof}
  We have \labelcref{item:topological} $\iff$
  \labelcref{item:owntopological} because the underlying topological
  space of the germs comonad of a topological space is the same
  topological space (that is, $\cmdtotop{(\cX)} \cong \tX$).

  \labelcref{item:topological} $\iff$ \labelcref{item:subsets} $\iff$
  \labelcref{item:subsetsterm} follows from \cref{prop:subbases}.
  
  \labelcref{item:topological} $\iff$ \labelcref{item:densesubterms}
  $\iff$ \labelcref{item:subterms} $\iff$ \labelcref{item:alletale}
  follows from \cref{lem:bangetale}.

  \labelcref{item:subterms} $\iff$ \labelcref{item:localic}
  follows immediately, since localic toposes are precisely those
  generated under colimits by subterminal objects.
\end{proof}

\chaptertocspace
\chapter{Bases}\label[section]{sec:bases}

In \cref{sec:spaces}, we saw that functors $\pP \colon \bB \to \Set$
generalize \emph{sub}bases of topological spaces. This leads one to ask:
which functors should be called generalized \emph{bases}?

In \cref{sec:basessub} we introduce a definition of basis that applies
to arbitrary comonads on arbitrary categories $\bS$. In fact we give
two equivalent characterizations of bases $\pP \colon \bB \to \bS$,
one analogous to the condition that coalgebras may be interpreted as
sheaves on $\bB$ (\cref{def:basis}) and the other analogous to the
condition that the sum of stalks yields the points of the \'etal\'e
space (\cref{prop:basis}). These equivalent conditions coincide with
the usual notion of topological basis for diagrams of subsets
$\pP \colon \bB \to \Set$. However, these are far from the only
categorical bases of comonads on $\Set$. (Indeed, every comonad on
every category admits a canonical, possibly large, basis, given by the
comonadic functor itself.) We provide other examples as well.

In \cref{sec:subbases} we examine the relationship between abstract
bases and abstract subbases.  We observe that whereas bases of $\cC$
amount to functors $\bB \to \CS$ such that coalgebras are realized as
presheaves on $\bB$, (pointwise) subbases of $\cC$ amount to functors
$\bB \to \CS$ such that just \emph{cofree} coalgebras are realized as
presheaves on $\bB$ (\cref{prop:subbasis}).

In \cref{sec:coconfluence} we study a more concrete condition on
functors $\pP \colon \bB \to \Set$ resembling the definition of
topological basis (\cref{def:influent}). We show that all functors
satisfying this condition are categorical bases
(\cref{cor:influentbasis}), and that the corresponding comonads are
precisely those that preserve weak pullbacks (\cref{cor:wpb}) as
studied in the field of
coalgebra~\cite{rutten,johnstone-power-tsujishita-watanabe-worrell,jacobs,adamek-velebil,gumm-schroder,gumm:elements}. We
also relate these bases to the \emph{bases of ionads} studied
in~\cite{garner:ionads}, which are more specific.

\section{Bases}\label{sec:basessub}

Our definition of basis of a comonad will be motivated by the
following characterization of bases of topological spaces.

\begin{lemma}\label{lem:topbasis}
  A subbasis $\bB \hookrightarrow \opens{\tX}$ of a topological
  space $\tX$ is a basis if and only if the canonical inclusion
  $\bB \hookrightarrow \Sh(\tX)$ is dense.
\end{lemma}
\begin{proof}
  By definition, $\bB$ is a topological basis of $\tX$ if and only if
  every open subset $\xU \in\opens{\tX}$ is a union of those in
  $\bB$. If $\bB \hookrightarrow \Sh(\tX)$ is dense, then $\bB$ is a
  basis, since every $\xU \in \opens{\tX}$ is a canonical colimit within
  $\Sh(\tX)$ of basic opens, thus in particular is covered by
  them. Conversely, if $\bB$ is a topological basis, then
  $\bB \hookrightarrow \Sh(\tX)$ is dense, since bases of $\tX$ lift
  to bases of \'etal\'e spaces over $\tX$, and thus by
  \cref{cor:subcoalgebrabasis} every sheaf is the canonical colimit of
  basic opens in $\Sh(\tX)$.
\end{proof}

\begin{definition}\label{def:basis}
  A functor $\pP \colon \bB \to \bS$ with density comonad $\cP$ is a
  \emph{basis} if the induced functor $\lP \colon \bB \to \PS$ from
  \cref{def:lift}, sending $U\in \bB$ to the coalgebra 
  $P(U)\to \cP(P(U))$ determined by the construction of $\cP$ as a left
  Kan extension, is dense.
\end{definition}

Intuitively, this says $\pP$ is a basis if coalgebras may be regarded
as ``sheaves'' on $\bB$, in the sense that the embedding
$\PS \to \Ps{\bB}$ is fully faithful.

\begin{remark}\label{rem:easybases}
  Suppose $\cC$ is a comonad on $\bS$. Recall from
  \cref{cor:easybases} that given an arbitrary dense functor
  $F \colon \bB \to \CS$, we recover $\cC$ as the density comonad
  $\cP$ of \[\pP \colon \bB \xto{F} \CS \xto{\carC} \bS.\] Moreover,
  we noted in \cref{rem:liftiscomparison} that $F$ is recovered as
  the lift $\lP$ associated to $\pP$. Thus a basis for the comonad $\cC$
  simply amounts to a dense functor into $\CS$.
\end{remark}

\begin{remark}\label{rem:pointwisebasis}
  Density comonads of bases $\pP \colon \bB \to \bS$ are always
  pointwise: by definition the lift $\lP$ is dense, i.e.\ has identity
  pointwise density comonad, and composing with the left adjoint
  $\carP$ preserves the pointwiseness by \cref{rem:composepointwise}.
\end{remark}

\begin{example}[Canonical bases]\label{rem:indeedbases}
  The ``canonical basis'' $\topbasis{\tX} \colon \opens{\tX} \to \Set$
  of a topological space $\tX$ from \cref{ex:topbasis} (which is also
  just a particular example of an ordinary topological basis), the
  ``canonical basis'' $\catbasis{\bC} \colon \bC\op \to \Set$ of a
  category from \cref{ex:catbasis}, and the ``big basis''
  $\carC \colon \CS \to \bS$ of an arbitrary comonad $\cC$ from
  \cref{cor:densityofleft} are all obtained by composing a dense
  functor with a comonadic functor. Hence as per \cref{rem:easybases},
  all are indeed bases.
\end{example}

The following is useful in practice for establishing functors are
bases of comonads:

\begin{lemma}\label{lem:extcoref}
  If $\pP \colon \bB \to \bS$ with $\bB$ small and $\bS$ cocomplete is
  such that the extension $\extend{\pP} \colon \Ps{\bB} \to \bS$
  preserves coreflexive equalizers --- or equivalently by
  \cref{lem:functorpreserve}, preserves finite intersections --- then
  $\pP$ is a basis.
\end{lemma}
\begin{proof}
  \cref{lem:crepff} tells us that, under these assumptions, the right adjoint
  $R \colon \PS \to \Ps{\bB}$ to the comparison functor
  $K \colon \Ps{\bB} \to \PS$ is fully faithful. As observed in
  \cref{rem:liftiscomparison}, the comparison functor $K$ is the
  extension $\extend{\lP}$ of $\lP$, so $R$ is $\nerve{\lP}$. And by
  definition $\pP$ is a basis just when $\nerve{\lP}$ is fully
  faithful.
\end{proof}

\anoteNI{Can the assumptions be weakened? And does this also work for
  some notion of ``coreflexive-equalizer-flatness'', more general than
  coreflexive-equalizer-preservation?}

\begin{remark}\label{rem:simplergerms}
  The condition that $\pP \colon \bB \to \Set$ is a basis simplifies
  the description of coalgebras in terms of germs
  from~\cref{rem:densitygerms}. Recall that an element of $\cP(\xA)$
  is a germ $[f]_x$ of a function $f \colon \pP(\xU) \to \xA$, and
  that a coalgebra $\hH \colon \eE \to \cP(\eE)$ consists of an
  $\eE$-valued germ $\hH(e) = [f_e]_{x_e}$ for each $e$ in $\eE$
  satisfying $f_e(x_e) = e$ (the unit law) and
  $[\hH \circ f_e]_{x_e} = [[f_e]_{(\dash)}]_{x_e}$ (the associativity
  law). The associativity law is clearly the more complicated
  condition to check.

  However, assuming $\pP$ is a basis, i.e.\ $\lP \colon \bB \to \PSet$
  is dense, then by the canonical colimit formula, every element of
  every coalgebra $\eE$ is in the image of a coalgebra map from some
  $\lP(\xU)$. (This condition is strictly weaker than the condition
  that $\pP$ is a basis; see \cref{cex:colimitsubbasis}.)
  Equivalently, this is to say every element of $\eE$ is in the image
  of some function $f \colon \pP(\xU) \to \eE$ where for all
  $x \in \pP(\xU)$, the coalgebra structure map $\hH$ assigns the
  element $f(x)$ the germ $[f]_{x}$. This simpler condition implies
  both the unit law and the associativity law, and therefore it is all
  we need to check in order to verify $\eE$ is a coalgebra.

  In the case of a basis of a topological space $\tX$, this means that
  a set $\eE$ is equipped with the structure of \'etal\'e space over
  $\tX$ if and only if every element is assigned the $\eE$-valued germ
  of some map from a basic open set $f \colon U \to \eE$ where for all
  $x \in U$, the germ assigned to $f(x)$ is $[f]_x$. In short, an
  \'etal\'e space over $\tX$ is equivalently a set equipped with
  infinitesimal neighborhoods at points induced by some covering of it
  with maps from basic open sets of $\tX$.
\end{remark}

We now give a worked example, which is in particular an example of a
\emph{basis of an ionad}, a notion introduced in~\cite{garner:ionads}.

\begin{example}[Infinitesimal group actions]\label{ex:real_line}\label{ex:etendue}\label{ex:topgroup}
  A familiar topological basis of the real line $\rr$ is given by the open
  intervals
  \[
    \setof{(x_1,x_2)\mid x_1 \leq x_2}.
  \]
  These subsets form a poset under inclusion, and by
  \cref{prop:subbases} the density comonad of the forgetful functor
  into $\Set$ is the germs comonad $\toptocmd{\rr}$ of $\rr$. By
  \cref{lem:topbasis}, this is also a basis in the abstract sense of
  \cref{def:basis}: every sheaf, or \'etal\'e space, over $\rr$ is the
  canonical colimit of the diagram of all maps into it from open
  intervals.

  But we can also incorporate the additive group structure of
  $\rr$. Define $\bB \coloneqq \Int$ to be the category of open
  intervals in $\rr$ and all orientation-and-distance-preserving maps,
  i.e.\ translations,\footnote{The undirected version of this, with
    distance-preserving but not necessarily orientation-preserving
    maps, is also interesting to consider.}\footnote{$\Int$ is
    equivalently the twisted arrow category of the monoid of
    nonnegative reals.} and let $\pP \colon \bB \hookrightarrow \Set$
  be the forgetful functor sending an interval to its set of
  points.\footnote{The \emph{topos of temporal behavior types} from
    \cite{schultz-spivak} is the category of coalgebras on a different
    basis with the same category $\Int$ of opens, namely the functor
    sending $(x_1,x_2)$ to
    $\setof{(y_1,y_2)\mid x_1<y_1\leq y_2<x_2}$, the set of closed
    intervals within the open interval. For comparison, this is 
    sheaves on the topological action groupoid of $\mathbb{R}^\delta$ 
    acting on the interval site, rather than just $\mathbb{R}$, so 
    that our topos sits over that of \cite{schultz-spivak} as the instantaneously
    generated temporal behaviors.}  Or alternatively, we may
  take $\bB$ to be the equivalent full subcategory consisting of just
  those intervals centered about $0$, since any two intervals of the
  same size are isomorphic.

  The density comonad $\cP$ is then described in the germ terminology
  of \cref{rem:densitygerms} as follows: an $\xA$-valued germ
  $[f]_x \in \cP(\xA)$ consists of an element of a basic open set
  $x \in U$ and a function $f \colon U \to \xA$, modulo
  $[f]_x = [g]_y$ whenever there exists a translation
  $i \colon U \hookrightarrow V$ such that $i(x) = y$ and
  $f = g \circ i$. Any germ $[g]_y$ about $y$ is then equivalent to a
  canonical germ $[f]_0$ about $0$ by translation; therefore an element of
  $\cP(\xA)$ may simply be described as a germ about $0$ of an
  $\xA$-valued function $\rr \to \xA$. In particular, setting $\xA=1$,
  since there is only one germ of a function into $1$ about $0$, we see
  that the ``set of points'' $\pts{\cP}=\cP(1)$ of $\cP$ is a singleton.

  Moreover, the forgetful functor $\pP \colon \bB \to \Set$ is
  \emph{flat},\footnote{A \emph{basis of an ionad} in the sense
    of~\cite{garner:ionads} is by definition a flat functor
    $\pP \colon \bB \to \Set^\xX$ for some set $\xX$. We will return
    to this concept in \cref{sec:coconfluence}.} i.e.\ its category of
  elements is cofiltered. Indeed, there is at most one arrow between
  any two objects in the category of elements, and any pair of
  elements admit an element mapping to both. It is well-known that
  flatness of a functor $\pP \colon\bB \to \Set$ is equivalent to
  the condition that $\extend{\pP} \colon \Ps{\bB} \to \Set$ preserves
  finite limits~\cite[Proposition 6.3.8]{borceux:i}. Therefore $\pP$
  is a basis in the abstract sense by \cref{lem:extcoref}.

  Now that we have established $\pP$ is a basis, we are permitted the
  simplified description of coalgebras from \cref{rem:simplergerms}. A
  $\cP$-coalgebra $\eE \to \cP(\eE)$ assigns to each element $e$ the
  germ about $0$ of some function $f \colon U \to \eE$ with $f(0) = e$
  such that for all $x \in U$, the element $f(x)$ is assigned the germ
  about $x$ of the same function $f$, where as already established we
  identify germs about points up to translation. In short, a
  $\cP$-coalgebra is a set for which ``action by an infinitesimal
  neighborhood about $0$'' is defined at each element.

  More precisely, observe that from any coalgebra $\eE$ we can
  construct a sheaf $\bar \eE$ on $\rr$ by the formula
  $\bar \eE(U) = \PSet(U,\eE)$, where an open interval $U$ has the
  tautological coalgebra structure (sending a point to the germ of the
  identity around that point). The \'etal\'e space of $\bar \eE$
  inherits an action of the \emph{discrete} translation group
  $\rr^\delta$ since the sections over $U$ are invariant under
  translations of $U$ by definition. (This action need not be
  continuous with respect to the usual topology on $\rr$; see
  \cref{cex:noncontinuous}.) Conversely, any \'etal\'e space over $\rr$
  equipped with an action of $\rr^\delta$ such that acting by $y$
  moves the fiber over $x$ to the fiber over $x + y$ gives rise to a
  $\pP$-coalgebra by taking the stalk about $0$. Thus $\PSet$ is
  equivalent to the topos of sheaves for the topological groupoid
  $\rr \times \rr^\delta \rightrightarrows \rr$, where the source map
  is the projection onto $\rr$ and the target map is the sum.


  The same construction evidently carries through more generally for
  any topological group $\tG$ equipped with an open neighborhood basis
  of the identity: take the category $\bB$ whose objects are basic
  open neighborhoods and whose arrows are arbitrary inclusions given
  by translations under the group action of $\tG$ on
  itself,\footnote{One can here decide to pick either the left or
    right action.} and define $\pP \colon \bB \to \Set$ to be the
  forgetful functor. (This $\bB$ is also equivalent to the larger
  category whose objects are images of the basic open neighborhoods of
  identity under translations under the group action, i.e.\ basic open
  sets of the space $\tG$, like $\Int$ above.) By the same calculation
  we again have that elements of $\cP$ are germs about the identity,
  and that coalgebras are ``infinitesimal $\tG$-sets''. Later in
  \cref{ex:topcat} we will further generalize this construction from
  groups to categories equipped with a suitable ``topology''.
  %
\end{example}

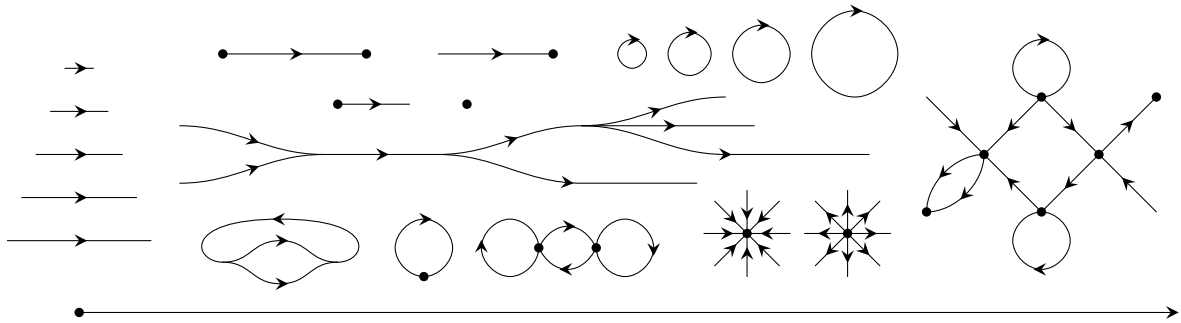
\usetikzlibrary{decorations.markings}
\tikzset{closed/.style={fill,circle,inner sep=1.25}}
\tikzset{flow/.style={draw,line cap=round,postaction={decorate},decoration={markings, mark=at position 0.5 with {\arrow[scale=1.5,xshift=2pt]{stealth}}}}}
\tikzset{ray/.style={postaction={decorate},decoration={markings, mark=at position 1 with {\arrow[scale=1.5,xshift=2pt]{stealth}}}}}

\begin{center}
  \vspace{-15pt}
  \begin{tikzpicture}[scale=.38]
    \draw [flow] (36.5,3.5) .. controls +(0.5,0) and +(0,0.5) .. (37.5,2.5) .. controls +(0,-0.5) and +(0.5,0) .. (36.5,1.5) .. controls +(-0.5,0) and +(0,-0.5) .. (35.5,2.5) .. controls +(0,0.5) and +(-0.5,0) .. cycle;
    \draw [flow] (29.75,2.75) -- (30.875,3.875);
    \draw [flow] (29.75,2.75) -- (31.25,2.75);
    \draw [flow] (29.75,2.75) -- (30.875,1.625);
    \draw [flow] (29.75,2.75) -- (29.75,1.25);
    \draw [flow] (29.75,2.75) -- (28.625,1.625);
    \draw [flow] (29.75,2.75) -- (28.25,2.75);
    \draw [flow] (29.75,2.75) -- (28.625,3.875);
    \draw [flow] (29.75,2.75) -- (29.75,4.25);
    \draw [flow] (21,2.25) .. controls +(0,0.5) and +(-0.5,0) .. (22,3.25) .. controls +(0.5,0) and +(0,0.5) .. (23,2.25) .. controls +(0,-0.5) and +(0.5,0) .. (22,1.25) .. controls +(-0.5,0) and +(0,-0.5) .. cycle;
    \draw [flow] (38.5,5.5) -- (40.5,7.5);
    \draw [flow] (40.5,3.5) -- (38.5,5.5);
    \draw [flow] (15.5,9) -- (19.5,9);
    \draw [flow] (8,9) -- (13,9);
    \draw [flow] (12,7.25) -- (14.5,7.25);
    \draw [flow] (38.5,5.5) -- (36.5,3.5);
    \draw [flow] (36.5,7.5) -- (38.5,5.5);
    \draw [flow] (32.5,7.5) -- (34.5,5.5);
    \draw [flow] (36.5,7.5) -- (34.5,5.5);
    \draw [flow] (36.5,3.5) -- (34.5,5.5);
    \draw [flow] (34.5,5.5) .. controls +(-1,0) and +(0,1) .. (32.5,3.5);
    \draw [flow] (34.5,5.5) .. controls +(0,-1) and +(1,0) .. (32.5,3.5);
    \node [closed] at (34.5,5.5) {};
    \node [closed] at (32.5,3.5) {};
    \draw [flow] (36.5,7.5) .. controls +(-0.5,0) and +(0,-0.5) .. (35.5,8.5) .. controls +(0,0.5) and +(-0.5,0) .. (36.5,9.5) .. controls +(0.5,0) and +(0,0.5) .. (37.5,8.5) .. controls +(0,-0.5) and +(0.5,0) .. cycle;
    \node [closed] at (36.5,7.5) {};
    \node [closed] at (36.5,3.5) {};
    \node [closed] at (16.5,7.25) {};
    \node [closed] at (19.5,9) {};
    \node [closed] at (8,9) {};
    \node [closed] at (13,9) {};
    \node [closed] at (12,7.25) {};
    \draw [flow] (0.5,2.5) -- (5.5,2.5);
    \draw [flow] (1,4) -- (5,4);
    \draw [flow] (1.5,5.5) -- (4.5,5.5);
    \draw [flow] (2,7) -- (4,7);
    \draw [flow] (2.5,8.5) -- (3.5,8.5);
    \draw [flow] (6.5,6.5) .. controls +(2,0) and +(-2,0) .. (11.5,5.5);
    \draw [flow] (6.5,4.5) .. controls +(2,0) and +(-2,0) .. (11.5,5.5);
    \draw [flow] (11.5,5.5) .. controls +(2,0) and +(-2,0) .. (15.5,5.5);
    \draw [flow] (15.5,5.5) .. controls +(2,0) and +(-2,0) .. (20.5,6.5);
    \draw [flow] (20.5,6.5) .. controls +(2,0) and +(-2,0) .. (25.5,7.5);
    \draw [flow] (20.5,6.5) -- (26.5,6.5);
    \draw [flow] (20.5,6.5) .. controls +(2,0) and +(-2,0) .. (25.5,5.5) .. controls +(2,0) and +(0,0) .. (30.5,5.5);
    \draw [flow] (15.5,5.5) .. controls +(2,0) and +(-2,0) .. (20.5,4.5) .. controls +(2,0) and +(0,0) .. (24.5,4.5);
    \node [closed] at (38.5,5.5) {};
    \draw [flow] (26.75,8) .. controls +(-0.5,0) and +(0,-0.5) .. (25.75,9) .. controls +(0,0.5) and +(-0.5,0) .. (26.75,10) .. controls +(0.5,0) and +(0,0.5) .. (27.75,9) .. controls +(0,-0.5) and +(0.5,0) .. cycle;
    \draw [flow] (22.25,8.5) .. controls +(-0.25,0) and +(0,-0.25) .. (21.75,9) .. controls +(0,0.25) and +(-0.25,0) .. (22.25,9.5) .. controls +(0.25,0) and +(0,0.25) .. (22.75,9) .. controls +(0,-0.25) and +(0.25,0) .. cycle;
    \draw [flow] (24.25,8.25) .. controls +(-0.375,0) and +(0,-0.375) .. (23.5,9) .. controls +(0,0.375) and +(-0.375,0) .. (24.25,9.75) .. controls +(0.375,0) and +(0,0.375) .. (25,9) .. controls +(0,-0.375) and +(0.375,0) .. cycle;
    \draw [flow] (15,1.25) .. controls +(-0.5,0) and +(0,-0.5) .. (14,2.25) .. controls +(0,0.5) and +(-0.5,0) .. (15,3.25) .. controls +(0.5,0) and +(0,0.5) .. (16,2.25) .. controls +(0,-0.5) and +(0.5,0) .. cycle;
    \draw [flow] (30,7.5) .. controls +(-0.75,0) and +(0,-0.75) .. (28.5,9) .. controls +(0,0.75) and +(-0.75,0) .. (30,10.5) .. controls +(0.75,0) and +(0,0.75) .. (31.5,9) .. controls +(0,-0.75) and +(0.75,0) .. cycle;
    \draw [flow] (19,2.25) .. controls +(0,-0.5) and +(0.5,0) .. (18,1.25) .. controls +(-0.5,0) and +(0,-0.5) .. (17,2.25) .. controls +(0,0.5) and +(-0.5,0) .. (18,3.25) .. controls +(0.5,0) and +(0,0.5) .. cycle;
    \draw [flow] (19,2.25) .. controls +(0.5,1) and +(-0.5,1) .. (21,2.25);
    \draw [flow] (21,2.25) .. controls +(-0.5,-1) and +(0.5,-1) .. (19,2.25);
    \node [closed] at (19,2.25) {};
    \node [closed] at (21,2.25) {};
    \draw [flow] (26.25,4.25) -- (26.25,2.75);
    \draw [flow] (24.75,2.75) -- (26.25,2.75);
    \draw [flow] (27.75,2.75) -- (26.25,2.75);
    \draw [flow] (26.25,1.25) -- (26.25,2.75);
    \draw [flow] (25.125,3.875) -- (26.25,2.75);
    \draw [flow] (27.375,3.875) -- (26.25,2.75);
    \draw [flow] (27.375,1.625) -- (26.25,2.75);
    \draw [flow] (25.125,1.625) -- (26.25,2.75);
    \node [closed] at (26.25,2.75) {};
    \node [closed] at (15,1.25) {};
    \node [closed] at (29.75,2.75) {};
    \node [closed] at (40.5,7.5) {};
    \draw [flow] (8,1.75) .. controls +(1,0) and +(-1,0) .. (10,2.5) .. controls +(1,0) and +(-1,0) .. (12,1.75);
    \draw [flow] (8,1.75) .. controls +(1,0) and +(-1,0) .. (10,1) .. controls +(1,0) and +(-1,0) .. (12,1.75);
    \draw [flow] (12,1.75) .. controls +(1,0) and +(3.5,0) .. (10,3.25) .. controls +(-3.5,0) and +(-1,0) .. (8,1.75);
    \draw [ray] (3,0) -- (41,0);
    \path (0.25,8.5) -- (0.25,2.5);
    \node [closed] at (3,0) {};
  \end{tikzpicture}
  \vspace{-10pt} \nopagebreak \captionof{figure}{Illustration of a set
    with an ``infinitesimal $\rr$-action'' at each point. The filled
    black dots represent points at which the action is
    constant. Unmarked ends of paths locally resemble an end of an
    open interval. Wherever a path is shown merging into another path
    or point, the gluing is along an open interval. Note that
    everything shown here is particularly simple; as in the case of
    ordinary $\rr$-actions, there are possibilities that are not
    easily visualizable, such as $\rr$ acting on the additive group
    $\rr/\qq$.}
\end{center}

We also have another equivalent characterization of bases.

\begin{proposition}\label{prop:basis}
  Let $\pP \colon \bB \to \bS$ with density comonad $\cP$. Then $\pP$
  is a basis (i.e.\ $\lP$ is dense) if and only if the
  comonadic functor $\carP \colon \PS \to \bS$ is the
  pointwise left Kan extension $\La{\lP}{P}$.
  \[
    \begin{tikzcd}[column sep=15pt, row sep=15pt]
      {\bB} &{}& {\bS} \\
      & {\PS}
      \arrow["\pP", from=1-1, to=1-3]
      \arrow["\lP"', from=1-1, to=2-2]
      \arrow["\carP"', dashed, from=2-2, to=1-3]
      \arrow["\Downarrow"{description, pos=0.3}, draw=none, from=1-2, to=2-2]
    \end{tikzcd}
  \]
\end{proposition}

For small $\bB$ and cocomplete $\bS$, this means the following
triangle commutes (up to isomorphism):
\[\begin{tikzcd}[column sep=15pt, row sep=15pt]
    {\Ps{\bB}} && \bS \\
    & \PS
    \arrow["{\extend{\pP}}", from=1-1, to=1-3]
    \arrow["\bsheaf{\pP}", from=2-2, to=1-1]
    \arrow["\carP"', from=2-2, to=1-3]
  \end{tikzcd}
  \qquad
\]
More intuitively: taking a coalgebra's underlying ``sheaf'', given by
$\bsheaf{\pP}$, then taking the ``sum of stalks'', given by
$\extend{\pP}$, yields the ``carrier of the \'etal\'e space'', given
by $\carP$.

\begin{proof}
  Recall that in general a functor is dense if and only if its density
  comonad is identity and is pointwise. We have the
  left-to-right direction because left adjoints preserve (pointwise)
  left Kan extensions:
  \[\La{\lP}{\pP} = \carP \circ \La{\lP}{\lP} =
    \carP \circ \gen{\lP} = \carP.\] Conversely, given the pointwise
  left Kan extension $\La{\lP}{\pP} = \carP$, we have the pointwise
  left Kan extension $\La{\lP}{\lP} = \id_{\PS}$ since $\carP$
  creates (in particular reflects) colimits.
\end{proof}

\begin{remark}\label{rem:densityformula}
  In the case $\pP \colon \bB \to \bS$ is a basis, we may describe the
  correspondence from \cref{prop:density} between the two types of
  commutative triangles shown below left (i.e.\ comonad maps from
  $\cP$ to a comonad $\cC$) and below right more explicitly.
  \[
    \begin{tikzcd}[row sep=10pt,column sep=15pt]
      {\PS} && {\CS} \\
      & {\bS}
      \arrow[dashed, from=1-1, to=1-3]
      \arrow["\carP"', from=1-1, to=2-2]
      \arrow["\carC", from=1-3, to=2-2]
    \end{tikzcd}
    \qquad\leftrightarrows\qquad
    \begin{tikzcd}[row sep=10pt,column sep=15pt]
      {\bB} && {\CS} \\
      & {\bS}
      \arrow[dashed, from=1-1, to=1-3]
      \arrow["\pP"', from=1-1, to=2-2]
      \arrow["\carC", from=1-3, to=2-2]
    \end{tikzcd}\]

  Regardless of whether $\pP$ is a basis, a commutative triangle shown
  left induces a commutative triangle shown right by composing with
  $\lP \colon \bB \to \PS$, and this correspondence is a natural
  bijection. However, assuming $\pP$ is a basis, we now also have an
  explicit construction of the inverse bijection. Namely, given
  $F \colon \bB \to \CS$ we take the pointwise left Kan extension
  $\La{\lP}{F} \colon \PS \to \CS$. Indeed, this gives a commutative
  triangle (i.e. comonad map):
  \[\carP = \La{\lP}{\pP} = \La{\lP}{(\carC \circ F)} = \carC \circ
    \La{\lP}{F} \] because comonadic functors create
  colimits. Moreover, this translation process
  $F \mapsto \La{\lP}{F}$ is indeed the inverse of the natural
  bijection $\com{\bS}{\phi} \mapsto (\com{\bS}{\phi} \circ \lP)$:
  \[\La{\lP}{(\com{\bS}{\phi} \circ \lP)} = \com{\bS}{\phi} \circ 
    \La{\lP}{\lP} = \com{\bS}{\phi}.\]

  These inverse translation processes, given by pasting triangles, are
  shown below, assuming $\bB$ is small and $\bS$ is cocomplete:
  \[\begin{tikzcd}
      \bB & \PS & \CS \\
      & \bS
      \arrow["{\pP}"', from=1-1, to=2-2]
      \arrow["{\lP}", from=1-1, to=1-2]
      \arrow["{\carP}", pos=.4, from=1-2, to=2-2]
      \arrow["{\com{\bS}{\phi}}", from=1-2, to=1-3]
      \arrow["{\carC}", from=1-3, to=2-2]
    \end{tikzcd}
    \qquad\quad
    \begin{tikzcd}
      \PS & \Ps{\bB} & \CS \\
      & \bS
      \arrow["{\carP}"', from=1-1, to=2-2]
      \arrow["{\bsheaf{\pP}}", from=1-1, to=1-2]
      \arrow["{\extend{\pP}}", pos=.4, from=1-2, to=2-2]
      \arrow["{\extend{F}}", from=1-2, to=1-3]
      \arrow["{\carC}", from=1-3, to=2-2]
    \end{tikzcd}\]
\end{remark}

Every basis of a comonad on $\Set$ also induces a topological basis:

\begin{proposition}
  If $\pP \colon \bB\to\Set$ is a basis, then the images of the maps
  to the colimit $\pP(\xU) \to \pts{\pP}$ determine a basis of the
  underlying topological space.
\end{proposition}
\begin{proof}
  Assuming $\pP$ is a basis, every coalgebra is a colimit of basic
  coalgebras (those in the image of $\lP \colon \bB \to \PSet$). In
  particular, there is a jointly epimorphic family of maps from basic
  coalgebras to any regular subobject (open subset) $\xU$ of the
  terminal coalgebra. Since by \cref{lem:regularmono} $\PSet$ is
  coregular, we may take image factorizations of these maps, whereby
  $\xU$ is covered by regular subobjects of the terminal coalgebra that
  arise as images of basic coalgebras.
\end{proof}

\begin{remark}
  A comonad on $\Set$ is small (i.e.\ accessible) if and only if it
  admits a small basis. Indeed, by \cref{lem:accessibility}, any
  category of coalgebras of an accessible comonad on an accessible
  category is accessible, and therefore has a small dense
  subcategory. Conversely, pointwise density comonads of small
  diagrams in accessible categories are always
  accessible~\cite[Proposition 2.4.3]{makkai-pare}. More generally by
  the same reasoning, we have that a comonad on any accessible
  category is accessible if and only if it admits a small basis.
\end{remark}

\section{Subbases}\label{sec:subbases}

A basis for a comonad $\cC$ on $\bS$ amounts to a dense functor into
$\CS$. This prompts the question of which functors into $\CS$
determine subbases --- for which functors $F \colon \bB \to \CS$ do we
recover $\cC$ as the density comonad $\cP$ of
$\pP \coloneqq \carC \circ F$?  Or more specifically, which functors
are of the form $\lP$? We will find the answer (in the pointwise case)
is precisely those which are ``dense relative to cofree coalgebras'',
in the sense that the restriction of the nerve functor
$\nerve{\lP} \colon \CS \to \Ps{\bB}$ to the Kleisli category $\KlCS$
\[\KlCS \hookrightarrow \CS
  \xto{\nerve{\lP}} \Ps{\bB}\] is fully faithful.

\begin{definition}\label{def:subbasis}
  Let $\cC$ be a comonad on $\bS$. A functor $F \colon \bB \to \CS$ is
  \emph{subdense} (with respect to $\carC \colon \CS \to \bS$) if
  $\cC$ is the pointwise density comonad $\cP$ of
  $\pP \coloneqq \carC \circ F$ and $F = \lP$ as in \cref{def:lift}.
  
  Explicitly, the condition $F = \lP$ means $F$ forms a universal
  triangle
  \[
    \begin{tikzcd}[row sep=10pt,column sep=15pt]
      {\bB} && {\CS} 
      \\
        & {\bS}
        \arrow["F", dashed, from=1-1, to=1-3]
        \arrow["\pP"', from=1-1, to=2-2]
        \arrow["\carC", from=1-3, to=2-2]
      \end{tikzcd}
    \]
    exhibiting $\carC$ as the free comonadic functor extending $\pP$
    as in \cref{prop:density}.
\end{definition}

\begin{remark}\label{rem:subbasismodule}
  By \cref{rem:lift}, the condition $F = \lP$ is also equivalent to
  saying that the
  natural transformation induced by the coalgebra structures in the
  image of $F$
  \[
    \begin{tikzcd}[column sep=15pt, row sep=15pt]
      {\bB} &{}& {\bS} \\
      & {\bS}
      \arrow["\pP", from=1-1, to=1-3]
      \arrow["\pP"', from=1-1, to=2-2]
      \arrow["\cC"', from=2-2, to=1-3]
      \arrow["\Downarrow"{description, pos=0.3}, draw=none, from=1-2, to=2-2]
    \end{tikzcd}
  \]
  witnesses $\cC$ as the density comonad of
  $\pP \coloneqq \carC \circ F$. More explicitly, this natural
  transformation may be written as $\carC \circ \eta \circ F$ where
  $\eta$ is the unit of the adjunction $\carC \dashv \rcarC$.

  Moreover, since we are assuming the density comonad $\cP$ is
  pointwise, this is equivalent to saying that the coalgebra
  structures are given by the colimit inclusion maps
  $\pP(\xU) \to \colim_{\pP(\xV) \to \pP(\xU)} \pP(\xV)$ indexed by
  the identities.
\end{remark}


\begin{lemma}\label{lem:cofreecolimit}
  Let $\pP \colon \bB \to \bS$ with pointwise density comonad
  $\cP$. The cofree coalgebra $\cP(\xX)$ on an object $\xX$ in $\bS$ is
  given by the canonical colimit of
  \[(\pP / \xX) \cong (\lP / \cP(\xX)) \to \bB \xto{\lP} \PS\] i.e.\
  the $\xX^\pP \cong \cP(\xX)^{\lP}$-weighted colimit of $\lP$.
\end{lemma}

\noindent
In short, $\lP$ is ``dense relative to'' cofree coalgebras.

\begin{proof}
  First note that we indeed have an isomorphism of $\bB$-presheaves
  $\xX^\pP \cong \cP(\xX)^{\lP}$, since maps $\pP(\xU) \to \xX$ in
  $\bS$ correspond to maps $\lP(\xU) \to \cP(\xX)$ in $\PS$ under the
  adjunction $\carP \dashv \rcarP$; thus we indeed also have an
  isomorphism between the categories of elements
  $(\pP / \xX) \cong (\lP / \cP(\xX))$. The functor
  $\carP \colon \PS \to \bS$ creates colimits, and the canonical
  cocone is colimiting in $\bS$ by definition of pointwise density
  comonad.
\end{proof}


\begin{lemma}\label{lem:littledense}
  Let $\pP \colon \bB \to \bS$. If an object $\xX$ in $\bS$ is given
  by the canonical colimit of the diagram
  \[(\pP / \xX) \to \bB \xto{\pP} \bS\] i.e.\ the
  $\xX^\pP$-weighted colimit of $\pP$, then maps from any weighted
  colimit of $\pP$ into $\xX$ are identified with maps between weights
  in $\Ps{\bB}$.
\end{lemma}
\begin{proof}
  %
  The universal property of the $F$-weighted colimit of $\pP$ is that
  maps from it to $\xX$ are identified with natural transformations
  from $F$ to $\xX^\pP \coloneqq \Ho{\bS}(\pP(\dash), \xX)$.
\end{proof}

Combining \cref{lem:cofreecolimit} and \cref{lem:littledense}, we see
that for each $\xX$ in $\bS$, the cofree coalgebra $\cP(\xX)$ may be
realized as the presheaf 
$\xX^\pP$, in the sense that coalgebra maps $\cP(\xX) \to \cP(\xY)$
are equivalently maps of presheaves
$\xX^\pP \Rightarrow \xY^\pP$. Intuitively, this means that the cofree
coalgebra on $\xX$ is realized as a ``sheaf of functions into $\xX$'',
as is in particular so in the case of comonads on $\Set$ corresponding
to topological spaces.

\anoteNI{I have a feeling this would be clarified by the theory of
  relative (co)monads. Is there a notion of relative density comonad?}

\begin{theorem}\label{prop:subbasis}
  Let $\cC$ be a comonad on $\bS$ and let
  $F \colon \bB \to \CS$. The following are equivalent.
  \begin{enumerate}[label=(\roman*)]
  \item\label{item:subbasis} $F$ is subdense (with respect to $\carC$).
  \item\label{item:subdensity} Each cofree coalgebra $\eE$ is the
    canonical colimit of the diagram
    \[(F / \eE) \to \bB \xto{F} \CS.\]
  \item\label{item:subdense} The restriction of the nerve functor
    $\nerve{F} \colon \CS \to \Ps{\bB}$ to the cofree coalgebras
    \[\KlCS \hookrightarrow \CS
      \xto{\nerve{F}} \Ps{\bB}\] is fully faithful.
  \end{enumerate} 
\end{theorem}

In particular when $\bB$ is small and $\bS$ is cocomplete (thus so is
$\CS$), the pointwise density comonad of $F$ exists, in
which case \labelcref{item:subdensity} equivalently says that this
density comonad $\gen{F}$ is identity when restricted to cofree
coalgebras, i.e.\ the counit induces an isomorphism
$\gen{F} \circ \rcarr \cong \rcarr \colon \bS \to \CS$.

\begin{proof}
  \labelcref{item:subbasis} $\imp$ \labelcref{item:subdensity}
  by \cref{lem:cofreecolimit}.
  \labelcref{item:subdensity} $\imp$ \labelcref{item:subdense}
  by \cref{lem:littledense}.

  For \labelcref{item:subdense} $\imp$ \labelcref{item:subbasis},
  suppose the restriction of the nerve $\nerve{F}$ to the Kleisli
  category is fully faithful. That is, maps between cofree coalgebras
  $\cC(\xX)$ in $\CS$ are identified with natural transformations
  between the $\bB$-presheaves 
  $\cC(\xX)^{F}$, which in turn are isomorphic to the $\bB$-presheaves
  $\xX^\pP$ where $\pP \coloneqq \carP \circ F$. Such maps
  $\cC(\xX) \to \cC(\xY)$ in $\CS$ are also identified with maps
  $\cC(\xX) \to \xY$ in $\bS$. Hence $\cC(\xX)$ is the
  $\xX^\pP$-weighted colimit of $\pP$ by definition. This is to say
  that $\cC$ is the pointwise density comonad $\cP$, at least at the
  level of underlying endofunctors.
  One verifies that the comonad structure is likewise as expected
  for the pointwise density comonad, with the counit and
  comultiplication given by evident maps from these colimits.

  Finally, it remains to check that $F = \lP$, meaning the coalgebra
  structure, i.e.\ component of the unit in $\CS$, at each object
  $F(\xU)$ is given by the colimit inclusion map
  $\pP(\xU) \to \cC(\pP(\xU)) \cong \colim_{\pP(\xV) \to \pP(\xU)}
  \pP(\xV)$ indexed by the identity $\pP(\xU) \to \pP(\xU)$ as in
  \cref{rem:subbasismodule}. Indeed, taking mates of arrows, this
  colimit formula in $\bS$ translates to the canonical colimit formula
  \labelcref{item:subdensity} in $\CS$ for the cofree coalgebra
  $\cC(\pP(\xU))$; in particular the identity $\pP(\xU) \to \pP(\xU)$
  in $\bS$ is mate to the desired unit map $F(\xU) \to \cC(\pP(\xU))$
  in $\CS$, and any map into a canonical colimit is the colimit
  inclusion map indexed by itself in the canonical colimit formula.
\end{proof}


\begin{remark}
  It follows from \cref{prop:subbasis} that if $\cC$ is a comonad on
  $\bS$ all of whose coalgebras are cofree, then every subdense
  functor into $\CS$ is dense, or in other words, all pointwise
  subbases of $\cC$ are bases. For instance, this property holds for
  idempotent comonads, which correspond to coreflective subcategories,
  and therefore a functor into a coreflective subcategory of $\bS$ is
  subdense if and only if it is dense. In particular, all comonads on
  preorders are idempotent, and therefore any functor into a preorder
  admitting a pointwise density comonad is a basis.
\end{remark}

\begin{example}[Infinitesimal group actions]\label{ex:infinitesimalsubbasis}
  A topological subbasis of $\rr$ that is not a basis is given by the
  infinite open intervals:
  \[
    \setof{(-\infty,x) \mid x \in \rr}
    \cup
    \setof{(x,\infty) \mid x \in \rr}
    \cup
    \setof{(-\infty,\infty)}.
  \]
  These subsets form a poset under inclusion, and by
  \cref{prop:subbases} (and \cref{rem:subbasisgerms}) the density
  comonad of the forgetful functor into $\Set$ is the germs comonad
  $\toptocmd{\rr}$ of $\rr$. By \cref{lem:topbasis}, it therefore
  fails to be a basis in the abstract sense of \cref{def:basis};
  indeed, not every sheaf, or \'etal\'e space, over $\rr$ is a
  colimit of infinite open intervals. In other words, the inclusion of
  infinite open intervals in sheaves on $\rr$ is subdense but not
  dense.

  In \cref{ex:real_line} we tweaked the usual basis of $\rr$ to yield
  a basis of the comonad whose coalgebras are ``infinitesimal
  $\rr$-sets''. And we can tweak this subbasis in just the same way to
  obtain a subbasis of the same comonad. Namely, let $\bB$ be the the
  category whose objects are the same infinite open intervals as
  above, and whose arrows are orientation-and-distance-preserving
  maps, i.e.\ translations, and let
  $\pP \colon \bB \hookrightarrow \Set$ be the forgetful functor.
  
  Note that the poset of subbasic neighborhoods of $0$ is connected.
  (This is why above we included $\rr = (-\infty,\infty)$ itself in
  the subbasis of $\rr$.) As noted in \cref{rem:subbasisgerms}, this
  guarantees that an $\xA$-valued germ about $0$ calculated with
  respect to our chosen neighborhood subbasis is the same as a germ
  calculated with respect to the entire system of neighborhoods. Note
  also that, as in \cref{ex:real_line}, any point in any subbasic open
  set is the image under a translation of an identity. (Together these
  properties ensure the category of elements $\El(\pP)$ is connected,
  i.e.\ $\pts{\pP}$ is $1$.) The calculation of the density comonad
  via germs is then exactly the same as in \cref{ex:real_line}, and we
  obtain the same comonad, whose coalgebras are sets for which
  ``action by an infinitesimal neighborhood about $0$'' is defined at
  each element. In particular, the canonical inclusion of subbasic
  open sets into this category of coalgebras is subdense.

  Here we can see the correspondence of \cref{prop:subbasis}: the
  cofree coalgebra on $\xA$ is given by the set of $\xA$-valued germs,
  and in order to verify that the density comonad is correct we
  exhibited the germs as the expected canonical colimit.

  Like \cref{ex:real_line}, the construction can be generalized to
  arbitrary topological groups equipped with a neighborhood subbasis
  of identity satisfying suitable connectedness properties. And we
  generalize this further to categories in \cref{ex:topcat}.
\end{example}

\section{Concrete bases}\label{sec:coconfluence}

We next establish a sufficient condition for a $\Set$-valued functor
to be a basis, which resembles the more familiar definition of basis
from topology.

A diagram $\pP \colon \bB \to \Set$ of subsets of $\pts{\pP}$
constitutes a basis for a topology on $\pts{\pP}$ if and only if for
every $x \in \pts{\pP}$ and $\xU_1, \xU_2 \in \bB$ with
$x \in \xU_1 \cap \xU_2$, there exists $\xU \in \bB$ with $x \in \xU$
and $\xU \subseteq \xU_1 \cap \xU_2$. This is equivalent to saying that
the category of elements $\El(\pP)$ satisfies the following condition.

\begin{definition}\label{def:influent}
  A category is \emph{co-confluent} (a.k.a. satisfies the \emph{right
    Ore condition}) if every cospan extends to a commutative square:
  
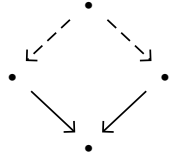
\[\begin{tikzcd}[column sep=15pt, row sep=15pt]
      & \smbullet \\
      \smbullet && \smbullet \\
      & \smbullet
      \arrow[dashed, from=1-2, to=2-1]
      \arrow[dashed, from=1-2, to=2-3]
      \arrow[from=2-1, to=3-2]
      \arrow[from=2-3, to=3-2]
    \end{tikzcd}\]

  We call a functor $\pP \colon \bB \to \Set$ a \emph{concrete basis}
  if its density comonad exists and its category of elements
  $\El(\pP)$ is co-confluent.
\end{definition}

Observe that co-confluence implies, by an inductive argument, that
every zig-zag of arrows between two objects extends to a span. Note
also that given $\pP \colon \bB \to \Set$, two elements
$x_1, x_2 \in \El(\pP)$ lie over the same point in the colimit
$\pts{\pP}$ just when they are connected by some zig-zag. Thus
co-confluence of $\El(\pP)$ implies that for all $x_1 \in \pP(\xU_1)$
and $x_2 \in \pP(\xU_2)$ lying over the same $x \in \pts{\pP}$,
there is a span $\xU_1 \ot \xU' \to \xU_2$ and $x' \in \pP(\xU')$
mapped to respectively $x_1$ and $x_2$. This agrees with the above
usual notion of basis for topological spaces.

We will show that every concrete basis is a categorical basis in the
sense of \cref{def:basis}. Not every comonad on $\Set$ admits such a
well-behaved basis. Nevertheless, those that do form a particular
natural class of comonads on $\Set$, more general than the ionads, as
mentioned in the introduction. These are the comonads on $\Set$ that
preserve weak pullbacks.

\begin{lemma}\label{lem:weakel}
  Let $\pP \colon \bB \to \Set$ and let $\bJ$ be a category. If
  $\El(\pP)$ has a cone to every diagram $\bJ \to \El(\pP)$
  lifting $D \colon \bJ \to \bB$, then $\pP$ preserves weak limits
  of $D$.

  Conversely, if $D \colon \bJ \to \bB$ is a functor with a weak
  limit preserved by $\pP$, then $\El(\pP)$ has a cone to every
  diagram lifting $D$.
\end{lemma}
\begin{proof}
  Let $L$ be a weak limit of $D \colon \bJ \to \bB$.  For the
  first statement, assume $\El(\pP)$ has a cone to every diagram
  $\bJ \to \El(\pP)$ lifting $D \colon \bJ \to \bB$. We must
  show that $\pP(L)$ is a weak limit of $\pP \circ D$, or equivalently
  that the induced map $\pP(L) \to \lim \pP \circ D$ is
  surjective. Giving an element $v \in \lim (\pP \circ D)$ is
  equivalent to giving a diagram $D' \colon \bJ \to \El(\pP)$
  lifting $D$. By assumption we then obtain a cone from an element
  $u \in \pP(\xU)$ to $D'$ in $\El(\pP)$, yielding a cone from $\xU$
  to $D$ in $\bB$, which in turn factors through the weak limit via
  $f \colon U \to L$. This yields an element $\pP(f)(u) \in \pP(L)$
  lying over $v$ in $\lim \pP \circ D$.

  Conversely, assume $\pP$ preserves $L$. Let
  $D' \colon \bJ \to \El(\pP)$, equivalent to a functor
  $D \colon \bJ \to \bB$ and an element of $\lim \pP \circ D$. By
  assumption the comparison map $\pP(L) \to \lim \pP \circ D$ is
  surjective, so there exists $u \in \pP(L)$ forming a cone to $D'$ in
  $\El(\pP)$.
\end{proof}

In particular, if $\bB$ has weak limits of shape $\bJ$, then
$\pP \colon \bB \to \Set$ preserves weak limits of shape $\bJ$ if
and only if $\El(\pP)$ has a cone to every diagram of shape
$\bJ$. By \cref{lem:complete}, it follows that a comonadic functor
into $\Set$ preserves weak pullbacks if and only if it is a concrete
basis.

\begin{lemma}\label{lem:tensorbypullback}
  Let $\pP \colon \bB \to \Set$ and $W \colon \bB\op \to \Set$. The
  tensor product of functors (a.k.a.\ profunctor composition)
  $\pP \otimes W$ is the set of connected components of the pullback
  of the projections $\El(\pP) \to \bB \ot \El(W)$.
\end{lemma}
\begin{proof}

  The tensor product of functors $\pP \otimes W$ is the $W$-weighted
  colimit of $\pP$ (or equivalently vice versa), which can be
  expressed as the colimit of $\El(W) \to \bB \xto{\pP} \Set$ (or
  equivalently the colimit of $\El(\pP)\op \to \bB\op \xto{W}
  \Set$). Reindexing one discrete (op)fibration over $\bB$ along the
  other, we obtain that the category of elements of this composite is
  the desired pullback.
\end{proof}

\begin{proposition}\label{prop:influent}
  Let $\pP \colon \bB \to \Set$ with $\bB$ small. Then it is a
  concrete basis, i.e.\ its category of elements $\El(\pP)$ is
  co-confluent, if and only if its Yoneda extension
  $\extend{\pP} \colon \Ps{\bB} \to \Set$ preserves weak pullbacks.
\end{proposition}
\begin{proof}
  By \cref{lem:weakel}, it suffices to show that $\El(\pP)$ is
  co-confluent if and only if $\El(\extend{\pP})$ is co-confluent.
  Suppose $\El(\extend{\pP})$ is co-confluent. Then every cospan
  within $\El(\pP) \hookrightarrow \El(\extend{\pP})$ admits a cone
  (i.e.\ square) from an element of $\extend{\pP}(W) = \pP \otimes W$
  for some presheaf $W$. Every element of $\pP \otimes W$ originates
  from some $\xU$ in $\bB$ and pair of elements in $x \in \pP(\xU)$ and
  $u \in W(\xU)$, so the given cone from $W$ yields a cone from $x$ by
  pre-composing the map corresponding to $u$ from (the representable
  presheaf of) $\xU$ into $W$.

  Conversely, suppose $\El(\pP)$ is co-confluent, and suppose there is
  a cospan between $a_1 \in \extend{\pP}(W_1)$ and
  $a_2 \in \extend{\pP}(W_2)$ in $\El(\extend{\pP})$, lying over a
  cospan $W_1 \xto{f} W_0 \xot{g} W_2$ in $\Ps{\bB}$. Here the element
  $a_1$ originates from some $x_1 \in \pP(\xU_1)$ and
  $u_1 \in W_1(\xU_1)$ and the element $a_2$ originates from some
  $x_2 \in \pP(\xU_2)$ and $u_2 \in W_2(\xU_2)$. By the assumption
  that $f$ and $g$ underlie a cospan between $a_1$ and $a_2$, the
  element of $\extend{\pP}(W_0)$ originating from $x_1$ and
  $f_{\xU_1}(u_1)$ and the element originating from $x_2$ and
  $g_{\xU_2}(u_2)$ are identified. By \cref{lem:tensorbypullback} this
  means there is a sequence of intermediate pairs related by a zig-zag
  of arrows in $\bB$ corresponding to arrows in both $\El(\pP)$ and
  $\El(W_0)$. Co-confluence of $\El(\pP)$ implies the zig-zag between
  $x_1$ and $x_2$ admits a cone from some element $x \in \pP(\xU)$ for
  some $\xU \in \bB$. We hence obtain a cone in $\Ps{\bB}$ from (the
  representable presheaf of) $\xU$ to the zig-zag between (the
  representable presheaves of) $\xU_1$ and $\xU_2$, and we have a
  canonical cocone from this zig-zag to $W_0$ as well by its
  construction. The outer legs of this cone and cocone form a square,
  whose cospan legs factor through $f$ and $g$. The induced span
  $W_1 \ot \xU \to W_2$ underlies a cone in $\El(\extend{\pP})$ from
  $x$ to the given cospan between $a_1$ and $a_2$.
\end{proof}

\begin{corollary}\label{cor:influentwp}
  Let $\pP \colon \bB \to \Set$ be a concrete basis. Then $\cP$
  preserves weak pullbacks.
\end{corollary}
\begin{proof}
  For $\bB$ small, the density comonad is given by
  $\extend{\pP} \circ \nerve{\pP}$, as in \cref{prop:densityadj}.

  It remains to show that the statement holds for large $\bB$ as
  well. If the density comonad of a functor from a large category into
  $\Set$ exists, then it is pointwise by \cref{lem:setpointwise} and
  so agrees with the density comonad calculated in any larger universe
  $\SET$. (The inclusion $\Set \hookrightarrow \SET$ preserves all,
  even large, limits and colimits that happen to exist in $\Set$, as
  can be seen by probing with maps from $1$ and maps into $2$.) Now if
  the category of elements of a functor into $\SET$ is co-confluent,
  its density comonad on $\SET$ preserves weak pullbacks by the same
  reasoning within the larger universe, in which case the original
  density comonad on $\Set$ also preserves weak pullbacks.
\end{proof}

\begin{corollary}\label{cor:influentbasis}
  Let $\pP \colon \bB \to \Set$ with density comonad $\cP$. If $\pP$
  is a concrete basis, then $\pP$ is a basis (in the sense of
  \cref{def:basis}).
\end{corollary}
\begin{proof}
  For $\bB$ small, by \cref{lem:weaktaut} and \cref{prop:influent}, $\extend{\pP}$ is
  taut, i.e.\ preserves pullbacks of monomorphisms along arbitrary
  maps. In particular $\extend{\pP}$ preserves finite intersections,
  so the result follows from \cref{lem:extcoref}.

  It remains to show that the statement holds for large $\bB$ as
  well. But should we find $\pP$ is a basis for the density comonad
  calculated in any larger universe $\SET$, then $\pP$ is already a
  basis in the smaller universe $\Set$: if the category of coalgebras
  in the larger universe $\CoalgOn{\cP}{\SET}$ embeds fully faithfully
  in (large) presheaves on $\bB$, then so does the full subcategory of
  coalgebras in the smaller universe $\CoalgOn{\cP}{\Set}$.
\end{proof}

\anoteNI{Does an analogue work in the generality of $\pP \colon \bB \to \bS$,
  using general ``weak-pullback-flatness''?}

\begin{theorem}\label{cor:wpb}
  A comonad on $\Set$ preserves weak pullbacks if and only if it
  admits a (possibly large) concrete basis.
\end{theorem}
\begin{proof}
  If $\cC$ preserves weak pullbacks, then so does
  $\carC \colon \CSet \to \Set$ by~\cref{lem:weakpreserve}, which then
  itself is a concrete basis. The right-to-left direction is given by
  \cref{cor:influentwp} (and \cref{cor:influentbasis}).
\end{proof}

\begin{example}
  The diagrams from \cref{ex:part}, \cref{ex:topbasis},
  \cref{ex:catbasis}, and \cref{ex:topgroup} are all concrete bases,
  so all are bases of weak-pullback-preserving comonads.
  
  On the other hand, the comonadic forgetful functor $\Shear \to \Set$
  from \cref{ex:shear} does not preserve weak pullbacks. Indeed, its
  category of elements is not co-confluent: the elements given by the
  defining diagram $2 \to 1$ constitute a cospan extending to no
  square. The comonad is therefore not the density comonad of any
  concrete basis.
\end{example}

The following recovers the comonad whose coalgebras are undirected
multigraphs from \cref{ex:undgraphs}.

\begin{example}[Symmetrized categories]\label{ex:symcat}
  Let $\bC$ be a small category, and let $\bG$ be a subgroup of the
  group of automorphisms of $\bC$ that fix the objects
  $\obs{\bC}$.\footnote{Here the assumption that $\bG$ fixes the
    objects of $\bC$ is not necessary, except to ensure that the
    points of the comonad we generate are given by the objects
    $\pts{\cC}$ of $\bC$, as well as to simplify the
    presentation.}  (For instance, to recover
  \cref{ex:undgraphs}, take $\bC$ to be the walking parallel pair of
  arrows $\rightrightarrows$, and take $\bG$ to be the two-element
  group of automorphisms that fix objects, with the generating
  involution swapping the two arrows.) Let $\groth \bG$ be the
  Grothendieck construction of the inclusion
  $\bG \hookrightarrow \Cat$ sending the unique object of $\bG$ to
  $\bC$.
  
  Explicitly, the objects of $\groth \bG$ are the same as those of $\bC$, and
  an arrow $c \to d$ is a pair $(g, x)$ of an element $g$ in $\bG$ and
  an arrow $x \colon c \to d$ in $\bC$, with composition given by
  \[(g', x') \circ (g, x) = (g' \circ g, x' \circ
  g'(x)).\] Equivalently, it is the category obtained from $\bC$ by
  freely adjoining a copy of $\bG$ at the automorphisms of each object
  and imposing the equations $g \circ x = g(x) \circ g$ for all $g$ in
  $\bG$ and $x \colon c \to d$ in $\bC$.

  \begin{center}
    \begin{tikzcd}
      && e \\
      & d & d \\
      c & c & c
      \arrow["{(g', x')}", from=2-2, to=1-3]
      \arrow["{g'}"', from=2-2, to=2-3]
      \arrow["{x'}"', from=2-3, to=1-3]
      \arrow["{(g, x)}", from=3-1, to=2-2]
      \arrow["g"', from=3-1, to=3-2]
      \arrow["x"', from=3-2, to=2-2]
      \arrow["{g'}"', dashed, from=3-2, to=3-3]
      \arrow["{g'(x)}"', dashed, from=3-3, to=2-3]
    \end{tikzcd}
    \nopagebreak \captionof{figure}{Composition in $\groth \bG$. Here we
      are using $g$ and $x$ as short for $(g, \id)$ and $(\id, x)$.}
  \end{center}

  In \cref{ex:catbasis}, we defined the ``canonical basis'' of the
  comonad $\cbC$ corresponding to a category $\bC$ by letting
  $\bB \coloneqq \bC\op$ and letting
  $\catbasis{\bC} \colon \bB \to \Set$ be the functor sending each
  object $c$ to the set $\obs{c / \bC}$ of arrows out of $c$ and
  sending each arrow $x \colon c \to d$ to pre-composition by $x$ in
  $\bC$. Now similarly let $\bB \coloneqq (\groth \bG)\op$, let
  $\pP \colon \bB \to \Set$ be the functor sending each object $c$ to
  the set $\obs{c / \bC}$ of arrows out of $c$, sending each arrow of
  the form $x = (\id, x) \colon c \to d$ to pre-composition by $x$ in
  $\bC$ (just like for the canonical basis of a category), and sending
  each arrow of the form $g = (g, \id)$ to the action of $g^{-1}$ on
  arrows of $\bC$. In short, this is just like the canonical basis of
  a category, except we have added in the actions by elements of the
  group $\bG$.

  This $\pP$ is a concrete basis. Indeed, given any two arrows
  $k \colon c \to d$ and $k' \colon c \to d$ in $\bC$ (i.e.\ objects of
  $\El(\pP)$) and two arrows $(g, x)$ and $(g', x')$ in $\bB$
  satisfying $g^{-1}(k \circ x) = {g'}^{-1}(k' \circ x')$ (i.e.\ a
  cospan in $\El(\pP)$), pick any $h$ and $h'$ in $\bG$ satisfying
  $h\circ g = h' \circ g'$. We find that $\id_d$ lies over $k$ and
  $k'$ via $(h, h(k))$ and $(h', h'(k'))$, forming a commutative
  square in $\El(\pP)$ as desired.

  Let us calculate the density comonad $\cP$. The set of points
  $\pts{\pP}$ of $\cP$ is the set $\obs{\bC}$ of objects of $\bC$. We
  find that $\cP(\xA)$ is the set of functions $\obs{c / \bC} \to \xA$
  modulo reindexing by the action of $\bG$ on $\obs{c / \bC}$:
  \[\cP(\xA) = \sum_{c\in \obs{\bC}} \colim_{g \in \bG} \xA^{\obs{c / \bC}}.\]
  Indeed, in terms of germs as in \cref{rem:densitygerms}, the
  elements of $\cP(\xA)$ are germs $[f]_x$ of functions
  $f \colon {c / \bC} \to \xA$ about arrows $x \colon c \to d$ in
  $\bC$, where we set $[f]_x = [f']_{x'}$ if there exists
  $y \colon c' \to c$ such that $x \circ y = x'$ and
  $f = f'(\dash \circ y)$, \emph{or} if there exists $g$ in $\bG$ such
  that $g^{-1}(x) = x'$ and $f = f' \circ g^{-1}$. By the first kind
  of equation (which also appears when calculating the simpler density
  comonad $\cbC$ corresponding to the category $\bC$), every germ
  $[f]_x$ reduces to one about identity $[f(\dash \circ x)]_{\id_d}$,
  and the second kind of equation (corresponding to the group of
  automorphisms $\bG$) identifies germs up to reindexing by the action
  of $\bG$ on arrows.
  
  As for coalgebras, since $\pP$ is a basis we are permitted to use
  the description from \cref{rem:simplergerms}. A $\cP$-coalgebra
  $\eE \to \cP(\eE)$ assigns to each element $e$ the germ of some
  function $f \colon \obs{c / \bC} \to \eE$ with $f(\id_c) = e$ and
  such that for all $x \colon c \to d$, the element $f(x)$ in $\eE$ is
  assigned the germ about $x$ of the same function $f$ --- which
  reduces to a germ about $\id_d$ of $f(\dash \circ x)$. That is, we
  have at each element $e$ an assigned object $c$ of $\bC$ and actions
  by arrows out of $c$ \emph{defined up to reindexing by the group
    $\bG$}, satisfying unit and associativity laws. To put it another
  way, a $\cP$-coalgebra is similar to a functor $\bC \to \Set$, but
  where at any element the collective actions by arrows in $\bC$ are
  only defined up to the group action on $\bC$.
  
  In particular, if we take $\bC$ to be $\rightrightarrows$, and we
  may take $\bG$ generated by the involution swapping the two arrows,
  we obtain the category of undirected multigraphs $\MultiGph$ from
  the introduction and \cref{ex:undgraphs}. By \cref{prop:regular},
  the category of coalgebras is regular, and by \cref{prop:classifier}
  it has a subobject classifier, as is the case for any example
  constructed in this way.
\end{example}

\begin{remark}
  Not only is the above \cref{ex:symcat} a concrete basis; in fact the
  category of elements $\El(\pP)$ satisfies a stronger infinitary
  version of co-confluence. Namely, every \emph{wide cospan} (a
  diagram indexed by a category obtained by freely adjoining a
  terminal object to a discrete category) admits a cone. Moreover, the
  proofs of \cref{lem:weaktaut},\footnote{To generalize
    \cref{lem:weaktaut}, the appropriate analogue of a pullback of a
    cospan one of whose legs is a monomorphism is a wide pullback of a
    wide cospan all but one of whose legs are monomorphisms, in which
    case the opposing leg of the resulting wide pullback wide span is
    a monomorphism.}  \cref{lem:weakpreserve}, \cref{prop:influent},
  \cref{cor:influentwp}, and \cref{cor:wpb} all generalize to
  analogues with respect to wide cospans up to any particular size.
  In particular, a comonad on $\Set$ preserves (small) \emph{weak wide
    pullbacks} if and only if it admits a basis
  $\pP \colon \bB \to \Set$ such that $\El(\pP)$ is infinitarily
  co-confluent (in the sense that every small wide cospan admits a
  cone).

  The weak-wide-pullback-preserving functors $\Set \to \Set$ that are
  small (a.k.a.\ accessible) are precisely the coproducts of
  \emph{symmetrized representables}~\cite[Theorem
  3.3]{adamek-velebil}. Such sums of symmetrized representables are
  also called \emph{generalized analytic functors}
  in~\cite{pare:taut}, since the finitary such endofunctors are the
  \emph{analytic functors}~\cite{joyal:analytic}. The comonad from
  \cref{ex:part} whose coalgebras are partitioned sets $\Part$ is also
  given by a basis whose category of elements satisfies the
  (small-sized) infinitary co-confluence condition, so it preserves
  (small) weak wide pullbacks as well, though it is not small.
\end{remark}

Co-confluence is not the only natural condition on $\El(\pP)$ to
generalize the notion of basis. Recall that a functor
$\pP\colon \bB \to \Set$ is called \emph{flat} if its category of
elements is \emph{cofiltered}, meaning every finite diagram admits a
cone. (Or equivalently: the empty diagram, every pair of objects, and
every parallel pair of morphisms admits a cone.)
In~\cite{garner:ionads}, Garner defines a \emph{basis of an ionad} to
be a functor $\bB \to \Set^\xX$ (specifically with $\bB$ small) that
is componentwise flat. Composing with
$\slicel \colon \Set^\xX \to \Set$, we see this is equivalent to the
data of a functor $\pP \colon \bB \to \Set$ that decomposes as a sum
of $\xX$ many flat functors. Such a functor $\pP$ into $\Set$ may also
be intrinsically characterized as one with category of elements
$\El(\pP)$ satisfying the following condition.

\begin{definition}\label{def:ionadbasis}
  We call a category \emph{connectedly-cofiltered} (a.k.a.\
  \emph{FINCL-cofiltered} in the terminology of
  \cite{adamek-borceux-lack-rosicky}) if each connected component is
  cofiltered. Equivalently, every finite connected diagram admits a
  cone.

  We call a functor $\pP \colon \bB \to \Set$ \emph{connectedly-flat}
  (a.k.a.\ \emph{FINCL-flat}) if its category of elements is
  connectedly-cofiltered.
\end{definition}

The only distinction between co-confluent and connectedly-cofiltered
is that the latter additionally requires every parallel pair of arrows
to admit an arrow equalizing them. For preorder-shaped diagrams, the
two conditions coincide, and thus both generalize the topological
definition of basis. For non-preorder-shaped diagrams, co-confluent is
strictly weaker than connectedly-cofiltered.

It is well-known that $\pP$ is flat (assuming $\bB$ small) if and only
if its Yoneda extension $\extend{\pP} \colon \Ps{\bB} \to \Set$
preserves finite limits. The following is similar.

\begin{lemma}[{\cite{adamek-borceux-lack-rosicky}}]\label{lem:connectedly}
  Let $\pP \colon \bB \to \Set$ with $\bB$ small. Then $\pP$ is
  connectedly-flat, i.e.\ its category of elements $\El(\pP)$ is
  connectedly-cofiltered, if and only if its Yoneda extension
  $\extend{\pP} \colon \Ps{\bB} \to \Set$ preserves pullbacks
  (equivalently, preserves finite connected limits).\qed
\end{lemma}

This yields analogues of the above \cref{cor:influentwp,cor:wpb} for
pullback-preserving comonads:

\begin{corollary}\label{cor:connpb}
  Let $\pP \colon \bB \to \Set$ be connectedly-flat with density
  comonad $\cP$. Then $\cP$ preserves pullbacks.\qed
\end{corollary}

\begin{corollary}\label{cor:pb}
  A comonad on $\Set$ preserves pullbacks if and only if it admits a
  connectedly-flat basis.\qed
\end{corollary}

\begin{remark}
  It was shown in~\cite[Proposition 3.4]{garner:ionads} that given any
  flat functor $\bB \to \Set^X$ with $\bB$ small (a.k.a.\ basis of an
  ionad), the category of coalgebras of the induced
  finite-limit-preserving density comonad on $\Set^X$ (a.k.a.\ ionad)
  is equivalent to the topos of sheaves on the site given by $\bB$
  where the covering families are those families sent to jointly
  epimorphic families in $\Set^X$. As we will explain in
  \cref{part:II}, the categories of coalgebras of such
  finite-limit-preserving comonads on $\Set^X$ are identified with
  those of corresponding pullback-preserving comonads on $\Set$, with
  the translation process given on comonadic functors by composing
  with the sum functor $\slicel \colon \Set^X \to \Set$. Moreover this
  translation process respects bases of comonads (see
  \cref{prop:liftbasis}), and it is easy to see the flat functors
  $\bB \to \Set^X$ correspond to the connectedly-flat functors
  $\bB \to \Set$ with colimit $X$. Hence Garner's result on bases
  applies equally well to pullback-preserving comonads on $\Set$: if
  $\pP \colon \bB \to \Set$ is a connectedly-flat functor with $\bB$
  small, the category of coalgebras $\SetP$ of the induced
  pullback-preserving density comonad $\cP$ on $\Set$ is equivalent to
  the topos of sheaves on the site given by $\bB$ where the covering
  families are those families sent to jointly surjective families in
  $\Set$.
\end{remark}

\begin{remark}
  There exist analogues of these results for other well-behaved
  classes of limits (those described by \emph{sound doctrines} in the
  sense of~\cite{adamek-borceux-lack-rosicky}, or more generally
  \emph{weakly sound doctrines} in the sense of~\cite{tendas}).
\end{remark}

\begin{example}
  The diagrams from \cref{ex:topbasis}, \cref{ex:catbasis}, and
  \cref{ex:topgroup} all have connectedly-cofiltered category of
  elements, so all are bases of pullback-preserving comonads.
  
  On the other hand, the comonadic forgetful functor $\Part \to \Set$
  from \cref{ex:part} does not preserve pullbacks.  (For example, it
  does not preserve the fixed points of a transposition permutation on
  a single part $3$, whereas a pullback-preserving functor would
  preserve all connected limits.)
  Neither does the comonadic forgetful functor $\MultiGph \to \Set$
  from \cref{ex:symcat}. (For example, it does not preserve the fixed
  points of the automorphism swapping the two vertices of the graph
  with two vertices and one edge between them.) These comonads are
  therefore not density comonads of any diagram whose category of
  elements is co-confluent.
\end{example}

For our final example, we generalize the construction of the comonad
from \cref{ex:topgroup} whose coalgebras are ``sets with an
infinitesimal $\tG$-action'' beyond topological groups $\tG$ to
categories $\bC$ equipped with a suitable ``topology''.\footnote{This
  is \emph{not} a Grothendieck topology. Here every ``open
  neighborhood'' of an object $c$ contains the identity $\id_c$,
  whereas if every covering sieve in a Grothendieck topology were to
  contain the identity then the topology would be trivial.}

\begin{example}[Infinitesimal categories]\label{ex:topcat}
  Let $\bC$ be a small category. Assign to each object $c$ a set
  $\neighbors{c}$ of subsets of $\pts{c/\bC}$ (the arrows out of $c$)
  called \emph{basic open neighborhoods} such that
  \begin{itemize}
  \item every member of $\neighbors{c}$ contains the identity $\id_c$;
  \item for any pair $U$ and $V$ in $\neighbors{c}$ and
    $x \colon c \to d$ in $U \cap V$, there exists some $W$ in
    $\neighbors{d}$ such that for every $y \colon d \to e$ in $W$, we
    have that $y \circ x$ is in $U \cap V$; and
  \item there is some member of $\neighbors{c}$.\footnote{This
      condition may be understood as the nullary intersection form of
      the previous binary intersection condition.}
  \end{itemize}
  If $\bC$ is indiscrete, then this is the same as endowing
  $\pts{\bC}$ with a topology and assigning an open neighborhood basis
  to each point. Any topological group (viewed as a one-object
  category) equipped with a neighborhood basis of identity is also an
  example.

  Now define a category $\bB$ in which an object is a pair $(c, U)$
  where $c \in \bC$ and $U \in \neighbors{c}$, and in which an arrow
  $(d, V) \to (c, U)$ is an arrow $x \colon c \to d$ such that
  pre-composing $x$ sends $U$ within $V$, i.e.\ for all $y \colon d \to e$
  in $V$, we have that $y \circ x$ is in $U$. There is a forgetful
  functor $\pP \colon \bB \to \Set$ sending $(c, U)$ to $U$. This
  $\pP$ is connectedly-cofiltered. Indeed, the first two conditions
  above precisely ensure that in $\El(\pP)$ any two elements with
  codomain $c$ admit arrows from an element of the form $\id_c$, and
  the last condition ensures the set of elements with codomain $c$ is
  inhabited; moreover any parallel pair in $\El(\pP)$ is equalized by
  any map from an identity element. Hence assuming $\bC$ is small,
  $\pP$ is a basis of a pullback-preserving comonad on $\Set$ by
  \cref{cor:influentbasis} and \cref{cor:connpb}.

  Quite similarly to \cref{ex:topgroup}, we find that $\cP(\xA)$ is
  the set of all germs about identity of $\xA$-valued functions
  \[\cP(\xA) = \sum_{c\in \obs{\bC}} \colim_{U \ni c} \xA^U\]
  where $U$ ranges over $\neighbors{c}$. And again using
  \cref{rem:simplergerms}, we find that a $\cP$-coalgebra
  $\eE \to \cP(\eE)$ assigns to each element $e$ the germ of some
  function $f \colon U \to \eE$ with $f(\id_c) = e$ and such that for
  all $x \colon c \to d$ in $U$, the element $f(x)$ in $\eE$ is
  assigned the germ \emph{about $x$} of the same function $f$ ---
  which reduces to a germ about $\id_d$ of $f(\dash \circ x)$. That
  is, we have at each element $e$ an assigned object $c$ of $\bC$ and an
  infinitesimal neighborhood of arrows out of $c$ acting on $e$, with
  the action satisfying unit and associativity laws. To put it another
  way, a $\cP$-coalgebra is similar to a functor $\bC \to \Set$, but
  where actions by arrows in $\bC$ are only defined infinitesimally.
  
  The set of points of $\cP$ is the set $\obs{\bC}$ of objects of
  $\bC$. If the category $\bC$ is indiscrete, then we obtain the
  comonad corresponding to the generated topology on $\pts{\bC}$. And
  regardless of what $\bC$ is, if we take each $\neighbors{c}$ to
  consist of only $\id_c$, then we obtain the discrete comonad on
  $\pts{C}$ points. On the other hand if we take each $\neighbors{c}$
  to consist only of the maximal subset, then we obtain the canonical
  basis of the category $\bC$ from \cref{ex:catbasis}.

  Any comonad induced by a topological group discussed in
  \cref{ex:topgroup} is an example of this example. For an interesting
  example that is not given by a topological group, consider $\nn$,
  the free monoid on an element viewed as a one-object category, and
  let the open neighborhoods about $0$ be sets of the form
  $\setof{0} \cup \setof{i \mid i > n}$. Then $\cP(\xA)$ sends a set
  to the set of pairs $(a, t)$, where $a$ is an element of $\xA$ and
  $t$ is a \emph{tail} of an $\xA$-valued sequence, i.e.\ an
  equivalence class of $\xA$-valued sequences modulo eventual
  agreement. A $\cP$-coalgebra is similar to a $\nn$-set, or discrete
  dynamical system, but where elements only have eventual infinite
  behaviors.

  There is also an analogous subbasis form, similar to
  \cref{ex:infinitesimalsubbasis}, replacing the last two conditions
  above with the following two conditions: for any $U$ in
  $\neighbors{c}$ and $x \colon c \to d$ in $U$, there exists some $V$
  in $\neighbors{d}$ such that for every $y \colon d \to e$ in $V$, we
  have that $y \circ x$ is in $U$; and $\neighbors{c}$ is connected as
  a poset under subset inclusion. (Together these ensure that
  $\pts{\pP}$ is $\obs{\bC}$ and the germs are calculated as
  expected.)
\end{example}

\begin{remark}\label{rem:preorders}
  Some comonads on $\Set$ correspond to both topological spaces
  (\cref{ex:top}) and categories (\cref{ex:cat}) at once. Such
  comonads may be identified with preorders, viewed as either
  Alexandroff topological spaces, whose open sets are upward-closed
  sets in the preorder, or as thin categories.

  Indeed, we may show this abstractly as follows. It is well-known
  that such Alexandroff topological spaces $\tX$ are those for which
  arbitrary intersections of open sets remain open. The category of
  elements $\El(\topbasis{\tX})$ of its diagram of open subsets
  $\topbasis{\tX} \colon \opens{\tX}\to \Set$ includes a cone not only to
  every \emph{finite} connected diagram, but to every connected
  diagram simpliciter. (In the framework of
  \cite{adamek-borceux-lack-rosicky}, $\topbasis{\tX}$ is flat with respect to
  connected diagrams up to any particular size.)  Analogously to
  \cref{lem:connectedly}, this is equivalent to the condition that the
  Yoneda extension $\extend{\topbasis{\tX}} \colon \Psh(\tX) \to \Set$
  preserves wide pullbacks (or equivalently connected
  limits). Therefore the density comonad
  $\gen{\topbasis{\tX}} \coloneqq \extend{\topbasis{\tX}} \circ
  \nerve{{\topbasis{\tX}}}$ preserves wide pullbacks, a.k.a.\ is
  polynomial. As mentioned in the introduction, and as we will further
  elaborate upon in \cref{part:II}, the polynomial comonads on $\Set$
  are precisely the comonads corresponding to
  categories~\cite{ahman-uustalu}. And conversely, Alexandroff spaces
  are the only topological spaces $\tX$ for which the comonadic
  forgetful functor $\Sh(\tX) \to \Set$ creates wide pullbacks, since
  this implies arbitrary intersections of opens are open. (In general,
  the underlying topological space of the comonad associated to any
  category is an Alexandroff topological space by
  \cref{prop:alextop}.)

  We can also view this correspondence through
  \cref{thm:tfaetop}. Categories yielding topological comonads are
  always thin: if, as in \cref{thm:tfaetop}, every object of
  $\fun{\bC}{\Set}$ is covered by subterminal objects, then in
  particular the subobject generated by the identity element in each
  representable presheaf is subterminal. Therefore, each object of
  $\bC$ has at most one arrow in from each other object. Conversely,
  thin categories yield topological comonads, since the ``canonical
  basis'' (\cref{ex:catbasis}) of a thin category can be viewed as a
  diagram of subsets of its set of objects.

  
  
  More explicitly, let us compare formulas for the comonad associated
  to a preorder $\pre{P}$ as an Alexandroff topological space $\tX$
  versus as a category $\bC$, which are supposedly the same:
  \[
    \cX(\xA) = \sum_{p\in \pts{\pre{P}}} \colim_{U \ni p} \xA^U
    \qqand
    \cbC(\xA) = \sum_{p\in \pts{\pre{P}}} \xA^{\pts{p / \pre{P}}}
  \]
  Here $U$ ranges over upward-closed sets of $\pre{P}$ containing the
  point $p$ (i.e.\ open neighborhoods of $p$ in the Alexandroff
  topology), and $\pts{p / \pre{P}}$ denotes the set of points above
  $p$, inclusive. The set $\pts{p / \pre{P}}$ is the minimal upward
  closed set containing $p$, and thus is a terminal object in the
  category that indexes the colimit $\colim_{U \ni p} \xA^U$, namely
  the opposite of the poset of upward-closed sets of $\pre{P}$. Hence
  this colimit simply evaluates to $\xA^{\pts{p / \pre{P}}}$ as
  claimed.
\end{remark}

\mypart{II}{Maps}

As mentioned in the introduction, comonads on $\Set$ subsume not only
topological spaces and categories, but also
ionads~\cite{garner:ionads}, an existing notion of space generalizing
both topological spaces and categories at once. In terms of topos
theory, ionads are ``toposes with points'': they generalize
topological spaces in the same way that toposes generalize
locales. There is also a notion of morphism between ionads called a
\emph{continuous map}, generalizing both continuous maps of
topological spaces and functors of categories at once. The study of
continuous maps of general comonads, as well as their interaction with
ordinary comonad maps, will be the main focus of the rest of the paper.

\section{Comonad functors}\label{sec:comonadfunctors}

More general than ordinary comonad maps (\cref{def:comonad}), there
are also standard definitions of morphism between comonads on
different categories, as well as 2-cells between these morphisms.

\tikzset{every picture/.style={baseline={([yshift=-3.5pt]current bounding box.center)},scale=.125}}
\begin{definition}\label{def:comfun}
  Let $\cC_1$ and $\cC_2$ be comonads on $\bS_1$ and $\bS_2$. A
  \emph{colax comonad functor} from $\cC_1$ to $\cC_2$ consists of a
  functor $F \colon \bS_1 \to \bS_2$ and a 2-cell
  $\chi \colon F \circ \cC_1 \Rightarrow \cC_2 \circ F$ satisfying
  
  \[
    \begin{tikzpicture}[scale=.825,yscale=-1]
      \node [ob] (s) at (-7,0) {$\bS_1$};
      \node [ob] (t) at (10.5,3.5) {$\bS_2$};
      \node [ob] (u) at (0,7) {$\bS_2$};
      \node [ob] (d) at (3.5,-3.5) {$\bS_1$};
      \node [ob] (tt) at (14,-7) {$\bS_2$};
      \node [ob] (dd) at (7,-14) {$\bS_1$};
      \draw [a] (s) -- node[arr,bl] {$F$} (u);
      \draw [a] (u) -- node[arr,below] {$\cC_2$} (t);
      \draw [a] (s) -- node[arr,above,pos=.3] {$\cC_1$} (d);
      \draw [a] (d) -- node[arr,bl] {$F$} (t);
      \draw [a] (t) -- node[arr,right] {$\cC_2$} (tt);
      \draw [a] (d) -- node[arr,left,pos=.7] {$\cC_1$} (dd);
      \draw [a] (dd) -- node[arr,ar] {$F$} (tt);
      \draw [a] (s) .. controls +(2,-8) and +(-8,2) .. node[arr,al] {$\cC_1$} (dd);
      \node [cell] at (1.75,1.75) {$\chi$};
      \node [cell] at (8.75,-5.25) {$\chi$};
      \node [cell] at (0,-7) {$\delta$};
    \end{tikzpicture}
    \eqq
    \begin{tikzpicture}[scale=.825,yscale=-1]
      \node [ob] (s) at (-7,0) {$\bS_1$};
      \node [ob] (t) at (10.5,3.5) {$\bS_2$};
      \node [ob] (u) at (0,7) {$\bS_2$};
      \node [ob] (tt) at (14,-7) {$\bS_2$};
      \node [ob] (dd) at (7,-14) {$\bS_1$};
      \draw [a] (s) -- node[arr,bl] {$F$} (u);
      \draw [a] (u) -- node[arr,below] {$\cC_2$} (t);
      \draw [a] (t) -- node[arr,right] {$\cC_2$} (tt);
      \draw [a] (dd) -- node[arr,ar] {$F$} (tt);
      \draw [a] (u) .. controls +(2,-8) and +(-8,2) .. node[arr,al,pos=.8] {$\cC_2$} (tt);
      \draw [a] (s) .. controls +(2,-8) and +(-8,2) .. node[arr,al] {$\cC_1$} (dd);
      \node [cell] at (.5,-5) {$\chi$};
      \node [cell] at (7,0) {$\delta$};
    \end{tikzpicture}
    \qquad\qquad
    \begin{tikzpicture}[scale=.825,yscale=-1]
      \node [ob] (s) at (-7,0) {$\bS_1$};
      \node [ob] (u) at (0,7) {$\bS_2$};
      \node [ob] at (0,7) {$\phantom{\bS_2}$};
      \node [ob] (tt) at (14,-7) {$\bS_2$};
      \node [ob] (dd) at (7,-14) {$\bS_1$};
      \draw [a] (s) -- node[arr,bl] {$F$} (u);
      \draw [a] (dd) -- node[arr,ar] {$F$} (tt);
      \draw [a] (s) .. controls +(2,-8) and +(-8,2) .. node[arr,al] {$\cC_1$} (dd);
      \draw [eq] (u) .. controls +(8,-2) and +(-2,8) .. (tt);
      \draw [eq] (s) .. controls +(8,-2) and +(-2,8) .. (dd);
      \node [cell] at (5.75,-1.25) {$=$};
      \node [cell] at (-.25,-7.25) {$\varepsilon$};
    \end{tikzpicture}
    \eqq
    \begin{tikzpicture}[scale=.825,yscale=-1]
      \node [ob] (s) at (-7,0) {$\bS_1$};
      \node [ob] (u) at (0,7) {$\bS_2$};
      \node [ob] (tt) at (14,-7) {$\bS_2$};
      \node [ob] (dd) at (7,-14) {$\bS_1$};
      \draw [a] (s) -- node[arr,bl] {$F$} (u);
      \draw [a] (dd) -- node[arr,ar] {$F$} (tt);
      \draw [a] (u) .. controls +(2,-8) and +(-8,2) .. node[arr,al,pos=.8] {$\cC_2$} (tt);
      \draw [eq] (u) .. controls +(8,-2) and +(-2,8) .. (tt);
      \draw [a] (s) .. controls +(2,-8) and +(-8,2) .. node[arr,al] {$\cC_1$} (dd);
      \node [cell] at (.5,-5) {$\chi$};
      \node [cell] at (6.75,-.25) {$\varepsilon$};
    \end{tikzpicture}
  \]
\end{definition}

In particular, colax comonad functors carried by identity are the same
as ordinary comonad maps. There are two types of 2-cell between colax
comonad functors. The first is from~\cite{street:monads}, and the
second is from~\cite{wood:ii} and \cite{lack-street}.

\begin{definition}\label{def:comtrans}
  Let $F_1, F_2\colon \bS_1 \to \bS_2$ be colax comonad functors. A
  \emph{comonad transformation} from $F_1$ to $F_2$ is a 2-cell
  $\gamma \colon F_1 \Rightarrow F_2$ satisfying
  \[
    \begin{tikzpicture}[xscale=.765625,yscale=-.875]
      \path (0,10) -- (0,-10);
      \node [ob] (s) at (-10,0) {$\bS_1$};
      \node [ob] (t) at (10,0) {$\bS_2$};
      \node [ob] (u) at (0,7) {$\bS_2$};
      \node [ob] (d) at (0,-7) {$\bS_1$};
      \draw [a] (s) -- node[arr,arcl,pos=.33] {$F_1$} (u);
      \draw [a] (u) -- node[arr,br] {$\cC_2$} (t);
      \draw [a] (s) -- node[arr,al] {$\cC_1$} (d);
      \draw [a] (d) -- node[arr,ar] {$F_1$} (t);
      \draw [a] (s) .. controls +(-1,6) and +(-6,2) .. node[arr,bl] {$F_2$} (u);
      \node [cell] at (.5,-.5) {$\chi$};
      \node [cell] at (-6.5,5) {$\gamma$};
    \end{tikzpicture}
    \eqq
    \begin{tikzpicture}[xscale=.765625,yscale=-.875]
      \path (0,10) -- (0,-10);
      \node [ob] (s) at (-10,0) {$\bS_1$};
      \node [ob] (t) at (10,0) {$\bS_2$};
      \node [ob] (u) at (0,7) {$\bS_2$};
      \node [ob] (d) at (0,-7) {$\bS_1$};
      \draw [a] (s) -- node[arr,bl] {$F_2$} (u);
      \draw [a] (u) -- node[arr,br] {$\cC_2$} (t);
      \draw [a] (s) -- node[arr,al] {$\cC_1$} (d);
      \draw [a] (d) -- node[arr,blcl,pos=.66] {$F_2$} (t);
      \draw [a] (d) .. controls +(6,-2) and +(1,-6) .. node[arr,ar] {$F_1$} (t);
      \node [cell] at (-.5,.5) {$\chi$};
      \node [cell] at (6.5,-5.25) {$\gamma$};
    \end{tikzpicture}
  \]
\end{definition}

The other kind of 2-cell was not explicitly given its own name
in~\cite{wood:ii} or~\cite{lack-street}, so here we call it by its
application in our paper.

\begin{definition}\label{def:comspec}
  A \emph{comonad specialization} from $F_1$ to $F_2$ is a 2-cell
  $\sigma \colon \cC_1 \then F_1 \Rightarrow F_2$ satisfying
  \[
    \begin{tikzpicture}[xscale=.765625,yscale=-.875]
      \node [ob] (ss) at (-20,7) {$\bS_1$};
      \node [ob] (s) at (-10,0) {$\bS_1$};
      \node [ob] (t) at (10,0) {$\bS_2$};
      \node [ob] (u) at (0,7) {$\bS_2$};
      \node [ob] (d) at (0,-7) {$\bS_1$};
      \draw [a] (ss) -- node[arr,alcl] {$\cC_1$} (s);
      \draw [a] (s) -- node[arr,arcl,pos=.15] {$F_1$} (u);
      \draw [a] (u) -- node[arr,br] {$\cC_2$} (t);
      \draw [a] (s) -- node[arr,alcl] {$\cC_1$} (d);
      \draw [a] (d) -- node[arr,ar] {$F_1$} (t);
      \draw [a] (ss) -- node[arr,below] {$F_2$} (u);
      \draw [a] (ss) .. controls +(-3,-11) and +(-11,-3) .. node[arr,al] {$\cC_1$} (d);
      \node [cell] at (.5,0) {$\chi$};
      \node [cell] at (-10,4.25) {$\sigma$};
      \node [cell] at (-13.75,-2.625) {$\delta$};
    \end{tikzpicture}
    \eqq
    \begin{tikzpicture}[xscale=.765625,yscale=-.875]
      \node [ob] (ss) at (-20,7) {$\bS_1$};
      \node [ob] (s) at (-10,0) {$\bS_1$};
      \node [ob] (t) at (10,0) {$\bS_2$};
      \node [ob] (u) at (0,7) {$\bS_2$};
      \node [ob] (d) at (0,-7) {$\bS_1$};
      \draw [a] (ss) -- node[arr,alcl] {$\cC_1$} (s);
      \draw [a] (s) -- node[arr,below] {$F_2$} (t);
      \draw [a] (u) -- node[arr,br] {$\cC_2$} (t);
      \draw [a] (s) -- node[arr,alcl] {$\cC_1$} (d);
      \draw [a] (d) -- node[arr,ar] {$F_1$} (t);
      \draw [a] (ss) -- node[arr,below] {$F_2$} (u);
      \draw [a] (ss) .. controls +(-3,-11) and +(-11,-3) .. node[arr,al] {$\cC_1$} (d);
      \node [cell] at (0,-2.75) {$\sigma$};
      \node [cell] at (-6.5,3.75) {$\chi$};
      \node [cell] at (-13.75,-2.625) {$\delta$};
    \end{tikzpicture}
  \]
\end{definition}
\tikzset{every picture/.style={}}

The following generalizes the correspondence between comonad maps and
commutative triangles from \cref{prop:respect} to a correspondence
between colax comonad functors and commutative squares.

\begin{proposition}[\cite{street:monads,lack-street}]\label{prop:formaltfae}
  Let $\cC_1$ and $\cC_2$ be comonads on $\bS_1$ and $\bS_2$, and
  let $F \colon \bS_1 \to \bS_2$ be an arbitrary
  functor.\footnote{These results hold more generally in an abstract
    2-category (assuming the relevant Eilenberg-Moore objects exist).}
  
  \begin{enumerate}[label=(\roman*)]
  \item Giving $F$ the structure of a colax comonad functor is
    equivalent to giving a functor between categories of coalgebras
    $G \colon \CoalgOn{\cC_1}{\bS_1} \to \CoalgOn{\cC_2}{\bS_2}$ such that
    the following square commutes:
    \[
      \begin{tikzcd}[column sep=15]
        \CoalgOn{\cC_1}{\bS_1} & \CoalgOn{\cC_2}{\bS_2} \\
        {\bS_1} & {\bS_2}
        \arrow["G",from=1-1, to=1-2]
        \arrow["\car{\cC_1}"', from=1-1, to=2-1]
        \arrow["\car{\cC_2}", from=1-2, to=2-2]
        \arrow["F"', from=2-1, to=2-2]
      \end{tikzcd}
    \]
  \end{enumerate}
  Now let $F_1$ and $F_2$ be colax comonad functors from $\cC_1$ to $\cC_2$,
  and let $\gamma \colon F_1 \Rightarrow F_2$ be an arbitrary natural
  transformation.
  \begin{enumerate}[label=(\roman*)]
    \setcounter{enumi}{1}
  \item Giving $\gamma$ the structure of a comonad transformation
    from $F_1$ to $F_2$ is equivalent to giving a natural
    transformation $\sigma' \colon G_1 \Rightarrow G_2$ forming a
    cylinder
    \[
      \begin{tikzcd}
        \CoalgOn{\cC_1}{\bS_1} & \CoalgOn{\cC_2}{\bS_2} \\
        \bS_1 & \bS_2
        \arrow["\car{\cC_1}"', from=1-1, to=2-1]
        \arrow["G_1", curve={height=-12pt}, from=1-1, to=1-2]
        \arrow["\car{\cC_2}", from=1-2, to=2-2]
        \arrow[""{name=0, anchor=center, inner sep=0}, "F_2"', curve={height=12pt}, from=2-1, to=2-2]
        \arrow[""{name=1, anchor=center, inner sep=0}, "F_1", curve={height=-12pt}, from=2-1, to=2-2]
        \arrow["\gamma"{description}, draw=none, from=1, to=0]
      \end{tikzcd}
      \quad=\quad
      \begin{tikzcd}
        \CoalgOn{\cC_1}{\bS_1} & \CoalgOn{\cC_2}{\bS_2} \\
        \bS_1 & \bS_2
        \arrow["\car{\cC_1}"', from=1-1, to=2-1]
        \arrow[""{name=0, anchor=center, inner sep=0}, "G_1", curve={height=-12pt}, from=1-1, to=1-2]
        \arrow[""{name=1, anchor=center, inner sep=0}, "G_2"', curve={height=12pt}, from=1-1, to=1-2]
        \arrow["\car{\cC_2}", from=1-2, to=2-2]
        \arrow["F_2"', curve={height=12pt}, from=2-1, to=2-2]
        \arrow["\sigma'"{description}, draw=none, from=0, to=1]
      \end{tikzcd}
    \]
  \item Giving a comonad specialization from $F_1$ to $F_2$ is
    equivalent to giving an arbitrary natural transformation
    $\sigma' \colon G_1 \Rightarrow G_2$.\qed
  \end{enumerate}
\end{proposition}

We refer the reader to~\cite{fairbanks:monads}, originally written as
an appendix for this paper, for more details and explicit proofs, as
well as illustrations of the definitions with string diagrams.

\section{Ionads}\label[bgchapter]{sec:ionads}

Ionads were introduced in~\cite{garner:ionads} as a
category-theoretic analogue of the order-theoretic notion of
topological space. We will give different proofs of some of Garner's
results below, more closely adapted to the themes of this paper.

\begin{definition}[{\cite[Definition 2.1]{garner:ionads}}]\label{def:ionadcontinuous}
  An \emph{ionad} $\iX$ consists of a set $\pts{\iX}$ and a
  finite-limit-preserving comonad 
  on $\SetpX$.
\end{definition}


We will denote\footnote{This notation $\Sh$ differs from Garner's
  notation $\mathbf{O}$.} the category of coalgebras of an ionad $\iX$
by $\ionopens{\iX}$. 

Morphisms of ionads generalize continuous maps of topological spaces.

\begin{definition}[{\cite[Definition 4.1]{garner:ionads}}]\label{def:ionadmap}
  A \emph{continuous map} of ionads $f \colon \iX \to \iY$ consists of
  a function $\one{f} \colon \pts{\iX} \to \pts{\iY}$ and a functor
  $\conion{f} \colon \ionopens{\iY} \to \ionopens{\iX}$ making the
  following diagram commute
  \[
    \begin{tikzcd}
      \ionopens{\iX} & \ionopens{\iY} \\
      {\SetpX} & {\SetpY}
      \arrow["\carion{\iX}"', from=1-1, to=2-1]
      \arrow["{\conion{f}}"', dashed, from=1-2, to=1-1]
      \arrow["\carion{\iY}", from=1-2, to=2-2]
      \arrow["{\reindex{\one{f}}}", from=2-2, to=2-1]
    \end{tikzcd}
  \]
  where $\reindex{\one{f}}$ is the reindexing functor given by a
  function $\one{f} \colon \pts{\iX} \to \pts{\iY}$.

  By \cref{prop:formaltfae}, this is the same as a colax comonad
  functor (\cref{def:comfun}) carried by $\reindex{\one{f}}$.
\end{definition}

\begin{remark}
  The definition of ionad is obtained from the order-theoretic
  definition of topological space as an \emph{interior operator}, or
  finite-limit-preserving comonad, on the poset of subsets
  $2^{\pts{\tX}}$ of a set $\pts{\tX}$, by replacing the subsets
  $2^{\pts{\tX}}$ with indexed sets $\SetpX$. The definition of
  continuous map of ionads is obtained similarly from the definition
  of continuous map of topological spaces.
  \[
    \begin{tikzcd}
      \opens{\tX} & \opens{\tY} \\
      {2^{\pts{\tX}}} & {2^{\pts{\tY}}}
      \arrow[hook, from=1-1, to=2-1]
      \arrow[dashed, from=1-2, to=1-1]
      \arrow[hook, from=1-2, to=2-2]
      \arrow["{\one{f}^{-1}}", from=2-2, to=2-1]
    \end{tikzcd}
  \]
\end{remark}

We write $\IonadO$ for the category of ionads and continuous maps.  An
ionad is equivalent to a topos equipped with a \emph{geometric
  surjection}~\cite[A4.2]{johnstone:elephant} from a discrete topos
$\Set^X$, a.k.a. a set $X$ of \emph{enough points} or a
\emph{separating set of points}~\cite[C2.2.12]{johnstone:elephant}.

\begin{proposition}\label{lem:ionadtopos}
  The category $\IonadO$ is equivalent to the category in which an
  object is a topos equipped with a separating set of points and an
  arrow is a geometric morphism equipped with a compatible map on
  points (i.e.\ a commutative square of geometric morphisms).
\end{proposition}
\begin{proof}
  Categories of coalgebras of finite-limit-preserving comonads on
  toposes (such as $\Set^\xX$) are again toposes~\cite[Theorem
  A4.2.1]{johnstone:elephant}.  A separating set of points $\xX$ of a
  topos, essentially by definition, amounts to a
  finite-limit-preserving comonadic functor into $\Set^\xX$ (the
  inverse image of a geometric surjection).

  Moreover, given a continuous map of ionads $f$, the functor
  $\conion{f}$ is always the inverse image of a geometric morphism.
  Indeed, since the comonads preserve finite limits and the comonadic
  functors create finite limits, $\reindex{\one{f}}$ preserves finite
  limits, so $\conion{f}$ preserves finite limits as well. Moreover
  $\conion{f}$ has a right adjoint by the adjoint triangle theorem
  (\cref{lem:lifting}). Thus the definition of continuous map of
  ionads is precisely a commutative square of geometric morphisms.
\end{proof}

Garner also introduced an appropriate notion of 2-cell between
continuous maps of ionads, so that ionads form not only a category but
a 2-category. Recall that if $x$ and $y$ are points in a topological
space $\tX$ and every open set containing $x$ contains $y$, then we
say $x$ is a \emph{specialization} of $y$, written $x \leq
y$. Moreover, if $f$ and $g$ are continuous maps of topological spaces
$\tX \to \tY$ and $f(x) \leq g(x)$ for each $x \in \tX$ we say $f$ is
a \emph{(pointwise) specialization} of $g$, written $f \leq g$. This
generalizes to the following definition.

\begin{definition}[{\cite[Definition 5.1]{garner:ionads}}]\label{def:ionadspecialization}
  A \emph{specialization} between continuous maps of ionads
  $f \Rightarrow g$ is a natural transformation
  $\conion{f} \Rightarrow \conion{g}$.

  By \cref{prop:formaltfae}, this is the same as a comonad
  specialization, as defined in \cref{def:comspec}, from $f$ to $g$.
\end{definition}



Ionads provide a common generalization of topological spaces --- as
sheaf toposes equipped with the topological space's set of points ---
and categories --- as presheaf toposes equipped with the category's
set of objects.


\begin{proposition}[{\cite[Example 4.8]{garner:ionads}}]\label{prop:topionad}
  The full sub-2-category of $\Ionad$ consisting of ionads
  corresponding to topological spaces is equivalent to the 2-category
  of topological spaces, continuous functions, and pointwise
  specializations.
  
  Here the comonadic functor corresponding to a topological space
  $\tX$ is the inverse image functor $\Sh(\tX) \to \Set^{\pts{\tX}}$
  given by the inclusion of points $\pts{\tX} \to \tX$.
\end{proposition}


\begin{proof}[Proof sketch]
  Suppose given a continuous map $f$ between ionads corresponding to
  topological spaces $\tX$ and $\tY$. The inverse image map $\conion{f}$
  preserves both colimits and finite limits, in particular subterminal
  objects. Thus $\conion{f}$ restricts to a functor
  $\opens{\tY} \to \opens{\tX}$ between posets of opens (subterminal
  sheaves), and in fact this determines $\conion{f}$ since all sheaves
  are canonical colimits of opens.
  Natural transformations between functors $\conion{f}$ are likewise
  determined by natural transformations between maps between these
  restrictions. The condition that the map
  $\opens{\tY} \to \opens{\tX}$ lifts $\reindex{\one{f}}$ is precisely
  the condition that preimages of open subsets are open, and natural
  transformations between such inverse image functors are precisely
  pointwise specializations.  Conversely, every continuous map of
  topological spaces extends to a geometric morphism between sheaf
  categories, respecting points as desired.
\end{proof}

Continuous maps of ionads corresponding to categories are functors,
and specializations are natural transformations.

\begin{proposition}[{\cite[Example 4.4]{garner:ionads}}]\label{lem:categories}\label{prop:cationad}
  The full sub-2-category of $\Ionad$ consisting of ionads whose
  comonad part preserves limits is equivalent to the 2-category $\Cat$
  of small categories, functors, and natural transformations.
  
  Here the comonadic functor corresponding to a category $\bC$ is the
  reindexing functor $\fun{\bC}{\Set} \to \Set^{\obs{\bC}}$ given by
  the inclusion of objects $\pts{\bC} \to \bC$.
\end{proposition}

This provides additional motivation for the definition of continuous
map by relating it to the definition of functor. The proof sketch we
give here is based on~\cite{trimble:functor}.

\begin{proof}[Proof sketch]
  A category $\bC$ may equivalently be described as a monad on
  $\xX = \obs{\bC}$ in the bicategory $\Span$ of sets and spans. This
  is equivalently the 2-category of profunctors between discrete
  categories, or the 2-category of left adjoint functors
  $\Set^\xX \to \Set^\xY$. Thus a category amounts to a left adjoint
  monad on a category $\Set^\xX$, or equivalently a right adjoint
  comonad.

  The bicategory $\Span$ is moreover the underlying loose bicategory
  of a (weak) double category $\SSpan$ (in fact, an equipment ---
  see~\cite{cruttwell-shulman}) in which the loose arrows are spans
  and the tight arrows are functions. Taking monads and bimodules
  (a.k.a.\ the $\MMod$ construction;
  see~\cite{wood:ii,leinster:enrichment,cruttwell-shulman}) yields the
  (weak) double category of categories, profunctors and functors, and
  natural transformations. Under the correspondence of
  \cref{lem:categories}, categories are then identified with left
  adjoint monads on slices of $\Set$, functors between categories are
  identified with colax monad functors of the form
  $\leindex{\one{f}}$, and natural transformations are identified with
  specializations (definitions dual to
  \cref{def:comfun,def:comspec}). Taking mates of everything in sight,
  categories are equivalently right adjoint comonads on slices of
  $\Set$, functors between categories are colax comonad functors of
  the form $\reindex{\one{f}}$ in the opposite direction, and natural
  transformations are appropriate 2-cells between colax comonad
  functors. These translate to the above definitions of continuous map
  and specialization by~\cref{prop:formaltfae}.

  Now for the second statement, a functor $\bC \to \Set$, or
  profunctor between $\point$ and $\bC$, may be equivalently described
  as a (bi)module, carried by a span between $\point$ and $\obs{\bC}$,
  acted upon by the monad corresponding to $\bC$ in
  $\Span$. This translates to a left adjoint
  functor $\Set \to \Set^\xX$ acted upon as a $\xT$-module by the
  corresponding left adjoint monad $\xT$ on $\Set^\xX$. Noting that
  $\Set$ is the free cocompletion of $\point$, this amounts to simply
  a $\xT$-algebra $\point \to \Set^\xX$. Taken across the same
  translation process, natural transformations between functors
  $\bC\to \Set$ correspond to algebra maps, and so $\fun{\bC}{\Set}$
  is equivalent to the category of $\xT$-algebras. Moreover, the
  monadic functor associated to any left adjoint monad is equivalent
  to the comonadic functor associated to its right adjoint comonad.
\end{proof}

\begin{remark}
  If a comonadic functor into $\Set^X$ admits a left adjoint, then it
  preserves limits; therefore by \cref{prop:cationad} it corresponds
  to a category with object set $X$. The Yoneda lemma states that a
  free functor on an element at an object is a representable functor.
  Thus the comonads on $\Set^X$ corresponding to categories are
  precisely those admitting ``representable'' coalgebras, in the sense
  that the comonadic functor has a left adjoint.
\end{remark}

\begin{remark}\label{rem:formallimit}
  Since the limit-preserving comonads on $\Set^\xX$ correspond to
  categories with object set $\xX$, in particular, the
  limit-preserving comonads on $\Set$ --- equivalently, comonads
  carried by a representable endofunctor --- correspond to one-object
  categories, i.e.\ monoids.

  There is another perspective that makes the correspondence between
  monoids and representable comonads on $\Set$ more apparent. Note
  that the opposite of the category of small endofunctors on $\Set$
  constitutes the \emph{free completion} $\widetilde{\Set}$ of
  $\Set$. Moreover, composition $\circ$ behaves as cartesian product
  on the representable functors, i.e.
  $((\dash)^\xX)^\xY \cong (\dash)^{\xX \times \xY}$. Hence monoids in
  this monoidal category $\widetilde{\Set}$ may be understood as
  monoids carried by arbitrary formal limits of sets rather than just
  sets, and those which are carried by sets (representable functors)
  are ordinary monoids. Such generalized monoids, monoids carried by
  formal limits of sets, are equally well small comonads on $\Set$.

  We will return again to formal limits in \cref{sec:halos}.
\end{remark}

\section{Comonads on slices}\label[bgchapter]{sec:slices}

A topological space, or comonad on $2^X$ preserving finite meets, is
evidently analogous to an ionad, or comonad on $\Set^X$ preserving
finite limits. But another more general analogue is also available,
namely a comonad on $\Set^X$ preserving the terminal object and finite
intersections. In this section we will show that such comonads are
simply equivalent to arbitrary comonads on $\Set$ itself. And in
particular, ionads are identified with the pullback-preserving
comonads on $\Set$.

To this end, we will be concerned with the relationship between
comonads on a category and comonads on its slice categories; in the
case of $\Set$, since $\Set^\xX \simeq \Set/\xX$, we may work with
whichever is most convenient at any given time.


\begin{definition}\label{def:slice}
  Any functor $F \colon \bC \to \bD$ with colimit $\pts{F}$ admits a
  unique factorization
  \[\begin{tikzcd}
      \bC & {\bD/\pts{F}} & \bD \arrow["{\ptlift{F}}", from=1-1,
      to=1-2] \arrow["{\slicel}", from=1-2, to=1-3]
    \end{tikzcd}\] where the colimit of $\ptlift{F}$ is the terminal
  object $\pts{F}$ in $\bD/\pts{F}$ and $\slicel$ is the projection
  functor from the slice category. In this case, $\ptlift{F}$ sends 
  each object of $\bC$ to the corresponding leg of the colimiting cocone for $|F|$. 
  We call this the \emph{{\slice} factorization} of $F$.
\end{definition}

\begin{remark}
  If $\bC$ has a terminal object $1$, then any functor $F$ out of
  $\bC$ has colimit $\pts{F} \cong F(1)$. Conversely by
  \cref{lem:pts}, a comonadic functor $\carC \colon \CS \to \bS$ has a
  colimit $\pts{\cC}$ if and only if its category of coalgebras $\CS$
  has a terminal object. In particular this is the case whenever $\bS$
  itself has a terminal object, by applying the right adjoint
  $\rcarC$. Note also that there is no distinction between the
  colimits $\qpts{\cC}$ and $\qpts{\carC}$ because the right adjoint
  $\rcarC$ is final (see \cref{def:final}), as all right
  adjoints are.
\end{remark}

\begin{definition}\label{def:sliceable}
  We call a comonad $\cC$ on $\bS$ with colimit $\pts{\cC}$
  \emph{sliceable} if in the {\slice} factorization of the comonadic
  functor $\carC \colon \CS \to \bS$
  \[\begin{tikzcd}
      {\CS} & {\qSC} & \bS
      \arrow["\ptlift{\carC}", from=1-1, to=1-2]
      \arrow["{\slicel}", from=1-2, to=1-3]
    \end{tikzcd}\] $\ptlift{\carC}$ is a left adjoint. We denote the
  induced comonad on $\qSC$ by $\clift{\cC}$.\footnote{This is an
    abuse of notation, since $\clift{\cC} \colon \qSC \to \qSC$ is
    \emph{not} the functor $\bS \to \qSC$ induced by the {\slice}
    factorization of $\cC \colon \bS \to \bS$. However, we will never
    in this paper apply the {\slice} factorization to the carrier of a
    comonad, so confusion should not arise.} We call comonads on
  $\qSC$ of the form $\clift{\cC}$ \emph{sliced}.
\end{definition}

\begin{convention}\label{rem:indexedsliced}
  To minimize clutter, in the case of $\Set$ we will abusively use the
  same notation for constructions involving the equivalent but
  non-isomorphic categories $\Set/\xX$ and $\Set^\xX$. Given a comonad
  $\cC$ on $\Set$, we still denote the corresponding comonad
  (determined up to isomorphism) on $\SetpC$ as in
  \cref{def:sliceable} by $\clift{\cC}$; we also still call this
  comonad $\clift{\cC}$ a ``sliced comonad''. We still denote its
  category of coalgebras
  $\CoalgOn{\clift{\cC}}{(\SetpC)} \simeq \CSet$ by simply $\CSet$.
  Likewise given a comonad map $\phi \colon \cC \to \cD$ we still
  denote the corresponding functor between these categories by
  $\com{\Set}{\phi} \colon \CSet \to \DSet$.
\end{convention}

\begin{example}[Topological spaces and categories]\label{ex:topslice}\label{ex:catslice}
  Recall the comonad $\cX$ on $\Set$ corresponding to the topological
  space $\tX$ from \cref{ex:top}, whose coalgebras are sheaves
  $\Sh(\tX)$. The slice factorization of the comonadic sum of stalks
  functor $\Sh(\tX) \to \Set$ yields
  \[\Sh(\tX) \xto{i^*} \Set^{\pts{\tX}} \xto{\slicel} \Set\]
  where the first factor $i^*$ sends a sheaf to its
  $\pts{\tX}$-indexed set of stalks --- or equivalently, sends an
  \'etal\'e space to its $\pts{\tX}$-indexed set of fibers. The
  induced sliced comonad $\clift{\cX}$ is precisely the comonad on
  $\Set^{\pts{\tX}}$ defining the ionad corresponding to $\tX$ under
  \cref{prop:topionad}. That is, $i^*$ is the inverse image functor of
  the surjective geometric morphism equipping the topos $\Sh(\tX)$
  with the set of points $\pts{\tX}$.

  Whereas the set $\cX(\xA)$ consists of germs of functions
  $f \colon \pts{\tX} \to \xA$, the $\pts{\tX}$-indexed set
  $\clift{\cX}((\xA_x)_{x \in \pts{\tX}})$ consists of germs of
  sections $f \colon \pts{\tX} \to \sum_{x \in\pts{\tX}}\xA_x$ or
  families of elements $(a_x \in \xA_x)_{x \in \pts{\tX}}$:
  \[
    \clift{\cX}\left((\xA_x)_{x \in \pts{\tX}}\right) = \left(\colim_{U \ni x} \prod_{y \in U} \xA_y\right)_{x \in \pts{\tX}}
  \]
  where $U$ ranges over neighborhoods of $x$.
  
  Likewise, recall the comonad $\cbC$ on $\Set$ corresponding
  to the category $\bC$ from \cref{ex:cat}, whose coalgebras are given
  by $\fun{\bC}{\Set}$. The slice factorization of the comonadic
  functor $\obs{\El(\dash)} \colon \fun{\bC}{\Set} \to \Set$ yields
  \[\fun{\bC}{\Set} \xto{\reindex{i}} \Set^{\obs{\bC}} \xto{\slicel} \Set\]
  where the first factor $\reindex{i}$ reindexes a functor
  $\bC \to \Set$ along the inclusion of the set of objects
  $i \colon \obs{\bC} \to \bC$. The induced sliced comonad
  $\clift{\cbC}$ is precisely the comonad on $\Set^{\obs{\bC}}$
  defining the ionad corresponding to $\bC$ under
  \cref{prop:cationad}.

  Whereas the set $\cbC(\xA)$ consists of functions
  $\obs{c/\bC} \to \xA$ where $c$ is an object of $\bC$, the
  $\obs{\bC}$-indexed set $\clift{\cbC}((\xA_c)_{c \in \obs{\bC}})$
  is given by
  \[
    \clift{\cbC}\left((\xA_c)_{c \in \obs{\bC}}\right) = \left(\prod_{d \in \obs{\bC}} {\xA_d}^{\bC(c, d)}\right)_{c \in \obs{\bC}}.
  \]
\end{example}

\begin{lemma}\label{lem:slicecomonadic}
  If $\cC$ is sliceable, then $\ptlift{\carC} \colon \CS \to \SC$ is
  strictly comonadic, and $\clift{\cC}$ is {\final}.\footnote{More
    generally by the same argument, if $\cC$ is an arbitrary comonad,
    and the first factor in the \emph{comprehensive factorization} of
    $\carC$ is a left adjoint, then it is strictly comonadic, and the
    induced comonad is final. \cref{app:comprehensive} discusses the
    relationship between the comprehensive factorization system and
    slice factorizations.}
\end{lemma}
\begin{proof}
  The functor $\carC$ strictly creates limits that are absolute in
  $\bS$ by~\cref{prop:comonadicity}, and the functor $\slicel$
  strictly creates connected limits~\cite[Proposition
  3.4.8]{riehl}. All absolute (co)limits are connected, as can be seen
  from the fact that any constantly $1$ functor into $\Set$ preserves
  only trivial coproducts. Thus $\ptlift{\carC}$ strictly creates
  limits that are absolute $\qSC$, hence is strictly comonadic by
  \cref{prop:comonadicity}. The comonad $\clift{\cC}$ is {\final}
  because $\ptlift{\carC}$ and its right adjoint are.
\end{proof}

In fact, the sliced comonad $\clift{\cC}$ on $\bS/\pts{\cC}$ always
admits a more explicit description:

\begin{proposition}\label{rem:sliceableformula}
  A comonad $\cC$ is sliceable if and only if its colimit $\pts{\cC}$
  exists and for each object $q \colon \xA \to \pts{\cC}$ of
  $\bS/\pts{\cC}$, the category $\CS$ admits the equalizer
  \[
    \begin{tikzcd}[row sep=5pt,column sep=15pt]
      \eE && {\cC(\xA)} && {\cC(\pts{\cC})} \\
      &&& {\pts{\cC}}
      \arrow["e", dashed, hook, from=1-1, to=1-3]
      \arrow["{{\cC(q)}}", from=1-3, to=1-5]
      \arrow["{\bang_{\cC(\xA)}}"', from=1-3, to=2-4]
      \arrow["\hH_{\pts{\cC}}"', from=2-4, to=1-5]
    \end{tikzcd}
  \]
  where $\bang_{\cC(\xA)}$ denotes the unique coalgebra map to the
  terminal $\cC$-coalgebra and $\hH_{\pts{\cC}}$ denotes the terminal
  $\cC$-coalgebra structure. In this case
  $\clift{\cC}(q) \colon \eE \to \pts{\cC}$ is given by
  $\bang_{\eE} \coloneqq \bang_{\cC(\xA)} \circ e$.
\end{proposition}

\anoteNI{The right-to-left direction of this is mysteriously similar
  to an application of the adjoint triangle theorem.}

Note that if ($\bS$ has and) both $\cC$ and $\cC \circ \cC$ preserve
such equalizers, then the comonadic functor (preserves and) creates
them by \cref{lem:create}. Therefore in this case we may calculate
$\clift{\cC}$ entirely in terms of equalizers in $\bS$.


\begin{proof}
  For the left-to-right direction, suppose $\cC$ is sliceable. Let
  $\eE$ denote the cofree coalgebra on $q \colon \xA \to \pts{\cC}$,
  and denote by $e \colon \eE \to \cC(\xA)$ the unique coalgebra map
  such that $\varepsilon_\xA \circ e = \varepsilon_q$, where
  $\varepsilon_\xA$ is the $\cC$-counit at $\xA$ and $\varepsilon_q$
  is given by the $\clift{\cC}$-counit at $q$. We will show that $e$
  is the desired equalizer of $\cC(f)$ and
  $\hH_{\pts{\cC}}\circ \bang_{\cC(\xA)}$. We form the commutative
  diagram
  \[
    \begin{tikzcd}[column sep=35,row sep=25]
      & {\cC(\xA)} && {\cC(\pts{\cC})} & \\
      \eE && \xA && {\pts{\cC}} \\
      & {\cC(\xA)} & {\pts{\cC}} & {\cC(\pts{\cC})}
      \arrow[equals, from=3-3, to=2-5]
      \arrow["{{\cC(f)}}", from=1-2, to=1-4]
      \arrow["{{\varepsilon_\xA}}", from=1-2, to=2-3]
      \arrow["{{\varepsilon_{\pts{\cC}}}}", from=1-4, to=2-5]
      \arrow["e", from=2-1, to=1-2]
      \arrow["{{\varepsilon_q}}", from=2-1, to=2-3]
      \arrow["{{\bang_\eE}}"', curve={height=15pt}, from=2-1, to=2-5]
      \arrow["e"', from=2-1, to=3-2]
      \arrow["f", from=2-3, to=2-5]
      \arrow["{{\bang_{\cC(\xA)}}}"', from=3-2, to=3-3]
      \arrow["{{\hH_{\pts{\cC}}}}"', from=3-3, to=3-4]
      \arrow["{{\varepsilon_{\pts{\cC}}}}"', from=3-4, to=2-5]
    \end{tikzcd}
  \]
  using naturality of the $\cC$-counit $\varepsilon$, the defining
  property of $e$, the fact that $\varepsilon_q$ is a morphism over
  $\pts{\cC}$, uniqueness of the coalgebra map $\bang_\eE$, and the
  unit law for the terminal coalgebra $\pts{\cC}$. Thus $e$ indeed
  equalizes $\cC(f)$ and $\hH_{\pts{\cC}}\circ \bang_{\cC(\xA)}$,
  since composing with the counit $\varepsilon_{\pts{\cC}}$ induces a
  bijection between coalgebra maps $\eE \to \cC(\pts{\cC})$ and
  arbitrary maps $\eE \to \pts{\cC}$.

  Now suppose $e' \colon \eE' \to \cC(\xA)$ is another map equalizing
  $\cC(f)$ and $\hH_{\pts{\cC}}\circ \bang_{\cC(A)}$. We form the
  commutative diagram
  \[
    \begin{tikzcd}
      {\eE'} & {\cC(\xA)} & {\pts{\cC}}\\
      & \cC(\xA) && {\cC(\pts{\cC})} \\
      && \xA && {\pts{\cC}}
      \arrow[curve={height=-40pt}, equals, from=1-3, to=3-5]
      \arrow["{e'}"', from=1-1, to=1-2]
      \arrow["{\bang_{\eE'}}", curve={height=-25pt}, from=1-1, to=1-3]
      \arrow["{e'}"', from=1-1, to=2-2]
      \arrow["{\bang_{\cC(\xA)}}"', from=1-2, to=1-3]
      \arrow["{\hH_\pts{\cC}}"', from=1-3, to=2-4]
      \arrow["{\cC(f)}"', from=2-2, to=2-4]
      \arrow["{\varepsilon_\xA}"', from=2-2, to=3-3]
      \arrow["{\varepsilon_{\pts{\cC}}}"', from=2-4, to=3-5]
      \arrow["f"', from=3-3, to=3-5]
    \end{tikzcd}
  \]
  using uniqueness of the coalgebra map $\bang_{\eE'}$, the unit law for
  the terminal coalgebra $\pts{\cC}$, the assumed property of $e'$,
  and naturality of the $\cC$-counit $\varepsilon$. Thus
  $\varepsilon_\xA \circ e'$ defines a morphism over $\pts{\cC}$,
  therefore by the adjunction we obtain a unique coalgebra map
  $\hK \colon \eE\to \eE$ such that
  $\varepsilon_q\circ \hK = \varepsilon_\xA \circ e'$. But this corresponds
  by composing with $\varepsilon_\xA$ to the condition that
  $e \circ \hK = \varepsilon_\xA$. Therefore $e$ is the desired
  equalizer.
  
  For the right-to-left direction, assume $\CS$ has such
  equalizers $e \colon \eE \hookrightarrow \cC(\xA)$. We will show
  that $\eE$ defines the right adjoint $R \colon \SC \to \CS$ of
  $\ptlift{\carC} \colon \CS \to \SC$ at $q$. Then the comonad
  $\clift{\cC} \coloneqq \ptlift{\carC} \circ R$ sends $q$ to
  $\bang_\eE = \bang_{\cC(\xA)} \circ e$ as desired, since by
  definition $\ptlift{\carC}$ sends $\eE$ to the unique coalgebra map
  $\bang_\eE \colon \eE \to \pts{\cC}$.

  We check the adjunction condition, that is, the correspondence
  between maps $f\colon \ptlift{\carC}(\eX)\to q$ in $\bS/\pts{\cC}$
  and coalgebra maps $\hK\colon \eX\to \eE$, where $\eX$ is any
  $\cC$-coalgebra. Given such $f$, we form the commutative diagram
  
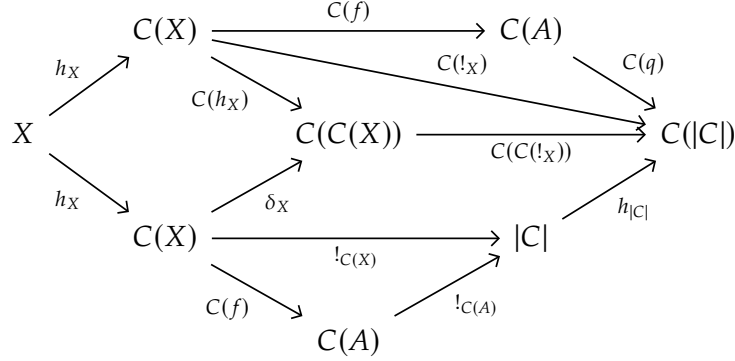
\[
    \begin{tikzcd}
      & {\cC(\eX)} && {\cC(\xA)} & \\
      \eX && {\cC (\cC(\eX))} && {\cC(\pts{\cC})} \\
      & {\cC(\eX)} && {\pts{\cC}} \\
      && {\cC(\xA)}
      \arrow["{{\cC(f)}}", from=1-2, to=1-4]
      \arrow["{{\cC(\hH_\eX)}}"', from=1-2, to=2-3]
      \arrow["{{\cC(\bang_\eX)}}", from=1-2, to=2-5]
      \arrow["{{\cC(q)}}", from=1-4, to=2-5]
      \arrow["{{\hH_\eX}}", from=2-1, to=1-2]
      \arrow["{{\hH_\eX}}"', from=2-1, to=3-2]
      \arrow["{{\cC (\cC(\bang_\eX))}}"', from=2-3, to=2-5]
      \arrow["{{\delta_\eX}}"', from=3-2, to=2-3]
      \arrow["{{\bang_{\cC(\eX)}}}"', from=3-2, to=3-4]
      \arrow["{{\cC(f)}}"', from=3-2, to=4-3]
      \arrow["{{\hH_{\pts{\cC}}}}"', from=3-4, to=2-5]
      \arrow["{{\bang_{\cC(\xA)}}}"', from=4-3, to=3-4]
    \end{tikzcd}
  \]
  using the defining triangle of $f$, uniqueness of the coalgebra map
  $\bang_\eX$, the associativity law for the
  coalgebra $\eX$, the
  coalgebra map law for $\bang_\eX$, and uniqueness of the coalgebra map $\bang_{\cC(\eX)}$.
  Therefore by the universal property of
  the equalizer $e$ in the category $\CS$ we get a unique coalgebra
  map $\hK\colon \eX\to \eE$ satisfying
  $e \circ \hK =\cC(f)\circ \hH_\eX$.

  In the other direction, given the coalgebra map
  $\hK\colon \eX\to \eE$, we define $f\colon \eX\to \xA$ as
  $\varepsilon_\xA\circ e\circ \hK$. To see that this gives a morphism
  over $\pts{\cC}$ we form the commutative diagram
  
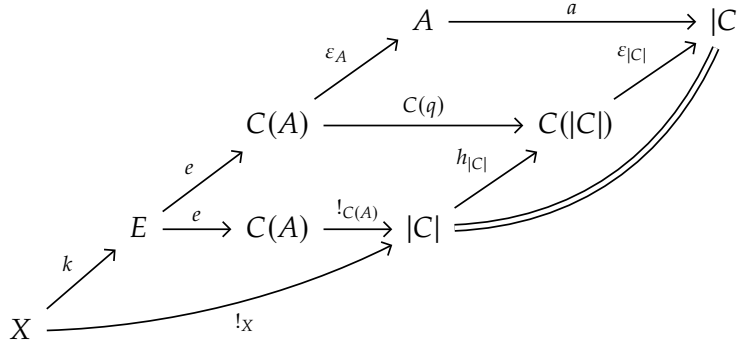
\[
    \begin{tikzcd}
      &&& \eA && {\pts{\cC}} \\
      && {\cC(\eA)} && {\cC(\pts{\cC})} \\
      & \eE & {\cC(\eA)} & {\pts{\cC}} \\
      \xX
      \arrow[curve={height=30pt}, equals, from=3-4, to=1-6]
      \arrow["a", from=1-4, to=1-6]
      \arrow["{{{\varepsilon_\xA}}}", from=2-3, to=1-4]
      \arrow["{{{\cC(q)}}}", from=2-3, to=2-5]
      \arrow["{{{\varepsilon_{\pts{\cC}}}}}", from=2-5, to=1-6]
      \arrow["e", from=3-2, to=2-3]
      \arrow["e", from=3-2, to=3-3]
      \arrow["{{{\bang_{\cC(\eA)}}}}", from=3-3, to=3-4]
      \arrow["{{{\hH_{\pts{\cC}}}}}", from=3-4, to=2-5]
      \arrow["k", from=4-1, to=3-2]
      \arrow["{{{\bang_\eX}}}"', curve={height=12pt}, from=4-1, to=3-4]
    \end{tikzcd}
  \]
  using naturality of the counit $\varepsilon$, the fact that $e$
  equalizes its fork, the unit law for the terminal coalgebra
  $\pts{\cC}$, and uniqueness of the coalgebra map $\bang_\eX$.
\end{proof}

\begin{remark}\label{rem:setsliceformula}
  Given a comonad $\cC$ on $\Set$ and a function
  $q \colon \xA \to \pts{\cC}$, we may explicitly describe the two
  coalgebra maps $\cC(q)$ and $\hH_{\pts{\cC}}\circ \bang_{\cC(\xA)}$
  of type $\cC(\xA) \to \cC(\pts{\cC})$ from
  \cref{rem:sliceableformula}, whose equalizer defines the domain of
  $\clift{\cC}(q)$, in terms of \emph{germs} of $\xA$-valued functions
  (as introduced in \cref{rem:densitygerms}). Namely, $\cC(q)$ is the map
  given by $[f]_x \mapsto [q \circ f]_x$, i.e.\ given by composing
  with $q$, and $\hH_{\pts{\cC}}\circ \bang_{\cC(\xA)}$ is the map
  sending any germ whose center lies over the point $x \in \pts{\cC}$
  to $\delta_1(x)$, the identity germ $[\id_{\pts{\cC}}]_x$ in the
  underlying topological space of $\cC$. Thus the equalizer
  $\clift{\cC}(q)$ consists of the $\xA$-valued germs lying over the
  identity germ via $q \colon \xA \to \pts{\cC}$ --- or rather the
  largest subcoalgebra of $\cC(\xA)$ included in this subset of germs,
  if the equalizer in $\CSet$ is not calculated as in $\Set$.

  This description indeed agrees with the description of the comonad
  $\clift{\cX}$ in \cref{ex:topslice}: a germ of an $\xA$-valued
  function that is sent to an identity germ by composing with
  $q \colon \xA \to \pts{\xX}$ is the germ of a section of $q$.
\end{remark}

In passing, we also note that sliceable comonads interact harmoniously
with the bases from \cref{sec:bases}.


\begin{lemma}\label{lem:liftcompatibility}
  For any $\pP \colon \bB \to \bS$ with colimit $\pts{\pP}$ and pointwise
  density comonad $\cP$, the lift
  $\ptlift{\pP} \colon \bB \to \bS/\pts{\pP}$ is the composite
  $\bB \xto{\lP} \PS \xto{\ptlift{\carP}}
  \bS/\pts{\pP}$.
\end{lemma}
\begin{proof}
  By \cref{lem:pts}, $\lP$ has terminal colimit. The functor
  $\ptlift{\carP}$ preserves the terminal object and all
  colimits. Thus
  \[\bB \xto{\lP} \CoalgS{\pP}
    \xto{\ptlift{\carP}} \bS/\pts{\pP} \xto{\slicel} \bS\]
  exhibits the slice factorization of $\pP$.
\end{proof}

\begin{corollary}\label{prop:liftbasis}
  Let $\pP \colon \bB \to \bS$ with sliceable density comonad $\cP$. If
  $\pP$ is a basis, then $\ptlift{\pP}$ provides a basis for
  $\clift{\cP}$, i.e.\ $\gen{\ptlift{\pP}} \cong \clift{\cP}$.
\end{corollary}
\begin{proof}
  By \cref{lem:liftcompatibility},
  $\ptlift{\pP} = \ptlift{\carP} \circ \lP$. Assuming $\pP$ is a
  basis, by \cref{cor:easybases}, $\gen{\ptlift{\pP}}$ is the comonad
  induced by $\ptlift{\carP}$, i.e.\ $\clift{\cP}$.
\end{proof}

Recall that we call a comonad \emph{crude} if it preserves coreflexive
equalizers, \cref{def:crep}. (Also recall from
\cref{lem:functorpreserve} that crudeness is the same as preservation
of regular finite intersections, so long as we merely assume that the
domain category has pushouts.)  Now given a category with finite
limits $\bS$, we identify crude comonads on $\bS$ with crude,
terminal-object-preserving comonads on slice categories of $\bS$.
  
\begin{lemma}\label{prop:factor}
  If $\cC$ is a crude comonad on a category $\bS$ with finite limits,
  then $\cC$ is sliceable, and $\clift{\cC}$ is crude.
\end{lemma}
\begin{proof}
  The slice category projection $\slicel \colon \SC\to \bS$ has a
  right adjoint given by $\pts{\cC} \times \dash$ and it strictly
  creates connected colimits. Thus it is strictly comonadic by
  \cref{prop:crude}. The category $\CS$ has coreflexive equalizers by
  \cref{lem:create}, and thus $\ptlift{\carC}$ is a left adjoint by
  the adjoint triangle theorem (\cref{lem:lifting}). Given that
  $\carC$ and $\carD$ create (and preserve) coreflexive equalizers, so
  does $\ptlift{\carC}$. Therefore the induced comonad $\clift{\cC}$
  preserves coreflexive equalizers as well.
\end{proof}

\begin{corollary}\label{cor:allslice}
  If $\cC$ is a comonad on $\Set$, then $\cC$ is sliceable, and
  $\clift{\cC}$ is crude.
\end{corollary}
\begin{proof}
  Follows from \cref{lem:corefl} and \cref{prop:factor}.
\end{proof}


Just as comonads may be lifted to slices, so too may comonad maps.

\begin{proposition}\label{prop:comtfae}
  The category consisting of sliceable comonads on $\bS$ and comonad
  maps is equivalent to the category consisting of sliced comonads on
  slice categories of $\bS$ and colax comonad functors
  (\cref{def:comfun}) of the form $\leindex{f}$ (where
  $f \colon \xX \to \xY$ is an arbitrary map in $\bS$ and
  $\leindex{f} \colon \bS/\xX \to \bS/\xY$ is the functor induced by
  composition with $f$).
\end{proposition}

More specifically, let $\cC$ and $\cD$ be sliceable comonads on $\bS$
with colimits $\pts{\cC}$ and $\pts{\cD}$, and let
$f \colon \pts{\cC} \to \pts{\cD}$ be an arbitrary arrow in
$\bS$. Then we claim that giving a comonad map
$\phi \colon \cC \to \cD$ with $\one{\phi} = f$ is the same as giving
a map $\one{\phi} \colon \pts{\cC} \to \pts{\cD}$ in $\bS$ and a colax
comonad functor from $\clift{\cC}$ to
$\clift{\cD}$ carried by $\leindex{f} \colon \SC \to \SD$.

Here $\one{\phi}$ is the map between colimits induced by the natural
transformation $\phi \colon \cC \to \cD$. In the case $\bS$ has a
terminal object $1$, this is the component of $\phi$ at $1$.

\begin{proof}
  By \cref{prop:respect} and \cref{prop:formaltfae}, a comonad map is
  equivalently specified by a functor $\com{\bS}{\phi}$ making the
  triangle below left commute, whereas a colax comonad functor
  structure on $\leindex{\one{\phi}}$ is equivalently specified by a
  functor $\com{\bS}{\phi}$ making the square below right commute.
  \[
    \begin{tikzcd}[column sep=15pt, row sep=15pt]
      \CS\ar[rr, "\com{\bS}{\phi}"]\ar[dr, "\carC"']&&\DS\ar[dl, "\carD"]\\
      &\bS
    \end{tikzcd}
    \qquad\qquad\qquad
    \begin{tikzcd}[column sep=20]
      \CS & \DS \\
      {\SC} & {\SD}
      \arrow["{\com{\bS}{\phi}}", from=1-1, to=1-2]
      \arrow["{\ptlift{\carC}}"', from=1-1, to=2-1]
      \arrow["{\ptlift{\carD}}", from=1-2, to=2-2]
      \arrow["{\leindex{\one{\phi}}}"', from=2-1, to=2-2]
    \end{tikzcd}
  \]
  
  Not just for comonadic functors, but for any functors into $\bS$
  with colimits, there is a correspondence between such commutative
  triangles and commutative squares.  Indeed, we may translate from
  squares to triangles like so:
  \[\begin{tikzcd}[column sep=10pt, row sep=10pt]
      \CS && \DS \\
      \\
      \SC && \SD \\
      & \bS
      \arrow["{\com{\bS}{\phi}}", from=1-1, to=1-3]
      \arrow["{\ptlift{\carC}}"', from=1-1, to=3-1]
      \arrow["{\ptlift{\carD}}", from=1-3, to=3-3]
      \arrow["{\leindex{\one{\phi}}}", from=3-1, to=3-3]
      \arrow["\slicel"', from=3-1, to=4-2]
      \arrow["\slicel", from=3-3, to=4-2]
    \end{tikzcd}
  \]
  Conversely, such a triangle induces a map
  $\one{\phi} \colon \pts{\cC} \to \pts{\cD}$ between the colimits, in
  fact a cocone in $\SD$ from $\ptlift{\carD}\circ \com{\bS}{\phi}$ to
  $\one{\phi}$ lifting the colimiting cocone from $\carC$ to
  $\pts{\cC}$. Hence the colimit of
  $\ptlift{\carD}\circ \com{\bS}{\phi}$ is $\one{\phi}$ because
  $\slicel \colon \SD \to \bS$ creates colimits. On the other hand the
  functor $\leindex{\one{\phi}}$ sends the terminal object
  ${\id}_{\pts{\cC}}$ in $\SC$ to $\one{\phi}$ in $\SD$. Therefore
  since $\slicel \colon \SC \to \bS$ creates and
  $\leindex{\one{\phi}}$ preserves colimits, the colimit of
  $\leindex{\one{\phi}}$ in $\SD$ is also $\one{\phi}$. The square
  then commutes because liftings of functors into $\bS$ along
  $\slicel \colon \SD \to \bS$ are determined by just lifting the
  colimit, as can be shown using that $\slicel$ is a discrete fibration
  (\cref{def:df}).
\end{proof}

\begin{corollary}\label{lem:liftedcategory}
  The full subcategory of $\Com(\bS/\xX)$ consisting of sliced
  comonads $\clift{\cC}$ is equivalent to the subcategory of
  $\Com(\bS)$ consisting of comonads with colimit $\xX$ and comonad
  maps $\phi$ for which $\one{\phi} \colon \one{\cC} \to \one{\cD}$ is
  identity on $\xX$.
\end{corollary}
\begin{proof}
  When $\one{\phi}$ is identity, the commutative squares of
  \cref{prop:comtfae} reduce to commutative triangles as in
  \cref{prop:respect}, i.e.\ comonad maps.
\end{proof}




We now come to the main takeaway of the section:

\begin{proposition}\label{prop:objective}
  Let $\bS$ be a category with finite limits. The category consisting
  of crude comonads on $\bS$ and comonad maps is equivalent to the
  category consisting of crude, terminal-object-preserving comonads on
  slice categories of $\bS$ and colax comonad functors of the form
  $\leindex{f}$.
\end{proposition}
\begin{proof}
  We have a one-to-one correspondence between crudely and strictly
  comonadic functors into $\bS$ and terminal-object-preserving,
  crudely and strictly comonadic functors into slice categories of
  $\bS$. The left-to-right direction of this correspondence is given
  by composing the comonadic functors, where the result is comonadic
  by the crude comonadicity theorem (\cref{lem:crudecompose}), and the
  right-to-left direction is given by \cref{prop:factor}. The
  equivalence of categories is by \cref{prop:comtfae}.
\end{proof}

\begin{corollary}\label{cor:setslices}
  The category $\Com(\Set)$ of comonads on $\Set$ and comonad maps is
  equivalent to the category of sliced comonads on categories of the
  form $\Set^\xX$ (or equivalently categories of the form $\Set/\xX$),
  and colax comonad functors of the form $\leindex{f}$.

  Moreover, here the sliced comonads are precisely the comonads that
  preserve the terminal object and finite intersections.
\end{corollary}
\begin{proof}
  Follows from \cref{prop:objective}, since by
  \cref{lem:corefl} all comonads on $\Set$ are crude, and by
  \cref{lem:functorpreserve} crude comonads on $\Set$ and its slices
  are equivalent to those preserving finite intersections.
\end{proof}

Now we see that ionads are identified with pullback-preserving
comonads on $\Set$ itself.

\begin{proposition}\label{prop:ionads}
  Let $\bS$ be a category with finite limits. The equivalence of
  \cref{prop:objective} restricts to an equivalence between
  pullback-preserving comonads on $\bS$ and finite-limit-preserving
  comonads on slice categories of $\bS$.
\end{proposition}
\begin{proof}
  In light of \cref{prop:objective}, we just need to show that $\cC$
  preserves pullbacks if and only if $\clift{\cC}$ preserves finite
  limits. Indeed, $\clift{\cC}$ preserves the terminal object, and
  since $\slicel$ creates connected limits $\clift{\cC}$ preserves
  pullbacks if and only if $\cC$ does, since then
  $\carC \colon \CS \to \bS$ does by \cref{lem:create}.
\end{proof}

\begin{corollary}\label{cor:ionads}
  Pullback-preserving comonads on $\Set$ are identified with
  finite-limit-preserving comonads on slice categories of $\Set$,
  i.e.\ ionads.
\end{corollary}

\begin{proof}
  Follows from \cref{prop:ionads}.
\end{proof}

\begin{warning}\label{rem:notcontinuous}
  Comonad maps between pullback-preserving comonads on $\Set$ are
  \emph{not} the same as continuous maps of ionads
  (\cref{def:ionadcontinuous}); \cref{cor:ionads} does not refer to
  continuous maps at all. Nevertheless, it is clear upon inspection
  that the isomorphisms in the category of comonads on $\Set$ and
  comonad maps do contravariantly correspond to the isomorphisms in
  the category of ionads and continuous maps, since if $f$ is
  bijective then $\leindex{f} \cong \reindex{f^{-1}}$. Hence this does
  establish an equivalence of \emph{groupoids} between
  pullback-preserving comonads on $\Set$ with invertible comonad maps
  and ionads with invertible continuous maps.  In \cref{sec:set} we
  will further dissolve the discrepancy by studying continuous maps
  between arbitrary comonads on $\Set$ as well.
\end{warning}

As noted in \cref{ex:topslice}, the translation process from comonads
on categories $\Set^X$ to comonads on $\Set$, given by composing with
the adjunction $\slicel \dashv \slicer$, takes the ionad corresponding
to a topological space $\tX$ to the comonad $\cX$ on $\Set$
corresponding to $\tX$ from \cref{ex:top}, and takes the ionad
corresponding to a small category $\bC$ to the comonad $\cbC$ on
$\Set$ corresponding to $\bC$ from \cref{ex:cat}.

By \cref{prop:cationad}, the ionads corresponding to categories are
the limit-preserving comonads on $\Set^X$. Now we see that the
corresponding comonads $\cbC$ on $\Set$ are precisely the polynomial,
or equivalently wide-pullback-preserving, or equivalently
connected-limit-preserving comonads, recovering the result
of~\cite{ahman-uustalu} that categories are identified with polynomial
comonads on $\Set$.

\begin{proposition}\label{prop:alex}
  Let $\bS$ be a category with finite limits. The equivalence of
  \cref{prop:objective} restricts to an equivalence between
  wide-pullback-preserving comonads on $\bS$ and limit-preserving
  comonads on slice categories of $\bS$.
\end{proposition}
\begin{proof}
  Entirely analogous to \cref{prop:ionads}.
\end{proof}

\begin{corollary}[{\cite{ahman-uustalu}}]\label{cor:au}
  Polynomial comonads on $\Set$ are identified with small categories.
\end{corollary}
\begin{proof}
  Follows from \cref{lem:categories} and \cref{prop:alex}.
\end{proof}

\begin{warning}\label{rem:notfunctors}
  Just as in \cref{rem:notcontinuous}, comonad maps between polynomial
  comonads on $\Set$ are \emph{not} the same as functors of
  categories. But again, we do have an equivalence of \emph{groupoids}
  between the groupoid of polynomial comonads on $\Set$ and invertible
  comonad maps and the groupoid of categories and isomorphisms.
\end{warning}

\begin{remark}\label{rem:yesgo}
  Unlike those corresponding to categories and ionads, the comonads
  $\cX$ on $\Set$ corresponding to topological spaces $\tX$
  \emph{cannot} be characterized solely in terms of some property on
  the underlying endofunctor on $\Set$; see \cref{cex:nogo}.

  However, we \emph{can} characterize the corresponding comonads
  $\clift{\cX}$ in terms of the underlying endofunctor on
  $\Set^X$. Specifically, they are those sliced comonads given
  componentwise by colimits of representable functors of the form
  $\Ho{\Set^X}(\xA, \dash)$ where $\xA$ is a subterminal object in
  $\Set^X$. Indeed, such functors $F$ have the property that for every
  object $\xS$, the maps $F(f) \colon F(\xA) \to F(\xS)$ with
  subterminal domain are jointly surjective. If $\clift{\cC}$ has this
  property, then every cofree $\clift{\cC}$-coalgebra is covered by
  regular subterminal objects, since $\clift{\cC}$ is the composite of
  the cofree functor $\rcar{\clift{\cC}} \colon \SetpC \to \SetC$ and
  the comonadic functor $\car{\clift{\cC}} \colon \SetC \to \SetpC$,
  and (regular) subterminal objects are preserved by any sliced
  $\clift{\cC}$. The result then follows from \cref{thm:tfaetop},
  since as always all coalgebras are covered by cofree coalgebras.
  Conversely, for any topological space $\tX$, the sliced comonad
  $\clift{\tX}$ is the density comonad of the induced functor
  $\ptlift{\topbasis{\tX}} \colon \opens{\tX} \to \SetpX$ by
  \cref{prop:liftbasis}, which is a colimit of the desired form on
  inspection (giving the formula from \cref{ex:topslice}).
\end{remark}

\begin{convention}\label{rem:cats}
  Take note that in both~\cite{garner:ionads} and this paper, the
  category $\bC$ is identified with a comonad whose coalgebras are
  given by $\fun{\bC}{\Set}$, not
  $\Ps{\bC} \coloneqq \fun{\bC\op}{\Set}$ as might be expected by
  analogy to sheaves over topological spaces. To keep these
  distinctions clear, we denote the comonad on $\Set$ whose coalgebras
  are equivalent to functors $\bC \to \Set$ (copresheaves on $\bC$) by
  $\cbC$, whereas we denote the comonad on $\Set$ whose coalgebras are
  equivalent to sheaves on the topological space $\tX$ by $\cX$. This
  leaves room to denote the comonad on $\Set$ whose coalgebras are
  equivalent to presheaves on the category $\bC$ by $\cattocmdp{\bC}$.
  Thus we have
  \[
    \CoalgSet{\cbC} \simeq \fun{\bC}{\Set}
    \qqand
    \CoalgSet{\cattocmdp{\bC}} \simeq \Ps{\bC}
    \qqand
    \CoalgSet{\cX} \simeq \Sh(\tX).
  \]
  Identifying ionads with pullback-preserving comonads on $\Set$, the
  embedding $\Top \hookrightarrow \Ionad$ is given by
  $\tX \mapsto \cX$, whereas the embedding
  $\Cat \hookrightarrow \Ionad$ is given by
  $\bC \mapsto \cbC$.\footnote{If one wanted to take the opposite
    convention, identifying $\bC$ with the comonad $\cattocmdp{\bC}$
    whose coalgebras are $\bC$-presheaves, one would then also want to
    redefine specializations (\cref{def:ionadspecialization}) to go
    the other way to agree with natural transformations. This is not
    so unreasonable, since the direction of the specialization
    preorder is already inconsistent in the literature on topological
    spaces.}
\end{convention}


\chaptertocspace
\chapter{Continuous maps}\label{sec:set}

A continuous map between topological spaces is a function
$\one{f} \colon \pts{\tX} \to \pts{\tY}$ such that taking inverse
images sends opens to opens. Continuous maps of comonads are similar,
replacing the opens of the topological spaces with the coalgebras of
the comonads. In this section we consider continuous maps between
comonads on $\Set$, as well as their interaction with ordinary comonad
maps.

In \cref{sec:continuousset} we define the 2-category $\ConSet$ of
comonads on $\Set$, continuous maps, and specializations. We
generalize some results from \cite{garner:ionads} about continuous
maps of ionads. We show that the functor $\ConOSet \to \Set$ sending a
comonad to its set of points is a fibration; this induces a
factorization system on comonads and continuous maps that generalizes
the (bijective-on-objects, fully-faithful) factorization system on
$\Cat$ (\cref{rem:booff}). We show that topological spaces $\Top$ form
a reflective sub-2-category of $\ConO$ (\cref{thm:topreflection}). We
also show that every comonad on $\Set$ has an underlying small
category (\cref{cor:smallcategory}), and that small categories $\Cat$
form a coreflective sub-2-category of $\ConSet$
(\cref{thm:catcoreflection}).

In \cref{sec:doubleset} we introduce a double category whose objects
are comonads on $\Set$ and whose two types of arrows are comonad maps
and continuous maps. In fact there are two such double categories: the
cells of the first double category $\CmdSet$ are generalized
specializations; the second double category $\CmdIdSet$ is a canonical
thin sub-double-category of the first, whose cells we call
\emph{identifications}. An identification expresses a compatibility
condition between comonad maps and continuous maps. We show that the
full sub-double-categories of $\CmdSet$ and $\CmdIdSet$ consisting of
polynomial comonads recover the two known double categories of
categories, retrofunctors, and functors from~\cite{clarke}
and~\cite{clarke-dimeglio} (\cref{thm:doublecat}). We re-find the
category of coalgebras of each comonad on $\Set$ within the double
category $\CmdIdSet$ (\cref{rem:doublecoalgebra}). We show that the
conjoint pairs in $\CmdIdSet$ are precisely the bijective-on-objects
continuous maps and comonad maps, and the companion comonad maps in
$\CmdIdSet$ are precisely the cartesian comonad maps
(\cref{prop:companion}).
We show that the thin double category of ionads has strong tabulators
of comonad maps (\cref{thm:strongtabulators}).

In \cref{sec:limits} we study limits and colimits of comonads on
$\Set$. The category of small comonads on $\Set$ and comonad maps is
both complete and cocomplete, and the category of arbitrary comonads
on $\Set$ and comonad maps is cocomplete (\cref{prop:comlimits}), for
general well-known reasons. We furthermore show that the category of
small comonads on $\Set$ and continuous maps is both complete and
cocomplete, and the category of arbitrary comonads on $\Set$ and
continuous maps is complete
(\cref{prop:contcocomplete,thm:contcomplete}). We also briefly discuss
a monoidal product $\otimes$ on $\ConO$ generalizing the usual product
$\times$ of topological spaces (\cref{rem:dirichlet}).

\section{Continuous maps}\label{sec:continuousset}

We saw in \cref{sec:slices} that comonads on $\Set$ are equivalent to
comonads on categories $\Set^X$ preserving the terminal object and
finite intersections. The definition of continuous map of ionads
(\cref{def:ionadmap}) now readily generalizes to arbitrary comonads on
$\Set$, as does the definition of specialization
(\cref{def:ionadspecialization}).

\begin{definition}\label{def:continuousset}\label{def:specializationset}
  A \emph{continuous map} $f \colon \cC \to \cD$ of comonads on $\Set$
  consists of a function $\one{f} \colon \pts{\cC} \to \pts{\cD}$ and
  a functor $\con{\Set}{f} \colon \DSet \to \CSet$ making the following
  diagram commute
  \[
    \begin{tikzcd}
      \CSet & \DSet \\
      {\SetpC} & {\SetpD}
      \arrow["\car{\clift{\cC}}"', from=1-1, to=2-1]
      \arrow["{\con{\Set}{f}}"', dashed, from=1-2, to=1-1]
      \arrow["\car{\clift{\cD}}", from=1-2, to=2-2]
      \arrow["{\reindex{\one{f}}}", from=2-2, to=2-1]
    \end{tikzcd}
  \]
  where $\car{\clift{\cC}}$ and $\car{\clift{\cD}}$ denote the
  strictly\footnote{Recall from \cref{rem:indexedsliced} that we use
    the mildly abusive notation $\CSet$ for the category of coalgebras
    of the induced comonad $\clift{\cC}$ on $\SetpC$ --- which is
    equivalent, but not in general isomorphic, to the category
    of coalgebras of $\cC$ on $\Set$.} comonadic functors induced by $\cC$ and
  $\cD$ as in \cref{def:sliceable} and $\reindex{\one{f}}$ is the
  reindexing functor given by a function
  $\one{f} \colon \pts{\cC} \to \pts{\cD}$.

  A \emph{specialization} between continuous maps of comonads
  $f \Rightarrow g$ is a natural transformation
  $\con{\Set}{f} \Rightarrow \con{\Set}{g}$.
\end{definition}

\begin{remark}\label{rem:laxtriangle}
  There is also a way of defining continuous maps of comonads directly
  in terms of the comonadic functors into $\Set$ itself. Recall that a
  natural transformation is \emph{cartesian} if all naturality squares
  are pullback squares. A continuous map $f \colon \cC \to \cD$ as
  above may be equivalently defined as a terminal-object-preserving
  functor $\con{\Set}{f}$ equipped a cartesian natural transformation
  $\connat{f}$
  \[
    \begin{tikzcd}[column sep=15pt]
      \CSet && \DSet \\
      & \Set
      \arrow["{\con{\Set}{f}}"', from=1-3, to=1-1]
      \arrow["\carC"', ""{name=0, anchor=center, inner sep=0}, from=1-1, to=2-2]
      \arrow["\carD", ""{name=1, anchor=center, inner sep=0}, from=1-3, to=2-2]
      \arrow["\Rightarrow_{\connat{f}}", shift right=10pt, draw=none, from=1, to=0]
    \end{tikzcd}
  \]
  where such morphisms are identified up to natural isomorphism on
  the functors $\con{\Set}{f}$ commuting with the triangles. The
  terminal-object-preservation allows $\carC \circ \con{\Set}{f}$ to
  be lifted to $\SetpC$, and the cartesian condition expresses that
  $\con{\Set}{f}$ lifts the reindexing functor
  $\reindex{\one{f}} \colon \SetpD \to \SetpC$ as in
  \cref{def:continuousset}. It is this definition of continuous map
  that best generalizes to arbitrary comonads on arbitrary categories;
  see \cref{sec:maps} for details.
\end{remark}

Density comonads of functors into $\Set$ satisfy a universal property
with respect to continuous maps analogous to their universal property
from \cref{prop:density} with respect to comonad maps.

\begin{restatable}{theorem}{thmuniversal}\label{thm:universal}
  For any $\pP \colon \bB \to \Set$ with density comonad $\cP$, we
  have a bijection between continuous maps of comonads $\cC \to \cP$
  (squares shown left) and squares shown right:
  
  \[
    \begin{tikzcd}[column sep=20, ampersand replacement=\&]
      \CSet \& \PSet \\
      {\SetpC} \& {\SetpP}
      \arrow["\car{\clift{\cC}}"', from=1-1, to=2-1]
      \arrow[dashed, from=1-2, to=1-1]
      \arrow["\car{\clift{\cP}}", from=1-2, to=2-2]
      \arrow["{\reindex{\one{f}}}", from=2-2, to=2-1]
    \end{tikzcd}
    \qquad\qquad\qquad
    \begin{tikzcd}[column sep=17.5, ampersand replacement=\&]
      \CSet \& \bB \\
      {\SetpC} \& {\SetpP}
      \arrow["\car{\clift{\cC}}"', from=1-1, to=2-1]
      \arrow[dashed, from=1-2, to=1-1]
      \arrow["\ptlift{\pP}", from=1-2, to=2-2]
      \arrow["{\reindex{\one{f}}}", from=2-2, to=2-1]
    \end{tikzcd}
  \]
\end{restatable}

Intuitively, this generalizes the fact that for topological spaces,
when one verifies continuity, it is only necessary to check that
inverse images of subbasic opens are open.
In fact, there are also similar results for specializations, which we
give in \cref{app:universal}.

\begin{proof}
  Deferred to \cref{app:universal}.
\end{proof}


\begin{example}
  To better understand the general behavior of continuous maps between
  comonads on $\Set$, it is illustrative to consider not just
  continuous maps of topological spaces or continuous maps of
  categories, but continuous maps between the two. If $\tX$ is a space
  and $\bG$ is a group (viewed as a one-object category), continuous
  maps from $\cX$ to $\cattocmd{\bG}$ are precisely principal
  $\bG$-bundles.

  Indeed, \cref{thm:universal} tells us precisely this. According to
  \cref{ex:catbasis}, the ``canonical basis'' for $\bG$ is given by
  $\pP \colon \bG \to \Set$ specifying the action of $\bG$ on itself,
  i.e.\ a $\bG$-torsor. Now a square
  \[
    \begin{tikzcd}[column sep=17.5]
      \Sh(\tX) & \bG \\
      {\Set^{\pts{\tX}}} & {\Set} \arrow["\car{\cX}"',
      from=1-1, to=2-1] \arrow[dashed, from=1-2, to=1-1]
      \arrow["\pP", from=1-2, to=2-2]
      \arrow["{\slicer}", from=2-2, to=2-1]
    \end{tikzcd}
  \]
  picks out a sheaf over $\tX$, or equivalently \'etal\'e space over
  $\tX$, whose fibers are given by $\bG$, and with $\bG$-action given
  fiberwise by torsors.

  A continuous map from $\cX$ to the density comonad $\cP$ of an
  arbitrary functor $\pP \colon \bB \to \Set$ is similar. Again by
  \cref{thm:universal}, such a continuous map $f \colon \cX \to \cP$
  is equivalently a functor $F \colon \bB \to \Sh(\tX)$ given
  fiberwise by reindexing the functor $\pP$ along some
  $\pts{f} \colon \pts{\tX} \to \pts{\pP}$. Specifically, over each
  point $x$ of $\tX$, we require that the action on $x$-fibers
  $F_x \colon \bB \to \Set$ is some summand of $\pP$ (i.e. a connected
  component of $\El(\pP)$, lying over some point of
  $\pts{\cP}$). Still more generally, continuous maps out of arbitrary
  comonads on $\Set$ may be described in the same way, where general
  coalgebras replace the \'etal\'e spaces.
\end{example}

We note that continuous maps between comonads on $\Set$ retain some of
the flavor of geometric morphisms of toposes:


\begin{lemma}\label{lem:conleft}
  If $f$ is a continuous map between comonads on $\Set$, then
  $\con{\Set}{f}$ is a left adjoint.
\end{lemma}
\begin{proof}
  Follows from the adjoint triangle theorem (\cref{lem:lifting}).
\end{proof}

\begin{remark}\label{rem:concomonadic}
  Here $\con{\Set}{f}$ is itself comonadic if and only if it is
  conservative. Indeed, $\con{\bS}{f}$ preserves coreflexive
  equalizers, since the comonadic forgetful functors create them and
  $\reindex{\one{f}}$ preserves them, so the statement follows by the
  crude comonadicity theorem (\cref{prop:crude}).
  
  A sufficient condition for this is that $\one{f}$ be
  surjective. Indeed, then $\reindex{\one{f}}$ is conservative, so
  $\con{\Set}{f}$ is as well by cancellation. (However, surjectivity
  of $\one{f}$ is not necessary for $\con{\Set}{f}$ to be comonadic,
  since for example any equivalence between categories, whose
  corresponding comonads on $\Set$ have coalgebra categories
  $\Set^{\bC}$, induces an equivalence between the categories
  $\Set^{\bC}$.)
\end{remark}




We now observe that continuous maps between comonads on $\Set$ may
also be viewed as certain ordinary comonad maps.


\begin{lemma}\label{lem:coniscom}
  A continuous map $f \colon \cC \to \cD$ between comonads on $\Set$
  is equivalently specified by a function
  $\one{f} \colon \pts{\cC} \to \pts{\cD}$ and an identity-on-points
  comonad map from $\one{f}^*{\cD}$ to $\cC$, where $\one{f}^*{\cD}$ denotes
  the comonad $\slicel \circ \reindex{\one{f}} \circ \clift{\cD}
  \circ \Pi_{\one{f}} \circ \slicer$
  \[\begin{tikzcd}
      \SetpD & \SetpC & \Set
      \arrow["{\clift{\cD}}", from=1-1, to=1-1, loop, in=160, out=200, distance=10mm]
      \arrow["{\reindex{\one{f}}}", shift left=2, from=1-1, to=1-2]
      \arrow["{\Pi_{\one{f}}}", shift left=2, from=1-2, to=1-1]
      \arrow["\bot"{anchor=center}, draw=none, from=1-1, to=1-2]
      \arrow[""{name=0, anchor=center, inner sep=0}, "\slicel", shift left=2, from=1-2, to=1-3]
      \arrow["\bot"{anchor=center}, draw=none, from=1-2, to=1-3]
      \arrow[""{name=1, anchor=center, inner sep=0}, "\slicer", shift left=2, from=1-3, to=1-2]
      \arrow[draw=none, from=0, to=1]
    \end{tikzcd}
  \]
\end{lemma}
\begin{proof}
  The defining square of a continuous map may be viewed as a triangle:
  \[\begin{tikzcd}[column sep=15pt]
      \CSet && \DSet \\
      & \SetpC
      \arrow["\carC"', from=1-1, to=2-2]
      \arrow["{\con{\Set}{f}}"', from=1-3, to=1-1]
      \arrow["\reindex{\one{f}} \circ \carD", from=1-3, to=2-2]
    \end{tikzcd}\] By \cref{prop:density} and \cref{cor:densityofleft}
  this is equivalently a comonad map from the induced comonad
  $\reindex{\one{f}} \circ \clift{\cD} \circ \Pi_{\one{f}}$ to
  $\clift{\cC}$.  Both comonads preserve the terminal object and
  finite intersections, so by \cref{lem:liftedcategory} this is
  equivalently an identity-on-points comonad map between the
  corresponding comonads on $\Set$.
\end{proof}

\begin{proposition}\label{prop:fibration}
  The functor $\one{\dash} \colon \ConOSet \to \Set$ is a
  fibration.
\end{proposition}
\begin{proof}
  \cref{lem:coniscom} expresses the condition that the canonical
  continuous maps $\one{f}^*{\cD} \to \cD$ are weakly cartesian. They are
  also closed under composition, so they are cartesian.
\end{proof}

To put it another way, $\one{\dash}$ is the fibration given by
Grothendieck construction of the functor $\Set \to \CAT$ that sends a
set $\xY$ to the category of sliced (i.e.\
terminal-object-and-finite-intersection-preserving) comonads on
$\Set^\xY$, and sending $f \colon \xX \to \xY$ to
$\reindex{f} \circ \dash \circ \Pi_{f}$.


\begin{remark}\label{rem:booff}
  The induced factorization system on $\ConOSet$ specializes to the
  factorization system on $\Cat$ where the two classes of maps are the
  bijective-on-objects functors and fully faithful functors, and to the
  factorization system on $\Top$ where the two classes of maps are
  bijective continuous maps and continuous maps such that the domain is
  endowed with the coarsest possible topology making the map
  continuous.
\end{remark}

In topology, the terminal space (the point) $\point$ and the
Sierpinski space $\sier$ play conceptually dual roles: continuous maps
out of $\point$ classify points, whereas continuous maps into $\sier$
classify open subsets. Essentially the same is true for general comonads on
$\Set$.

For convenience let us denote the identity comonad on $\Set$, which
corresponds to the terminal topological space or category, by
$\pointcmd \coloneqq \idSet \cong \toptocmd{\point} \cong
\cattocmd{\point}$, and let us similarly denote the comonad on $\Set$
corresponding to the Sierpinski space or walking arrow category by
$\arrocmd \coloneqq \toptocmd{\sier} \cong \cattocmd{\arro}$.

\begin{lemma}
  Let $\cC$ be a comonad on $\Set$. The frame of opens
  $\opens{\cmdtotop{\cC}}$ of the underlying topological space
  $\cmdtotop{\cC}$ is the (thin) category
  $\Ho{\ConSet}(\cC, \arrocmd)$.
\end{lemma}
\begin{proof}
  Note that when $f$ is a continuous map $\cC \to \cD$, the inverse
  image functor $\con{\Set}{f} \colon \DSet \to \CSet$ preserves
  colimits, the terminal object, and its regular subobjects, i.e.\
  open subsets as in \cref{def:open}. Since $\CoalgSet{\arrocmd}$ is
  the free cocompletion of the walking arrow, a continuous map
  $\cC \to \arrocmd$ is therefore determined by an open subset
  (alternatively, this can also be seen with \cref{thm:universal},
  using that $\sier$ is the density comonad of a diagram in $\Set$
  picking out an single arrow $1 \to 2$); likewise specializations are
  given by natural transformations of diagrams of shape $\arro$, which
  are open subset inclusions.
\end{proof}

\begin{theorem}\label{thm:topreflection}
  $\Top$ is a reflective sub-2-category of $\ConSet$. The reflection
  sends a comonad to its underlying topological space.
\end{theorem}
\begin{proof}
  The unit of the reflection at $\cC$ is the bijective-on-points
  continuous map $\cC \to \cmdtotop{\cC}$ corresponding as in
  \cref{prop:boo} to the bijective-on-points comonad map
  $\cmdtotop{\cC} \to \cC$ that is given by the fully faithful
  embedding $\Sh(\cmdtotop{\cC}) \hookrightarrow \CSet$ from
  \cref{prop:sheavesff}. If $f \colon \cC \to \cX$ is a
  continuous map, then $\con{\Set}{f} \colon \Sh(\tX) \to \CSet$ preserves
  regular subterminal objects and their colimits. Thus by
  \cref{lem:bangetale} it factors through this fully faithful
  embedding. Hence
  $\Ho{\Top}(\cmdtotop{\cC}, \tX) \cong \Ho{\ConSet}(\cC,
  \cX)$ as desired.
\end{proof}

We already introduced the underlying topological space of a comonad on
$\Set$ in~\cref{def:underlyingspace}, but there is also a conceptually
dual notion of underlying category.

\begin{definition}\label{def:undercat}
  The \emph{underlying category} of a comonad $\cC$ on $\Set$ is the
  category $\Ho{\ConSet}(\pointcmd, \cC)$ of continuous maps from the
  identity comonad into $\cC$ and specializations.
\end{definition}

Note that a continuous maps from $\pointcmd$ to $\cC$ is simply a
point in $\pts{\cC}$. Hence the objects of the underlying category of
$\cC$ are given by the set $\pts{\cC}$.

We now show that the underlying category of an arbitrary comonad on
$\Set$ is small.

\begin{lemma}[{\cite[Corollary 2.3]{adamek-sousa}}]\label{lem:setdensity}
  If $\pP \colon \bB \to \Set$ admits a density comonad, then there is
  a small set of natural transformations
  $\pP \Rightarrow \pP$.\footnote{Here $\Set$ may be replaced by any
    category with a coseparating set.}\footnote{The reader is warned
    that there is an error in the formula for density comonads on
    $\Set$ from~\cite[Corollary 4.3]{adamek-sousa}:
    $CX = \mathbf{Nat}(X^F, 2^F)$ should instead read
    $\Set(CX, 2) = \mathbf{Nat}(X^F, 2^F)$.}\qed
\end{lemma}

\begin{corollary}\label{lem:setofnatural}
  If $F, G \colon \bB \to \Set$ admit density comonads, then there is
  a small set of natural transformations $F \Rightarrow G$.
\end{corollary}
\begin{proof}
  If $F$ and $G$ admit density comonads (which are necessarily
  pointwise by \cref{lem:setpointwise}), then the coproduct functor
  $F + G$ does also by the colimit formula for density comonads:
  specifically, $\gen{F + G} \cong \gen{F} + \gen{G}$. The natural
  transformations from $F$ to $G$ are identified with a subset of the
  natural transformations from $F + G$ to itself (namely those sending
  $F$ to $G$ and which are identity on $G$), so the result follows
  from \cref{lem:setdensity}.
\end{proof}

\begin{lemma}\label{prop:locallysmall}
  $\ConSet$ is locally locally small.
\end{lemma}
\begin{proof}
  Suppose $f, g \colon \cC \to \cD$ are continuous maps between
  comonads on $\Set$. By \cref{lem:conleft}, $\con{\Set}{f}$ and $\con{\Set}{g}$
  are left adjoint. Therefore $\carC \circ \con{\Set}{f}$ and
  $\carC \circ \con{\Set}{g}$ have density comonads by
  \cref{cor:densityofleft}, and thus there is a small set of natural
  transformations
  $\carC \circ \con{\Set}{f} \Rightarrow \carC \circ \con{\Set}{g}$ by
  \cref{lem:setofnatural}. Since $\carC$ is faithful we conclude there
  is a small set of natural transformations
  $\con{\Set}{f} \Rightarrow \con{\Set}{g}$, i.e.\ specializations from $f$ to
  $g$.
\end{proof}

\begin{corollary}\label{cor:smallcategory}
  Every comonad on $\Set$ has a small underlying category.\qed
\end{corollary}

\begin{theorem}\label{thm:catcoreflection}
  $\Cat$ is a coreflective sub-2-category of $\ConSet$. The
  coreflection sends a comonad to its underlying category.
\end{theorem}
\begin{proof}
  The claim is that $\cattocmd{(\dash)} \colon \Cat \to \ConSet$ is
  left 2-adjoint to
  $\Ho{\ConSet}(\pointcmd, \dash) \colon \ConSet \to \Cat$. Indeed, a
  functor $F \colon \bC \to \Ho{\ConSet}(\pointcmd, \cD)$ amounts to a
  function $\obs{F} \colon \obs{\bC} \to \pts{\cD}$ and a functor
  $F' \colon \bC \to \fun{(\DSet)}{\Set}$ such that each
  $F'(c)$ is the functor
  $\reindex{\obs{F}(c)} \circ \car{\clift{\cD}}$, which sends a
  $\cD$-coalgebra to the elements in it lying over $\obs{F}(c)$. This
  $F'$ curries to a functor $\DSet \to \fun{\bC}{\Set}$ lifting
  $\reindex{\obs{F}} \colon \SetpD \to \fun{\obs{\bC}}{\Set}$,
  corresponding to a continuous map $\cbC \to \cD$. Natural
  transformations also transport along currying.
\end{proof}

Another way of saying this is that $\cbC$ is the
$\bC$-copower of $\pointcmd$ in $\ConSet$. This is an instance of the
more general construction of copowers in $\ConSet$, discussed briefly
later in \cref{rem:copowers}.

\section{The double category of comonads}\label{sec:doubleset}

Comonad maps and continuous maps between comonads on $\Set$ amount to two kinds of commutative
squares, shown respectively below:
\[
  \begin{tikzcd}
    \CSet & \DSet \\
    {\SetpC} & {\SetpD}
    \arrow["{\com{\Set}{\psi}}", from=1-1, to=1-2]
    \arrow["{\car{\clift{\cC}}}"', from=1-1, to=2-1]
    \arrow["{\car{\clift{\cD}}}", from=1-2, to=2-2]
    \arrow["{\leindex{\one{\psi}}}"', from=2-1, to=2-2]
  \end{tikzcd}
  \qquad\qquad\qquad
  \begin{tikzcd}
    \CSet & \DSet \\
    {\SetpC} & {\SetpD}
    \arrow["{\car{\clift{\cC}}}"', from=1-1, to=2-1]
    \arrow["{\con{\Set}{f}}"', from=1-2, to=1-1]
    \arrow["{\car{\clift{\cD}}}", from=1-2, to=2-2]
    \arrow["{\reindex{\one{f}}}", from=2-2, to=2-1]
  \end{tikzcd}
\]

The two are related together in the structure of a double category.

\begin{definition}\label{def:doublespecset}
  A \emph{specialization} between comonad maps
  $\phi \colon \cA \to \cC$ and $\psi \colon \cB \to \cD$ and
  continuous maps $f \colon \cA \to \cB$ and $g \colon \cC \to \cD$ of
  comonads on $\Set$ is a natural transformation of the following
  form:
  \[
    \begin{tikzcd}[column sep=15pt]
      & \BSet & \\
      \ASet && \DSet \\
      & \CSet
      \arrow["{\con{\Set}{f}}"', from=1-2, to=2-1]
      \arrow["{\com{\Set}{\psi}}", from=1-2, to=2-3]
      \arrow["{\Rightarrow_\sigma}"{description}, shift left, draw=none, from=1-2, to=3-2]
      \arrow["{\com{\Set}{\phi}}"', from=2-1, to=3-2]
      \arrow["{\con{\Set}{g}}", from=2-3, to=3-2]
    \end{tikzcd}
  \]
\end{definition}

We denote by $\CmdSet$ the double category consisting of comonads on
$\Set$, comonad maps and continuous maps, and
specializations.\footnote{We have been deliberately vague about which
  direction in this strict double category is ``vertical'' versus
  ``horizontal'' (or ``tight'' versus ``loose''; this is a strict
  double category).} As may be worth noting, this double category has
a small set of 2-cells per boundary by the same argument as in
\cref{prop:locallysmall}.

Specializations evidently generalize the specializations between
continuous maps from \cref{def:specializationset}, and therefore one
of the underlying 2-categories of $\CmdSet$ is $\ConSet$. But we also
obtain a notion of specialization between parallel comonad maps,
endowing $\ConSet$ with the structure of a 2-category as well:
specifically, the 2-cells are the comonad specializations of
\cref{def:comspec}, corresponding to arbitrary natural transformations
between functors $\com{\bS}{\phi}$.

However, this is not the only definition of double-categorical 2-cell
available: it would appear natural to ask for compatibility with the
continuous map and comonad map structures.

\begin{definition}\label{def:doubleidset}
  An \emph{identification} between comonad maps and continuous maps as
  above consists of a specialization $\sigma$ as above and a natural
  transformation
  $\leindex{\one{\phi}} \circ \reindex{\one{f}} \Rightarrow
  \reindex{\one{g}} \circ \leindex{\one{\psi}}$ satisfying the
  following cube axiom:
  \[
    \begin{tikzcd}[column sep=17]
      & \BSet & \\
      \ASet & \SetpB & \DSet \\
      \SetpA && \SetpD \\
      & \SetpC
      \arrow["{\con{\Set}{f}}"'{pos=0.4}, from=1-2, to=2-1]
      \arrow["{\car{\clift{\cB}}}"{pos=0.4}, from=1-2, to=2-2]
      \arrow["{\com{\Set}{\psi}}"{pos=0.4}, from=1-2, to=2-3]
      \arrow["{\car{\clift{\cA}}}"'{pos=0.4}, from=2-1, to=3-1]
      \arrow["{\reindex{\one{f}}}"'{pos=0.4}, from=2-2, to=3-1]
      \arrow["{\leindex{\one{\psi}}}"{pos=0.4}, from=2-2, to=3-3]
      \arrow["\Rightarrow"{description}, draw=none, from=2-2, to=4-2]
      \arrow["{\car{\clift{\cD}}}"{pos=0.4}, from=2-3, to=3-3]
      \arrow["{\leindex{\one{\phi}}}"', from=3-1, to=4-2]
      \arrow["{\reindex{\one{g}}}", from=3-3, to=4-2]
    \end{tikzcd}
    \qquad = \qquad
    \begin{tikzcd}[column sep=17]
      & \BSet & \\
      \ASet && \DSet \\
      \SetpA & \CSet & \SetpD \\
      & \SetpC
      \arrow["{\con{\Set}{f}}"'{pos=0.4}, from=1-2, to=2-1]
      \arrow["{\com{\Set}{\psi}}"{pos=0.4}, from=1-2, to=2-3]
      \arrow["{\Rightarrow_\sigma}"{description}, shift left, draw=none, from=1-2, to=3-2]
      \arrow["{\car{\clift{\cA}}}"'{pos=0.4}, from=2-1, to=3-1]
      \arrow["{\com{\Set}{\phi}}"'{pos=0.6}, from=2-1, to=3-2]
      \arrow["{\con{\Set}{g}}"{pos=0.6}, from=2-3, to=3-2]
      \arrow["{\car{\clift{\cD}}}"{pos=0.4}, from=2-3, to=3-3]
      \arrow["{\leindex{\one{\phi}}}"', from=3-1, to=4-2]
      \arrow["{\car{\clift{\cC}}}"{pos=0.4}, from=3-2, to=4-2]
      \arrow["{\reindex{\one{g}}}", from=3-3, to=4-2]
    \end{tikzcd}
  \]
\end{definition}

In fact, there is at most one such natural transformation
$\leindex{\one{\phi}} \circ \reindex{\one{f}} \Rightarrow
\reindex{\one{g}} \circ \leindex{\one{\psi}}$, which exists if and
only if $\one{\psi} \circ \one{f} = \one{g} \circ \one{\phi}$. Indeed,
specifying such a natural transformation is equivalent to specifying a
natural transformation
$\reindex{\one{f}} \circ \reindex{\one{\psi}} = \reindex{\one{\phi}}
\circ \reindex{\one{g}}$ by taking mates, and the only such natural
transformations between reindexing functors are identities. Therefore
identifications are merely specializations satisfying a property.

We denote by $\CmdIdSet$ the sub-double-category of $\CmdSet$
consisting of identifications.

\begin{proposition}
  For any identification as above, we have
  $\one{\psi} \circ \one{f} = \one{g} \circ \one{\phi}$. Furthermore,
  the double category $\CmdIdSet$ is \emph{thin}, i.e.\ each boundary
  admits at most one 2-cell. Moreover, the only 2-cells in its
  underlying 2-categories are identities.
\end{proposition}
\begin{proof}
  
  We have the first statement by the discussion above. The second
  statement follows because the comonadic functor $\car{\clift{\cA}}$
  is faithful.  The third statement follows because the comonadic
  functor $\car{\clift{\cA}}$ reflects identities.
\end{proof}

Although the only identifications between parallel maps are
identities, nontrivial double-categorical identifications do exist,
expressing a compatibility condition between continuous maps and
comonad maps.

\begin{remark}\label{rem:doublecoalgebra}
  For any comonad $\cC$ on $\Set$, the category of coalgebras $\CSet$
  can be found within the double category $\CmdIdSet$. Note that by
  \cref{prop:density}, $\cC$-coalgebras carried by the set $S$ are
  identified with comonad maps from the density comonad $\gen{S}$ of
  $S \colon \point \to \Set$ to $\cC$ (as seen in
  \cref{ex:indiscrete}). It is also straightforward to see by
  \cref{thm:universal} that functions between sets
  $\one{f} \colon S \to T$ are identified with continuous maps
  $f \colon \gen{S} \to \gen{T}$.  Now we claim that a coalgebra map
  carried by $\one{f}$ is precisely an identification
  \[\begin{tikzcd}[column sep=2.5pt,row sep=10pt]
      & {\CoalgSet{\gen{T}}} & \\
      {\CoalgSet{\gen{S}}} \\
      && \CSet
      \arrow[""{name=0, anchor=center, inner sep=0}, "{\con{\Set}{f}}"', from=1-2, to=2-1]
      \arrow["{\com{\Set}{\psi}}", from=1-2, to=3-3]
      \arrow["{\com{\Set}{\phi}}"', from=2-1, to=3-3]
      \arrow["{\Rightarrow_\gamma}"{anchor=center,pos=.45}, shift right=1pt, draw=none, from=0, to=3-3]
    \end{tikzcd}\]
  


  Indeed, just as comonad maps and continuous maps are determined by
  their restriction along subbases (\cref{prop:density} and
  \cref{thm:universal}), so too are specializations and
  identifications, which we will establish alongside
  \cref{thm:universal} in \cref{app:universal}. A specialization of
  the above form then reduces to
  
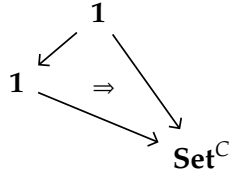
\[
    \begin{tikzcd}[column sep=15pt,row sep=10pt]
      & \point & \\
      \point \\
      && \CSet
      \arrow[""{name=0, anchor=center, inner sep=0}, from=1-2, to=2-1]
      \arrow[from=1-2, to=3-3]
      \arrow[from=2-1, to=3-3]
      \arrow["{\Rightarrow}"{anchor=center,pos=.45}, shift right=1pt, draw=none, from=0, to=3-3]
    \end{tikzcd}
  \]
  which simply amounts to a map of coalgebras, and the identification
  axiom ensures that it has underlying function $\one{f}$.
  %

  Abstractly, this all more or less amounts to observing that the
  \emph{cotabulator} (a standard kind of double-categorical colimit
  defined in~\cite{grandis-pare:limits}) in $\CmdIdSet$ of the
  continuous map $f \colon \gen{S} \to \gen{T}$ consists of the domain
  and codomain inclusion comonad maps from $\gen{S}$ and $\gen{T}$ to
  the density comonad of the diagram $\arro \to \Set$ given by the
  arrow $f \colon S \to T$.
\end{remark}

\begin{remark}
  Identifications may also be expressed in terms of the triangles
  corresponding to comonad maps as in \cref{prop:respect} and the
  triangles corresponding to continuous maps as in
  \cref{rem:laxtriangle}. Specifically, an identification may
  equivalently be defined as a specialization satisfying the following
  pyramid axiom:
  
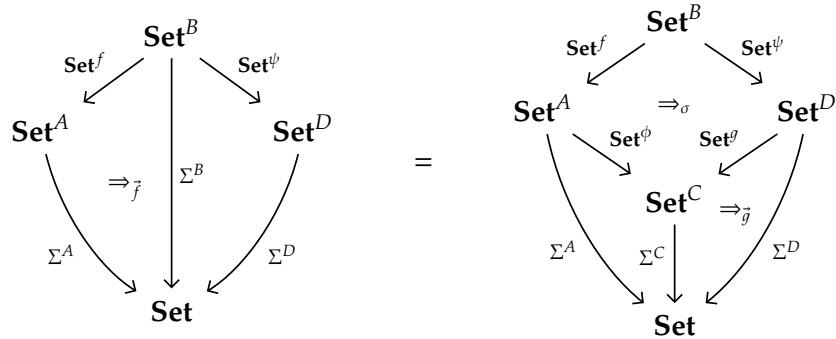
\[
    \begin{tikzcd}[row sep=17pt,column sep=17pt]
      & \BSet & \\
      \ASet && \DSet \\
      \\
      \\
      & \Set
      \arrow[""{name=0, anchor=center, inner sep=0}, "{\con{\Set}{f}}"', from=1-2, to=2-1]
      \arrow["{\com{\Set}{\psi}}", from=1-2, to=2-3]
      \arrow["\carB", from=1-2, to=5-2]
      \arrow["\carA"', curve={height=12pt}, from=2-1, to=5-2]
      \arrow["\carD", curve={height=-12pt}, from=2-3, to=5-2]
      \arrow["{\Rightarrow_{\connat{f}}}"{description}, shift left=2, draw=none, from=5-2, to=0]
    \end{tikzcd}
    \qquad = \qquad
    \begin{tikzcd}[row sep=15pt,column sep=17pt]
      & \BSet & \\
      \ASet && \DSet \\
      & \CSet \\
      \\
      & \Set
      \arrow["{\con{\Set}{f}}"', from=1-2, to=2-1]
      \arrow["{\com{\Set}{\psi}}", from=1-2, to=2-3]
      \arrow["{\Rightarrow_\sigma}"{description}, shift left, draw=none, from=2-1, to=2-3]
      \arrow["{\com{\Set}{\phi}}", from=2-1, to=3-2]
      \arrow["\carA"', curve={height=12pt}, from=2-1, to=5-2]
      \arrow["{\con{\Set}{g}}"', from=2-3, to=3-2]
      \arrow["\carD", curve={height=-12pt}, from=2-3, to=5-2]
      \arrow["\carC"'{pos=0.4}, from=3-2, to=5-2]
      \arrow["{\Rightarrow_{\connat{g}}}"{description}, shift left, draw=none, from=5-2, to=2-3]
    \end{tikzcd}
  \]
  Again, it is this style of definition that best generalizes to
  arbitrary comonads on arbitrary categories; see \cref{sec:maps} for
  details.
\end{remark}

\begin{remark}\label{rem:doubleformal}
  The double categories $\Cmd$ and $\CmdId$ can also be interpreted
  through the lens of the formal theory of (co)monads. We have seen
  that comonads on $\Set$ are identified with sliced (i.e.\
  terminal-object-and-finite-intersection-preserving) comonads on
  categories $\Set^X$, and comonad maps and continuous maps are
  identified with commutative squares carried by functors of the form
  $\leindex{\one{\phi}}$ and $\reindex{\one{\phi}}$ respectively. In
  turn, these are equivalently colax comonad functors
  (\cref{def:comfun}) of the form $\leindex{\one{\phi}}$ and
  $\reindex{\one{f}}$ respectively, by \cref{prop:formaltfae}. Also by
  \cref{prop:formaltfae}, a specialization as in
  \cref{def:doublespecset} is the same as a comonad specialization
  (\cref{def:comspec}) between the colax comonad functors
  $\leindex{\one{\phi}} \circ \reindex{\one{f}}$ and
  $\reindex{\one{g}} \circ \leindex{\one{\psi}}$. Likewise an
  identification as in \cref{def:doubleidset} is the same as a comonad
  transformation (\cref{def:comtrans}) between the colax comonad
  functors $\leindex{\one{\phi}} \circ \reindex{\one{f}}$ and
  $\reindex{\one{g}} \circ \leindex{\one{\psi}}$, since another view
  of the cube axiom is
  \[
    \begin{tikzcd}[column sep=17]
      & \ASet & \\
      \CSet & \SetpA & \BSet \\
      \SetpC && \SetpB \\
      & \SetpD
      \arrow["{\com{\Set}{\phi}}"', pos=.4, from=1-2, to=2-1]
      \arrow["{\car{\clift{\cA}}}", from=1-2, to=2-2]
      \arrow["{\car{\clift{\cC}}}"', from=2-1, to=3-1]
      \arrow["{\leindex{\one{\phi}}}"', pos=.4, from=2-2, to=3-1]
      \arrow["\Downarrow"{description}, draw=none, from=2-2, to=4-2]
      \arrow["{\con{\Set}{f}}"', from=2-3, to=1-2]
      \arrow["{\car{\clift{\cB}}}", from=2-3, to=3-3]
      \arrow["{\reindex{\one{f}}}"', from=3-3, to=2-2]
      \arrow["{\leindex{\one{\psi}}}", from=3-3, to=4-2]
      \arrow["{\reindex{\one{g}}}", pos=.4, from=4-2, to=3-1]
    \end{tikzcd}
    \quad=\quad
    \begin{tikzcd}[column sep=17]
      & \ASet & \\
      \CSet && \BSet \\
      \SetpC & \DSet & \SetpB \\
      & \SetpD
      \arrow["{\com{\Set}{\phi}}"', pos=.4, from=1-2, to=2-1]
      \arrow["{\Downarrow_\sigma}"{description}, shift left, draw=none, from=1-2, to=3-2]
      \arrow["{\car{\clift{\cC}}}"', from=2-1, to=3-1]
      \arrow["{\con{\Set}{f}}"', from=2-3, to=1-2]
      \arrow["{\com{\Set}{\psi}}", pos=.55, from=2-3, to=3-2]
      \arrow["{\car{\clift{\cB}}}", from=2-3, to=3-3]
      \arrow["{\con{\Set}{g}}", pos=.4, from=3-2, to=2-1]
      \arrow["{\car{\clift{\cD}}}", pos=.4, from=3-2, to=4-2]
      \arrow["{\leindex{\one{\psi}}}", from=3-3, to=4-2]
      \arrow["{\reindex{\one{g}}}", pos=.4, from=4-2, to=3-1]
    \end{tikzcd}
  \]
  which by \cref{prop:formaltfae} corresponds to a comonad
  transformation.

  Not only can comonad maps be viewed as colax comonad functors
  $\leindex{\one{\phi}}$; by taking mates, they can also be viewed as
  \emph{lax} comonad functors (which are defined dually)
  $\reindex{\one{\phi}}$. Indeed, it is a general fact that colax
  comonad functor structures on a left adjoint are in bijection with
  lax comonad functor structures on the corresponding right adjoint.
  Thus comonad maps and continuous maps between comonads on $\Set$ are
  precisely the lax and colax comonad functors carried by reindexing
  functors between sliced comonads on categories $\Set^X$. Thus the
  two types of 1-cells in the double categories $\CmdSet$ and
  $\CmdIdSet$ are respectively lax and colax comonad functors of the
  form $\reindex{f}$.
  
  In general, there are two double categories of comonads in any
  2-category $\tcat{C}$, defined in~\cite[Definitions 4.1 and
  4.2]{fairbanks:monads} (there for the dual case of monads), both of
  which have the same 1-cells, namely lax and colax comonad functors
  internal to $\tcat{C}$, but the 2-cells are different: they
  respectively generalize comonad transformations and comonad
  specializations (\cref{def:comtrans,def:comspec}). Moreover, left
  and right adjoint pairs of colax and lax comonad functors in
  $\tcat{C}$ yield companion pairs in both of these double categories,
  i.e.\ the double-categorical comonad transformation and comonad
  specialization cells of type shown below left are in correspondence
  with the simpler 2-categorical comonad transformation and comonad
  specialization cells of type shown below right
  \[
    \begin{tikzcd}[column sep=12.5,row sep=15]
      & B & \\
      A && D \\
      & C
      \arrow["{{F_1}}"', from=1-2, to=2-1]
      \arrow["{{R_2}}"', from=2-3, to=1-2]
      \arrow["\Downarrow"{description}, draw=none, from=2-3, to=2-1]
      \arrow["{{F_2}}", from=2-3, to=3-2]
      \arrow["{{R_1}}", from=3-2, to=2-1]
    \end{tikzcd}
    \qquad
    \leftrightarrows
    \qquad
    \begin{tikzcd}[row sep=20]
      B \\
      \\
      C
      \arrow[""{name=0, anchor=center, inner sep=0}, "{F_1 \circ L_1}"', curve={height=12pt}, from=1-1, to=3-1]
      \arrow[""{name=1, anchor=center, inner sep=0}, "{L_2 \circ F_2}", curve={height=-12pt}, from=1-1, to=3-1]
      \arrow["\Rightarrow"{description}, draw=none, from=1, to=0]
    \end{tikzcd}
  \]
  where $L_1 \dashv R_1$ and $L_2 \dashv R_2$ at the level of carrying 1-cells.

  Unravelling details, the categories $\CmdIdSet$ and $\CmdSet$ then
  arise --- up to reversing the directions of 1-cells --- as locally
  full sub-double-categories of the two double categories of comonads
  in $\CAT$, in both cases consisting of precisely the
  terminal-object-and-finite-intersection-preserving comonads on
  categories $\Set^\xX$ and 1-cells carried by reindexing functors.
\end{remark}

As a way of justifying our double-categorical definitions of
specialization and identification, we observe that they recover the
two known double categories of categories, retrofunctors, and
functors.

\begin{theorem}\label{thm:doublecat}
  The full sub-double-categories of $\CmdIdSet$ and $\CmdSet$
  consisting of polynomial comonads are respectively equivalent to the
  thin and non-thin double categories of small categories, retrofunctors
  (a.k.a.\ cofunctors), and functors introduced
  in~\cite[Definition 3.1]{clarke}
  and~\cite[Definition 5]{clarke-dimeglio}.
\end{theorem}

\begin{proof}
  The two double categories of $\bS$-internal categories,
  retrofunctors, and functors arise (up to equivalence) as the locally
  full sub-double-categories of the two (doubly weak) double
  categories of monads in the (weak) 2-category $\Span(\bS)$ --- where
  the 1-cells in both cases are lax and colax monad 1-cells, and the
  two types of double-categorical cells are generalized monad
  transformations and monad specializations as discussed in
  \cref{rem:doubleformal} --- consisting of those 1-cells carried by
  maps in $\bS$~\cite[Example 4.8]{fairbanks:monads}.
  We are interested in the case $\bS \coloneqq \Set$.

  Here spans from $\xX$ to $\xY$ correspond to left adjoint functors
  from $\Set^\xX$ to $\Set^\xY$; in particular, spans given by
  functions $\one{f}$ correspond to functors of the form
  $\leindex{\one{f}}$. Indeed, we already used this in
  \cref{lem:categories} to identify small categories with right
  adjoint comonads on categories $\Set^\xX$ and to identify functors
  with continuous maps; the result to be shown is a continuation of
  the same correspondence.

  Taking mates of everything in sight, thus turning left adjoints into
  right adjoints in the opposite direction, the two double categories
  of interest consisting of categories, retrofunctors, and functors
  are identified (suitably dualized) with the two double categories
  consisting of right adjoint comonads on categories $\Set^\xX$ and
  lax and colax comonad functors carried by reindexing functors
  $\reindex{\one{f}}$. And by \cref{rem:doubleformal}, these are precisely
  the desired full sub-double-categories of $\CmdIdSet$ and $\CmdSet$.
\end{proof}


We now characterize companions and conjoints in $\CmdIdSet$.

\begin{proposition}
  A comonad map or continuous map is a conjoint in $\CmdIdSet$ if and
  only if it is bijective-on-points.
\end{proposition}
\begin{proof}
  This is a special case of a more general description of conjoints in
  double categories of 2-functors between 2-categories. As described
  in \cref{rem:doubleformal}, comonad maps and continuous maps
  correspond to lax and colax comonad functors between sliced comonads
  carried by reindexing functors, and identifications correspond to
  (double-categorically generalized) comonad transformations. The
  double category of comonads, lax and colax comonad functors, and
  comonad transformations in a 2-category is an example of the double
  category of 2-functors, lax and colax transformations, and
  modifications between two 2-categories defined in~\cite[Proposition
  A.8]{fairbanks-shulman} where the domain is taken to be the free
  2-category containing a monad on an object. As noted in~\cite[Remark
  A.9]{fairbanks-shulman}, the conjoint pairs in this double category
  are the mate pair lax and colax 1-cells
  $\clift{\cC} \circ \reindex{\one{\phi}} \Rightarrow
  \reindex{\one{\phi}} \circ \clift{\cD}$ and
  $\reindex{\one{f}} \circ \clift{\cD} \Rightarrow \clift{\cC} \circ
  \reindex{\one{f}}$ where
  $\reindex{\one{f}} \dashv \reindex{\one{\phi}}$. But the only case
  in which $\reindex{\one{f}} \dashv \reindex{\one{\phi}}$ is if
  $\one{f} = \one{\phi}^{-1}$.
\end{proof}

\begin{proposition}\label{prop:companion}
  A comonad map $\phi$ between comonads on $\Set$ is a companion in
  $\CmdIdSet$ if and only if $\phi$ is a cartesian natural
  transformation.
\end{proposition}
\begin{proof}
  Also noted in~\cite[Remark A.9]{fairbanks-shulman}, the companion
  pairs here are, up to isomorphism, the inverse pair lax and colax
  1-cells
  $\clift{\cC} \circ \reindex{\one{\phi}} \Rightarrow
  \reindex{\one{\phi}} \circ \clift{\cD}$ and
  $\reindex{\one{f}} \circ \clift{\cD} \Rightarrow \clift{\cC} \circ
  \reindex{\one{f}}$ where $\one{\phi} = \one{f}$.

  Note that a natural transformation between
  terminal-object-preserving functors is invertible if and only if it
  is cartesian. In particular a comonad map $\phi \colon \cC \to \cD$
  is a companion in $\Cmd(\bS)$ if and only if its underlying
  $\clift{\cC} \circ \reindex{\one{\phi}} \Rightarrow
  \reindex{\one{\phi}} \circ \clift{\cD}$ is cartesian. In turn,
  $\clift{\cC} \circ \reindex{\one{\phi}} \Rightarrow
  \reindex{\one{\phi}} \circ \clift{\cD}$ is cartesian if and only if
  its mate
  $\clift{\phi} \colon \leindex{\one{\phi}} \circ \clift{\cC} \circ
  \reindex{\one{\phi}} \to \clift{\cD}$ is cartesian, since both
  $\leindex{\one{\phi}}$ and $\reindex{\one{\phi}}$ preserve pullbacks
  and both the unit and counit of the adjunction
  $\leindex{\one{\phi}} \dashv \reindex{\one{\phi}}$ are cartesian. And
  now $\clift{\phi}$ is cartesian if and only if
  $\phi = \slicel \circ \clift{\phi} \circ \slicer$ itself is cartesian,
  since $\slicer$ creates pullbacks.
\end{proof}

In the following two propositions, we will show that for
pullback-preserving comonads on $\Set$, \'etale comonad maps
(\cref{def:etale}) and companion (equivalently, cartesian) comonad
maps coincide. This tells us that any \'etal\'e comonad map between
pullback-preserving comonads on $\Set$ corresponds to a continuous map
in the same direction as well. This reflects the fact that an \'etale
map of topological spaces is in particular also a continuous map, and
a discrete opfibration of categories is in particular also a
functor. (Indeed, it was already observed that companions in the thin
double category of categories are discrete opfibrations
in~\cite{clarke}.)

\begin{proposition}\label{prop:booetale}
  The bijective-on-points and \'etale comonad maps form a
  factorization system on $\ComSet$.

  Moreover, this factorization system restricts to the full
  subcategory of pullback-preserving comonads.
\end{proposition}
\begin{proof}
  This is the {\slice} factorization (or equivalently, comprehensive
  factorization) applied to the maps $\com{\Set}{\phi}$.

  The comonad of elements of a pullback-preserving comonad is
  pullback-preserving, since slice category projections preserve
  all connected limits.
\end{proof}


\begin{proposition}\label{prop:cartesianetale}
  A comonad map into a pullback-preserving comonad on $\Set$ is
  \'etale if and only if it is cartesian, i.e.\ a companion in
  $\CmdSet$.
\end{proposition}
\begin{proof}
  First, note that a comonad admitting either an \'etale or cartesian
  comonad map into a pullback-preserving comonad is itself
  pullback-preserving.
  
  By \cref{prop:booetale} the bijective-on-points and \'etale form a
  factorization system on the category of pullback-preserving comonads
  and comonad maps. It is therefore sufficient to show that the
  bijective-on-points and cartesian comonad maps also form a
  factorization system on the same category, since the right class of
  a factorization system determines the left class.

  Such a factorization system already exists on $\SetSet$, the
  (vertical, cartesian) factorization system arising from the
  fibration $\one{\dash} \colon \SetSet \to \Set$ given by
  $F \mapsto F(1)$. We show that the factorization in $\SetSet$ of
  a comonad map between pullback-preserving comonads lifts to
  $\ComSet$.

  Suppose given a comonad map $\phi \colon \cC \coto \cD$. We obtain a
  factorization of $\phi$ in $\SetSet$ as $\cC \xto{\lambda} F \xto{\rho} \cD$
  where $\lambda$ is vertical and $\rho$ is cartesian. We equip $F$ with the
  unique comonad structure that makes $\lambda$ and $\rho$ comonad maps, as
  follows. The counit $F \Rightarrow \idSet$ is induced by the counit
  for $\cD$ as $\varepsilon \circ \rho$:
  \[\begin{tikzcd}
      \cC & F & \cD \\
      & \idSet
      \arrow["\lambda", from=1-1, to=1-2]
      \arrow["\varepsilon"', from=1-1, to=2-2]
      \arrow["\rho", from=1-2, to=1-3]
      \arrow[dashed, from=1-2, to=2-2]
      \arrow["\varepsilon", from=1-3, to=2-2]
    \end{tikzcd}\]
  The comultiplication is obtained by
  orthogonality of vertical and cartesian maps
  \[\begin{tikzcd}
      \cC & F & \cD \\
      {\cC \circ \cC} & {F \circ F} & {\cD \circ \cD}
      \arrow["\lambda", from=1-1, to=1-2]
      \arrow["\delta"', from=1-1, to=2-1]
      \arrow["\rho", from=1-2, to=1-3]
      \arrow[dashed, from=1-2, to=2-2]
      \arrow["\delta", from=1-3, to=2-3]
      \arrow["{\lambda \circ \lambda}"', from=2-1, to=2-2]
      \arrow["{\rho \circ \rho}"', from=2-2, to=2-3]
    \end{tikzcd}\] where we have used that $\cD$ is
  pullback-preserving to deduce that $\rho \circ \rho$ is cartesian. The
  unit and associativity laws for the above counit and
  comultiplication on $F$ are also deduced from orthogonality via
  uniqueness of diagonal fillers.
  %
  %
\end{proof}


Combining \cref{prop:booetale} and \cref{prop:cartesianetale}, we have
that every comonad map between pullback-preserving comonads can be
written, uniquely up to isomorphism, as the composite of a conjoint
and a companion.  With this, we are able to explicitly answer the
question of what comonad maps between topological spaces are:

\begin{corollary}\label{cor:topcom}
  For comonads $\cC$ and $\cD$ on $\Set$ corresponding to topological
  spaces $\tX$ and $\tY$, a comonad map $\phi \colon \cC \to \cD$
  consists of a function $\one{\phi} \colon \pts{\tX} \to \pts{\tY}$
  and a refinement of the topology on $\tX$ making $\one{\phi}$
  \'etale.
\end{corollary}
\begin{proof}
  The comonad map $\phi$ factors uniquely as a bijective-on-points
  comonad map followed by an \'etale comonad map. Being a conjoint,
  the bijective-on-points comonad map constitutes a
  bijective-on-points continuous map in the opposite direction, i.e.\
  a refinement of the topology $\tX$. Being a companion, the \'etale
  comonad map constitutes a continuous map in the same direction, with
  the same \'etale map on points, i.e.\ an \'etale map of topological
  spaces.
\end{proof}

We will now show that the thin double category of ionads has
tabulators of comonad maps. In fact these are \emph{strong
  tabulators}, in the sense of \cite{aleiferi,clarke}. This
generalizes the same result about the thin double category of
categories, retrofunctors, and functors, which was shown
in~\cite[Proposition 3.5]{clarke}.

\begin{theorem}\label{thm:strongtabulators}
  The full sub-double-category of $\CmdIdSet$ consisting of
  pullback-preserving comonads --- i.e.\ ionads --- has (strong)
  tabulators of comonad maps.
\end{theorem}

Our proof of this fact will be based upon the notion of \emph{halo},
which we will introduce in \cref{sec:halos}, a generalization of the
notion of infinitesimal neighborhood of a point in a topological
space. We give our proof here, but in order to understand it, the
reader is instructed to first read \cref{sec:halos}, specifically
\cref{rem:doublehalospec}.

\begin{proof}
  We will show that for a comonad map
  $\psi \colon \cB \to \cD$ and continuous maps $f \colon \cA \to \cB$
  and $g \colon \cA \to \cD$ where all comonads preserve pullbacks,
  any identification as shown below left factors uniquely as shown
  below right
  \[
    \begin{tikzcd}[column sep=2.5pt,row sep=10pt]
      & \BSet & \\
      && \DSet \\
      \ASet
      \arrow[""{name=0, anchor=center, inner sep=0}, "{\com{\Set}{\psi}}", from=1-2, to=2-3]
      \arrow["{\con{\Set}{f}}"', from=1-2, to=3-1]
      \arrow["{\con{\Set}{g}}", from=2-3, to=3-1]
      \arrow["{\Rightarrow_\gamma}"{anchor=center,pos=.55}, shift right=4pt, draw=none, from=3-1, to=0]
    \end{tikzcd}
    \qquad
    =
    \qquad
    \begin{tikzcd}[column sep=2.5pt,row sep=10pt]
      & \BSet && \\
      & {} & \CSet \\
      & \CSet & {} & \DSet \\
      \ASet
      \arrow["{\com{\Set}{\lambda}}", from=1-2, to=2-3]
      \arrow[from=1-2, to=3-2]
      \arrow["\Rightarrow"{description}, draw=none, from=2-3, to=2-2]
      \arrow[equals, from=2-3, to=3-2]
      \arrow["{\com{\Set}{\rho}}", from=2-3, to=3-4]
      \arrow["{\exists\bang}"', "{\con{\Set}{h}}", dashed, from=3-2, to=4-1]
      \arrow["\Rightarrow"{description}, draw=none, from=3-3, to=2-3]
      \arrow[from=3-4, to=3-2]
    \end{tikzcd}
  \]
  where $\lambda$ and $\rho$ are given by the (bijective-on-points,
  cartesian)-factorization of the comonad map $\psi$, the unmarked
  arrows and cells are given by their respective conjoint and
  companion continuous maps, and the dashed arrow is given by a
  continuous map $h \colon \cA \to \cC$.

  As in \cref{rem:doublehalospec}, such an identification asserts that
  $\one{\psi} \circ \one{f} = \one{g}$, and that for each point $a$ in
  $\pts{\cA}$, the two induced maps from the halo of $a$ to the halo
  of $\one{f}(a)$ are equal. Since conjoints are bijections-on-points
  and companions are isomorphisms-on-halos (\cref{rem:companionsviahalos}), 
  our continuous map
  $h \colon \cA \to \cC$ must be defined on points in the same way as
  $f$ and defined on halos in the same way as $g$. 
  
  We must still check
  that these $\pts{A}$-many maps of halos moreover constitute maps of
  $\pts{\cB}$-labelled halos,
  making up a natural transformation
  \[\begin{tikzcd}[column sep=15,row sep=15]
      & \SetpB & \\
      \SetpA && \SetpB \\
      & \SetpA
      \arrow["{\reindex{\one{f}}}"', from=1-2, to=2-1]
      \arrow["\Downarrow"{description}, draw=none, from=1-2, to=3-2]
      \arrow["{\clift{\cC}}"', from=2-3, to=1-2]
      \arrow["{\reindex{\one{f}}}", from=2-3, to=3-2]
      \arrow["{\clift{\cA}}", from=3-2, to=2-1]
    \end{tikzcd}
  \]

  Since halos correspond to a comonad's summands, halos of
  pullback-preserving comonads on $\Set$ preserve both pullbacks and
  the terminal object: that is, they are encoded by
  finite-limit-preserving functors $\Set \to \Set$.  By
  \cref{lem:cancellabel}, since the maps of halos given by $g$ and by
  $\rho$ constitute maps of labelled halos, and the maps of halos
  given by $g$ are obtained as composites of the maps of halos given
  by $h$ and by $\rho$, we have that the maps of halos given by $h$
  also constitute maps of labelled halos. It remains to verify that
  the above natural transformation defines a continuous map
  $h \colon \cA \to \cC$, i.e.\ satisfies the laws of a colax comonad
  functor $\clift{\cC} \to \clift{\cA}$.

  By construction, our putative colax comonad functor
  $\clift{\cC} \to \clift{\cA}$ factors the colax comonad functor
  $\clift{\cD} \to \clift{\cA}$ encoding the continuous map $g$
  through the colax comonad functor $\clift{\cD} \to \clift{\cC}$
  encoding the companion to $\rho$.  Using the colax comonad functor
  laws for both of these continuous maps, we then obtain
  \[
    \begin{tikzcd}
      \SetpD && \SetpD \\
      \SetpB && \SetpB \\
      \SetpA & \SetpB & \SetpA \\
      & \SetpA
      \arrow["{\reindex{\one{\psi}}}"', from=1-1, to=2-1]
      \arrow[""{name=0, anchor=center, inner sep=0}, "{\clift{\cD}}"', from=1-3, to=1-1]
      \arrow["{\reindex{\one{\psi}}}", from=1-3, to=2-3]
      \arrow["{\reindex{\one{f}}}"', from=2-1, to=3-1]
      \arrow[""{name=1, anchor=center, inner sep=0}, "{\clift{\cC}}"', from=2-3, to=2-1]
      \arrow[""{name=2, anchor=center, inner sep=0}, "{\clift{\cC}}", from=2-3, to=3-2]
      \arrow["{\reindex{\one{f}}}", from=2-3, to=3-3]
      \arrow[""{name=3, anchor=center, inner sep=0}, "{\clift{\cC}}", from=3-2, to=2-1]
      \arrow["{\reindex{\one{f}}}", from=3-2, to=4-2]
      \arrow[""{name=4, anchor=center, inner sep=0}, "{\clift{\cA}}", from=3-3, to=4-2]
      \arrow[""{name=5, anchor=center, inner sep=0}, "{\clift{\cA}}", from=4-2, to=3-1]
      \arrow["\cong"{description}, draw=none, from=0, to=1]
      \arrow["{\Downarrow_\delta}"{description}, draw=none, from=1, to=3-2]
      \arrow["\Downarrow"{description}, draw=none, from=2, to=4]
      \arrow["\Downarrow"{description}, draw=none, from=3, to=5]
    \end{tikzcd}
    \qquad=\qquad
    \begin{tikzcd}
      \SetpD && \SetpD \\
      \SetpB && \SetpB \\
      \SetpA && \SetpA \\
      & \SetpA
      \arrow["{\reindex{\one{\psi}}}"', from=1-1, to=2-1]
      \arrow[""{name=0, anchor=center, inner sep=0}, "{\clift{\cD}}"', from=1-3, to=1-1]
      \arrow["{\reindex{\one{\psi}}}", from=1-3, to=2-3]
      \arrow["{\reindex{\one{f}}}"', from=2-1, to=3-1]
      \arrow[""{name=1, anchor=center, inner sep=0}, "{\clift{\cC}}"', from=2-3, to=2-1]
      \arrow["{\reindex{\one{f}}}", from=2-3, to=3-3]
      \arrow[""{name=2, anchor=center, inner sep=0}, "{\clift{\cA}}"', from=3-3, to=3-1]
      \arrow["{\clift{\cA}}", from=3-3, to=4-2]
      \arrow["{\clift{\cA}}", from=4-2, to=3-1]
      \arrow["\cong"{description}, draw=none, from=0, to=1]
      \arrow["\Downarrow"{description}, draw=none, from=1, to=2]
      \arrow["{\Downarrow_\delta}"{description}, draw=none, from=2, to=4-2]
    \end{tikzcd}
  \]
  and similarly for the unit law. (Note that
  $\one{\rho} = \one{\psi}$, and the structure natural transformations
  for $\rho$ and its companion are inverse isomorphisms.)
  We may cancel the isomorphisms, and we may apply
  \cref{prop:faithfulhalo} to cancel $\reindex{\one{\psi}}$ as well,
  yielding the colax comonad functor laws as desired.
\end{proof}

\section{Limits and colimits}\label{sec:limits}

Similar to \cite[Section 6]{garner:ionads}, we conclude with results
on limits and colimits of diagrams of continuous maps. But let us
first recall some standard general results about limits and colimits
of diagrams of comonad maps.


\begin{proposition}\label{prop:comoncolim}
  For any category $\bS$, the forgetful functor $\Com(\bS) \to \End{\bS}$
  strictly creates colimits. In particular, if $\bS$ is cocomplete then
  $\Com(\bS)$ is cocomplete.
\end{proposition}
\begin{proof}
  It is a general fact that for a category of comonoids $\Comon{\bC}$
  the forgetful functor $\Comon{\bC} \to \bC$ strictly creates
  colimits (``colimits of comonoids are pointwise''). Likewise for a
  functor category $\bD^\bC$ the forgetful functor
  $\bD^\bC \to \bD^{\obs{\bC}}$ strictly creates colimits (``colimits
  of functors are pointwise'').
\end{proof}

\begin{example}
  In \cref{ex:symcat} we constructed a comonad on $\Set$ from a
  category $\bC$ and a subgroup $\bG$ of automorphisms of $\bC$, which
  we called a ``symmetrized category''. One can in fact carry out the
  same construction more generally using an arbitrary groupoid $\bG$
  and functor $F \colon \bG \to \Cat$ to yield a comonad $\cC$ on
  $\Set$. Since we have a functor from $\Cat$ to the category
  $\ConSet$ of comonads on $\Set$ and continuous maps, and since
  bijective-on-objects comonad maps are dual to continuous maps, we
  therefore obtain from $F$ a functor
  $F' \colon \bG\op \cong \bG \to \ComSet$. The comonad $\cC$ is then
  the colimit of $F'$ in $\ComSet$. Thus a ``symmetrized category''
  (in the more general sense) is precisely a groupoid-indexed colimit
  of comonads corresponding to categories.
\end{example}

Limits are less straightforward than colimits. At least the following
is straightforward:

\begin{proposition}
  For any category $\bS$, the identity comonad $\id_{\bS}$ is terminal in
  $\Com(\bS)$.
\end{proposition}
\begin{proof}
  It is a general fact that in a monoidal category, the unit object
  equipped with its unique comonoid structure is the terminal
  comonoid.
\end{proof}

The comprehensive reference for results on limits in categories of
comonads (or colimits in categories of monads)
is~\cite{kelly:transfinite}. Rather than simply citing the results, we
provide proofs in order to set up the techniques for the case of
continuous maps afterward.

\begin{proposition}[{\cite[Proposition 26.3]{kelly:transfinite}}]\label{prop:multicoalgebra}
  Let $\bS$ be a (locally small) cocomplete category. The fully
  faithful functor $\Com(\bS) \to \CAT/\bS$, sending a comonad to its
  associated strictly comonadic functor, preserves limits.
  
  In particular, a coalgebra of a product comonad
  $\prod_{i \in I} \cC_i$, if it exists, is a set equipped with a
  $\cC_i$-coalgebra structure for each $i \in I$.
\end{proposition}
\begin{proof}
  The density comonad $\gen{\xX}$ of an object
  $\xX \colon \point \to \bS$ (given by the colimit formula
  $\xA \mapsto \Ho{\bS}(\xX, \xA)\cdot \xA$ as in
  \cref{def:pointwise}) has the universal property
  (\cref{prop:density}) that comonad maps out of $\gen{\xX}$
  correspond to coalgebra structures on $\xX$. Thus for any limit of
  comonads $\lim_{i \in I} \cC_i$ that exists, a coalgebra structure
  on $\xX$ amounts to a cone of comonad maps $\gen{\xX} \to \cC_i$, or
  equivalently $\cC_i$-coalgebra structures on $\xX$ related by the
  diagram of comonad maps. This is precisely an object in the desired
  limit taken in $\CAT/\bS$.
  
  
  Maps of coalgebras can be similarly obtained by probing with density
  comonads of shape $\arro$. Moreover, we can extract source and
  target objects from arrows by composing with appropriate morphisms
  between the diagrams of shape $\point$ and $\arro$. This determines
  the entire structure of the category of coalgebras since the
  forgetful functor into $\bS$ is faithful.
\end{proof}


Not all small limits of comonads on $\Set$ need exist. However, we do
have that small limits of \emph{small} comonads exist and are
small. This comes from the general fact that small limits of
accessible comonads on a locally presentable category exist and are
accessible, as we show next. Unsurprisingly, the proof involves the
theory of accessible categories; we will also make use of
bicategorical limits, a.k.a.\ bilimits. For preliminaries on these
topics we refer the reader to~\cite{makkai-pare}. We will need a few
lemmas.

\begin{lemma}[{\cite[Theorem 5.1.6]{makkai-pare}}]\label{lem:accessiblebilimits}
  The 2-category of accessible categories and accessible functors has
  small bilimits, preserved by the forgetful 2-functor into
  $\CAT$.\qed
\end{lemma}

For simplicity our presentation will be focused on ordinary (i.e.\
conical) limits. (In fact several arguments in this section generalize
from ordinary limits or colimits to 2-categorical limits or colimits
essentially without changing the argument; see \cref{rem:copowers}.)

\begin{definition}\label{def:limbilim}
  Given a diagram in a 2-category, we will say that its \emph{bilimit
    and strict limit coincide} if every pseudo-cone is isomorphic to a
  strict cone. Assuming the strict 2-categorical limit exists, this is
  equivalent to the condition that the strict limit cone is also a
  bilimit pseudo-cone.
\end{definition}

\begin{lemma}\label{lem:limbilim}
  Suppose $D \colon \bB \to \CAT$ is a diagram whose bilimit and
  strict limit coincide. Now suppose $D' \colon \bB \to \CAT$ is
  another diagram with a natural transformation
  $\alpha \colon D' \Rightarrow D$ whose components are amnestic
  isofibrations (\cref{def:amnesticisofibration}). Then the bilimit
  and strict limit of $D'$ coincide.
\end{lemma}
\begin{proof}
  By composing with the natural transformation $\alpha$, any
  pseudocone to $D$ induces a pseudocone to $D'$, which in turn is
  isomorphic to a strict cone. Per leg, the isomorphisms lift along
  the isofibrations; transporting the original pseudocone along these
  isomorphisms yields a new pseudo-cone to $D$ lifting the strict cone
  to $D'$, which is itself a strict cone since the functors are
  amnestic.
\end{proof}

\begin{lemma}\label{lem:inducedcreation}
  Suppose $D \colon \bB \to \CAT$ is a diagram with a limit cone in
  which all functors (the functors in the diagram $D$ as well as the
  limit cone's legs) preserve limits or colimits of a certain class of
  diagrams.

  (More precisely, here the ``class of diagrams'' may be any
  collection of diagrams closed under composition with functors in the
  aforementioned cone. For example these could be all diagrams of
  shape $\bJ$ for fixed $\bJ$, or for another example they could be
  all diagrams admitting absolute limits.)

  Now suppose $D' \colon \bB \to \CAT$ is another diagram admitting
  a limit, with a natural transformation
  $\alpha \colon D' \Rightarrow D$ whose components strictly create
  limits or colimits of the same class. Then the induced functor
  $\lim D' \to \lim D$ strictly creates limits or colimits of the
  same class.
\end{lemma}
\begin{proof}
  Straightforward diagram chase.
\end{proof}

\begin{proposition}[{\cite[Section VIII]{kelly:transfinite},~\cite{blackwell}}]\label{prop:comonadlimit}
  Small limits of accessible comonads on a locally presentable
  category exist and are accessible.
\end{proposition}
\begin{proof}
  The limit of a diagram in a slice category
  $D \colon \bJ \to \bC/\xX$ is calculated as the limit of the diagram
  $D' \colon \bJ' \to \bC$ where $\bJ'$ is obtained from $\bJ$ by
  freely adjoining a terminal object, and $D'$ is obtained from $D$ by
  sending the terminal object to $\xX$. Note that such a diagram is
  connected and automatically comes with a natural transformation
  (cocone) to the constant diagram at $\xX$, whose limit is $\xX$; the
  limit in $\bC/\xX$ is the induced map $\lim D' \to \xX$.

  Suppose given a small diagram in the category $\Com(\bS)$ of
  comonads on a category $\bS$. Since the forgetful functor
  $\Com(\bS) \to \CAT/\bS$ is fully faithful, we need only show that
  the limit taken in $\CAT/\bS$ is a comonadic functor. As noted
  above, the limit of the diagram $D$ in $\CAT/\bS$ is calculated as
  the limit of a connected diagram $D'$ in $\CAT$, which comes with a
  natural transformation from $D'$ to the constant diagram at $\bS$
  whose components are strictly comonadic functors. Since the diagram
  shape has a terminal object, the bilimit and strict limit of the
  latter constant diagram at $\bS$ coincide as $\bS$. Since the
  components of the natural transformation are amnestic isofibrations,
  the bilmit and strict limit of $D'$ also coincide by
  \cref{lem:limbilim}. Also, since the latter constant diagram and its
  limit cone (being identities) preserve colimits and absolute limits,
  and the components of the natural transformation (being strictly
  comonadic) strictly create them, the induced map $\lim D' \to \bS$
  creates them as well by \cref{lem:inducedcreation}. Thus it remains
  to show that $\lim D' \to \bS$ is a left adjoint, hence strictly
  comonadic by \cref{prop:comonadicity}.

  Now assuming $\bS$ is accessible, the categories of coalgebras are
  also accessible by \cref{lem:accessibility}, as are the induced
  functors between them since they preserve colimits. Thus the limit
  $\lim D'$ is also accessible by \cref{lem:accessiblebilimits}.
  Further assuming $\bS$ is locally presentable (i.e.\ has colimits),
  $\lim D'$ is also locally presentable. Thus $\lim D' \to \bS$ is a
  left adjoint, as is any colimit-preserving functor from a total
  category to a locally small category.
\end{proof}


We now have the following facts about colimits and limits of comonads
on $\Set$ for the above entirely general reasons.

\begin{proposition}\label{prop:comlimits}
  The categories $\ComSet$ and $\AccComSet$ are both cocomplete.
  The category $\AccComSet$ is complete.\qed
\end{proposition}

Now we move on to limits and colimits in the category $\ConSet$ of
comonads on $\Set$ and \emph{continuous maps}.  First we address
colimits. The following brings to mind that the category of sheaves on
a coproduct of topological spaces is the product of the sheaf
categories. It also generalizes the fact that the functor
$\fun{(\dash)}{\Set} \colon \CAT\op \to \CAT$ sends colimits to
limits.

\begin{proposition}\label{prop:contcolimits}
  The functor $(\ConOSet)\op \to \CAT^\arro$ sending a comonad $\cC$ to the
  strictly comonadic functor
  $\car{\clift{\cC}} \colon \CSet \to \SetpC$ sends colimits (that
  exist) to limits.
\end{proposition}


\begin{proof}
  The proof will be similar in spirit to that of
  \cref{prop:multicoalgebra}: we determine a category of coalgebras by
  probing with maps. This time we classify objects and arrows in a
  category of coalgebras via continuous maps into certain special
  comonads, rather than via comonad maps out of certain special
  comonads.

  First of all, note that for any comonad $\cC$, functions from
  $\pts{\cC}$ to a fixed set $I \colon \point \to \Set$ are identified
  with continuous maps from $\cC$ to the density comonad $\gen{I}$
  (corresponding to the indiscrete topological space with points $I$
  or indiscrete category with objects $I$). Thus the set of points
  $\pts{\cC}$ can be determined up to isomorphism by probing with
  continuous maps from $\cC$ to such comonads $\gen{I}$.
  
  The $\cC$-coalgebras carried by objects of $\SetpC$ taking values in
  a fixed family of sets $(S_i)_{i \in I}$ are identified with
  continuous maps from $\cC$ to $\cP$ where
  $\pP \colon \arro \to \Set$ is the diagram given by the arrow
  $\sum_{i \in I} S_i \to I$ in $\Set$. Indeed, by
  \cref{thm:universal}, such a continuous map is specified by a
  diagram
  \[
    \begin{tikzcd}[column sep=22.5]
      \CSet & \arro \\
      {\SetpC} & {\Set^I}
      \arrow["\car{\clift{\cC}}"', from=1-1, to=2-1]
      \arrow[dashed, from=1-2, to=1-1]
      \arrow["\ptlift{\pP}", from=1-2, to=2-2]
      \arrow["{\reindex{\one{f}}}", from=2-2, to=2-1]
    \end{tikzcd}
  \]
  where $\one{f}$ selects the sets assigned to elements of
  $\pts{\cC}$. The target object in $\arro$ is automatically sent to
  the terminal object of $\CSet$, and the source object may be sent to
  any particular coalgebra structure on the chosen family of sets.
  
  Arrows between such $\cC$-coalgebras are then identified with maps
  from $\cC$ to the density comonad of the diagram of shape $\compo$
  given by
  \[\sum_{f \in I'}\mathrm{domain}(f) \to \sum_{f \in
      I'}\mathrm{codomain}(f) \to I'\] where $I'$ is the set of all
  functions $f \colon S \to T$ for $S, T \in I$.


  We can also extract underlying families of sets from objects, as well as
  source and target objects from arrows, by composing with appropriate
  morphisms the between diagrams of shape $\point$, $\arro$, and
  $\compo$.
  The category of coalgebras $\CSet$, as well as the faithful
  forgetful functor $\car{\clift{\cC}} \colon \CSet \to \SetpC$, may
  thus be explicitly obtained by probing with continuous maps out of
  $\cC$. This determines the desired description of the colimit.
\end{proof}

Similar to limits of diagrams of comonad maps, colimits of diagrams of
continuous maps need not exist in general, but in the case of small
comonads they do always exist.

\begin{theorem}\label{prop:contcocomplete}
  Colimits of diagrams of small comonads and continuous maps exist and
  are small. In particular, the category $\AccConOSet$ is cocomplete.
\end{theorem}

\begin{proof}
  Let $D \colon \bJ \to \ConOSet$ be a diagram of comonads and
  continuous maps, and let $\xX$ denote the colimit of the functor
  $\pts{D(\dash)}$
  \[\bJ \xto{D} \ConOSet \xto{\pts{\dash}} \Set.\] We have that
  $\Set^{\xX}$ is a strict limit as well as a bilimit of the diagram of
  reindexing functors $\Set^{\pts{D(\dash)}}$
  \[\bJ \xto{D} \ConOSet \xto{\pts{\dash}} \Set
    \xto{\Set^{(\dash)}} \CAT\] since the representable functor
  $\Set^{(\dash)}$ sends strict colimits and bicolimits (which
  coincide in $\Set$, viewed as a full sub-2-category of $\Cat$) to
  strict limits and bilimits respectively.

  Now consider the diagram of functors between categories of
  coalgebras $\Set^{D(\dash)}$
  \[\bJ \xto{D} \ConOSet \xto{\Set^{(\dash)}} \CAT.\] The strictly
  comonadic functors
  $\car{\clift{D(j)}} \colon \Set^{D(j)} \to \Set^{\pts{D(j)}}$
  constitute a natural transformation whose components are amnestic
  isofibrations. Thus the strict limit $\bL$ of $\Set^{D(\dash)}$ is
  itself a bilimit by \cref{lem:limbilim}.


  Because $\Set^{D(\dash)}$ and its limiting cone consist entirely of
  functors preserving all colimits, coreflexive equalizers, and
  terminal objects (indeed all limits), and the components of the
  natural transformation
  $\car{\clift{D(j)}} \colon \Set^{D(j)} \to \Set^{\pts{D(j)}}$
  strictly create such (co)limits, the induced functor between limits
  $\bL \to \Set^\xX$ strictly creates such (co)limits as well by
  \cref{lem:inducedcreation}. Now assume all comonads in
  the image of $D$ are small, so their categories of coalgebras are
  accessible. Bilimits in the 2-category of accessible categories and
  accessible functors are as in $\CAT$
  by~\cref{lem:accessiblebilimits}, so the limit $\bL$ is accessible
  as well. It is therefore locally presentable since it has colimits,
  and it follows that the colimit-preserving functor
  $\bL \to \Set^\xX$ is a left adjoint. It thus corresponds to a small
  (i.e.\ accessible) comonad on $\Set$, satisfying the universal
  property of the colimit in $\ConOSet$ by construction.
\end{proof}



In particular, comparing with \cref{prop:multicoalgebra} we can see
that colimits of connected diagrams of bijective-on-points continuous
maps are the same as limits of the corresponding comonad maps in the
opposite direction.
Interestingly, it also turns out that coproducts in
$\ComSet$ are just the same as coproducts in $\ConSet$. \anoteNI{Are
  they double-categorical coproducts?}


\begin{proposition}\label{prop:coproduct}
  We have
  $\CoalgSet{\sum_{i \in I}\cC_i} \simeq \prod_{i \in I}
  \CoalgSet{\cC_i}$.
\end{proposition}

Again, this is a bit cleaner in terms of the corresponding comonads on
categories $\Set^\xX$ rather than comonads on $\Set$: the above
equivalence of categories of coalgebras becomes an isomorphism.

\begin{proof}
  Since taking products is 2-functorial, we have an adjunction
  \[
    \begin{tikzcd}[column sep=40pt]
      \prod_{i \in I}\CoalgSet{\cC_i}\ar[r, shift left=5pt, "\prod_{i \in I} \car{\cC_i}"]\ar["\bot"{anchor=center},r,draw=none]&
      \Set^I\ar[l, shift left=5pt, "\prod_{i \in I} \rcar{\cC_i}"]
    \end{tikzcd}
  \]
  with induced comonad $\prod_{i \in I}\cC_i$ on $\Set^I$. Also, the
  functor $\prod_{i \in I} \car{\cC_i}$ strictly creates coreflexive
  equalizers just like its factors $\car{\cC_i}$, so it is strictly
  crudely comonadic by the crude comonadicity theorem
  (\cref{prop:crude}). The functor
  $\sum_{\bang} \colon \Set^I \to \Set$ is also (crudely)
  comonadic. Thus since the composite of a crudely comonadic functor
  with a (crudely) comonadic functor is (crudely) comonadic
  (\cref{lem:crudecompose}), the composite
  $\prod_{i \in I}\CoalgSet{\cC_i} \to \Set$ is (crudely)
  comonadic. The induced comonad on $\Set$ is the desired coproduct,
  calculated as in \cref{prop:comoncolim}.
\end{proof}


\begin{corollary}\label{prop:contcoprod}
  $\ConOSet$ has all small coproducts.
\end{corollary}
\begin{proof}
  The only place the assumption of small comonads was used in the proof of
  \cref{prop:contcocomplete} was in showing that the putative
  comonadic functor is a left adjoint. In the case of a coproduct, the
  left adjoint exists automatically as we saw in the proof of
  \cref{prop:coproduct}.
  %
\end{proof}

\begin{remark}\label{rem:copowers}
  Although we have been focused on ordinary (i.e.\ conical) limits and
  colimits rather than more general 2-categorical limits and colimits,
  it is worth at least remarking on 2-categorical generalizations
  of the above results.

  Assuming $\bS$ is cocomplete, then just like for ordinary limits in
  $\Com(\bS)$, precisely the same probing method used in the proof of
  \cref{prop:multicoalgebra} explicitly determines all flavors of
  2-categorical limits that exist in the 2-category $\Com(\bS)$ where
  2-cells are comonad specializations (\cref{def:comspec}), i.e.\
  arbitrary natural transformations between functors
  $\com{\bS}{\phi}$. In particular, we find that for a 2-limit in
  $\Com(\bS)$ the category of coalgebras is a subcategory of the
  corresponding 2-limit in $\CAT$, consisting of those objects whose
  underlying coalgebras are all carried by the same object in $\bS$
  and those arrows whose underlying coalgebra maps are all carried by
  the same arrow in $\bS$.

  \anoteNI{Does the existence result for limits of accessible comonads
    \cref{prop:comonadlimit} also generalizes to 2-limits?}


  And just like for ordinary colimits in $\ConOSet$, precisely the
  same probing method used in \cref{prop:contcolimits} explicitly
  determines all flavors of 2-categorical colimits that exist in
  $\ConSet$. In particular, we find that for a 2-colimit in $\ConSet$
  the category of coalgebras is a subcategory of the corresponding
  2-limit in $\CAT$, consisting of those objects and arrows whose
  underlying coalgebras and coalgebra maps are carried by sets and
  functions related by the relevant $\reindex{\one{f}}$ reindexing
  functors.

  Moreover, $\ConOSet$ has copowers (as does $\AccConSet$),
  constructed the same way as in \cite[Theorem 6.1]{garner:ionads} ---
  which are sent to powers by
  $\Set^{(\dash)} \colon (\ConOSet)\op \to \CAT$ --- and it follows
  that $\AccConSet$ is cocomplete as a 2-category.
\end{remark}

Next we investigate limits. First, just as coproducts with respect to
continuous maps are the same as coproducts with respect to comonad
maps, so too is the terminal object with respect to continuous maps
the same as the terminal object with respect to comonad maps.

\begin{proposition}\label{prop:conterminal}
  The identity $\idSet$ is terminal in $\ConOSet$.
\end{proposition}
\begin{proof}
  There is a unique functor $\Set \to \CSet$ lifting
  $\slicer \colon \Set \to \SetpC$ along the comonadic functor
  $\CSet \to \SetpC$, since such a functor is necessarily left adjoint
  and terminal-object-preserving. Namely, this is the functor sending
  a set $\xA$ to the $\xA$-fold coproduct of the terminal object in
  $\CSet$.
\end{proof}

Finally, we show completeness of the category $\ConSet$. Our approach
will be quite similar to \cite{garner:ionads}. We require a few
lemmas.

\begin{lemma}\label{lem:reindexlimits}
  The reindexing functors of the fibration
  $\one{\dash} \colon \ConOSet \to \Set$ preserve limits.
\end{lemma}
\begin{proof}
  Recall that the fiber of this fibration over a set $X$ is the
  opposite of the category of sliced comonads ---- i.e.\
  terminal-object-and-finite-intersection-preserving comonads --- on
  $\Set^\xX$, and the fibration reindexing functor along
  $f \colon X \to Y$ is given by (the opposite of) the functor
  $\reindex{f} \circ \dash \circ \Pi_{f}$ from sliced comonads on
  $\Set^Y$ to sliced comonads on $\Set^\xX$. Thus the claim is that
  these functors preserve colimits. Since colimits of comonoids are
  pointwise, it suffices to show that the functors
  $\reindex{f} \circ \dash \circ \Pi_{f}$ between categories of
  terminal-object-and-finite-intersection-preserving
  \emph{endofunctors} preserve colimits. Indeed, they have right
  adjoints given by the functors
  $\Pi_{f} \circ \dash \circ \reindex{f}$.
\end{proof}

Let $\Com(\bS)_\xX$ denote the subcategory of $\Com(\bS)$ (the
category of comonads on $\bS$) consisting of comonads $\cC$ with
$\qpts{\cC} = \xX$ and comonad maps $\phi$ with
$\one{\phi} = \id_\xX$.

\begin{lemma}\label{lem:objectivereflective}
  The following is an adjunction, exhibiting $\Com(\Set)_\xX$ as a
  reflective subcategory of $\Com(\Set/\xX)_\xX$, the category of
  terminal-object-preserving comonads on $\Set/\xX$.
  \[
    \begin{tikzcd}[column sep=40pt]
      \Com(\Set/\xX)_\xX\ar[r, shift left=5pt, "\slicel\circ\, \dash\, \circ \slicer"]\ar["\bot"{anchor=center},r,draw=none]&
      \Com(\Set)_\xX\ar[l, shift left=5pt, "\clift{(\dash)}"]
    \end{tikzcd}
  \]
\end{lemma}

(Here we use the slice categories $\Set/X$ for convenience instead of
our usual $\Set^X$, since this makes the construction of sliced
comonads $\clift{(\dash)}$ more direct.)

\begin{proof}
  Let $\cC$ be comonad on $\Set/\xX$ with $\pts{\cC} = \xX$ and let
  $\cD$ be a comonad on $\Set$ with $\pts{\cD} = \xX$. We have a
  correspondence between commutative triangles
  \[
    \begin{tikzcd}[column sep=15pt, row sep=15pt]
      \CoalgOn{\cC}{(\Set/\xX)}\ar[rr, "F"]\ar[dr, "\carC"']&&\DSet\ar[dl, "\ptlift{\carD}"]\\
      &\Set/\xX
    \end{tikzcd}
    \qquad\qquad
    \begin{tikzcd}[column sep=15pt, row sep=15pt]
      \CoalgOn{\cC}{(\Set/\xX)}\ar[rr, "F"]\ar[dr, "\slicel\circ \carC"']&&\DSet\ar[dl, "\carD"]\\
      &\Set
    \end{tikzcd}
  \]
  such that $F$ preserves the terminal coalgebra $\xX$, since {\slice}
  factorizations are unique. A triangle as shown left amounts to a
  comonad map $\cC \to \clift{\cD}$, whereas a triangle as shown right
  amounts to a comonad map $\slicel\circ \cC \circ \slicer \to \cD$ by
  \cref{prop:density} and \cref{cor:densityofleft}.
\end{proof}

\begin{lemma}\label{lem:fiberlimits}
  Each fiber of the fibration $\one{\dash} \colon \ConOSet \to \Set$
  is complete.
\end{lemma}
\begin{proof}
  The fiber over the set $\xX$ is the opposite of the category of
  comonads of the form $\clift{\cC}$ on $\Set^\xX$. Hence we must show
  that the category of sliced comonads on $\Set^\xX$ is cocomplete.

  Since such sliced comonads are reflective in
  terminal-object-preserving comonads on $\Set^\xX$ by
  \cref{lem:objectivereflective}, it suffices to show that the
  category of terminal-object-preserving comonads on $\Set^\xX$ is
  cocomplete. Since colimits of comonoids are pointwise
  (\cref{prop:comoncolim}), it suffices to show that the category of
  terminal-object-preserving endofunctors on $\Set^\xX$ is
  cocomplete. We have
  $\fun{\Set^\xX}{\Set^\xX} \simeq \fun{\xX}{(\fun{\Set^\xX}{\Set})}$
  whereby the terminal-object-preserving endofunctors correspond to
  families of terminal-object-preserving functors, so in fact it
  suffices to show that the category of terminal-object-preserving
  functors from $\Set^\xX$ to $\Set$ is cocomplete. Indeed, this is a
  reflective subcategory of $\fun{\Set^\xX}{\Set}$. Explicitly, the
  free terminal-object-preserving functor on a functor
  $F \colon \bC \to \Set$ is obtained as the colimit of a diagram
  \[\begin{tikzcd}
      {\Ho{\bC}(1, \dash)} & F
      \arrow["\vdots" {xshift=-1pt}, shift right=2.75, from=1-1, to=1-2]
      \arrow[shift left=2.75, from=1-1, to=1-2]
    \end{tikzcd}\] identifying all elements of
  $F(1)$.
\end{proof}

\begin{theorem}\label{thm:contcomplete}
  The categories $\ConSet$ and $\AccConSet$ are complete.
\end{theorem}
\begin{proof}
  To show that a category fibered over a complete category is complete
  it suffices to show the fibers are complete and the reindexing
  functors preserve limits; these are shown in \cref{lem:fiberlimits}
  and \cref{lem:reindexlimits}.

  We thus obtain the completeness of the category $\ConOSet$. We also
  obtain the same result for $\AccConSet$, noting all of the
  constructions involved in \cref{lem:fiberlimits} and
  \cref{lem:reindexlimits} respect accessibility.
\end{proof}

\anoteNI{Is the above fibration a 2-fibration over the 2-category $\Set$
  with the locally indiscrete enrichment? Do the fibers have and the
  reindexing functors preserve 2-limits? Does the analogue of the
  result on limits for fibrations hold for 2-limits for 2-fibrations?
  If all true, we would have 2-completeness.}

\begin{remark}\label{rem:dirichlet}
  Let $\tX_1$ and $\tX_2$ be topological spaces generated by the
  subbases $\pP_1 \colon \bB_1 \hookrightarrow \Set$ and
  $\pP_2 \colon \bB_2 \hookrightarrow \Set$ (and let us assume these
  subbases cover $\tX_1$ and $\tX_2$). The product space
  $\tX_1 \times \tX_2$ is given by $\pts{\tX_1} \times \pts{\tX_2}$
  with the \emph{box topology}, generated by the subbasis
  $\bB_1 \times \bB_2 \hookrightarrow \Set$ defined by
  $(U, V) \mapsto U \times V$. It was shown in \cite{garner:ionads}
  (specifically for the more restricted bases considered there) that a
  similar formula describes the cartesian product in the category of
  bounded (a.k.a.\ small) ionads and continuous maps, which also
  recovers the usual cartesian product of categories. However, in the
  wider context of comonads on $\Set$, this product ceases to coincide
  with the cartesian product (see \cref{cex:nonproductproduct}),
  instead defining a different monoidal product.

  If $F, G \colon \Set \to \Set$ are functors, let $F \otimes G$
  denote their \emph{Day convolution} with respect to $\times$ (if it
  exists).  That is, $F \otimes G$ is the pointwise left Kan extension
  of the external product
  \[\cat{\Set} \times \cat{\Set} \xto{F \times G} \Set \times \Set
    \xto{\times} \Set.\] along
  $\times \colon \Set \times \Set \to \Set$.  As discussed
  in~\cite{mcdermott-rivas-uustalu}, the product $\otimes$ lifts to
  comonads, giving a symmetric monoidal closed structure on the
  category of small comonads and comonad maps
  $\AccComSet$.\footnote{Analogous statements hold replacing $\Set$ by
    any symmetric monoidal closed locally presentable category.}

  This recovers the usual product of spaces.  In
  fact, it is given by the following generalized box topology
  formula. If $\pP_1 \colon \bB_1 \to \Set$ and
  $\pP_2 \colon \bB_2 \to \Set$ are functors with density comonads
  such that $\gen{\pP_1} \otimes \gen{\pP_2}$
  exists, then it is the density comonad of the external product
  \[\bB_1 \times \bB_2 \xto{\pP_1 \times \pP_2} \Set \times \Set
    \xto{\times} \Set.\] Indeed, writing the Kan extensions
  $\gen{\pP_1}$ and $\gen{\pP_2}$ (which are pointwise by
  \cref{lem:setpointwise}) as weighted colimits, we calculate their
  external product as
  \begin{align*}
    \gen{\pP_1}(c) \times \gen{\pP_2}(d)
    &\cong (\La{\pP_1}{\pP_1})(c) \times (\La{\pP_2}{\pP_2})(d)\\
    &\cong (\colimw^{\Ho{\Set}(\pP_1 (\dash), c)}\pP_1) \times (\colimw^{\Ho{\Set}(\pP_2 (\dash), d)} \pP_2)\\
    &\cong \colimw^{\Ho{\Set}(\pP_1 (\dash), c)}\colimw^{\Ho{\Set}(\pP_2 (\dash), d)}(\pP_1(\dash) \times \pP_2(\dash))\\
    &\cong \colimw^{\Ho{(\Set \times \Set)}((\pP_1(\dash), \pP_2(\dash)), (c,d))} (\pP_1(\dash) \times \pP_2(\dash))\\
    &\cong (\La{\pP_1 \times \pP_2}{(\pP_1(\dash) \times \pP_2(\dash))})(c, d)
  \end{align*} where we have used preservation of colimits by $\times$ in each
  variable, as $\Set$ is cartesian closed.
  Then by the general property of Kan extensions that
  $\La{f}{(\La{g}{h})} \cong \La{f \circ g}{h}$, we have as desired
  \[
    \gen{\pP_1} \otimes \gen{\pP_2} \cong \La{\times}{(\La{\pP_1
        \times \pP_2}{(\pP_1(\dash) \times \pP_2(\dash))})} \cong
    \La{\pP_1(\dash) \times \pP_2(\dash)}{(\pP_1(\dash) \times
      \pP_2(\dash))} \cong \gen{\pP_1(\dash) \times \pP_2(\dash)}.\]
  

  This product $\otimes$ not only gives a symmetric monoidal structure
  on the category of small comonads and comonad maps $\AccComSet$, but
  also on the category of small comonads and continuous maps
  $\AccConSet$. This is because we may take external products of
  continuous maps: given two continuous maps between small comonads,
  one can define a continuous map between products using
  \cref{thm:universal}. It can further be verified that $\otimes$
  defines symmetric monoidal structure on both of the double
  categories $\AccCmdSet$ and $\AccCmdIdSet$ of small comonads on
  $\Set$ using the further double-categorical forms of
  \cref{thm:universal} in \cref{app:universal}.
\end{remark}

\chaptertocspace
\chapter{Infinitesimal neighborhoods}\label{sec:halos}

In this section, we introduce the concept of infinitesimal
neighborhood, or \emph{halo},\footnote{Here we are borrowing the term
  for an infinitesimal neighborhood from nonstandard
  analysis~\cite{goldblatt}.} of a point of a comonad on
$\Set$. Through this lens we can view comonads on $\Set$ as spaces
rather more literally than before and import intuitions from
topology. We will find that various concepts encountered throughout
the paper have interpretations in terms of halos, which agree with the
expected behavior of infinitesimal neighborhoods in topological spaces.

\section{Halos}

Let us return to the comonad $\cX$ corresponding to a topological
space $\tX$, whose coalgebras are sheaves on $\tX$. Recall from
\cref{ex:top} that $\cX$ is given by the following formula (where $U$
ranges over the neighborhoods of a point $x$):

\begin{equation}\label{eq:comonadformula}
  \cX(\xA) = \sum_{x\in \pts{\tX}} \colim_{U \ni x} \xA^\xU
\end{equation}

\noindent
One can read off the set of points $\pts{\tX}$ as the value of the
comonad at $1$. Interestingly, it is also possible to read off the
``infinitesimal neighborhoods'' of points from the formula of a
comonad.

Let us formalize the notion of ``infinitesimal neighborhood'' of a
point as a \emph{formal limit of sets}. By definition, a (small)
formal limit of sets is an object of $\widetilde{\Set}$, the free
limit completion of $\Set$. Equivalently, $\widetilde{\Set}$ is the
opposite of the category of (small) functors $\Set \to \Set$. The
``infinitesimal neighborhood'' of a point $x$ in a topological space
$\tX$ is then the formal limit of all neighborhoods of $x$, i.e.\ the
limit of the functor
\begin{equation}\label{eq:comonadlimit}
  {\neighbors{x}} \hookrightarrow \Set \hookrightarrow \widetilde{\Set}
\end{equation}
where $\neighbors{x}$ denotes the poset of neighborhoods of $x$, the
inclusion $\neighbors{x} \hookrightarrow \Set$ sends a neighborhood to
its set of points, and $\Set \hookrightarrow \widetilde{\Set}$ is the embedding of
$\Set$ in $\widetilde{\Set}$ (i.e.\ sending a set $\xA$ to the
representable functor $(\dash)^\xA \colon \Set \to \Set$).

We call this formal limit of all neighborhoods of $x$ the \emph{halo}
of $x$, and we denote it by $\halo{x}$. That is,
\[\halo{x} = \flim_{U \ni x} U\]
where $U$ ranges over the neighborhoods of $x$, and using $\flim$ to denote the
formal limit, i.e.\ the limit taken in $\widetilde{\Set}$.

As for why we should think of an infinitesimal neighborhood
specifically as a formal limit of sets, consider some properties that
one would expect such an object to have. First, mapping into the halo
of $x$ ought to be the same as mapping into all of the individual
neighborhoods $\xU$ of $x$ at once --- for instance, we should be able
to pick out the point $x$ itself, or for that matter any other point
infinitesimally close to it. This is the defining property of the
limit of all neighborhoods of $x$. But unlike the ordinary limit taken
in $\Set$, we also expect that mapping \emph{out} of the halo of $x$
into a set $\xA$ should be the same as giving the \emph{germ} of an
$\xA$-valued function about $x$ --- a function only defined
infinitesimally, up to restricting to ever-smaller neighborhoods.

This is precisely what formal limits do: we have that
$\Ho{\widetilde{\Set}}(\xX, \dash)$ sends limits to limits as always,
but additionally, when $\xA$ is a set (itself viewed as a formal limit
of sets via the embedding $\Set \to \widetilde{\Set}$), we have that
$\Ho{\widetilde{\Set}}(\dash, \xA)$ sends limits to
\emph{colimits}.\footnote{An object for which homming in sends limits
  to colimits is also known as a \emph{co-tiny object}. The ordinary
  objects of a category are all co-tiny within its free limit
  completion.}

\[
  \Ho{\widetilde{\Set}}(\xX, \lim_{\xU \ni x} \xU) \cong \lim_{\xU \ni x} \Ho{\widetilde{\Set}}(\xX, \xU) \qquad\quad\qquad \Ho{\widetilde{\Set}}(\lim_{\xU \ni x} \xU, \xA) \cong \colim_{\xU \ni x} \Ho{\widetilde{\Set}}(\xU, \xA)
\]
\nopagebreak
\vspace{-1.5em}
\[
  \hspace{17.5em}\text{(where $\xA$ is a set)}
\]

In particular, maps from the halo $\halo{x}$ (the formal limit of
neighborhoods) of $x$ to a set $\xA$ are given by
\[\xA^{\halo{x}} = \colim_{U \ni x} \xA^\xU\]
i.e.\ germs about $x$ of $\xA$-valued functions.

\begin{remark}\label{rem:continuous}
  The notion of continuous map of topological spaces can be seen to
  arise from the notion of halo. A function
  $f \colon \pts{\tX} \to \pts{\tY}$ between topological spaces is
  continuous precisely when it induces a map of halos, i.e.\ a map of
  formal limits of sets, from the halo of $x$ to the halo of $f(x)$
  for each point $x$ of $\tX$. Indeed, as just described, formal
  limits of sets are such that mapping into the formal limit of a
  diagram by definition means mapping compatibly into every set in the
  diagram, whereas mapping out of a formal limit into a set just
  requires mapping out of some set in the diagram.  For $f$ to induce
  a map of halos for each $x$ then means that for every neighborhood
  $\epsilon \ni f(x)$, there exists a neighborhood $\delta \ni x$
  mapping into it, which is the usual epsilon-delta definition of
  continuity.
\end{remark}

The limit of the diagram \eqref{eq:comonadlimit} of all neighborhoods
of $x$ in $\widetilde{\Set}$, which recall is the opposite of the
category of (small) functors $\Set \to \Set$, is dually viewed as a
colimit in $\SetSet$, specifically:
\[\colim_{U \ni x} \, (\dash)^\xU.\] We see the above formula \eqref{eq:comonadformula}
for the comonad $\cX$ is simply the sum over all points of the
functors $\Set \to \Set$ encoding their halos. Intuitively, the
comonad structure then specifies how the halos cohere together as a
space.\footnote{Unravelling details, the counit specifies the central
  point in each halo, and the comultiplication encodes that every
  neighborhood $\xU$ includes some neighborhood $\xV$ in which every
  point $x$ has a neighborhood $\xW_x$ included in $\xU$, as in the
  system-of-neighborhoods axiomatization of topological spaces given
  in~\cite{willard}.}

Now we generalize from topological spaces to comonads on $\Set$.



\begin{definition}\label{def:halo}
  Let $\cC$ be a comonad on $\Set$ and let $x$ be a point in
  $\pts{\cC}$. The \emph{halo} of $x$ is the object of $(\SetSet)\op$,
  or formal limit of sets,\footnote{For large comonads on $\Set$, halos
    of points are objects of $(\SetSet)\op$, which is larger than the
    free small limit completion $\widetilde{\Set}$. We can still
    interpret such objects as large formal limits.} corresponding to
  the functor $\cC_x \colon \Set \to \Set$ defined as the composite
  \[
    \Set \xto{\slicer} \SetpC \xto{\clift{\cC}} \SetpC
    \xto{\reindex{x}} \Set.
  \]
  The functors $\cC_x$ are characterized by the fact that
  $\cC \cong \sum_{x \in \pts{\cC}} \cC_x$ and each $\cC_x$ preserves the
  terminal object (so is indecomposable as a summand). Indeed, we have
  \[\cC \cong \slicel \circ \clift{\cC} \circ \slicer \cong \sum_{x \in \pts{\cC}} \reindex{x} \circ \clift{\cC}
    \circ \slicer = \sum_{x \in \pts{\cC}} \cC_x\] since the functor
  $\slicel \colon \SetpC \to \Set$ is given by
  $\sum_{x \in \pts{\cC}} \reindex{x}$.

  The summand $\cC_x$ may also be defined as the pullback
  \[
    \begin{tikzcd}
      {\cC_x(\xA)} & {\cC(A)} \\
      1 & {\cC(1)}
      \arrow[hook, from=1-1, to=1-2]
      \arrow[from=1-1, to=2-1]
      \arrow["\lrcorner"{anchor=center, pos=0.125}, draw=none, from=1-1, to=2-2]
      \arrow["{\cC(\bang)}", from=1-2, to=2-2]
      \arrow["x"', hook, from=2-1, to=2-2]
    \end{tikzcd}
  \]

  We denote by $\Halo$ the full subcategory of $(\SetSet)\op$
  consisting of functors that preserve the terminal object (thus,
  the formal \emph{connected} limits, specifically). We also denote
  the halo of a point $x \in \pts{\cC}$ by $\halo{x}$. That is,
  $\halo{x}$ is the summand $\cC_x$ viewed as an object of $\Halo$
  (rather than dually as an object of $\SetSet$).
\end{definition}

\begin{remark}[Germs]\label{rem:halogerm}
  If $\cC$ is a comonad on $\Set$, then in the terminology of
  \cref{rem:densitygerms}, elements of $\cC(\xA)$ are $\xA$-valued
  \emph{germs}. An $\xA$-valued germ is simply a map from a halo of
  $\cC$ to the set $\xA$.
  
  Indeed, if $\cC$ is a comonad on $\Set$, then elements of $\cC(\xA)$
  are identified with natural transformations
  $\Set(\xA, \dash) \Rightarrow \cC$ by the Yoneda lemma. Since
  $\Set(\xA, \dash)$ preserves the terminal object, such a natural
  transformation factors through a unique indecomposable summand of
  $\cC$. Dually in terms of formal limits, this is to say that an
  element of $\cC(\xA)$ is equivalently specified by a map from one of
  the halos of $\cC$ to $\xA$ (itself viewed as a formal limit of sets
  via the Yoneda embedding). Thus the formula of any comonad on $\Set$
  (indeed, any endofunctor on $\Set$) can be written in the form
  \[
    \cC(\xA) = \sum_{x \in \pts{\cC}} \xA^{\halo{x}}
  \]
  which is like a polynomial, but where the exponent $\halo{x}$ is not
  merely a set but a formal connected limit of sets.

  This is as would be intuitively expected from the interpretation of
  halos as infinitesimal neighborhoods: a germ of an $\xA$-valued
  function is a mapping from some infinitesimal neighborhood of some
  $x$ to $\xA$.
\end{remark}

\begin{remark}[Subbases]\label{rem:subbasishalo}
  The halos of the density comonad $\cP$ of $\pP \colon \bB \to \Set$
  can be calculated explicitly from $\pP$, just as the neighborhoods
  of a point in a topological space can be calculated from a
  subbasis. By the formula for density comonads
  \cref{prop:setformula}~\labelcref{item:pcolimit}, the endofunctor
  $\cP$ is given by the colimit of
  \[
    (\El(\pP))\op\to\bB\op\xto{\pP\op}\Set\op\hookrightarrow\SetSet
  \]
  

  First, the above colimit decomposes as a coproduct of colimits
  indexed by the connected components of $(\El(\pP))\op$. Colimits in
  $\SetSet$ are limits in $(\SetSet)\op$, encoding formal limits of
  sets. The halo $\halo{x}$ of a point $x$ in $\pts{\pP}$ is then the
  formal limit indexed by the associated connected component
  $\neighbors{x}$ of $\El(\pP)$. (That is, $\neighbors{x}$ is the
  category of elements $\El(\pP_x)$ of the component $\pP_x$ of $\pP$,
  the fiber lying over $x$ in the colimit $\pts{\pP}$.) We can then
  write

  \[
    \halo{x} \;\; \cong \; \flim_{(U, u)\, \in \,\neighbors{x}} \; \pP(U)
  \]
  where $\flim$ denotes the limit taken in $(\SetSet)\op$.

  In particular, in the case of a topological space $\tX$ with
  subbasis $\pP \colon \bB \hookrightarrow \Set$ generating the
  comonad $\gen{\topbasis{\tX}} \cong \cX$,\footnote{As in
    \cref{rem:subbasiswarning}, we here need that the poset of
    subbasic open neighborhoods of each point is connected, so that
    the density comonad is $\cX$ as expected.} each connected
  component in the category of elements $\El(\pP)$ is the poset of
  subbasic open neighborhoods around a particular point in $\tX$,
  since maps in $\El(\pP)$ from $(U,x)$ to $(V,y)$ can exist only if
  $x=y$. That is, the halo of each point can be calculated as the
  formal limit of its subbasic open neighborhoods.

  In general, we may interpret the connected component $\neighbors{x}$
  of $\El(\pP)$ corresponding to the point $x$ in $\cP$ as its
  ``category of subbasic open neighborhoods''.
\end{remark}

\begin{remark}[Categories]\label{rem:categoryhalos}
  Recall from \cref{cor:au} that the comonads on $\Set$ corresponding
  to categories are precisely those carried by polynomial functors
  $\Set \to \Set$, i.e.\ sums of representable functors. The
  representable functor $(\dash)^\xX \colon \Set \to \Set$, as an
  object of $(\SetSet)\op$, represents the set $\xX$ itself viewed as
  a formal limit. Hence categories are the ``spaces'' (comonads) in
  which the infinitesimal neighborhoods (halos) of points are simply
  sets, rather than formal limits of sets. (Or to put it another way,
  categories are the ``spaces'' in which every object has a
  \emph{minimal} neighborhood, since the formal limit of a diagram
  over a category with an initial object is isomorphic to the
  diagram's value at the initial object.)

  The halo of each object $c$ is given by the set $\obs{c / \bC}$ of
  arrows out of $c$, since the formula from \cref{ex:cat} for the
  comonad $\cbC$ corresponding to a category $\bC$ is
  \[
    \cbC(\xA) = \sum_{c \in \pts{\bC}} \xA^{\obs{c / \bC}}.
  \]
  This halo is not merely a subset of the points (objects)
  $\obs{\bC}$. Rather, its elements (arrows in $\bC$) may be
  interpreted as \emph{witnesses} (or \emph{proofs}) of nearness to
  another point. Indeed, results throughout the paper have hinted at
  the principle that to pass from topological spaces to arbitrary
  comonads on $\Set$ is to replace truth values with sets: we are
  passing from density comonads of diagrams of subsets to density
  comonads of arbitrary diagrams (\cref{thm:tfaetop}); or we are
  passing from comonads on $\Pow{\xX}$ that preserve the terminal
  object and finite intersections to such comonads on $\Set^X$
  (\cref{cor:setslices}); or for polynomial comonads we are passing
  from preorders to categories (\cref{rem:preorders}).
  


  As noted in \cref{rem:preorders}, polynomial comonads on $\Set$
  (corresponding to categories) and comonads arising from topological
  spaces intersect at the comonads corresponding to Alexandroff
  topological spaces (or corresponding to preorders when viewed as
  categories). In terms of halos, this translates to the claim that
  Alexandroff spaces are the topological spaces in which halos of
  points are simply sets. And indeed, as is well-known, such an
  Alexandroff topological space (whose open sets correspond to
  upward-closed sets in the preorder) is equivalently one in which
  every point $x$ has a minimal neighborhood (the principal
  upward-closed set generated by $x$, consisting of just $x$ and the
  points above it). Since the limit of a diagram with an initial
  object is simply that initial object, we have that the halo of $x$,
  calculated as the formal limit of its diagram of neighborhoods, is
  then simply given by its minimal neighborhood. Thus the halos in
  Alexandroff topological spaces are sets; and here the halos are mere
  subsets of the set of points (no multiple witnesses of nearness).
\end{remark}

\begin{remark}[Continuous maps]\label{rem:halomap}
  A continuous map $f \colon \cC \to \cD$ of comonads on $\Set$
  induces a map of halos from the halo of each $x$ in $\pts{\cC}$ to
  the halo of $\one{f}(x)$ in $\pts{\cD}$. Thus we find there is a
  forgetful functor
  \[\ConSet \to \Fam(\Halo).\]

  To see this, observe that the colax comonad functor
  (\cref{def:comfun}) structure
  \[\reindex{\one{f}} \circ \clift{\cD} \Rightarrow \clift{\cC} \circ
  \reindex{\one{f}}\] defining a continuous map $f$ induces the
  following natural transformation of functors $\Set \to \Set$.
  \[\begin{tikzcd}[row sep=15pt]
      && \SetpD && \\
      \Set & \SetpC && \SetpD & \Set \\
      && \SetpC
      \arrow["{\reindex{\one{f}(x)}}"', curve={height=12pt}, from=1-3, to=2-1]
      \arrow["{\reindex{\one{f}}}", from=1-3, to=2-2]
      \arrow["\Downarrow"{description}, draw=none, from=1-3, to=3-3]
      \arrow["{\reindex{x}}", from=2-2, to=2-1]
      \arrow["{\clift{\cD}}"', from=2-4, to=1-3]
      \arrow["{\reindex{\one{f}}}"', from=2-4, to=3-3]
      \arrow["\slicer"', from=2-5, to=2-4]
      \arrow["\slicer", curve={height=-12pt}, from=2-5, to=3-3]
      \arrow["{\clift{\cC}}", from=3-3, to=2-2]
    \end{tikzcd}
  \]

  \noindent
  Here the bottom composite is the functor
  $\cC_{x} \colon \Set \to \Set$ encoding the halo $\halo{x}$ of $x$,
  and the top composite is the functor
  $\cD_{\one{f}(x)} \colon \Set \to \Set$ encoding the halo of $\one{f}(x)$. Maps of halos (formal limits) are by
  definition natural transformations in the opposite direction between
  the $\Set$ functors encoding them. Thus this is indeed a map of the
  expected halos. In fact we will see shortly
  (\cref{prop:continuoushalo}) that this constitutes a faithful and
  conservative functor $\ConOSet \to \Fam(\Halo)$, so continuous maps
  can be viewed as structure-preserving maps of families of halos.

  Halos are dually viewed as endofunctors on $\Set$ sending $\xA$ to
  the set of germs at a particular point of the comonad, so the dual
  interpretation of this is that a continuous map $\cC \to \cD$
  induces a natural mapping from $\xA$-valued germs at $f(x)$ to
  $\xA$-valued germs at $x$, and moreover $f$ can be entirely
  recovered from this information.
  
  A comonad map $\phi \colon \cC \to \cD$ by contrast induces for each
  $x$ in $\pts{\cC}$ a map of halos from the halo of $\one{\phi}(x)$
  in $\pts{\cD}$ to the halo of $x$, in accordance with the forgetful
  functor \[\ComSet \to \SetSet \cong \Fam(\Halo\op).\]
  
  This can also be seen using the presentation of a comonad map as a
  \emph{lax} comonad functor
  \[\clift{\cC} \circ \reindex{\one{f}} \Rightarrow \reindex{\one{f}}
    \circ \clift{\cD}\] (mate to the colax comonad functor
  $\leindex{\one{\phi}} \circ \clift{\cC} \Rightarrow \clift{\cD}
  \circ \leindex{\one{\phi}}$ given by \cref{prop:comtfae}).
  \[\begin{tikzcd}[row sep=15pt]
      && \SetpC && \\
      \Set & \SetpC && \SetpD & \Set \\
      && \SetpD
      \arrow["{\clift{\cC}}"', from=1-3, to=2-2]
      \arrow["\Downarrow"{description}, draw=none, from=1-3, to=3-3]
      \arrow["{\reindex{x}}"', from=2-2, to=2-1]
      \arrow["{\reindex{\one{\phi}}}", from=2-4, to=1-3]
      \arrow["{\clift{\cD}}", from=2-4, to=3-3]
      \arrow["\slicer"', curve={height=12pt}, from=2-5, to=1-3]
      \arrow["\slicer", from=2-5, to=2-4]
      \arrow["{\reindex{\one{\phi}(x)}}", curve={height=-12pt}, from=3-3, to=2-1]
      \arrow["{\reindex{\one{\phi}}}"', from=3-3, to=2-2]
    \end{tikzcd}\] Here the bottom composite is the functor
  $\cD_{\one{\phi}(x)} \colon \Set \to \Set$ encoding the halo of
  $\one{\phi}(x)$, and the top composite is the functor
  $\cC_{x} \colon \Set \to \Set$ encoding the halo of $x$. Such a map
  between families of terminal-object-preserving objects in $\SetSet$
  is just the same as a map between the sum endofunctors on $\Set$, as
  expected for a comonad map.

  In particular, recall from \cref{ex:indiscrete} that a
  $\cD$-coalgebra carried by the set $\eE$ is the same as a comonad
  map $\gen{\eE} \to \cD$ where $\gen{\eE}$ is the indiscrete comonad
  on $\eE$ points. The halos of $\eE$ are all given by the set
  $\eE$. Thus a $\cD$-coalgebra may be seen as picking for each
  element $e$ in $\eE$ a point $x$ in $\cD$ and a map from the halo of
  $x$ to $\eE$ --- or by \cref{rem:halogerm}, a \emph{germ about $x$}
  --- in a suitably compatible way. This matches what one finds in
  $\cX$-coalgebras as \'etal\'e spaces or $\cbC$-coalgebras as
  $\cC$-sets, and more generally the germ description of coalgebras of
  comonads from \cref{rem:densitygerms}.
\end{remark}

\begin{remark}[Specializations]\label{rem:halospec}
  If $f, g \colon \tX \to \tY$ are continuous maps between topological
  spaces, a pointwise specialization $f \Rightarrow g$ asserts, for
  each point $x$ of $\tX$, that the point $g(x)$ is in all
  neighborhoods of the point $f(x)$. Similarly, if
  $f, g \colon \cC \to \cD$ are continuous maps between comonads on
  $\Set$, then a specialization $f \Rightarrow g$ induces, for each
  point $x$ in $\pts{\cC}$, a choice of point in the halo of $f(x)$ in
  $\pts{\cD}$ --- i.e.\ a map in $\Halo$ from the set $1$ to the halo
  of $f(x)$.

  To see this, observe that the comonad specialization
  (\cref{def:comspec}) structure
  \[\reindex{\one{f}} \circ \clift{\cD} \Rightarrow \reindex{\one{g}}\]
  induces the following natural transformation of functors
  $\Set \to \Set$.
  \[
    \begin{tikzcd}
      && \SetpD && \\
      \Set & \SetpC && \SetpD & \Set
      \arrow["{\reindex{\one{f}(x)}}"', curve={height=12pt}, from=1-3, to=2-1]
      \arrow["{\reindex{\one{f}}}"', from=1-3, to=2-2]
      \arrow["{\reindex{x}}"', from=2-2, to=2-1]
      \arrow["{\clift{\cD}}"', from=2-4, to=1-3]
      \arrow[""{name=0, anchor=center, inner sep=0}, "{\reindex{\one{g}}}", from=2-4, to=2-2]
      \arrow[curve={height=-30pt}, equals, from=2-5, to=2-1]
      \arrow["\slicer"', from=2-5, to=2-4]
      \arrow["\Downarrow"{description}, draw=none, from=1-3, to=0]
    \end{tikzcd}
  \]
  
  \noindent
  Here the bottom composite is the identity functor on $\Set$, which
  encodes the one-element set $1$ as a formal limit of sets, and the
  top composite is the functor $\cD_{\one{f}(x)} \colon \Set \to \Set$
  encoding the halo of $\one{f}(x)$. Again, maps of formal limits are
  by definition natural transformations in the opposite direction
  between the $\Set$ functors encoding them, so this indeed specifies
  a point in the halo of $\one{f}(x)$.  In fact, in the same way that a
  continuous map $f \colon \cC \to \cD$ is entirely determined by the
  maps of halos it induces (i.e.\ the forgetful functor
  $\ConOSet \to \Fam(\Halo)$ is faithful), we will see in
  \cref{prop:spechalo} that a specialization $f \Rightarrow g$ is
  entirely determined by the points of halos it induces.\anoteNI{Does
    this give another way of seeing $\ConSet$ is locally locally
    small? Is there a small set of natural transformations from any
    endofunctor $F$ on $\Set$ to the identity? Or for that matter is
    there a small set of natural transformations of type
    $\reindex{\one{f}} \circ \clift{\cD} \Rightarrow
    \reindex{\one{g}}$ in the first place?}

  What is more, we will in particular see that a specialization from a
  \emph{point} $x \colon \point \to \cD$ is simply a point in the halo
  $\halo{x}$ of $x$. And by definition, the underlying category of a
  comonad (\cref{def:undercat}) is the category of specializations
  between its points. Thus, continuing the thread of
  \cref{rem:categoryhalos}, we can interpret the underlying category
  of a comonad as the category of \emph{witnesses of nearness} between
  points.
  %
  %
\end{remark}

Although we do have a forgetful functor from continuous maps to
families of halos, this forgets more than we may wish. For instance,
when we examine the halos of a category, which are simply the sets
$\pts{c/\bC}$ (see \cref{rem:categoryhalos}), we have no way of extracting
from these plain sets the codomains of arrows. But one of the defining
features of a functor, or equivalently continuous map of comonads
corresponding to categories, is that it preserves codomains (as well
as domains) of arrows. Fortunately, there is a simple variation on the
notion of halo that is able to retain the information of codomains of
arrows.

Instead of a formal limit of sets, we may consider formal limits of
\emph{$\xX$-indexed sets}. By definition, a (small) formal limit of
$\xX$-indexed sets is an object of the free limit completion
$\widetilde{\Set^X}$ of $\Set^X$. Equivalently, $\widetilde{\Set^X}$
is the opposite of the category of (small) functors $\Set^X \to \Set$.

Recall that in \cref{sec:slices}, we identified comonads $\cC$ on
$\Set$ having points $\pts{\cC} = X$ with those comonads $\clift{\cC}$ on
$\Set^X$ that preserve the terminal object and finite
intersections. Just as an endofunctor on $\Set$ decomposes as a set of
formal (connected) limits of sets, an endofunctor on $\Set^\xX$
decomposes into an $\xX$-indexed set of formal (connected) limits of
$\xX$-indexed sets. Indeed, we have
\[\fun{(\Set^\xX)}{(\Set^\xX)} \cong (\fun{(\Set^\xX)}{\Set})^X\]
by currying. This leads to the following notion of ``infinitesimal
neighborhood'' consisting of $\xX$-indexed sets rather than sets.

\begin{definition}\label{rem:labelledhalo}
  Let $\cC$ be a comonad on $\Set$, and let $x$ be a point in
  $\pts{\cC}$. The \emph{$\pts{\cC}$-labelled halo} of $x$ (or simply
  the \emph{labelled halo} of $x$) is the object of
  $(\fun{(\Set^\xX)}{\Set})\op$, or formal limit of $\xX$-indexed
  sets, defined by the composite
  \[
    \SetpC \xto{\clift{\cC}} \SetpC \xto{\reindex{x}} \Set.
  \]
  i.e.\ the component $\clift{\cC}_x$ of $\clift{\cC}$ at $x$.

  We denote by $\LHalo{X}$ the full subcategory of
  $(\fun{(\Set^\xX)}{\Set})\op$ consisting of functors that preserve the
  terminal object (again encoding formal \emph{connected} limits,
  specifically). We denote the labelled halo of $x$ by $\lhalo{x}$.
\end{definition}

There is a canonical forgetful functor $\LHalo{\xX}$, given by
pre-composing the functor $\slicer \colon \Set \to \SetpC$, which
``forgets labels''.  This is also the functor induced by freely
extending the sum functor $\slicel \colon \Set^X \to \Set$ to formal
limits, since dually the functor
$\dash \circ \slicer \colon \fun{(\Set^X)}{\Set} \to \SetSet$ is the
left adjoint of the functor
$\dash \circ \slicel \colon \SetSet \to \fun{(\Set^X)}{\Set}$
reindexing along $\slicel$.

\begin{remark}[Functors and natural transformations]\label{rem:halonat}
  Labelled halos are indeed more expressive than halos. Whereas the
  halo of an object $c$ in a category $\bC$ is given by just the set
  $\obs{c / \bC}$ (\cref{rem:categoryhalos}), the $\obs{\bC}$-labelled halo
  is the $\obs{\bC}$-indexed set $c' \mapsto \Ho{\bC}(c, c')$
  (\cref{ex:catslice}).
  
  Continuous maps are also more directly interpreted in terms of
  labelled halos than halos. The natural transformation defining a
  continuous map
  $\reindex{\one{f}} \circ \clift{\cD} \Rightarrow \clift{\cC} \circ
  \reindex{\one{f}}$ precisely encodes, for each $c$ in $\pts{\cC}$, a
  map of formal limits from the $\one{f}$-relabelling of the
  $\pts{\cC}$-labelled halo of $c$ to the $\pts{\cD}$-labelled halo of
  $f(c)$. In the case of comonads corresponding to categories,
  upgrading a map of unlabelled halos (i.e.\ sets of arrows) to a map
  of labelled halos (i.e.\ sets of arrows labelled by codomain)
  amounts to requiring preservation of codomains. The colax comonad
  functor laws may be understood as specifying preservation of
  ``identities'' and ``composition''.

  Likewise, the natural transformation defining a specialization
  $\reindex{\one{f}} \circ \clift{\cD} \Rightarrow \reindex{\one{g}}$
  precisely encodes, for each $c$ in $\pts{\cC}$, a
  $\one{g}(c)$-labelled element in the $\pts{D}$-labelled halo of
  $\one{f}(c)$. In the case of comonads corresponding to categories,
  these are the arrows $F(c) \to G(c)$ that make up a natural
  transformation $F \Rightarrow G$ between the functors $F$ and $G$
  corresponding to the continuous maps $f$ and $g$.

  Moreover, some of the intuition from natural transformations carries
  over to arbitrary specializations: a specialization between
  continuous maps is an identity if and only if each halo's induced
  point is its ``center'' (or ``identity''), the point induced by the
  counit $\varepsilon \colon \cC \rightarrow \idSet$, generalizing the
  situation for natural transformations.  Indeed, it is an easy
  exercise to show that, in general, a natural transformation between
  colax comonad functors constitutes a comonad transformation if and
  only if its composition with a counit constitutes a comonad
  specialization~\cite[Remark 3.3]{fairbanks:monads}. Thus a
  specialization between continuous maps (which is equivalent to a
  comonad specialization between the colax comonad functors
  $\reindex{\one{f}}$ and $\reindex{\one{g}}$) comes from an
  identification (which is equivalent to a comonad transformation
  between these colax comonad functors, as observed in
  \cref{rem:doubleformal}) if and only if the defining natural
  transformation factors through the counit.
\end{remark}

\begin{remark}[Topological spaces]
  The comonads $\cX$ corresponding to topological spaces $\tX$ are
  precisely those whose labelled halos are formal limits of
  \emph{subterminal} indexed sets, i.e.\ formal limits of diagrams of
  subsets of the set of points $X$. (We noted this in
  \cref{rem:yesgo}.) This captures the intuition that topological
  spaces are precisely the comonads whose points can be near one
  another in at most one way.

  Without labels, we lose the ability to characterize topological
  spaces in terms of halos; indeed, we cannot separate the halos of
  topological spaces from those of categories (\cref{cex:nogo}). Halos
  of a topological space as well as halos of categories (and more
  generally halos of a comonad given by any \emph{infinitesimal
    category} as defined in \cref{ex:topcat}) are both of the form
  \[\flim_{U \in \mathcal{F}} U\]
  where $U$ ranges over the subsets in a chosen filter $\mathcal{F}$
  on a fixed set $\xX$. Viewed dually as functors $\Set \to \Set$,
  these are the \emph{reduced powers}
  \[\colim_{U \in \mathcal{F}} \, (\dash)^U\]
  Reduced powers are equivalently the finite-limit-preserving
  endofunctors of $\Set$ that arise as quotients of
  representables~\cite[Theorem 1]{blass}.
  %
  %
  %
\end{remark}

We now establish that labelled halos of comonads on $\Set$ may be
viewed simply as halos equipped with additional structure. We
use the following lemma.

\begin{lemma}\label{lem:freeintslice}
  $\Set^\xX$ is the free completion under finite intersections of its
  full subcategory consisting of globally inhabited objects.
\end{lemma}
\begin{proof}
  Deferred to \cref{app:universal} (\cref{cor:freeintslice}).
\end{proof}

\begin{lemma}\label{prop:faithfulhalo}
  The functor from the category of finite-intersection-preserving
  functors $\Set^\xX \to \Set^Y$ to the category of
  finite-intersection-preserving functors $\Set \to \Set^Y$ given by
  pre-composing $\slicer \colon \Set \to \Set^\xX$
  is faithful and conservative.
\end{lemma}

Dually, this in particular tells us that the forgetful functor
$\Halo_X \to \Halo$ from $X$-labelled halos to halos (given by
forgetting labels) becomes faithful and conservative when restricted
to those labelled halos that arise from comonads on $\Set$. Indeed,
for any comonad $\cC$ on $\Set$, the comonad $\clift{\cC}$ preserves
finite intersections (\cref{cor:setslices}).

\begin{proof}
  By \cref{lem:freeintslice}, if $F, G \colon \Set^\xX \to \Set^Y$
  preserve finite intersections, then a natural transformation
  $\alpha \colon F \Rightarrow G$, as well as its invertibility, is
  determined by its restriction to the full subcategory $\sS$ of
  globally inhabited objects of $\Set^\xX$. Moreover, the functor
  $\slicer' \colon \Set_{>0} \to \sS$, induced by restricting
  $\slicer \colon \Set \to \Set^\xX$ to nonempty sets, is
  \emph{co-faithful}, meaning natural transformations between functors
  out of $\sS$ are determined by their restrictions along
  $\slicer'$. Indeed, a functor is co-faithful if and only if every
  object in the codomain is a retract of an object in the
  image~\cite[Proposition
  3.1]{carboni-johnson-street-verity}\cite[Proposition
  4.1]{adamek-elbashir-sobral-velebil}, which is clearly true in this
  case. Also, a functor is co-faithful if and only if it is
  \emph{co-conservative} (a.k.a. \emph{liberal}~\cite[Definition
  2.2]{carboni-johnson-street-verity}), meaning a natural transformation
  between functors out of $\sS$ is invertible if its restriction along
  $\slicer'$ is.
\end{proof}

From this we are able to see that continuous maps are identified with
certain maps of families of halos as claimed in \cref{rem:halomap}, justifying the view of comonads on
$\Set$ as families of halos equipped with additional structure.

\begin{proposition}\label{prop:continuoushalo}
  The forgetful functor $\ConSet \to \Fam(\Halo)$ from the
  category of comonads on $\Set$ and continuous maps to the category
  of families of halos is faithful and conservative.
\end{proposition}
\begin{proof}
  The forgetful functor, as described in \cref{rem:halomap}, sends a
  continuous map $f$ with colax comonad structure
  $\reindex{\one{f}} \circ \clift{\cD} \Rightarrow \clift{\cC} \circ
  \reindex{\one{f}}$ to the $\pts{\cC}$-indexed family of halo maps
  given by whiskering with functors $\slicer \colon \Set \to \SetpD$
  and $\reindex{c} \colon \SetpC \to \Set$.
  
  Since the functors $\reindex{c}$ are the product projections of
  $\SetpC$, they are jointly faithful. Thus faithfulness of
  $\ConSet \to \Fam((\SetSet)\op)$ follows from
  \cref{prop:faithfulhalo}. As for conservativity, note that the
  functors $\reindex{c}$ are also jointly conservative, and note that
  a map in the category $\Fam(\bC)$ is an isomorphism if and only if
  it is a family of isomorphisms indexed by a bijection. Thus if a
  continuous map $f \colon \cC \to \cD$ has an underlying invertible
  map in $\Fam((\SetSet)\op)$, then its colax comonad functor
  structure natural transformation
  $\reindex{\one{f}} \circ \clift{\cD} \Rightarrow \clift{\cC} \circ
  \reindex{\one{f}}$ is invertible and its map on points
  $\one{f} \colon \pts{\cC} \to \pts{\cD}$ is a bijection. The
  companion of the conjoint (or equivalently the conjoint of the
  companion) of this continuous map, given by the mate of the inverse
  (or equivalently the inverse of the mate)
  $\reindex{\one{f^{-1}}} \circ \clift{\cD} \Rightarrow \clift{\cC}
  \circ \reindex{\one{f^{-1}}}$ of its colax comonad functor structure
  natural isomorphism, is then its inverse.
\end{proof}

\begin{warning}\label{rem:slicefunctorcoslice}
  Naively, one might expect that an $\xX$-labelled halo is the same as
  a halo equipped with a map to $\xX$. In other words, one might
  expect that $\LHalo{\xX}$ is the slice category $\Halo/\xX$ (where
  we view the set $\xX$ itself as a formal limit of sets via the
  Yoneda embedding). But this is not so: as discussed in
  \cref{rem:halogerm}, a map from a formal limit of sets to $\xX$ is
  an $\xX$-valued \emph{germ}, or element of the corresponding
  endofunctor $\Set \to \Set$ at $\xX$. Not all germs necessarily
  correspond to liftings to $\LHalo{\xX}$. The issue is essentially
  that given a diagram of sets (defining a formal limit of sets) a map
  into $\xX$ from \emph{some} set in the diagram (defining a germ) may
  not consistently determine a map into $\xX$ from \emph{every} set in
  the diagram (defining a formal limit of $\xX$-indexed sets); see
  \cref{cex:haloslice}.

  We do have that this is so for particularly well-behaved halos: the
  \emph{full subcategory} of $\LHalo{\xX}$ that is dual to the full
  subcategory of \emph{finite-limit-preserving} functors
  $\Set^X \to \Set$ is in fact the slice category over $\xX$ of the
  \emph{full subcategory} of $\Halo$ that is dual to the full
  subcategory of \emph{finite-limit-preserving} functors
  $\Set \to \Set$. Indeed, we recall from~\cite[Proposition
  1.10.15]{jacobs:logic} that if $\bC$ is a category with finite
  limits and $\xX$ is an object of $\bC$, then the slice category
  $\bC/\xX$ is the free category with finite limits extending $\bC$
  with a global element of $\xX$.  That is, if $\bD$ is also a
  category with finite limits, then giving a finite-limit-preserving
  functor $F' \colon \bC/\xX \to \bD$ is equivalent to giving a
  finite-limit-preserving functor $F \colon \bC \to \bD$ and an
  element $x \colon 1 \to F(\xX)$.
  In particular, giving a finite-limit-preserving functor
  $F' \colon \Set/\xX \simeq \Set^\xX \to \Set$ is equivalent to
  giving a finite-limit-preserving functor $F \colon \Set \to \Set$
  and an element $x \colon 1 \to F(\xX)$, or equivalently by the
  Yoneda lemma, a natural transformation
  $\Ho{\Set}(\xX, \dash) \Rightarrow \slicer$, which dualizes to the
  claimed statement about halos.
\end{warning}

According to the above \cref{rem:slicefunctorcoslice}, the forgetful
functor $\LHalo{\xX} \to \Halo$ is not a discrete fibration
(equivalently a slice category projection, since $\LHalo{\xX}$ has a
terminal object). Nevertheless, we have the following weaker property.

\begin{proposition}\label{lem:cancellabel}
  If $\overline{f} \colon \lhalo{x} \to \lhalo{z}$ and
  $\overline{g} \colon \lhalo{y} \to \lhalo{z}$ are maps between
  $\xX$-labelled halos with underlying maps of halos
  $f \colon \halo{x} \to \halo{z}$ and
  $g \colon \halo{y} \to \halo{z}$, and
  $h \colon \halo{x} \to \halo{y}$ is a map of halos such that
  $f = g \circ h$, then there exists a map of $\xX$-labelled
  halos $\overline{h} \colon \lhalo{x} \to \lhalo{y}$ whose underlying map of halos
  is $h$.
  \[
    \begin{tikzcd}
      {\halo{x}} && {\halo{z}} \\
      & {\halo{y}}
      \arrow["f", from=1-1, to=1-3]
      \arrow["h"', from=1-1, to=2-2]
      \arrow["g"', from=2-2, to=1-3]
    \end{tikzcd}
    \qquad\qquad\implies\qquad\qquad
    \begin{tikzcd}
      {\lhalo{x}} && {\lhalo{z}} \\
      & {\lhalo{y}}
      \arrow["{\overline{f}}", from=1-1, to=1-3]
      \arrow["{\overline{h}}"', dashed, from=1-1, to=2-2]
      \arrow["{\overline{g}}"', from=2-2, to=1-3]
    \end{tikzcd}
  \]
\end{proposition}

Moreover, by \cref{prop:faithfulhalo}, assuming these halos are in
fact halos of comonads on $\Set$ (so that in particular they
correspond to finite-intersection-preserving functors
$\Set^X \to \Set$), then this $\overline{h}$ is unique.

\begin{proof}
  Each $\xX$-labelled halo is a formal limit of some connected diagram
  of $\xX$-indexed sets. Let us without loss of generality freely
  adjoin a terminal object to each of these diagrams, sending it to
  the terminal $\xX$-indexed set; this does not change the formal
  limit, since the original diagram shape is initial (dual to final,
  \cref{def:final}) within the extended diagram shape.

  Now recall (as discussed in the beginning of the section) that a map
  of formal limits specifies, for each object in the codomain diagram,
  a map into it from some object in the domain diagram, modulo the
  equivalence relation generated by pre-composing maps in the domain
  diagram, and compatibly with respect to maps in the codomain
  diagram. Hence a map of halos $\halo{x} \to \halo{y}$ lifts to a map
  of $\xX$-labelled halos $\lhalo{x} \to \lhalo{y}$ if and only if the
  map at each object of the codomain diagram can be chosen to respect
  $\xX$-labels, which is evidently so if and only if the map at the
  terminal object can be chosen to be the identity.

  Assuming $f$ and $g$ satisfy this property, then $h$ satisfies this
  property as well. Indeed, the map at the terminal object is such
  that its composite with an identity is identity modulo the
  equivalence relation generated by pre-composing maps in the domain
  diagram, so it is as well.
\end{proof}

We show that specializations are determined by points as claimed in
\cref{rem:halospec}.

\begin{proposition}\label{prop:spechalo}\label{prop:spechalopoint}
  Let $f, g \colon \cC \to \cD$ be continuous maps between comonads on
  $\Set$. The canonical map from the set of specializations
  $f \Rightarrow g$ to the set of $\pts{\cC}$-indexed families
  consisting of a point in the halo of $f(x)$ for each
  $x \in \pts{\cC}$ is injective.

  Moreover, the canonical map from the set of all specializations out
  of a point $x \colon \point \to \cC$ to the set of all points in its
  halo $1 \to \halo{x}$ is bijective.
\end{proposition}
\begin{proof}

  The first statement follows from \cref{prop:faithfulhalo}, where the
  canonical assignment is as described in \cref{rem:halospec}.

  For the second statement, first note that a specialization between
  points $x, y \colon \point \to \cD$ is precisely a $y$-labelled
  point (i.e. a map from the element $y \in \pts{\cD}$, viewed as an
  $\pts{\cC}$-indexed set with one element at $y$ and no elements
  elsewhere) in the labelled halo $\lhalo{x}$ of $x$.  Indeed, colax
  comonad functors into identity comonads and comonad specializations
  between them are respectively 1-cells and 2-cells of the appropriate
  type satisfying no conditions at all, so such a specialization is
  simply a natural transformation
  $\reindex{\one{x}} \circ \clift{\cD} \Rightarrow \reindex{\one{y}}$.

  Therefore it suffices to show that such labelled points in
  $\xX$-labelled halos are in correspondence with points in the
  underlying (non-labelled) halos. Indeed, every $\xX$-labelled halo
  admits a unique map into the terminal $\xX$-indexed set $\xX$; thus
  by composing, any point in its underlying halo induces some point
  $y$ in $\xX$, which then lifts to a $y$-labelled point, as desired,
  by \cref{lem:cancellabel}.
  %
\end{proof}

Finally, we explain the double categories $\CmdSet$ and $\CmdIdSet$
from \cref{sec:doubleset} in terms of halos.

\begin{remark}[The double category of comonads]\label{rem:doublehalospec}
  Consider the boundary of a specialization or identification:
  
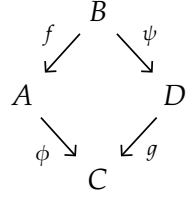
\[\begin{tikzcd}[column sep=10pt,row sep=15pt]
      & \cB & \\
      \cA && \cD \\
      & \cC
      \arrow["{f}"', from=1-2, to=2-1]
      \arrow["{\psi}", from=1-2, to=2-3]
      \arrow["{\phi}"', from=2-1, to=3-2]
      \arrow["{g}", from=2-3, to=3-2]
    \end{tikzcd}
  \]
  Here each comonad map or continuous map induces maps of halos as in
  \cref{rem:halomap}. Explicitly, for each point $a$ in $\pts{\cA}$,
  the comonad map $\phi$ induces a map from the halo of $\phi(a)$ to
  the halo of $a$, and the continuous map $f$ induces a map from the
  halo of $a$ to the halo of $\one{f}(a)$. Likewise for each point $c$
  in $\pts{\cC}$, the continuous map $g$ induces a map of halos from
  the halo of $c$ to the halo of $\one{g}(c)$, and for each point $b$
  in $\pts{\cB}$, the comonad map $\psi$ induces a map from the halo
  of $\one{\psi}(b)$ to the halo of $b$.

  This boundary is filled by an identification if and only if
  $\one{\psi} \circ \one{f} = \one{g} \circ \one{\phi}$ and, for each
  point $a$ in $\pts{\cA}$, the two induced maps from the halo of
  $\phi(a)$ to the halo of $f(a)$ are equal. Indeed, by
  \cref{rem:doubleformal} an identification is the same as a comonad
  transformation of colax comonad functors
  $\leindex{\one{\phi}} \circ \reindex{\one{f}} \Rightarrow
  \reindex{\one{g}} \circ \leindex{\one{\psi}}$; this is mate to an
  equality
  $\reindex{\one{f}} \circ \reindex{\one{\psi}} = \reindex{\one{\phi}}
  \circ \reindex{\one{g}}$
  which satisfies
  \[
    \begin{tikzcd}[column sep=10pt]
      & \SetpB & \SetpB & \\
      \SetpA & \SetpA && \SetpD \\
      & \SetpC & \SetpC
      \arrow["{\reindex{\one{f}}}"', from=1-2, to=2-1]
      \arrow["\Downarrow"{description}, draw=none, from=1-2, to=2-2]
      \arrow["{\clift{\cB}}"', from=1-3, to=1-2]
      \arrow["{\reindex{\one{f}}}", from=1-3, to=2-2]
      \arrow["{\clift{\cA}}"', from=2-2, to=2-1]
      \arrow["\Downarrow"{description}, draw=none, from=2-2, to=3-2]
      \arrow["{\reindex{\one{\psi}}}"', from=2-4, to=1-3]
      \arrow["{\reindex{\one{g}}}", from=2-4, to=3-3]
      \arrow["{\reindex{\one{\phi}}}", from=3-2, to=2-1]
      \arrow["{\reindex{\one{\phi}}}"', from=3-3, to=2-2]
      \arrow["{\clift{\cC}}", from=3-3, to=3-2]
    \end{tikzcd}
    \quad
    =
    \quad
    \begin{tikzcd}[column sep=10pt]
      & \SetpB & \SetpB & \\
      \SetpA && \SetpD & \SetpD \\
      & \SetpC & \SetpC
      \arrow["{\reindex{\one{f}}}"', from=1-2, to=2-1]
      \arrow["{\clift{\cB}}"', from=1-3, to=1-2]
      \arrow["\Downarrow"{description}, draw=none, from=1-3, to=2-3]
      \arrow["{\reindex{\one{\psi}}}", from=2-3, to=1-2]
      \arrow["{\reindex{\one{g}}}"', from=2-3, to=3-2]
      \arrow["\Downarrow"{description}, draw=none, from=2-3, to=3-3]
      \arrow["{\reindex{\one{\psi}}}"', from=2-4, to=1-3]
      \arrow["{\clift{\cD}}"', from=2-4, to=2-3]
      \arrow["{\reindex{\one{g}}}", from=2-4, to=3-3]
      \arrow["{\reindex{\one{\phi}}}", from=3-2, to=2-1]
      \arrow["{\clift{\cC}}", from=3-3, to=3-2]
    \end{tikzcd}
  \]
  Whiskering both sides of this equation with the functors
  $\slicer \colon \Set \to \SetpD$ and
  $\reindex{a} \colon \SetpA \to \Set$ yields the desired equation of
  halo maps. Conversely, the original equation can be deduced from the
  equations of halo maps by \cref{prop:faithfulhalo}.

  Double-categorical specializations are similar to
  \cref{rem:halospec}: a specialization filling the above boundary
  induces a choice of point in the halo of $\one{f}(a)$ in $\pts{\cB}$
  for each point $a$ in $\pts{\cA}$. The comonad specialization
  structure
  $\leindex{\one{\phi}} \circ \reindex{\one{f}} \circ \clift{\cB}
  \Rightarrow \reindex{\one{g}} \circ \leindex{\one{\psi}}$ is mate to
  a natural transformation
  $\reindex{\one{f}} \circ \clift{\cB} \circ \reindex{\one{\psi}}
  \Rightarrow \reindex{\one{\phi}} \circ \reindex{\one{g}}$,
  inducing the following natural transformation of functors
  $\Set \to \Set$.
  \[\begin{tikzcd}[column sep=5pt, row sep=15pt]
      &&&& \SetpB && \SetpB &&&& \\
      \Set & {} & \SetpA & {} &&&& {} & \SetpD & {} & \Set \\
      &&&&& \SetpC
      \arrow["{\reindex{\one{f}(a)}}"', curve={height=12pt}, from=1-5, to=2-1]
      \arrow["{\reindex{\one{f}}}", from=1-5, to=2-3]
      \arrow[""{name=0, anchor=center, inner sep=0}, "{\clift{\cB}}"', from=1-7, to=1-5]
      \arrow["{\reindex{a}}"', from=2-3, to=2-1]
      \arrow["{\reindex{\one{\psi}}}", from=2-9, to=1-7]
      \arrow["{\reindex{\one{g}}}", from=2-9, to=3-6]
      \arrow["\slicer"', curve={height=12pt}, from=2-11, to=1-7]
      \arrow[curve={height=-60pt}, equals, from=2-11, to=2-1]
      \arrow["\slicer"', from=2-11, to=2-9]
      \arrow["{\reindex{\one{\phi}}}", from=3-6, to=2-3]
      \arrow["\Downarrow"{description}, draw=none, from=0, to=3-6]
    \end{tikzcd}\]
  
  \noindent
  Again, the bottom composite is the identity functor on $\Set$, which
  encodes the one-element set $1$ as a formal limit of sets, and the
  top composite is the functor $\cB_{\one{f}(a)} \colon \Set \to \Set$
  encoding the halo of $\one{f}(a)$.  Again, a specialization is
  entirely determined by the points it induces by
  \cref{prop:faithfulhalo}. And again, a specialization is an
  identification if and only if every induced point of a halo is the
  one induced by the counit
  $\varepsilon \colon \cB \Rightarrow \idSet$ (an ``identity''), for
  the same reason as in the simpler case of specializations between
  continuous maps (\cref{rem:halonat}).
\end{remark}

\begin{remark}[Companions and conjoints]\label{rem:companionsviahalos}
  A comonad map or continuous map is a companion if and only if each
  induced map of halos, as in \cref{rem:halomap}, is an
  isomorphism.\footnote{More generally for arbitrary functors into
    $\Set$, a (uniformly) cartesian natural transformation is one that
    is an isomorphism per indecomposable summand (meaning a summand
    whose colimit is $1$). For functors with domain lacking a terminal
    object, it is appropriate to here use uniformly cartesian natural
    transformation as defined in \cref{app:comprehensive} instead of
    cartesian natural transformations.} Indeed, the invertibility of
  $\clift{\cC} \circ \reindex{\one{\phi}} \Rightarrow
  \reindex{\one{\phi}} \circ \clift{\cD}$ or
  $\reindex{\one{f}} \circ \clift{\cD} \Rightarrow \clift{\cC} \circ
  \reindex{\one{f}}$ in the proof of \cref{prop:companion} expresses
  that the induced maps of labelled halos are isomorphisms; and a map
  of labelled halos is an isomorphism if and only if the underlying
  map of halos is an isomorphism by \cref{prop:faithfulhalo}. The
  companion is then given by the same map on points and the inverse
  isomorphisms of halos. This also agrees with the explicit
  description of identifications in terms of halo maps from
  \cref{rem:doublehalospec}.

  On the other hand, as noted earlier, a comonad map is a conjoint in
  $\CmdIdSet$ if and only if the map on points is a bijection. The
  conjoint continuous map is then given by the same maps of halos and
  the inverse bijection on points.
\end{remark}

\clearpage
\phantomsection
\pdfbookmark[-1]{Appendices}{appendices}
\part*{Appendices}

\appendix

\addtocontents{toc}{\protect\setcounter{tocdepth}{0}}

\newcommand{\appendixtocspace}{\addtocontents{toc}{\vspace{-.5em}}}

\chapter{Comonadicity}\label{sec:comonadicity}
\appendixtocspace

A functor is called comonadic if it is essentially of the form
$\carC \colon \CS \to \bS$, the forgetful functor from a category of
coalgebras. More precisely:

\begin{definition}
  A functor $\bB \to \bS$ is \emph{comonadic} (resp.\ \emph{strictly
    comonadic}) if it is a left adjoint with induced comonad $\cC$,
  and the comparison functor $K \colon \bB \to \CS$ from
  \cref{lem:comparison} is an equivalence (resp.\ isomorphism).
\end{definition}

There are several useful criteria for determining whether a functor is
comonadic, called \emph{comonadicity theorems}. The purpose of this
section is to collect in one place assorted technical results used in
various proofs throughout the paper.

\begin{definition}\label{def:amnesticisofibration}
  A functor $F \colon \bB \to \bS$ is called \emph{amnestic} if for
  any isomorphism $f$ in $\bB$ such that $F(f)$ is an identity, we
  have that $f$ is an identity.

  A functor $F \colon \bB \to \bS$ is called an \emph{isofibration} if it lifts isomorphisms
  along objects: for any object $x$ in $\bB$ and isomorphism
  $f \colon F(x) \cong y$ in $\bS$ there exists an isomorphism
  $f' \colon x \cong y'$ in $\bB$ with $F(f') = f$.
\end{definition}

We note that a functor $F \colon \bB \to \bS$ is an amnestic isofibration
if and only if it lifts isomorphisms along objects uniquely: for any
object $x$ in $\bB$ and isomorphism $f \colon F(x) \cong y$ in $\bS$
there exists a unique isomorphism $f' \colon x \cong y'$ in $\bB$ with
$F(f') = f$.

\begin{lemma}
  A comonadic functor is strictly comonadic if and only if it is an
  amnestic isofibration.
\end{lemma}
\begin{proof}
  It is easy to verify that forgetful functors from categories of
  coalgebras are amnestic isofibrations. It is also easily shown that
  an equivalence is an isomorphism if and only if it is an amnestic
  isofibration. Amnestic isofibrations satisfy the same cancellative
  property as discrete fibrations, as well as being closed under
  composition. Therefore a comonadic functor is an amnestic
  isofibration if and only if the induced comparison equivalence is an
  isomorphism.
\end{proof}

\begin{definition}
  A functor $F \colon \bB \to \bS$ \emph{strictly creates} limits of a
  particular diagram $D \colon \bJ \to \bS$ if for any limiting
  cone $\alpha$ of $D$ and any diagram $D' \colon \bJ \to \bB$ with
  $F \circ D' = D$, there exists a unique cone $\alpha'$ of $D'$ with
  $F(\alpha') = \alpha$, and this $\alpha'$ is limiting.

  A functor $F$ \emph{creates} limits of $D$ if the above property holds up to
  isomorphism: for any limiting cone $\alpha$ of $D$ and any diagram
  $D' \colon \bJ \to \bB$ with $F \circ D' = D$, there exists a
  unique-up-to-isomorphism cone $\alpha'$ of $D'$ with $F(\alpha')$ isomorphic as a
  cone to $\alpha$, and this $\alpha'$ is limiting.%
  \footnote{%
    Our definition of (strict) creation is equivalent to reflection and (strict and
    unique) lifting. Some authors add assumptions of
    preservation and/or existence.
    In particular, when we say $F$ creates limits of a certain shape, we
    do not assume $\bB$ or $\bS$ has such limits, or even that $F$
    preserves such limits. If $\bS$ does have limits of a certain shape,
    then $F$ creating them does imply $F$ preserving them.%
  }
  Creation of colimits is defined dually.
\end{definition}

We note that a functor strictly creates (co)limits of shape $\point$
if and only if it is a conservative amnestic isofibration. This
property also follows from strict creation of any particular shape of
connected (co)limit. We also note that an amnestic isofibration
strictly creates any (co)limits it creates.


Strict versions of various results in this section follow from the
non-strict versions (as well as vice versa), using the above
observations about amnestic isofibrations.

\begin{proposition}[{\cite[Theorem 5.6.5]{riehl}}]\label{lem:create}\ \vspace{-1.5em}\\
  \begin{enumerate}[label=(\roman*)]
  \item Every (strictly) comonadic functor (strictly) creates, as well
    as preserves, colimits.
  \item Every (strictly) comonadic functor (strictly) creates limits
    preserved by $\cC$ and $\cC \circ \cC$, where $\cC$ is the induced
    comonad.\qed
  \end{enumerate}
\end{proposition}


\begin{proposition}[{\cite[Theorem 7.3]{pare:absolute}}]\label{prop:comonadicity}
  A left adjoint functor $\bB \to \bS$ is (strictly) comonadic if and
  only if it (strictly) creates limits that are absolute in $\bS$,
  i.e.\ limits preserved by all functors out of $\bS$.\qed
\end{proposition}

Another criterion, sufficient but not necessary, is the \emph{crude
  comonadicity theorem}, which makes use of the definition of
coreflexive equalizer (\cref{def:crep}).
A coreflexive equalizer may also be described as the limit of a
diagram whose shape is the free category containing a coreflexive
pair. When we say a functor \emph{creates coreflexive equalizers}, we
mean creation of limits of diagrams of this shape. (We do not require
creation of equalizers that are coreflexive in the codomain but not
the domain.)

\begin{proposition}[{\cite[Section 3.5]{barr-wells}}]\label{prop:crude}
  If a left adjoint functor (strictly) creates coreflexive equalizers,
  then it is (strictly) comonadic.\footnote{Again, assumptions of
    existence and preservation are often added, but they are not
    necessary.}\qed
\end{proposition}

In the case that $F \colon \bB \to \bS$ preserves and $\bB$ has limits
of a given shape, it is sufficient to check that $F$ is conservative in
order to deduce that $F$ reflects, and therefore creates, limits of that
shape.
Hence if a conservative left adjoint $F \colon \bB \to \bS$ preserves
and $\bB$ has coreflexive equalizers, then $F$ is comonadic. This is
the form of the result that is usually called the ``crude comonadicity
theorem''.
We note that by \cref{lem:coreflpull} coreflexive equalizers are
special cases of pullbacks, so it is here sufficient to check
preservation of pullbacks, or more specifically pullback squares
consisting of regular monomorphisms, i.e.\ regular finite
intersections.

Note also that if a comonad preserves coreflexive equalizers, then its
comonadic functor satisfies the conditions of \cref{prop:crude} by
\cref{lem:create}. The converse holds if $\bS$ has coreflexive
equalizers, and this is so in all our examples of interest.

In general, comonadic functors are not closed under composition, but
crudely comonadic functors are better behaved: the conditions of
\cref{prop:crude} compose. In fact, the composite of a functor
satisfying the conditions of \cref{prop:crude} with an arbitrary
comonadic functor is comonadic:

\begin{lemma}[{\cite[Proposition 3.5.1]{barr-wells}}]\label{lem:crudecompose}
  If $L_1 \colon \bB \to \bC$ is a left adjoint that creates
  coreflexive equalizers (hence is comonadic by \cref{prop:crude}) and
  $L_2 \colon \bC \to \bS$ is comonadic, then the composite
  $L_2 \circ L_1$ is comonadic.\qed
\end{lemma}

Comonadic functors obey the following \emph{cancellation law}:

\begin{lemma}[{\cite{bourn}}]\label{lem:cancel}
  Suppose the following is a commutative triangle of left adjoint
  functors:
  \[\begin{tikzcd}[column sep=15pt, row sep=15pt]
      {\bB} && {\bC} \\
      & {\bS}
      \arrow["F_1", from=1-1, to=1-3]
      \arrow["{F}"', from=1-1, to=2-2]
      \arrow["{F_2}", from=1-3, to=2-2]
    \end{tikzcd}\] If $F$ and $F_2$ are (strictly) comonadic, then
  $F_1$ is (strictly) comonadic.
\end{lemma}
\begin{proof}
  Follows from \cref{prop:comonadicity}.
\end{proof}

The following is the \emph{adjoint triangle theorem}.

\begin{lemma}[{\cite{dubuc:triangle}}]\label{lem:lifting}
  Let $\cC$ be a comonad on $\bS$, and suppose the following triangle
  commutes:
  \[\begin{tikzcd}[column sep=15pt, row sep=15pt]
      {\bB} && \CS \\
      & \bS
      \arrow["{L'}", from=1-1, to=1-3]
      \arrow["L"', from=1-1, to=2-2]
      \arrow["\carC", from=1-3, to=2-2]
    \end{tikzcd}\]
  If $\bB$ has coreflexive equalizers and $L$ is a left adjoint,
  then $L'$ is a left adjoint.\qed
\end{lemma}


\begin{corollary}\label{lem:adjointlift}
  Suppose given a comonad map $\phi \colon \cC \to \cD$, or
  equivalently a commutative triangle:
  \[
    \begin{tikzcd}[column sep=15pt, row sep=15pt]
      \CS\ar[rr, "\com{\bS}{\phi}"]\ar[dr, "\carC"']&&\DS\ar[dl, "\carD"]\\
      &\bS
    \end{tikzcd}
  \]
  If $\CS$ has coreflexive equalizers, then $\com{\bS}{\phi}$ is
  strictly comonadic.
\end{corollary}
\begin{proof}
  Follows from \cref{lem:lifting} and \cref{lem:cancel}.
\end{proof}

The following lemma is also useful.

\begin{lemma}[{\cite[Theorem IV.4.2]{maclane-moerdijk}}]\label{lem:crepff}
  Suppose $L \colon \bB \to \bS$ is a left adjoint functor inducing the
  comonad $\cC$. If $\bB$ has and $L$ preserves coreflexive
  equalizers, then the right adjoint $\CS \to \bB$ of the comparison
  functor $\bB \to \CS$ is fully faithful.\qed
\end{lemma}

\chapter{Continuous maps of comonads on general categories}\label{sec:maps}
\appendixsectionnumbering
\appendixtocspace

In this section we generalize continuous maps to comonads on
categories besides $\Set$. The only background from the paper needed
to read this is \cref{bg:II}; in particular it is not necessary to
first read \cref{sec:set}, which is about continuous maps in the case of
$\Set$.

In \cref{sec:continuous} we introduce the abstract definition of
continuous map between comonads (\cref{def:continuous}). We also study
two kinds of 2-cell between continuous maps: specializations
(\cref{def:specialization}) generalizing specializations between
continuous maps of topological spaces, defined just like
in~\cite{garner:ionads}; and identifications
(\cref{def:identification}), which are to be viewed as equating
different representations of the same continuous map.

In \cref{sec:double} we extend these two 2-categories to double
categories, combining comonad maps and continuous maps together. In
particular the identifications express a nontrivial compatibility
condition between comonad maps and continuous maps.

\begin{convention}\label{rem:noterminal}
  For simplicity, in this section $\bS$ always denotes a category with
  a terminal object. However, all definitions and results stated
  without qualification do generalize to the case that $\bS$ lacks a
  terminal object, by replacing \emph{terminal-object-preserving} with
  \emph{final} and \emph{cartesian} with \emph{uniformly cartesian} as
  defined in \cref{app:comprehensive}.\footnote{One simply modifies
    the proofs by passing to (large if necessary) cocompletions
    whenever a terminal object is needed: as described in
    \cref{app:comprehensive}, a functor $F$ is final if and only if
    the corresponding functor $\Coco{F}$ between free cocompletions is
    terminal-object-preserving, and a natural transformation $\alpha$
    is uniformly cartesian if and only if the corresponding natural
    transformation $\Coco{\alpha}$ is cartesian.}
\end{convention}

\section{Continuous maps}\label{sec:continuous}

Recall that a natural transformation is \emph{cartesian} if all
naturality squares are pullback squares.

\begin{definition}\label{def:continuous}
  Let $\cC$ and $\cD$ be comonads on $\bS$. A \emph{continuous map}
  $f \colon \cC \conto \cD$ consists of a {\final} functor
  $\con{\bS}{f} \colon \DS \to \CS$ and a {\cartesian} natural transformation
  \[
    \begin{tikzcd}[column sep=15pt]
      \CS && \DS \\
      & \bS
      \arrow["{\con{\bS}{f}}"', from=1-3, to=1-1]
      \arrow["\carC"', ""{name=0, anchor=center, inner sep=0}, from=1-1, to=2-2]
      \arrow["\carD", ""{name=1, anchor=center, inner sep=0}, from=1-3, to=2-2]
      \arrow["\Rightarrow_{\connat{f}}", shift right=10pt, draw=none, from=1, to=0]
    \end{tikzcd}
  \]
  
  We write $\one{f} \colon \qpts{\cC} \to \qpts{\cD}$ for the map induced
  by the component of $\connat{f}$ at the terminal object. We call
  $\con{\bS}{f}$ the \emph{inverse image} map, and we call $\one{f}$
  the \emph{points} map.
\end{definition}

As we will make more explicit in \cref{prop:contfae}, this definition
of continuous map is essentially the same as the definition for ionads
(\cref{def:ionadmap}) formulated in greater generality.

Intuitively, for $\con{\bS}{f}$ to preserve the terminal object means that
$\one{f}$ is a totally defined function: the inverse image of all
points $\pts{\cD}$ is all points $\pts{\cC}$. To make sense of the
rest of the definition, consider a similar diagram describing open
sets of topological spaces:
\[\begin{tikzcd}[column sep=5pt]
    \opens{\tX} && \opens{\tY} \\
    & \Set
    \arrow["{\contop{f}}"', from=1-3, to=1-1]
    \arrow["\pP_\tX"', ""{name=0, anchor=center, inner sep=0}, from=1-1, to=2-2]
    \arrow["\pP_\tY", ""{name=1, anchor=center, inner sep=0}, from=1-3, to=2-2]
    \arrow["\Rightarrow_{\connat{f}}", shift right=10pt, draw=none, from=1, to=0]
  \end{tikzcd}\]

For $\one{f}$ to extend to a cartesian natural transformation
$\connat{f}\colon \pP_\tX \circ \contop{f} \Rightarrow \pP_\tY$ means
precisely that $\contop{f}$ acts as the inverse image of $\one{f}$:
for any open set $U$ of $\tY$ we have a pullback square
\[
  \begin{tikzcd}
    \contop{f}(U)\ar[r, "\connat{f}_U"]\ar[d,hook]&U\ar[d,hook]\\
    \pts{\tX}\ar[r, "\one{f}"']&\pts{\tY}\ar[ul,phantom,very near end, "\lrcorner"]
  \end{tikzcd}
\]




\begin{lemma}
  Continuous maps form a category.\vspace{-.5em}
  \[\begin{tikzcd}[column sep=50pt]
      \CS & \DS & {\CoalgS{\cE}} \\
      & \bS
      \arrow[""{name=0, anchor=center, inner sep=0}, "\carC"', from=1-1, to=2-2]
      \arrow["\con{\bS}{f}"', from=1-2, to=1-1]
      \arrow[""{name=1, anchor=center, inner sep=0}, "\carD", from=1-2, to=2-2]
      \arrow["\con{\bS}{g}"', from=1-3, to=1-2]
      \arrow[""{name=2, anchor=center, inner sep=0}, "{\car{\cE}}", from=1-3, to=2-2]
      \arrow["{\Rightarrow_{\connat{f}}}"{pos=0.3}, shift right=3pt, draw=none, from=0, to=1]
      \arrow["{\Rightarrow_{\connat{g}}}"{pos=0.7}, shift right=3pt, draw=none, from=1, to=2]
    \end{tikzcd}\]
\end{lemma}
\begin{proof}
  {\Final} functors compose, whiskering an arrow with a {\cartesian}
  transformation on the domain side yields a {\cartesian} natural
  transformation, and {\cartesian} natural transformations compose
  vertically.
\end{proof}

The definition of specialization of continuous map between ionads
(\cref{def:ionadspecialization}) readily generalizes to continuous
maps between arbitrary comonads.

\begin{definition}\label{def:specialization}
  A \emph{specialization} between continuous maps
  $f \Rightarrow g \colon \cC \to \cD$ is a natural transformation
  $\con{\bS}{f} \Rightarrow \con{\bS}{g}$.
\end{definition}





\begin{remark}
  Continuous maps often retain some of the flavor of geometric
  morphisms of toposes. If $\bS$ is a category with finite limits and
  $f$ is a continuous map of crude (i.e.\
  coreflexive-equalizer-preserving) comonads on $\bS$ such that the
  pullback functor $\reindex{\one{f}}$ is a left adjoint (as is always
  the case when $\bS$ is locally cartesian closed), then
  $\con{\bS}{f}$ is a left adjoint by the adjoint triangle theorem
  (\cref{lem:lifting}).

  We here also have that $\con{\bS}{f}$ preserves coreflexive
  equalizers, since the comonadic forgetful functors create them and
  $\reindex{\one{f}}$ preserves them. Now $\con{\bS}{f}$ is itself
  comonadic if and only if it is conservative by the crude
  comonadicity theorem (\cref{prop:crude}).
\end{remark}

Given a continuous map $f \colon \tX \to \tY$ between topological
spaces, there may be several functors $\Sh(\tY) \to \Sh(\tX)$
representative of the same continuous map of topological spaces. For
this reason, continuous maps should be identified up to a suitable
notion of sameness.

\begin{definition}\label{def:identification}
  An \emph{identification} between continuous maps
  $f \Rightarrow g \colon \cC \to \cD$ is a specialization
  $\sigma \colon \con{\bS}{f} \Rightarrow \con{\bS}{g}$ satisfying
  \[
    \begin{tikzcd}[column sep=15pt, row sep=30pt]
      \CS && \DS \\
      & \bS
      \arrow["\carC"', ""{name=0, anchor=center, inner sep=0}, from=1-1, to=2-2]
      \arrow["{\con{\bS}{f}}"', from=1-3, to=1-1]
      \arrow["\carD", ""{name=1, anchor=center, inner sep=0}, from=1-3, to=2-2]
      \arrow["{\Rightarrow_{\connat{f}}}", shift right, draw=none, from=0, to=1]
    \end{tikzcd}
    \quad=\quad
    \begin{tikzcd}[column sep=15pt, row sep=30pt]
      \CS && \DS \\
      & \bS
      \arrow["\carC"', ""{name=0, anchor=center, inner sep=0}, from=1-1, to=2-2]
      \arrow[""{name=1, anchor=center, inner sep=0}, "{\con{\bS}{f}}"', curve={height=18pt}, from=1-3, to=1-1]
      \arrow[""{name=2, anchor=center, inner sep=0}, "{\con{\bS}{g}}", from=1-3, to=1-1]
      \arrow["\carD",""{name=3, anchor=center, inner sep=0}, from=1-3, to=2-2]
      \arrow["{\Rightarrow_{\connat{g}}}", shift right, draw=none, from=0, to=3]
      \arrow["{\Downarrow_{\sigma}}", shift right=2, draw=none, from=1, to=2]
    \end{tikzcd}
  \]
\end{definition}

\begin{lemma}\label{lem:one}
  If an identification exists between parallel continuous maps $f$ and
  $g$, then $\one{f} = \one{g}$.
\end{lemma}
\begin{proof}
  Given any specialization $\sigma$ as in \cref{def:specialization},
  $\sigma_{\qpts{\cD}}$ is an isomorphism since $f$ and $g$ are
  {\final}. Now by the identification law,
  $\connat{f}_{\qpts{\cD}} = \connat{g}_{\qpts{\cD}} \circ
  \carC(\sigma_{\qpts{\cD}})$, so
  \[\one{f} = \connat{f}_{\qpts{\cD}} \circ \carC(\bang^{-1}) =
    \connat{g}_{\qpts{\cD}} \circ \carC(\sigma_{\qpts{\cD}}) \circ
    \carC(\bang^{-1}) = \connat{g}_{\qpts{\cD}} \circ \carC(\bang^{-1})
    = \one{g}.\qedhere\]
\end{proof}

\begin{proposition}\label{lem:conid}
  The 2-category of comonads on $\bS$, continuous maps, and
  identifications is equivalent to a 1-category.
\end{proposition}
\begin{proof}
  Whiskering an identification 2-cell $\sigma$ as in
  \cref{def:identification} with the comonadic functor
  $\carC \colon \CS \to \bS$ yields a natural transformation which is
  {\cartesian} by pullback cancellation, and its component at
  $\qpts{\cD}$ is an isomorphism. Thus it is the unique natural
  transformation of its type whose components are the canonical
  isomorphisms between pullbacks.  Strictly comonadic functors reflect and
  uniquely lift natural isomorphisms because they are conservative
  amnestic isofibrations (i.e.\ reflect and uniquely lift isomorphisms
  along objects).
\end{proof}

\begin{definition}
  We write $\Con(\bS)$ for the above category (the category of
  comonads on $\bS$ and identification classes of continuous maps).
  
  Just as we use the same notation for the 2-category $\Cat$ and the
  1-category $\Cat$, we will also sometimes use $\Con(\bS)$ to refer
  to the 2-category of comonads on $\bS$, continuous maps, and
  specializations --- or rather the 2-category obtained from this by
  quotienting out identifications, just so that the underlying
  1-category agrees with $\Con(\bS)$.
\end{definition}

\begin{remark}
  Given a chosen terminal object of $\bS$, thus inducing chosen
  terminal objects of the categories $\CS$, we get a forgetful functor
  of 1-categories $\one{\dash} \colon \ConO(\bS) \to \bS$ by
  \cref{lem:one}.
\end{remark}


When viewed up to identification, continuous maps that are
``isomorphisms on points'' are just the same as ordinary comonad maps
that are ``isomorphisms on points'' in the opposite direction:

\begin{proposition}\label{prop:boo}
  The subcategory of $\Com(\bS)$ consisting of
  comonad maps $\phi$ such that $\one{\phi}$ is an isomorphism is isomorphic
  to the opposite of the subcategory of $\ConO(\bS)$ consisting of
  continuous maps $f$ such that $\one{f}$ is an isomorphism.
  
\end{proposition}
\begin{proof}
  Both kinds of map may be identified with commutative
  triangles
  \[
    \begin{tikzcd}[column sep=15pt]
      \CS && \DS \\
      & \bS
      \arrow[from=1-1, to=1-3]
      \arrow["\carC"', from=1-1, to=2-2]
      \arrow["\carD", from=1-3, to=2-2]
    \end{tikzcd}
  \]
  where the top arrow is {\final}.
  %
\end{proof}

Just as comonads and comonad maps may be lifted to slices
(\cref{def:sliceable,prop:comtfae}), so too may continuous maps,
specializations, and identifications. To fix terminology, let us say
that an \emph{underlying category structure} on a 2-category
$\tcat{C}$ consists of an identity-on-1-cells sub-2-category
$\tcat{I} \hookrightarrow \tcat{C}$, where $\tcat{I}$ is equivalent to
a 1-category (that is, to a locally discrete 2-category).  We will
refer to 2-cells in $\tcat{I}$ as \emph{identifying 2-cells}, and we
will refer to the quotient of $\tcat{I}$ by identifying 2-cells as the
\emph{underlying category} of $\tcat{C}$.

\begin{proposition}\label{prop:contfae}
  Let $\bS$ be a category with finite limits. The following
  2-categories equipped with underlying category structure are
  equivalent, inducing an isomorphism between underlying categories.
  \begin{enumerate}[label=(\roman*)]
  \item\label{item:continuousdef}
    \begin{itemize}
    \item 0-cells: comonads on $\bS$.
    \item 1-cells: continuous maps (\cref{def:continuous}).
    \item 2-cells: specializations (\cref{def:specialization}).
    \item Identifying 2-cells: identifications (\cref{def:identification}).
    \end{itemize}
    \vspace{1em}
    
  \item\label{item:continuouslift}
    \begin{itemize}
    \item 0-cells: comonads on $\bS$.
    \item 1-cells: commutative squares
      \[
        \begin{tikzcd}[column sep=25]
          \CS & \DS \\
          {\SC} & {\SD}
          \arrow["\ptlift{\carC}"', from=1-1, to=2-1]
          \arrow["{\con{\bS}{f}}"', dashed, from=1-2, to=1-1]
          \arrow["\ptlift{\carD}", from=1-2, to=2-2]
          \arrow["{\reindex{\one{f}}}", from=2-2, to=2-1]
        \end{tikzcd}
      \]
      where $\one{f} \colon \pts{\cC} \to \pts{\cD}$ is a map in
      $\bS$, the functor $\reindex{\one{f}} \colon \SD \to \SC$ is a
      choice of pullback functor along $\one{f}$, and
      $\con{\bS}{f} \colon \DS \to \CS$ is an arbitrary functor.

      (If $\cC$ and $\cD$ are sliceable, this is equivalently a
      \emph{colax comonad functor} $\clift{\cD} \to \clift{\cC}$
      carried by $\reindex{\one{f}}$; see \cref{prop:formaltfae}.)
    \item 2-cells: natural transformations, again just as in
      \cref{def:specialization}
      \[\begin{tikzcd}
          \CS & \DS
          \arrow[""{name=0, anchor=center, inner sep=0}, "{\con{\bS}{f}}"', curve={height=12pt}, from=1-2, to=1-1]
          \arrow[""{name=1, anchor=center, inner sep=0}, "{\con{\bS}{g}}", curve={height=-12pt}, from=1-2, to=1-1]
          \arrow["\Downarrow"{description}, draw=none, from=0, to=1]
        \end{tikzcd}\]

      (If $\cC$ and $\cD$ are sliceable, this is equivalently
      a \emph{comonad specialization}; see
      \cref{prop:formaltfae}.)
    \item Identifying 2-cells: cylinders
      \[
        \begin{tikzcd}
          \CS & \DS \\
          \SC & \SD
          \arrow["{\ptlift{\carC}}"', from=1-1, to=2-1]
          \arrow["{\con{\bS}{f}}"', curve={height=12pt}, from=1-2, to=1-1]
          \arrow["{\ptlift{\carD}}", from=1-2, to=2-2]
          \arrow[""{name=0, anchor=center, inner sep=0}, "{\reindex{\one{g}}}", curve={height=-12pt}, from=2-2, to=2-1]
          \arrow[""{name=1, anchor=center, inner sep=0}, "{\reindex{\one{f}}}"', curve={height=12pt}, from=2-2, to=2-1]
          \arrow["\cong"{description}, draw=none, from=1, to=0]
        \end{tikzcd}
        \qquad = \qquad
        \begin{tikzcd}
          \CS & \DS \\
          \SC & \SD
          \arrow["{\ptlift{\carC}}"', from=1-1, to=2-1]
          \arrow[""{name=0, anchor=center, inner sep=0}, "{\con{\bS}{f}}"', curve={height=12pt}, from=1-2, to=1-1]
          \arrow[""{name=1, anchor=center, inner sep=0}, "{\con{\bS}{g}}", curve={height=-12pt}, from=1-2, to=1-1]
          \arrow["{\ptlift{\carD}}", from=1-2, to=2-2]
          \arrow["{\reindex{\one{g}}}", curve={height=-12pt}, from=2-2, to=2-1]
          \arrow["\Downarrow"{description}, draw=none, from=0, to=1]
        \end{tikzcd}
      \]
      where $\one{f} = \one{g}$ and the bottom 2-cell $\cong$ in the left
      diagram is the canonical isomorphism between two choices of pullback
      functor.

      (If $\cC$ and $\cD$ are sliceable, this is equivalently
      a \emph{comonad transformation}, carried by $\cong$; see
      \cref{prop:formaltfae}.)
    \end{itemize}
  \end{enumerate}
\end{proposition}
\begin{proof}
  Note first that given a commutative square as in
  \labelcref{item:continuouslift}, $\con{\bS}{f}$ preserves the terminal
  object, since the terminal object is created by both
  $\ptlift{\carC}$ and $\ptlift{\carD}$ and preserved by
  $\reindex{\one{f}}$.
  
  We show such commutative squares are in correspondence, up to
  identification, with the triangles from \cref{def:continuous}.
  %
  %
  We may translate from squares to triangles like so:
  \[\begin{tikzcd}[column sep=10pt, row sep=10pt]
      \CS && \DS \\
      \\
      \SC && \SD \\
      & \bS
      \arrow["{\ptlift{\carC}}"', from=1-1, to=3-1]
      \arrow["{\con{\bS}{f}}"', dashed, from=1-3, to=1-1]
      \arrow["{\ptlift{\carD}}", from=1-3, to=3-3]
      \arrow["\slicel"', from=3-1, to=4-2]
      \arrow[""{name=0, anchor=center, inner sep=0}, "{\reindex{\one{f}}}"', from=3-3, to=3-1]
      \arrow["\slicel", from=3-3, to=4-2]
      \arrow["\Rightarrow"{description}, draw=none, from=0, to=4-2]
    \end{tikzcd}
  \]
  where the triangular 2-cell marked $\Rightarrow$ is the cartesian
  natural transformation mate to the identity
  $\slicel \circ \leindex{\one{f}} = \slicel$ given by the chosen
  pullback squares along $\one{f}$.

  Conversely, given a triangle defining a continuous map, we would
  like to obtain a square inducing it. Although such a square cannot
  always be found on the nose (the cartesian natural transformation
  may involve multiple choices of pullbacks in $\bS$, so that no
  particular pullback functor $\reindex{\one{f}}$ will do), we will be able to
  find one by working up to identification. Note that all functors
  $\DS \to \bS$ admitting a cartesian natural transformation $\alpha$
  into $\carD$ with component $\one{f}$ at the terminal object are
  canonically naturally isomorphic. In particular we may alter the
  domain of any such $\alpha$ by a canonical natural isomorphism to
  obtain a cartesian natural transformation given by any particular
  chosen pullback squares along $\one{f}$. Moreover, recall that
  $\carC$ is an amnestic isofibration, i.e.\ uniquely lifts
  isomorphisms along objects. We may thus alter $\con{\bS}{f}$ to lift
  any specified pullback functor $\reindex{\one{f}}$, via a unique
  identification. Hence for any particular pullback functor, the
  previous assignment from squares to triangles is essentially
  surjective.
  %

  The 2-cells of \labelcref{item:continuouslift} are then precisely
  specializations between underlying continuous maps.  Likewise, as we
  show next, the identifying 2-cells of
  \labelcref{item:continuouslift} correspond to identifications
  between the underlying continuous maps. Composing the triangular
  2-cell mate to identity (marked $\Rightarrow$) with the canonical
  natural isomorphism between pullback functors yields another
  triangular 2-cell mate to identity:
  \[
    \begin{tikzcd}[column sep=10pt, row sep=10pt]
      \CS && \DS \\
      \\
      \SC && \SD \\
      \\
      \\
      & \bS
      \arrow["{\ptlift{\carC}}"', from=1-1, to=3-1]
      \arrow["{\con{\bS}{f}}"', curve={height=12pt}, from=1-3, to=1-1]
      \arrow["{\ptlift{\carD}}", from=1-3, to=3-3]
      \arrow[""{name=0, anchor=center, inner sep=0}, "\slicel"', from=3-1, to=6-2]
      \arrow[""{name=1, anchor=center, inner sep=0}, "{\reindex{\one{g}}}", curve={height=-12pt}, from=3-3, to=3-1]
      \arrow[""{name=2, anchor=center, inner sep=0}, "{\reindex{\one{f}}}"', curve={height=12pt}, from=3-3, to=3-1]
      \arrow[""{name=3, anchor=center, inner sep=0}, "\slicel", from=3-3, to=6-2]
      \arrow["\Rightarrow"{description}, shift right=1, draw=none, from=0, to=3]
      \arrow["\cong"{description}, draw=none, from=2, to=1]
    \end{tikzcd}
    \qquad = \qquad
    \begin{tikzcd}[column sep=10pt, row sep=10pt]
      \CS && \DS \\
      \\
      \SC && \SD \\
      \\
      \\
      & \bS
      \arrow["{\ptlift{\carC}}"', from=1-1, to=3-1]
      \arrow["{\con{\bS}{f}}"', curve={height=12pt}, from=1-3, to=1-1]
      \arrow["{\ptlift{\carD}}", from=1-3, to=3-3]
      \arrow["\slicel"', from=3-1, to=6-2]
      \arrow[""{name=0, anchor=center, inner sep=0}, "{\reindex{\one{f}}}"', curve={height=12pt}, from=3-3, to=3-1]
      \arrow["\slicel", from=3-3, to=6-2]
      \arrow["\Rightarrow"{description}, draw=none, from=0, to=6-2]
    \end{tikzcd}
  \]
  Hence the cylinder axiom of \labelcref{item:continuouslift} yields
  the specialization axiom of \cref{def:specialization}. Conversely,
  from the specialization axiom we obtain (as in the proof of
  \cref{lem:conid}) that the natural transformation whose components
  are canonical isomorphisms between pullbacks is equal to the
  whiskering $\carC \circ \sigma$. The same equation is given by
  whiskering both sides of the cylinder axiom with
  $\slicel \colon \SC \to \bS$, which is faithful and so can be
  cancelled to obtain the cylinder axiom.
\end{proof}



\section{The double category of comonads}\label{sec:double}

Comonad maps and continuous maps amount to two kinds of triangles
between comonadic functors, shown respectively below:
\[
  \begin{tikzcd}[column sep=15pt]
    \CS\ar["\com{\bS}{\phi}", rr, dashed]\ar[dr, "\carC"']&&\DS\ar[dl, "\carD"]\\
    &\bS
  \end{tikzcd}
  \qquad\qquad\qquad
  \begin{tikzcd}[column sep=15pt]
    \CS && \DS \\
    & \bS
    \arrow["\con{\bS}{f}"', dashed, from=1-3, to=1-1]
    \arrow["\carC"', ""{name=0, anchor=center, inner sep=0}, from=1-1, to=2-2]
    \arrow["\carD", ""{name=1, anchor=center, inner sep=0}, from=1-3, to=2-2]
    \arrow["\Rightarrow", shift right=10pt, draw=none, from=1, to=0]
  \end{tikzcd}
\]
The two are related together in the structure of a double category. In
fact, just as there are specializations and identifications between
continuous maps of comonads on $\bS$, there are likewise two kinds of
double-categorical 2-cell relating comonad maps and continuous
maps. This gives us not just one but two notable double categories
$\Cmd(\bS)$ and $\CmdId(\bS)$ of comonads on $\bS$.

\begin{definition}
  A \emph{specialization} between comonad maps
  $\phi \colon \cA \to \cC$ and $\psi \colon \cB \to \cD$ and continuous maps
  $f \colon \cA \to \cB$ and $g \colon \cC \to \cD$ of comonads
  on $\bS$ is a natural transformation of the following form:
  
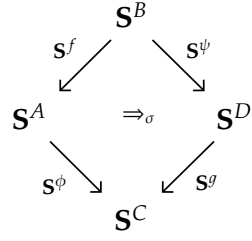
\[
    \begin{tikzcd}[column sep=15pt]
      & \BS & \\
      \AS && \DS \\
      & \CS
      \arrow["{\con{\bS}{f}}"', from=1-2, to=2-1]
      \arrow["{\com{\bS}{\psi}}", from=1-2, to=2-3]
      \arrow["{\Rightarrow_\sigma}"{description}, shift left, draw=none, from=1-2, to=3-2]
      \arrow["{\com{\bS}{\phi}}"', from=2-1, to=3-2]
      \arrow["{\con{\bS}{g}}", from=2-3, to=3-2]
    \end{tikzcd}
  \]
\end{definition}

\begin{definition}\label{def:doubleid}
  An \emph{identification} between comonad maps and continuous maps as
  above is a specialization $\sigma$ satisfying the following pyramid
  axiom:
  
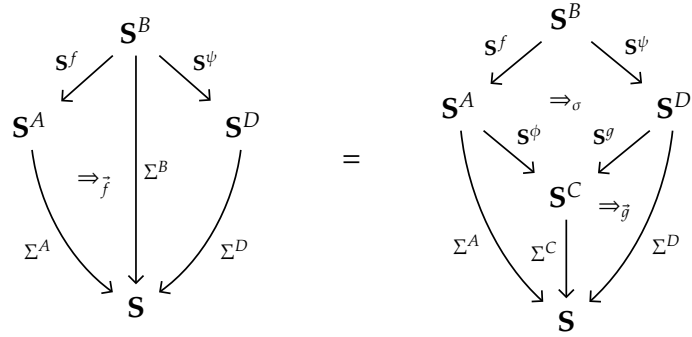
\[
    \begin{tikzcd}[row sep=17pt,column sep=17pt]
      & \BS & \\
      \AS && \DS \\
      \\
      \\
      & \bS
      \arrow[""{name=0, anchor=center, inner sep=0}, "{\con{\bS}{f}}"', from=1-2, to=2-1]
      \arrow["{\com{\bS}{\psi}}", from=1-2, to=2-3]
      \arrow["\carB", from=1-2, to=5-2]
      \arrow["\carA"', curve={height=12pt}, from=2-1, to=5-2]
      \arrow["\carD", curve={height=-12pt}, from=2-3, to=5-2]
      \arrow["{\Rightarrow_{\connat{f}}}"{description}, shift left=2, draw=none, from=5-2, to=0]
    \end{tikzcd}
    \qquad = \qquad
    \begin{tikzcd}[row sep=15pt,column sep=17pt]
      & \BS & \\
      \AS && \DS \\
      & \CS \\
      \\
      & \bS
      \arrow["{\con{\bS}{f}}"', from=1-2, to=2-1]
      \arrow["{\com{\bS}{\psi}}", from=1-2, to=2-3]
      \arrow["{\Rightarrow_\sigma}"{description}, shift left, draw=none, from=2-1, to=2-3]
      \arrow["{\com{\bS}{\phi}}", from=2-1, to=3-2]
      \arrow["\carA"', curve={height=12pt}, from=2-1, to=5-2]
      \arrow["{\con{\bS}{g}}"', from=2-3, to=3-2]
      \arrow["\carD", curve={height=-12pt}, from=2-3, to=5-2]
      \arrow["\carC"'{pos=0.4}, from=3-2, to=5-2]
      \arrow["{\Rightarrow_{\connat{g}}}"{description}, shift left, draw=none, from=5-2, to=2-3]
    \end{tikzcd}
  \]
\end{definition}

\begin{lemma}\label{lem:doubleone}
  For any identification as above, we have
  $\one{\psi} \circ \one{f} = \one{g} \circ \one{\phi}$.
\end{lemma}
\begin{proof}[Proof sketch]
  In fact, for any category $\bS$ there is a double
  category in which objects are functors into $\bS$, 1-cells are lax and
  colax triangles, and 2-cells are similar to identifications:
  
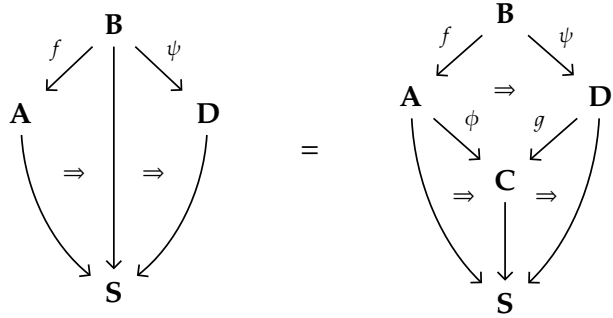
\[
    \begin{tikzcd}[row sep=17pt,column sep=17pt]
      & \bB & \\
      \bA && \bD \\
      \\
      \\
      & \bS
      \arrow[""{name=0, anchor=center, inner sep=0}, "f"', from=1-2, to=2-1]
      \arrow[""{name=1, anchor=center, inner sep=0}, "\psi", from=1-2, to=2-3]
      \arrow[from=1-2, to=5-2]
      \arrow[curve={height=12pt}, from=2-1, to=5-2]
      \arrow[curve={height=-12pt}, from=2-3, to=5-2]
      \arrow["\Rightarrow"{description}, shift left=2, draw=none, from=5-2, to=0]
      \arrow["\Rightarrow"{description}, shift right=2, draw=none, from=5-2, to=1]
    \end{tikzcd}
    \qquad = \qquad
    \begin{tikzcd}[row sep=15pt,column sep=17pt]
      & \bB & \\
      \bA && \bD \\
      & \bC \\
      \\
      & \bS
      \arrow["f"', from=1-2, to=2-1]
      \arrow["\psi", from=1-2, to=2-3]
      \arrow["\Rightarrow"{description}, shift left, draw=none, from=2-1, to=2-3]
      \arrow["\phi", from=2-1, to=3-2]
      \arrow[curve={height=12pt}, from=2-1, to=5-2]
      \arrow["g"', from=2-3, to=3-2]
      \arrow[curve={height=-12pt}, from=2-3, to=5-2]
      \arrow[from=3-2, to=5-2]
      \arrow["\Rightarrow"{description}, shift right, draw=none, from=5-2, to=2-1]
      \arrow["\Rightarrow"{description}, shift left, draw=none, from=5-2, to=2-3]
    \end{tikzcd}
  \]
  Within this double category we may consider the sub-double-category
  consisting of just those arrows into $\bS$ with colimits and just
  those lax triangles carried by final functors (see
  \cref{def:final}). Here the existence of a cell ensures the induced
  maps between colimits satisfy the desired equation.
  %
\end{proof}

\begin{proposition}\label{lem:thin}
  The double category of comonads, comonad maps and continuous maps,
  and identifications is thin, i.e.\ each boundary admits at most one
  2-cell. Moreover, both its underlying 2-categories are equivalent to
  1-categories (all globular 2-cells in this thin double category are
  invertible).
\end{proposition}
\begin{proof}
  As witnessed by the commutative diagram in the proof of
  \cref{lem:doubleone}, $\sigma_{\qpts{\cB}}$ is carried by the
  canonical map induced by the commutative square
  $\one{\psi} \circ \one{f} = \one{g} \circ \one{\phi}$ from
  $\qpts{A}$ into the pullback $\con{\bS}{g}(\qpts{G})$ of $\one{g}$
  and $\one{\psi}$. By pullback cancellation, the whiskering of
  $\sigma$ with $\carC$ is cartesian and determined by
  $\sigma_{\qpts{\cB}}$. Thus $\sigma$ itself is determined because
  $\carC$ is faithful.
  
  We have already seen that the 2-category of comonads on $\bS$,
  continuous maps, and identifications is equivalent to a
  1-category in \cref{lem:conid}. The 2-category of comonad maps and
  identifications is actually isomorphic to a 1-category: every
  identification of comonad maps is an identity, since comonadic
  functors are amnestic.
\end{proof}

\begin{definition}
  The \emph{double category of comonads on $\bS$} is the double
  category $\Cmd(\bS)$ obtained from the double category whose objects
  are comonads on $\bS$, 1-cells are comonad maps and continuous maps,
  and 2-cells are specializations, by quotienting out by
  identifications between parallel 1-cells (so that the underlying
  categories are correct).

  We write $\CmdId(\bS)$ for the identity-on-1-cells
  sub-double-category of $\Cmd(\bS)$ consisting of identifications.
\end{definition}

\begin{remark}
  Given a chosen terminal object of $\bS$, thus inducing chosen
  terminal objects of the categories $\CS$, we obtain a forgetful
  double functor from $\CmdId(\bS)$ to the double category of
  commutative squares in $\bS$ by \cref{lem:doubleone}.
\end{remark}

\begin{proposition}\label{prop:conjoint}
  The conjoint pairs in $\CmdId(\bS)$ are the comonad maps $\phi$ and
  continuous maps $f$ such that $\one{\phi} = \one{f}^{-1}$ (the
  isomorphism-on-points maps as in \cref{prop:boo}).

  A comonad map is a conjoint in $\Cmd(\bS)$ if and only if it is a
  conjoint in $\CmdId(\bS)$, and a continuous map is a conjoint in $\Cmd(\bS)$ if
  and only if it is isomorphic in $\Cmd(\bS)$ to a conjoint in $\CmdId(\bS)$.
\end{proposition}
\begin{proof}
  The conjoint laws for a comonad map $\phi \colon \cC \to \cD$
  and continuous map $f \colon \cD \to \cC$ in $\Cmd(\bS)$ specify an
  inverse pair of natural isomorphisms between $\com{\bS}{\phi}$ and
  $\con{\bS}{f}$.
  \[
    \begin{tikzcd}[column sep=0, row sep=30]
      \CS && \CS & \\
      & \DS && \DS
      \arrow[""{name=0, anchor=center, inner sep=0}, equals, from=1-1, to=1-3]
      \arrow["{\com{\bS}{\phi}}"', from=1-1, to=2-2]
      \arrow["{\con{\bS}{f}}"', pos=.3, from=1-3, to=2-2]
      \arrow["{\com{\bS}{\phi}}", from=1-3, to=2-4]
      \arrow[""{name=1, anchor=center, inner sep=0}, equals, from=2-4, to=2-2]
      \arrow["\Rightarrow"{description}, draw=none, from=0, to=2-2]
      \arrow["\Rightarrow"{description}, draw=none, from=1-3, to=1]
    \end{tikzcd}
    \!\!\!\!=
    \begin{tikzcd}[column sep=10, row sep=30]
      \CS & \CS \\
      \DS & \DS
      \arrow[shorten=-3pt, equals, from=1-1, to=1-2]
      \arrow[""{name=0, anchor=center, inner sep=0}, "{\com{\bS}{\phi}}"', from=1-1, to=2-1]
      \arrow[""{name=1, anchor=center, inner sep=0}, "{\com{\bS}{\phi}}", from=1-2, to=2-2]
      \arrow[shorten=-3pt, equals, from=2-2, to=2-1]
      \arrow["{=}"{description}, draw=none, from=0, to=1]
    \end{tikzcd}
    \qquad\qquad
    \begin{tikzcd}[column sep=0, row sep=30]
      & \CS && \CS \\
      \DS && \DS
      \arrow[""{name=0, anchor=center, inner sep=0}, equals, from=1-2, to=1-4]
      \arrow["{\con{\bS}{f}}"', from=1-2, to=2-1]
      \arrow["{\com{\bS}{\phi}}", pos=.3, from=1-2, to=2-3]
      \arrow["{\con{\bS}{f}}", from=1-4, to=2-3]
      \arrow[""{name=1, anchor=center, inner sep=0}, equals, from=2-1, to=2-3]
      \arrow["\Rightarrow"{description}, draw=none, from=0, to=2-3]
      \arrow["\Rightarrow"{description}, draw=none, from=1-2, to=1]
    \end{tikzcd}
    \!\!\!\!=
    \begin{tikzcd}[column sep=10, row sep=30]
      \CS & \CS \\
      \DS & \DS
      \arrow[shorten=-3pt, equals, from=1-1, to=1-2]
      \arrow[""{name=0, anchor=center, inner sep=0}, "{\con{\bS}{f}}"', from=1-1, to=2-1]
      \arrow[""{name=1, anchor=center, inner sep=0}, "{\con{\bS}{f}}", from=1-2, to=2-2]
      \arrow[shorten=-3pt, equals, from=2-2, to=2-1]
      \arrow["{=}"{description}, draw=none, from=0, to=1]
    \end{tikzcd}
  \]
  Thus $\phi$ is a conjoint in $\Cmd(\bS)$ if and only if
  $\com{\bS}{\phi}$ is naturally isomorphic to some $\con{\bS}{f}$,
  which is the case if and only if $\one{\phi}$ is {\final}, as in
  \cref{prop:boo}. In this case $\phi$ is conjoint in $\CmdId(\bS)$ to
  a unique-up-to-identification $f$ with $\one{\phi} =
  \one{f}^{-1}$. Therefore, since all conjoints of a given are 1-cell
  isomorphic, the conjoint continuous maps in $\Cmd(\bS)$ are those
  isomorphic in $\Cmd(\bS)$ to conjoints in $\CmdId(\bS)$.
\end{proof}

\begin{proposition}\label{prop:companionadjoint}
  The companion pairs in $\Cmd(\bS)$ are the comonad maps $\phi$ and
  continuous maps $f$ such that $\com{\bS}{\phi}$ is left adjoint to
  $\con{\bS}{f}$.
  
  The companion pairs in $\CmdId(\bS)$ are furthermore such that
  $\connat{f}$ (the {\cartesian} triangle for $f$) is mate to identity
  (the commutative triangle for $\phi$). In this case we have
  $\one{\phi} = \one{f}$.
\end{proposition}
\begin{proof}
  The companion laws for a comonad map $\phi \colon \cC \to \cD$
  and continuous map $f \colon \cD \to \cC$ in $\Cmd(\bS)$
  specify an adjunction between $\com{\bS}{\phi}$ and $\con{\bS}{f}$.
  \[
    \begin{tikzcd}[column sep=20, row sep=15]
      \CS & \\
      & \DS \\
      \CS \\
      & \DS
      \arrow["{\com{\bS}{\phi}}", from=1-1, to=2-2]
      \arrow[""{name=0, anchor=center, inner sep=0}, equals, from=1-1, to=3-1]
      \arrow["{\con{\bS}{f}}", from=2-2, to=3-1]
      \arrow["{\com{\bS}{\phi}}"', from=3-1, to=4-2]
      \arrow[""{name=1, anchor=center, inner sep=0}, equals, from=4-2, to=2-2]
      \arrow["\Rightarrow"{description}, draw=none, from=0, to=2-2]
      \arrow["\Rightarrow"{description}, shift right, draw=none, from=3-1, to=1]
    \end{tikzcd}
    =
    \begin{tikzcd}[column sep=20]
      \CS & \DS \\
      \CS & \DS
      \arrow[""{name=0, anchor=center, inner sep=0}, shorten=-3pt, "{\com{\bS}{\phi}}", from=1-1, to=1-2]
      \arrow[equals, from=1-2, to=2-2]
      \arrow[equals, from=2-1, to=1-1]
      \arrow[""{name=1, anchor=center, inner sep=0}, shorten=-3pt, "{\com{\bS}{\phi}}"', from=2-1, to=2-2]
      \arrow["{=}"{description}, draw=none, from=1, to=0]
    \end{tikzcd}
    \qquad\qquad
    \begin{tikzcd}[column sep=20, row sep=15]
      & \DS \\
      \CS \\
      & \DS \\
      \CS
      \arrow["{\con{\bS}{f}}"', from=1-2, to=2-1]
      \arrow["{\com{\bS}{\phi}}"', from=2-1, to=3-2]
      \arrow[""{name=0, anchor=center, inner sep=0}, equals, from=2-1, to=4-1]
      \arrow[""{name=1, anchor=center, inner sep=0}, equals, from=3-2, to=1-2]
      \arrow["{\con{\bS}{f}}", from=3-2, to=4-1]
      \arrow["\Rightarrow"{description}, shift right, draw=none, from=2-1, to=1]
      \arrow["\Rightarrow"{description}, draw=none, from=0, to=3-2]
    \end{tikzcd}
    =
    \begin{tikzcd}[column sep=20]
      \CS & \DS \\
      \CS & \DS
      \arrow[""{name=0, anchor=center, inner sep=0}, shorten=-3pt, "{\con{\bS}{f}}"', from=1-2, to=1-1]
      \arrow[equals, from=1-2, to=2-2]
      \arrow[equals, from=2-1, to=1-1]
      \arrow[""{name=1, anchor=center, inner sep=0}, shorten=-3pt, "{\con{\bS}{f}}", from=2-2, to=2-1]
      \arrow["{=}"{description}, draw=none, from=1, to=0]
    \end{tikzcd}
  \]
  The pyramid laws for identifications express that the triangles
  associated to $\com{\bS}{\phi}$ and $\con{\bS}{f}$ are
  mate. We have $\one{\phi} = \one{f}$ by \cref{lem:doubleone}.
\end{proof}

We can also describe the 2-cells of $\CmdId(\bS)$ in terms of the
slice factorizations of the comonadic functors, similarly to
\cref{prop:contfae}.

\begin{proposition}\label{prop:cmdtfae}
  Let $\bS$ be a category with finite limits.  Suppose given the
  following comonad maps $\phi$ and $\psi$ and continuous maps $f$ and
  $g$ between comonads on $\bS$ represented as in \cref{prop:comtfae}
  and \cref{prop:contfae}:
  \[
    \begin{tikzcd}[column sep=12.5]
      \AS & \CS \\
      {\SA} & {\SC}
      \arrow["{\com{\bS}{\phi}}", dashed, from=1-1, to=1-2]
      \arrow["{\ptlift{\carA}}"', from=1-1, to=2-1]
      \arrow["{\ptlift{\carC}}", from=1-2, to=2-2]
      \arrow["{\leindex{\one{\phi}}}"', from=2-1, to=2-2]
    \end{tikzcd}
    \hspace{6.5mm}
    \begin{tikzcd}[column sep=12.5]
      \BS & \DS \\
      {\SB} & {\SD}
      \arrow["{\com{\bS}{\psi}}", dashed, from=1-1, to=1-2]
      \arrow["{\ptlift{\carB}}"', from=1-1, to=2-1]
      \arrow["{\ptlift{\carD}}", from=1-2, to=2-2]
      \arrow["{\leindex{\one{\psi}}}"', from=2-1, to=2-2]
    \end{tikzcd}
    \hspace{6.5mm}
    \begin{tikzcd}[column sep=12.5]
      \AS & \BS \\
      {\SA} & {\SB}
      \arrow["\ptlift{\carA}"', from=1-1, to=2-1]
      \arrow["{\con{\bS}{f}}"', dashed, from=1-2, to=1-1]
      \arrow["\ptlift{\carB}", from=1-2, to=2-2]
      \arrow["{\reindex{\one{f}}}", from=2-2, to=2-1]
    \end{tikzcd}
    \hspace{6.5mm}
    \begin{tikzcd}[column sep=12.5]
      \CS & \DS \\
      {\SC} & {\SD}
      \arrow["\ptlift{\carC}"', from=1-1, to=2-1]
      \arrow["{\con{\bS}{g}}"', dashed, from=1-2, to=1-1]
      \arrow["\ptlift{\carD}", from=1-2, to=2-2]
      \arrow["{\reindex{\one{g}}}", from=2-2, to=2-1]
    \end{tikzcd}
  \]
  Then identifications between these comonad maps and continuous maps
  (\cref{def:doubleid}) are in one-to-one correspondence with
  cubes
  \[
    \begin{tikzcd}[column sep=17]
      & \BS & \\
      \AS & \SB & \DS \\
      \SA && \SD \\
      & \SC
      \arrow["{\con{\bS}{f}}"'{pos=0.4}, from=1-2, to=2-1]
      \arrow["{\ptlift{\carB}}"{pos=0.4}, from=1-2, to=2-2]
      \arrow["{\com{\bS}{\psi}}"{pos=0.4}, from=1-2, to=2-3]
      \arrow["{\ptlift{\carA}}"'{pos=0.4}, from=2-1, to=3-1]
      \arrow["{\reindex{\one{f}}}"'{pos=0.4}, from=2-2, to=3-1]
      \arrow["{\leindex{\one{\psi}}}"{pos=0.4}, from=2-2, to=3-3]
      \arrow["\Rightarrow"{description}, draw=none, from=2-2, to=4-2]
      \arrow["{\ptlift{\carD}}"{pos=0.4}, from=2-3, to=3-3]
      \arrow["{\leindex{\one{\phi}}}"', from=3-1, to=4-2]
      \arrow["{\reindex{\one{g}}}", from=3-3, to=4-2]
    \end{tikzcd}
    \qquad = \qquad
    \begin{tikzcd}[column sep=17]
      & \BS & \\
      \AS && \DS \\
      \SA & \CS & \SD \\
      & \SC
      \arrow["{\con{\bS}{f}}"'{pos=0.4}, from=1-2, to=2-1]
      \arrow["{\com{\bS}{\psi}}"{pos=0.4}, from=1-2, to=2-3]
      \arrow["{\Rightarrow_\sigma}"{description}, shift left, draw=none, from=1-2, to=3-2]
      \arrow["{\ptlift{\carA}}"'{pos=0.4}, from=2-1, to=3-1]
      \arrow["{\com{\bS}{\phi}}"'{pos=0.6}, from=2-1, to=3-2]
      \arrow["{\con{\bS}{g}}"{pos=0.6}, from=2-3, to=3-2]
      \arrow["{\ptlift{\carD}}"{pos=0.4}, from=2-3, to=3-3]
      \arrow["{\leindex{\one{\phi}}}"', from=3-1, to=4-2]
      \arrow["{\ptlift{\carC}}"{pos=0.4}, from=3-2, to=4-2]
      \arrow["{\reindex{\one{g}}}", from=3-3, to=4-2]
    \end{tikzcd}
  \]
  

  \noindent
  where $\one{\psi} \circ \one{f} = \one{g} \circ \one{\phi}$ and
  the bottom 2-cell
  $\leindex{\one{\phi}} \circ \reindex{\one{f}} \Rightarrow \reindex{\one{g}}
  \circ \leindex{\one{\psi}}$ in the left diagram is the mate to the
  identity
  $\leindex{\one{g}} \circ \leindex{\one{\phi}} = \leindex{\one{\psi}} \circ
  \leindex{\one{f}}$.
\end{proposition}
\begin{proof}
  We can translate from cubes to pyramids like so:
  
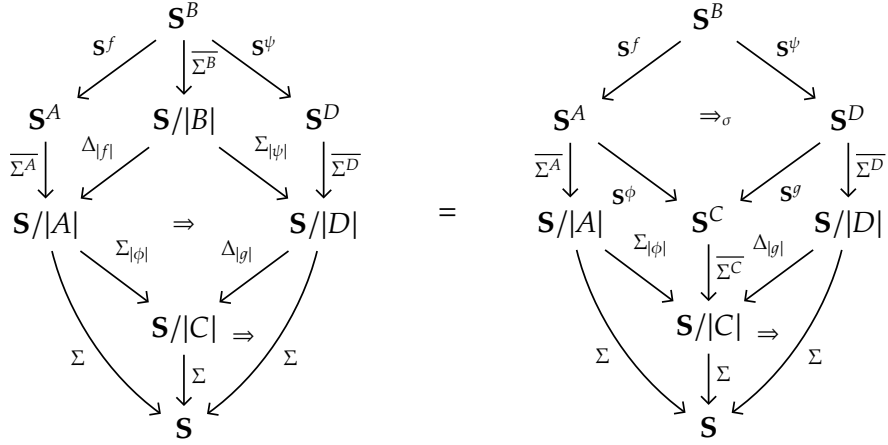
\[
    \begin{tikzcd}[column sep=17]
      & \BS & \\
      \AS & \SB & \DS \\
      \SA && \SD \\
      & \SC \\
      & \bS
      \arrow["{\con{\bS}{f}}"'{pos=0.4}, from=1-2, to=2-1]
      \arrow["{\ptlift{\carB}}"{pos=0.4}, from=1-2, to=2-2]
      \arrow["{\com{\bS}{\psi}}"{pos=0.4}, from=1-2, to=2-3]
      \arrow["{\ptlift{\carA}}"'{pos=0.4}, from=2-1, to=3-1]
      \arrow["{\reindex{\one{f}}}"'{pos=0.4}, from=2-2, to=3-1]
      \arrow["{\leindex{\one{\psi}}}"{pos=0.4}, from=2-2, to=3-3]
      \arrow["\Rightarrow"{description}, draw=none, from=2-2, to=4-2]
      \arrow["{\ptlift{\carD}}"{pos=0.4}, from=2-3, to=3-3]
      \arrow["{\leindex{\one{\phi}}}"{pos=0.4}, from=3-1, to=4-2]
      \arrow["\carr"', curve={height=12pt}, from=3-1, to=5-2]
      \arrow[""{name=0, anchor=center, inner sep=0}, "{\reindex{\one{g}}}"'{pos=0.4}, from=3-3, to=4-2]
      \arrow["\carr", curve={height=-12pt}, from=3-3, to=5-2]
      \arrow["\carr"{pos=0.4}, from=4-2, to=5-2]
      \arrow["{\Rightarrow}"{description, pos=0.6}, shift right=2, draw=none, from=5-2, to=0]
    \end{tikzcd}
    \qquad = \qquad
    \begin{tikzcd}[column sep=17]
      & \BS & \\
      \AS && \DS \\
      \SA & \CS & \SD \\
      & \SC \\
      & \bS
      \arrow["{\con{\bS}{f}}"'{pos=0.4}, from=1-2, to=2-1]
      \arrow["{\com{\bS}{\psi}}"{pos=0.4}, from=1-2, to=2-3]
      \arrow["{\Rightarrow_\sigma}"{description}, shift left, draw=none, from=1-2, to=3-2]
      \arrow["{\ptlift{\carA}}"'{pos=0.4}, from=2-1, to=3-1]
      \arrow["{\com{\bS}{\phi}}"'{pos=0.6}, from=2-1, to=3-2]
      \arrow["{\con{\bS}{g}}"{pos=0.6}, from=2-3, to=3-2]
      \arrow["{\ptlift{\carD}}"{pos=0.4}, from=2-3, to=3-3]
      \arrow["{\leindex{\one{\phi}}}"{pos=0.3}, from=3-1, to=4-2]
      \arrow["\carr"', curve={height=12pt}, from=3-1, to=5-2]
      \arrow["{\ptlift{\carC}}"{pos=0.4}, from=3-2, to=4-2]
      \arrow[""{name=0, anchor=center, inner sep=0}, "{\reindex{\one{g}}}"'{pos=0.3}, from=3-3, to=4-2]
      \arrow["\carr", curve={height=-12pt}, from=3-3, to=5-2]
      \arrow["\carr"{pos=0.4}, from=4-2, to=5-2]
      \arrow["{\Rightarrow}"{description, pos=0.6}, shift right=2, draw=none, from=5-2, to=0]
    \end{tikzcd}
  \]

  \noindent
  where the cells marked simply $\Rightarrow$ are mates of identities;
  they compose to yield more mates of identities. This allows us to
  rewrite the left diagram as follows:
  
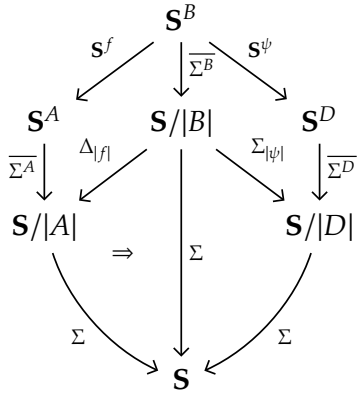
\[\begin{tikzcd}[column sep=17]
      & \BS & \\
      \AS & \SB & \DS \\
      \SA && \SD \\
      \\
      & \bS
      \arrow["{\con{\bS}{f}}"'{pos=0.4}, from=1-2, to=2-1]
      \arrow["{\ptlift{\carB}}"{pos=0.4}, from=1-2, to=2-2]
      \arrow["{\com{\bS}{\psi}}"{pos=0.4}, from=1-2, to=2-3]
      \arrow["{\ptlift{\carA}}"'{pos=0.4}, from=2-1, to=3-1]
      \arrow[""{name=0, anchor=center, inner sep=0}, "{\reindex{\one{f}}}"'{pos=0.4}, from=2-2, to=3-1]
      \arrow["{\leindex{\one{\psi}}}"{pos=0.4}, from=2-2, to=3-3]
      \arrow["\carr", from=2-2, to=5-2]
      \arrow["{\ptlift{\carD}}"{pos=0.4}, from=2-3, to=3-3]
      \arrow["\carr"', curve={height=12pt}, from=3-1, to=5-2]
      \arrow["\carr", curve={height=-12pt}, from=3-3, to=5-2]
      \arrow["\Rightarrow"{description, pos=0.6}, shift left=2, draw=none, from=5-2, to=0]
    \end{tikzcd}
  \]
  

  \noindent  
  Applying the correspondence of \cref{prop:contfae} for continuous
  maps, we now see that the original cube axiom gives the pyramid
  axiom from \cref{def:doubleid}. Conversely, from the pyramid axiom
  we obtain (as in the proof of \cref{lem:thin}) that the canonical
  cartesian natural transformation of its type induced by the
  commutative square
  $\one{\psi} \circ \one{f} = \one{g} \circ \one{\phi}$ is equal to
  the whiskering $\carC \circ \sigma$. The same equation is given by
  whiskering both sides of the cube axiom with
  $\slicel \colon \SC \to \bS$, which is faithful and so can be
  cancelled to obtain the cylinder axiom.
\end{proof}

\begin{remark}\label{rem:doubleslice}
  The double categories $\Cmd(\bS)$ and $\CmdId(\bS)$ --- at least,
  restricted to the sliceable comonads and assuming $\bS$ has finite
  limits --- can also be interpreted through the lens of the formal
  theory of (co)monads~\cite{street:monads,lack-street}. By
  \cref{prop:formaltfae}, a specialization between comonad maps and
  continuous maps between such sliceable comonads is a square in the
  2-category of colax comonad functors (\cref{def:comfun}) and comonad
  specializations (\cref{def:comspec}) involving the corresponding
  sliced comonads. An identification between comonad maps and
  continuous maps between sliceable comonads in the above cube form,
  which can also be rewritten as
  \[
    \begin{tikzcd}[column sep=17, row sep=15]
      & \AS & \\
      \CS & \SA & \BS \\
      \SC && \SB \\
      & \SD
      \arrow["{\com{\bS}{\phi}}"', pos=.4, from=1-2, to=2-1]
      \arrow["{\ptlift{\carA}}", from=1-2, to=2-2]
      \arrow["{\ptlift{\carC}}"', from=2-1, to=3-1]
      \arrow["{\leindex{\one{\phi}}}"', pos=.4, from=2-2, to=3-1]
      \arrow["\Downarrow"{description}, draw=none, from=2-2, to=4-2]
      \arrow["{\con{\bS}{f}}"', from=2-3, to=1-2]
      \arrow["{\ptlift{\carB}}", from=2-3, to=3-3]
      \arrow["{\reindex{\one{f}}}"', from=3-3, to=2-2]
      \arrow["{\leindex{\one{\psi}}}", from=3-3, to=4-2]
      \arrow["{\reindex{\one{g}}}", pos=.4, from=4-2, to=3-1]
    \end{tikzcd}
    \quad=\quad
    \begin{tikzcd}[column sep=17, row sep=15]
      & \AS & \\
      \CS && \BS \\
      \SC & \DS & \SB \\
      & \SD
      \arrow["{\com{\bS}{\phi}}"', pos=.4, from=1-2, to=2-1]
      \arrow["{\Downarrow_\sigma}"{description}, shift left, draw=none, from=1-2, to=3-2]
      \arrow["{\ptlift{\carC}}"', from=2-1, to=3-1]
      \arrow["{\con{\bS}{f}}"', from=2-3, to=1-2]
      \arrow["{\com{\bS}{\psi}}", pos=.55, from=2-3, to=3-2]
      \arrow["{\ptlift{\carB}}", from=2-3, to=3-3]
      \arrow["{\con{\bS}{g}}", pos=.4, from=3-2, to=2-1]
      \arrow["{\ptlift{\carD}}", pos=.4, from=3-2, to=4-2]
      \arrow["{\leindex{\one{\psi}}}", from=3-3, to=4-2]
      \arrow["{\reindex{\one{g}}}", pos=.4, from=4-2, to=3-1]
    \end{tikzcd}
  \]
  is a square in the 2-category of colax comonad functors and comonad
  transformations (\cref{def:comtrans}) carried by the canonical
  2-cell
  $\leindex{\one{\phi}} \circ \reindex{\one{f}} \Rightarrow
  \reindex{\one{g}} \circ \leindex{\one{\psi}}$.
  %
  %
\end{remark}

\chapter{Categories without 1}\label{app:comprehensive}
\appendixtocspace

Any functor $F \colon \bC \to \bD$ where $\bC$ has a terminal object
admits a unique factorization
\[\begin{tikzcd}
    \bC & {\bD/\pts{F}} & \bD \arrow["{\ptlift{F}}", from=1-1,
    to=1-2] \arrow["{\slicel}", from=1-2, to=1-3]
  \end{tikzcd}\] such that $\ptlift{F}$ preserves the terminal object
and $\slicel$ is the projection functor from the slice category
$\bD/\pts{F}$. This is an instance of the \emph{slice factorization}
(\cref{def:slice}).

In fact, we obtain a factorization system on the full subcategory of
$\Cat$ consisting of categories with terminal objects: the left class
consists of those functors preserving the terminal object, and the
right class consists of those functors given by slice category
projections (up to composing an isomorphism on the domain). Perhaps
surprisingly, this factorization system furthermore admits a canonical
extension to all of $\Cat$, known as the \emph{comprehensive
  factorization
  system}~\cite{street-walters:comprehensive}.\footnote{The reader is
  warned that there is an error in the proof of~\cite[Theorem
  3]{street-walters:comprehensive}, as noted in \cite[Section
  4.7]{kelly:enriched}.}

In the generality that $\bC$ lacks a terminal object, we may simply
apply the same process, not to $F$ itself, but to the induced functor
between free cocompletions\footnote{In the case that $\bC$ is large,
  we would here use the free large cocompletions instead of the usual
  free small cocompletions.}
$\Coco{F}\colon \Coco{\bC} \to \Coco{\bD}$.
\[\begin{tikzcd}
    \Coco{\bC} & {\Coco{\bD}/\fpts{F}} & \Coco{\bD}
    \arrow["{\ptlift{\Coco{F}}}", from=1-1, to=1-2]
    \arrow["{\slicel}", from=1-2, to=1-3]
  \end{tikzcd}\]

Here we have denoted the image of the terminal object by
$\fpts{F} \coloneqq \pts{\Coco{F}}$ for convenience.  Let us
furthermore denote by $\bD\sslash\fpts{F}$ the full subcategory of
$\Coco{\bD}/\fpts{F}$ consisting of objects with representable domain,
i.e.\ the comma category $(\yo/\fpts{F})$ where
$\yo \colon \bD \hookrightarrow \Ps{\bD}$ denotes the Yoneda embedding, or
equivalently the category of elements $\El(\fpts{F})$ of
$\fpts{F}$. (We use this notation $\bD\sslash\fpts{F}$ to evoke the
more familiar situation in which $\bC$ has a terminal object, in which
case $\bD\sslash\fpts{F} \cong \bD/\pts{F}$.) Now restricting
$\Coco{\bC}$ and $\Coco{\bD}$ back to their full subcategories of
representables $\bC$ and $\bD$, we have thus factored our original
functor $F$ as
\[
  \begin{tikzcd}
    \bC & \bD\sslash\fpts{F} & \bD \arrow["{\fptlift{F}}", from=1-1,
    to=1-2] \arrow["{\sliced}", from=1-2, to=1-3]
  \end{tikzcd}
\]
where $\fptlift{F}$ and $\sliced$ are the functors obtained by
restricting $\ptlift{\Coco{F}}$ and $\slicel$ above.

\begin{definition}\label{def:comprehensive}
  The above factorization of $F \colon \bC \to \bD$ is its
  \emph{comprehensive factorization}.
\end{definition}

It is evident that the comprehensive factorization agrees with the
slice factorization in the case that a terminal object of $\bC$ does
exist. Also, recall that any slice category of a presheaf category is
a presheaf category --- specifically
$\Coco{\bC}/X \simeq \Coco{\El(X)}$ --- thus the slice category
$\Coco{\bD}/\fpts{F}$ is itself (up to equivalence) the free
cocompletion of $\bD\sslash\fpts{F}$, the category of elements of
$\fpts{F}$. Since $\Coco{F}$ preserves and $\slicel$ creates colimits,
$\ptlift{\Coco{F}}$ preserves colimits as well. Thus
$\ptlift{\Coco{F}}$ is, up to equivalence, the functor
$\Coco{\fptlift{F}}$ between free cocompletions induced by
$\fptlift{F}$, and the functor $\slicel$ is, up to equivalence, the
functor $\Coco{\sliced}$ between free cocompletions induced by
$\sliced$. In this sense taking comprehensive factorizations commutes
with taking free cocompletions.

We will confirm that this determines a factorization system on $\Cat$
shortly. The left class consists of the \emph{final} functors, and the
right class consists of the \emph{discrete fibrations}.

\begin{definition}\label{def:final}
  A functor $F \colon \bC \to \bD$ is called \emph{final} if it
  satisfies any of the following equivalent conditions.
  \begin{enumerate}[label=(\roman*)]
  \item\label{item:finalcomp} The second factor
    $\sliced \colon \bD\sslash\fpts{F} \to \bD$ in the comprehensive
    factorization of $F$ is an isomorphism.
  \item\label{item:finalterm} The induced functor between free
    cocompletions $\Coco{F}\colon \Coco{\bC} \to \Coco{\bD}$ preserves
    the terminal object.
  \item \label{item:finalconn} The comma category $(d/F)$ is connected
    for each object $d$ of $\bD$.
  \item\label{item:finalcolim} For any functor $J \colon \bD \to \bE$,
    the induced mapping from cocones under the diagram $J$ to cocones
    under the diagram $J \circ F$ is bijective. (Thus, precomposing
    $F$ with any diagram $J$ leaves the colimit unchanged.)
  \end{enumerate}
\end{definition}

The equivalence \labelcref{item:finalcomp} $\iff$
\labelcref{item:finalterm} is evident; \labelcref{item:finalterm}
$\iff$ \labelcref{item:finalconn} is also relatively straightforward
(and in particular follows from the general theory of flatness for
limit doctrines~\cite{adamek-borceux-lack-rosicky}, whereby another
name for final would be \emph{representably TERM-flat});
\labelcref{item:finalconn} $\iff$ \labelcref{item:finalcolim} is given
in \cite[Proposition 2.5.2]{kashiwara-schapira}. If $\bC$ has a
terminal object, then $F$ is final if and only if it preserves the
terminal object. On the other hand, when $\bD$ has a terminal object,
$F$ is final if and only if its colimit is both terminal and absolute.


\begin{definition}\label{def:df}
  A functor $F \colon \bC \to \bD$ is called a \emph{discrete
    fibration} if it satisfies any of the following equivalent
  conditions.
  \begin{enumerate}[label=(\roman*)]
  \item\label{item:dfcomp} The first factor
    $\fptlift{F} \colon \bD\sslash\fpts{F} \to \bD$ in the
    comprehensive factorization of $F$ is an isomorphism.
  \item\label{item:dfel} $F$ is (up to composing with an isomorphism
    on the domain) the projection functor $\El(\xA) \to \bD$ from the
    category of elements $\bC \cong \El(\xA)$ of a presheaf $\xA$ over
    $\bD$.
  \item\label{item:dfconcrete} For each object $c$ in $\bC$ and
    arrow $f \colon d \to F(c)$ in $\bD$, there is a unique arrow
    $f' \colon d' \to c$ in $\bC$ such that $F(f') = f$.
  \item\label{item:dffib} $F$ is a fibration (in the standard sense of
    e.g.\ 
    \cite[Definition 1.1.3]{jacobs:logic}) such that each fiber is a
    discrete category.
  \end{enumerate}
\end{definition}

The equivalence \labelcref{item:dfcomp} $\iff$ \labelcref{item:dfel}
is evident; \labelcref{item:dfel} $\iff$ \labelcref{item:dfconcrete}
$\iff$ \labelcref{item:dffib} is given in \cite[Section
3.4]{vistoli}. If $\bC$ has a terminal object, then $F$ is a discrete
fibration if and only if it amounts to the projection from a slice
category $\bC \cong \bD/d \to \bD$. It is also worth noting that
whereas the property of being final is invariant under equivalence,
the property of being a discrete fibration is not; in fact, every
discrete fibration is an amnestic isofibration
(\cref{def:amnesticisofibration}).

\begin{proposition}[{\cite{street-walters:comprehensive},
    \cite[Section 4.7]{kelly:enriched}}]
  Final functors and discrete fibrations form a factorization system
  on $\Cat$ (as well as $\CAT$).
\end{proposition}
\begin{proof}
  Both classes are evidently closed under composition and include all
  isomorphisms. It remains to see that the comprehensive factorization
  of a functor $F \colon \bC \to \bD$ from \cref{def:comprehensive}
  is, up to isomorphism, the unique factorization into a final functor
  followed by a discrete fibration. Indeed, any functors
  $\bC \xto{F_1} \bE \xto{F_2} \bD$ factorizing $F$ where $F_1$ is
  final and $F_2$ is a discrete fibration extend to functors between
  free cocompletions
  $\Coco{\bC} \xto{\Coco{F_1}} \Coco{\bE} \xto{\Coco{F_2}} \Coco{\bD}$
  factorizing $\Coco{F}$. Since $F_1$ is final, $\Coco{F_1}$ preserves
  the terminal object. Since $F_2$ is a discrete fibration,
  $\Coco{F_2}$ has slice factorization whose first factor is an
  equivalence $\Coco{\bE} \simeq \Coco{\bD}/\xA$ that is moreover an
  isomorphism restricted to the full subcategory $\bE$ of
  representables, and whose second factor is the projection
  $\Coco{\bD}/\xA \to \Coco{\bD}$, where $\xA$ is the presheaf
  corresponding to the discrete fibration $F_2$. Composing
  $\Coco{F_1}$ with $\Coco{\bE} \simeq \Coco{\bD}/\xA$ yields the
  slice factorization of $\Coco{F}$, from which the comprehensive
  factorization is obtained by restricting to representables, by
  definition. Hence the given factorization is up to isomorphism the
  comprehensive factorization.
\end{proof}





Recall that a natural transformation is \emph{cartesian} if all
naturality squares are pullback squares. The following stronger
definition interacts more robustly with categories lacking terminal
objects.

\begin{definition}\label{def:cartesian}
  A natural transformation
  $\alpha \colon F \Rightarrow G \colon \bC \to \bD$ is
  \emph{uniformly cartesian} if the induced natural transformation
  $\Coco{\alpha} \colon \Coco{F} \Rightarrow \Coco{G} \colon
  \Coco{\bC} \to \Coco{\bD}$ of functors between free cocompletions is
  cartesian.\footnote{Again, in case $\bC$ is large we instead use
    free large cocompletions.}
\end{definition}

The definition of continuous map of comonads (\cref{def:continuous})
may be generalized to categories $\bS$ lacking a terminal object by
replacing the terminal-object-preserving functor $\con{\bS}{f}$ with a
final functor, and replacing the cartesian natural transformation
$\connat{f}$ with a uniformly cartesian natural transformation. (We
seem to need the uniformly cartesian condition in particular to
generalize our proofs of \cref{lem:conid} and
\cref{lem:thin}.\anoteNI{Are there counterexamples without this
  assumption?})

To establish a natural transformation $\alpha$ is uniformly cartesian,
it suffices to check that the naturality squares of $\Coco{\alpha}$ at
arrows from representables into the terminal object of $\Coco{\bC}$
are pullbacks. Indeed, this implies that $\Coco{\alpha}$ is cartesian
restricted to the representables, by pullback cancellation; moreover,
every object of $\Coco{\bC}$ is a canonical colimit of representables,
and colimits in $\Coco{\bD}$ are stable under pullback. Thus the
functor $\Coco{\bC} \to \Coco{\bD}$ obtained from $\Coco{G}$ by
forming pullbacks along the component
$\fpts{\alpha} \colon \fpts{F} \to \fpts{G}$ of $\Coco{\alpha}$ at the
terminal object is the colimit-preserving functor $\Coco{F}$ extending
$F$; and the thereby induced cartesian natural transformation
$\Coco{F} \Rightarrow \Coco{G}$ is $\Coco{\alpha}$, since any such
natural transformation is determined by its restriction
$F \Rightarrow G$, here $\alpha$.

When we assume the domain $\bC$ has a terminal object, just as the
notions of terminal-object-preserving and final become equivalent, the
notions of cartesian and uniformly cartesian become equivalent. This
follows from the previous paragraph, minding that the Yoneda embedding
preserves limits. And in general, uniformly cartesian implies
cartesian, minding that the Yoneda embedding reflects limits since it
is fully faithful.

Again in the case that $\bC$ has a terminal object, a cartesian
natural transformation $\alpha \colon F \Rightarrow G$ expresses
precisely that the diagram $F$ is obtained by reindexing the diagram
$G$ along the component of $\alpha$ at the terminal object
$\pts{\alpha} \colon \pts{F} \to \pts{G}$.
When we wish to generalize this description to domains that lack a
terminal object, the following even further strengthening of cartesian
is appropriate; as we note below, it ends up being the same as
uniformly cartesian in sufficiently well-behaved cases. We will refer
to this definition once in \cref{app:universal}.

\begin{definition}\label{def:strongcartesian}
  Let $F, G \colon C \to D$ be functors admitting colimits $\pts{F}$
  and $\pts{G}$, and let $\alpha \colon F \Rightarrow G$ be a natural
  transformation. Denote by $C^*$ the category obtained from $C$ by
  freely adjoining a terminal object, denote by $F^*$ and $G^*$ the
  functors given by the colimiting cocones of $F$ and $G$
  (respectively sending the terminal object to $\pts{F}$ and
  $\pts{G}$), and denote by $\alpha^* \colon F^* \Rightarrow G^*$ the
  natural transformation induced by $\alpha$. We say that $\alpha$ is
  \emph{strongly uniformly cartesian} if $\alpha^*$ is cartesian.

  Equivalently, this means that $\alpha$ is of the form
  \[\begin{tikzcd}[row sep=20,column sep=5pt]
      & \bC & \\
      {\bD/\pts{F}} && {\bD/\pts{G}} \\
      & \bD
      \arrow["{\ptlift{F}}"', from=1-2, to=2-1]
      \arrow["{\ptlift{G}}", from=1-2, to=2-3]
      \arrow["\slicel"', from=2-1, to=3-2]
      \arrow[""{name=0, anchor=center, inner sep=0}, "{\Delta_{\pts{\alpha}}}"', from=2-3, to=2-1]
      \arrow["\slicel", from=2-3, to=3-2]
      \arrow["\Rightarrow"{description}, draw=none, from=0, to=3-2]
    \end{tikzcd}\]
  up to altering $F$ by a natural isomorphism, where $\Rightarrow$
  denotes the canonical cartesian natural transformation mate to
  identity.
\end{definition}

Strongly uniformly cartesian implies uniformly cartesian. Indeed, if
$\alpha^*$ is cartesian, then it is uniformly cartesian since $\bC^*$
has a terminal object; it follows that $\alpha$ is also uniformly
cartesian, since restricting the domain of a (uniformly) cartesian
natural transformation evidently yields a (uniformly) cartesian
natural transformation.

Conversely, if we assume the codomain $\bD$ is a Grothendieck topos
such as $\Set$, then a natural transformation between
(small\footnote{If the domain $\bC$ is large, then we may apply the
  same argument in a larger universe $\SET$ instead, replacing the
  Grothendieck topos $\bD$ with the larger topos consisting of sheaves
  on the same small site valued in $\SET$. In fact the fully faithful
  inclusion from $\bD$ into such an enlargement preserves all limits
  and colimits that happen to exist in $\bD$, even large ones, because
  $\bD$ has a strong separator and a strong coseparator preserved by
  the inclusion. The universe enlargement likewise respects uniformly
  cartesian natural transformations. Therefore the original assumption
  that the diagrams are small is unnecessary, and we may replace it
  with the minimal assumption that the diagrams admit colimits.})
diagrams is uniformly cartesian if and only if it is strongly
uniformly cartesian. Indeed, in this case the Yoneda embedding
$\yo \colon \bD \to \Coco{\bD}$ admits a finite-limit-preserving left
adjoint. This left adjoint maps the pullback naturality squares in
$\Coco{\bD}$ witnessing that $\alpha$ is uniformly cartesian, shown
below left, to pullback naturality squares in $\bD$ witnessing that
$\alpha$ is strongly uniformly cartesian, shown below right.
\[
  \begin{tikzcd}[column sep=40pt]
    {\yo \circ F(c)} & {\yo \circ G(c)} \\
    {\fpts{F}} & {\fpts{G}}
    \arrow["{\Coco{\alpha}_{\yo_c}}", from=1-1, to=1-2]
    \arrow[from=1-1, to=2-1]
    \arrow["\lrcorner"{anchor=center, pos=0.125}, draw=none, from=1-1, to=2-2]
    \arrow[from=1-2, to=2-2]
    \arrow["{\fpts{\alpha}}"', from=2-1, to=2-2]
  \end{tikzcd}
  \qquad\qquad
  \begin{tikzcd}[column sep=40pt]
    {F(c)} & {G(c)} \\
    {\pts{F}} & {\pts{G}}
    \arrow["{\alpha_c}", from=1-1, to=1-2]
    \arrow[from=1-1, to=2-1]
    \arrow["\lrcorner"{anchor=center, pos=0.125}, draw=none, from=1-1, to=2-2]
    \arrow[from=1-2, to=2-2]
    \arrow["{\pts{\alpha}}"', from=2-1, to=2-2]
  \end{tikzcd}
\]
Here $\fpts{F}$ and $\fpts{G}$ denote $\Coco{F}(1)$ and $\Coco{G}(1)$,
equivalently the colimits of $\yo \circ F$ and $\yo \circ G$ in
$\Coco{\bD}$, and $\fpts{\alpha}$ denotes the component of
$\Coco{\alpha}$ at $1$, the induced map between colimits.

\begin{proposition}\label{prop:stronguniform}
  A natural transformation $F \Rightarrow G \colon \bC \to \Set$ is
  uniformly cartesian if and only if the induced functor between
  categories of elements $\El(F) \to \El(G)$ restricts to an
  isomorphism from each connected component of $\El(F)$ to some
  connected component of $\El(G)$. That is, $\alpha$ restricts to an
  isomorphism from each indecomposable summand of the functor $F$ to
  some indecomposable summand of the functor $G$.
\end{proposition}
\begin{proof}
  The indecomposable summands of a $\Set$-valued functor are the
  fibers lying over the elements of its colimit,\footnote{Again, even
    in the case that $\bC$ is large, we may enlarge the universe as
    necessary to ensure that $F$ and $G$ have colimits.} and pulling
  back the colimiting cocone of $G$ along a function
  $\pts{F} \to \pts{G}$ means reindexing these fibers lying over
  $\pts{F}$ to fibers lying over $\pts{G}$.
\end{proof}

This more concrete description of the $\Set$-valued uniformly
cartesian natural transformations also yields a more concrete
description of arbitrary uniformly cartesian natural transformations,
since a natural transformation $F \Rightarrow G \colon \bC \to \bD$ is
uniformly cartesian if and only if the induced $\Set$-valued natural
transformation
$\Ho{\bD}(d, F(\dash)) \Rightarrow \Ho{\bD}(d, G(\dash)) \colon \bC \to \Set$
is uniformly cartesian for all objects $d$ in $\bD$.

\chapter{Continuous universality}\label{app:universal}
\appendixtocspace

The purpose of this section is to prove \cref{thm:universal}:

\thmuniversal*


In fact, as we will see, we have a more general correspondence for not
just continuous maps but also the specialization 2-cells between them
(\cref{def:specializationset}), and, still more generally, the
double-categorical cells from \cref{sec:doubleset}.

We will make use of the following, which also holds in an abstract
2-category.

\begin{lemma}[{\cite[Proposition 7.6]{fairbanks:monads}}]\label{lem:densityformaltfae}
  Let $\pP \colon \bB \to \bS_1$ have density comonad $\cP$, let $\cC$
  be a comonad on $\bS_2$, and let $F \colon \bS_1 \to \bS_2$ be any
  functor preserving the left Kan extension $\cP$ (such as any left
  adjoint).
  
  \begin{enumerate}[label=(\roman*)]
  \item Giving a functor
    $G \colon \CoalgOn{\cP}{\bS_1} \to \CoalgOn{\cC}{\bS_2}$ forming a
    commutative square as shown below left --- i.e.\ a \emph{colax
      comonad functor} carried by $F$ (see \cref{prop:formaltfae}) ---
    is equivalent to giving a functor
    $G' \colon \bB \to \CoalgOn{\cC}{\bS_2}$ forming a commutative
    square as shown below right.
    \[
      \begin{tikzcd}[column sep=20]
        \CoalgOn{\cP}{\bS_1} & \CoalgOn{\cC}{\bS_2} \\
        {\bS_1} & {\bS_2}
        \arrow["G",dashed, from=1-1, to=1-2]
        \arrow["\carP"', from=1-1, to=2-1]
        \arrow["\carC", from=1-2, to=2-2]
        \arrow["F"', from=2-1, to=2-2]
      \end{tikzcd}
      \qquad\qquad\qquad\qquad
      \begin{tikzcd}[column sep=20]
        \bB & \CoalgOn{\cC}{\bS_2} \\
        {\bS_1} & {\bS_2}
        \arrow["G'", dashed, from=1-1, to=1-2]
        \arrow["\pP"', from=1-1, to=2-1]
        \arrow["\carC", from=1-2, to=2-2]
        \arrow["F"', from=2-1, to=2-2]
      \end{tikzcd}
    \]
    More specifically, the canonical map from left squares to right
    squares, given by pre-composing
    $\lP \colon \bB \to \CoalgOn{\cP}{\bS_1}$, is bijective.
  \end{enumerate}
  Now let $F_1, F_2\colon \bS_1 \to \bS_2$ with $F_1$ preserving the
  left Kan extension $\cP$ as above, let
  \[
    \begin{tikzcd}[column sep=20]
      \CoalgOn{\cP}{\bS_1} & \CoalgOn{\cC}{\bS_2} \\
      {\bS_1} & {\bS_2}
      \arrow["G_1", from=1-1, to=1-2]
      \arrow["\carP"', from=1-1, to=2-1]
      \arrow["\carC", from=1-2, to=2-2]
      \arrow["F_1"', from=2-1, to=2-2]
    \end{tikzcd}
    \quad
    \begin{tikzcd}[column sep=20]
      \bB & \CoalgOn{\cC}{\bS_2} \\
      {\bS_1} & {\bS_2}
      \arrow["G_1'", from=1-1, to=1-2]
      \arrow["\pP"', from=1-1, to=2-1]
      \arrow["\carC", from=1-2, to=2-2]
      \arrow["F_1"', from=2-1, to=2-2]
    \end{tikzcd}
    \qqand
    \begin{tikzcd}[column sep=20]
      \CoalgOn{\cP}{\bS_1} & \CoalgOn{\cC}{\bS_2} \\
      {\bS_1} & {\bS_2}
      \arrow["G_2", from=1-1, to=1-2]
      \arrow["\carP"', from=1-1, to=2-1]
      \arrow["\carC", from=1-2, to=2-2]
      \arrow["F_2"', from=2-1, to=2-2]
    \end{tikzcd}
    \quad
    \begin{tikzcd}[column sep=20]
      \bB & \CoalgOn{\cC}{\bS_2} \\
      {\bS_1} & {\bS_2}
      \arrow["G_2'", from=1-1, to=1-2]
      \arrow["\pP"', from=1-1, to=2-1]
      \arrow["\carC", from=1-2, to=2-2]
      \arrow["F_2"', from=2-1, to=2-2]
    \end{tikzcd}
  \]
  be two pairs of corresponding squares as above, and let
  $\gamma \colon F_1 \Rightarrow F_2$ be an arbitrary natural
  transformation.
  \begin{enumerate}[label=(\roman*)]
    \setcounter{enumi}{1}
  \item Giving an arbitrary natural transformation
    $\sigma \colon G_1 \Rightarrow G_2$ --- i.e.\ a \emph{comonad
      specialization} between the comonad 1-cells $F_1$ and $F_2$ (see
    \cref{prop:formaltfae}) --- is equivalent to giving an arbitrary
    natural transformation $\sigma' \colon G_1' \Rightarrow G_2'$.
  \item Giving a natural transformation
    $\sigma \colon G_1 \Rightarrow G_2$ forming a cylinder
    \[
      \begin{tikzcd}
        \CoalgOn{\cP}{\bS_1} & \CoalgOn{\cC}{\bS_2} \\
        \bS_1 & \bS_2
        \arrow["\carP"', from=1-1, to=2-1]
        \arrow["G_1", curve={height=-12pt}, from=1-1, to=1-2]
        \arrow["\carC", from=1-2, to=2-2]
        \arrow[""{name=0, anchor=center, inner sep=0}, "F_2"', curve={height=12pt}, from=2-1, to=2-2]
        \arrow[""{name=1, anchor=center, inner sep=0}, "F_1", curve={height=-12pt}, from=2-1, to=2-2]
        \arrow["\gamma"{description}, draw=none, from=1, to=0]
      \end{tikzcd}
      \quad=\quad
      \begin{tikzcd}
        \CoalgOn{\cP}{\bS_1} & \CoalgOn{\cC}{\bS_2} \\
        \bS_1 & \bS_2
        \arrow["\carP"', from=1-1, to=2-1]
        \arrow[""{name=0, anchor=center, inner sep=0}, "G_1", curve={height=-12pt}, from=1-1, to=1-2]
        \arrow[""{name=1, anchor=center, inner sep=0}, "G_2"', curve={height=12pt}, from=1-1, to=1-2]
        \arrow["\carC", from=1-2, to=2-2]
        \arrow["F_2"', curve={height=12pt}, from=2-1, to=2-2]
        \arrow["\sigma"{description}, draw=none, from=0, to=1]
      \end{tikzcd}
    \]
    --- i.e.\ a \emph{comonad transformation} between the comonad
    1-cells $F_1$ and $F_2$ carried by $\gamma$ (see
    \cref{prop:formaltfae}) --- is equivalent to giving a natural
    transformation $\sigma' \colon G_1' \Rightarrow G_2'$ forming a
    cylinder
    \[
      \begin{tikzcd}
        \;\bB\;\vphantom{\CoalgOn{\cP}{\bS_1}} & \CoalgOn{\cC}{\bS_2} \\
        \bS_1 & \bS_2
        \arrow["\pP"', from=1-1, to=2-1]
        \arrow["G_1'", curve={height=-12pt}, from=1-1, to=1-2]
        \arrow["\carC", from=1-2, to=2-2]
        \arrow[""{name=0, anchor=center, inner sep=0}, "F_2"', curve={height=12pt}, from=2-1, to=2-2]
        \arrow[""{name=1, anchor=center, inner sep=0}, "F_1", curve={height=-12pt}, from=2-1, to=2-2]
        \arrow["\gamma"{description}, draw=none, from=1, to=0]
      \end{tikzcd}
      \quad=\quad\;
      \begin{tikzcd}
        \;\bB\;\vphantom{\CoalgOn{\cP}{\bS_1}} & \CoalgOn{\cC}{\bS_2} \\
        \bS_1 & \bS_2
        \arrow["\pP"', from=1-1, to=2-1]
        \arrow[""{name=0, anchor=center, inner sep=0}, "G_1'", curve={height=-12pt}, from=1-1, to=1-2]
        \arrow[""{name=1, anchor=center, inner sep=0}, "G_2'"', curve={height=12pt}, from=1-1, to=1-2]
        \arrow["\carC", from=1-2, to=2-2]
        \arrow["F_2"', curve={height=12pt}, from=2-1, to=2-2]
        \arrow["\sigma'"{description}, draw=none, from=0, to=1]
      \end{tikzcd}
    \]
  \end{enumerate}
  
  More specifically, the canonical maps between these (given by
  pre-composing $\lP$) are bijective.\qed
\end{lemma}

Our desired result is easier to show when $\pP$ is a basis (in the
sense of \cref{def:basis}), and in this case the result holds in
greater generality than $\Set$. Let us address the full
double-categorical picture while we are at it.

\begin{proposition}\label{prop:lccbasis}
  Let $\bS$ be a locally cartesian closed category with a terminal
  object, and let $\pP \colon \bB \to \bS$ be a basis with sliceable
  density comonad $\cB \coloneqq \cP$. Comonad maps expressed as
  commutative squares, continuous maps expressed as commutative
  squares, specializations, and identifications as expressed as cubes
  as in \cref{prop:cmdtfae}
  \[
    \begin{tikzcd}[column sep=15]
      \BS & \DS \\
      {\SB} & {\SD}
      \arrow["{\com{\bS}{\psi}}", from=1-1, to=1-2]
      \arrow["{\ptlift{\carB}}"', from=1-1, to=2-1]
      \arrow["{\ptlift{\carD}}", from=1-2, to=2-2]
      \arrow["{\leindex{\one{\psi}}}"', from=2-1, to=2-2]
    \end{tikzcd}
    \qquad
    \begin{tikzcd}[column sep=15]
      \AS & \BS \\
      {\SA} & {\SB}
      \arrow["\ptlift{\carA}"', from=1-1, to=2-1]
      \arrow["{\con{\bS}{f}}"', from=1-2, to=1-1]
      \arrow["\ptlift{\carB}", from=1-2, to=2-2]
      \arrow["{\reindex{\one{f}}}", from=2-2, to=2-1]
    \end{tikzcd}
    \qquad
    \begin{tikzcd}[column sep=17, row sep=5]
      & \BS & \\
      \AS && \DS \\
      & \CS & \\
      \arrow["{\con{\bS}{f}}"'{pos=0.35}, from=1-2, to=2-1]
      \arrow["{\Rightarrow_\sigma}"{description}, shift left, draw=none, from=1-2, to=3-2]
      \arrow["{\com{\bS}{\psi}}"{pos=0.35}, from=1-2, to=2-3]
      \arrow["{\con{\bS}{g}}", from=2-3, to=3-2]
      \arrow["{\com{\bS}{\phi}}"', from=2-1, to=3-2]
    \end{tikzcd}
  \]
  \vspace{-2em}
  \[
    \begin{tikzcd}[column sep=17, row sep=15]
      & \BS & \\
      \AS & \SB & \DS \\
      \SA && \SD \\
      & \SC
      \arrow["{\con{\bS}{f}}"'{pos=0.4}, from=1-2, to=2-1]
      \arrow["{\ptlift{\carB}}"{pos=0.4}, from=1-2, to=2-2]
      \arrow["{\com{\bS}{\psi}}"{pos=0.4}, from=1-2, to=2-3]
      \arrow["{\ptlift{\carA}}"'{pos=0.4}, from=2-1, to=3-1]
      \arrow["{\reindex{\one{f}}}"'{pos=0.4}, from=2-2, to=3-1]
      \arrow["{\leindex{\one{\psi}}}"{pos=0.4}, from=2-2, to=3-3]
      \arrow["\Rightarrow"{description}, draw=none, from=2-2, to=4-2]
      \arrow["{\ptlift{\carD}}"{pos=0.4}, from=2-3, to=3-3]
      \arrow["{\leindex{\one{\phi}}}"', from=3-1, to=4-2]
      \arrow["{\reindex{\one{g}}}", from=3-3, to=4-2]
    \end{tikzcd}
    \;\;\; = \;\;\;
    \begin{tikzcd}[column sep=17, row sep=15]
      & \BS & \\
      \AS && \DS \\
      \SA & \CS & \SD \\
      & \SC
      \arrow["{\con{\bS}{f}}"'{pos=0.4}, from=1-2, to=2-1]
      \arrow["{\com{\bS}{\psi}}"{pos=0.4}, from=1-2, to=2-3]
      \arrow["{\Rightarrow_\sigma}"{description}, shift left, draw=none, from=1-2, to=3-2]
      \arrow["{\ptlift{\carA}}"'{pos=0.4}, from=2-1, to=3-1]
      \arrow["{\com{\bS}{\phi}}"'{pos=0.6}, from=2-1, to=3-2]
      \arrow["{\con{\bS}{g}}"{pos=0.6}, from=2-3, to=3-2]
      \arrow["{\ptlift{\carD}}"{pos=0.4}, from=2-3, to=3-3]
      \arrow["{\leindex{\one{\phi}}}"', from=3-1, to=4-2]
      \arrow["{\ptlift{\carC}}"{pos=0.4}, from=3-2, to=4-2]
      \arrow["{\reindex{\one{g}}}", from=3-3, to=4-2]
    \end{tikzcd}
  \]
  are each in bijection with the corresponding diagrams replacing
  $\BS$ with $\bB$ and $\ptlift{\carB}$ with $\ptlift{\pP}$
  \[
    \begin{tikzcd}[column sep=12.5]
      \bB & \DS \\
      {\SB} & {\SD}
      \arrow[dashed, from=1-1, to=1-2]
      \arrow["{\ptlift{\pP}}"', from=1-1, to=2-1]
      \arrow["{\ptlift{\carD}}", from=1-2, to=2-2]
      \arrow["{\leindex{\one{\psi}}}"', from=2-1, to=2-2]
    \end{tikzcd}
    \qquad
    \begin{tikzcd}[column sep=12.5]
      \AS & \bB \\
      {\SA} & {\SB}
      \arrow["\ptlift{\carA}"', from=1-1, to=2-1]
      \arrow[dashed, from=1-2, to=1-1]
      \arrow["\ptlift{\pP}", from=1-2, to=2-2]
      \arrow["{\reindex{\one{f}}}", from=2-2, to=2-1]
    \end{tikzcd}
    \qquad
    \begin{tikzcd}[column sep=17, row sep=5]
      & \bB & \\
      \AS && \DS \\
      & \CS & \\
      \arrow[dashed, from=1-2, to=2-1]
      \arrow["{\Rightarrow_{\sigma'}}"{description}, shift left, draw=none, from=1-2, to=3-2]
      \arrow[dashed, from=1-2, to=2-3]
      \arrow["{\con{\bS}{g}}", from=2-3, to=3-2]
      \arrow["{\com{\bS}{\phi}}"', from=2-1, to=3-2]
    \end{tikzcd}
  \]
  \vspace{-2em}
  \[
    \begin{tikzcd}[column sep=17, row sep=15]
      & \bB & \\
      \AS & \SB & \DS \\
      \SA && \SD \\
      & \SC
      \arrow[dashed, from=1-2, to=2-1]
      \arrow["{\ptlift{\pP}}"{pos=0.4}, from=1-2, to=2-2]
      \arrow[dashed, from=1-2, to=2-3]
      \arrow["{\ptlift{\carA}}"'{pos=0.4}, from=2-1, to=3-1]
      \arrow["{\reindex{\one{f}}}"'{pos=0.4}, from=2-2, to=3-1]
      \arrow["{\leindex{\one{\psi}}}"{pos=0.4}, from=2-2, to=3-3]
      \arrow["\Rightarrow"{description}, draw=none, from=2-2, to=4-2]
      \arrow["{\ptlift{\carD}}"{pos=0.4}, from=2-3, to=3-3]
      \arrow["{\leindex{\one{\phi}}}"', from=3-1, to=4-2]
      \arrow["{\reindex{\one{g}}}", from=3-3, to=4-2]
    \end{tikzcd}
    \;\;\; = \;\;\;
    \begin{tikzcd}[column sep=17, row sep=15]
      & \bB & \\
      \AS && \DS \\
      \SA & \CS & \SD \\
      & \SC
      \arrow[dashed, from=1-2, to=2-1]
      \arrow[dashed, from=1-2, to=2-3]
      \arrow["{\Rightarrow_{\sigma'}}"{description}, shift left, draw=none, from=1-2, to=3-2]
      \arrow["{\ptlift{\carA}}"'{pos=0.4}, from=2-1, to=3-1]
      \arrow["{\com{\bS}{\phi}}"'{pos=0.6}, from=2-1, to=3-2]
      \arrow["{\con{\bS}{g}}"{pos=0.6}, from=2-3, to=3-2]
      \arrow["{\ptlift{\carD}}"{pos=0.4}, from=2-3, to=3-3]
      \arrow["{\leindex{\one{\phi}}}"', from=3-1, to=4-2]
      \arrow["{\ptlift{\carC}}"{pos=0.4}, from=3-2, to=4-2]
      \arrow["{\reindex{\one{g}}}", from=3-3, to=4-2]
    \end{tikzcd}
  \]
  which in turn correspond --- up to choices of pullbacks --- to
  commutative triangles, strongly uniformly cartesian (\cref{def:strongcartesian})
  lax triangles, natural transformations (same as above), and pyramids
  \[
    \begin{tikzcd}[column sep=10pt]
      \bB\ar[rr, dashed]\ar[dr, "\pP"']&&\DS\ar[dl, "\carD"]\\
      &\bS
    \end{tikzcd}
    \qquad
    \begin{tikzcd}[column sep=10pt]
      \AS && \bB \\
      & \bS
      \arrow[dashed, from=1-3, to=1-1]
      \arrow["\carA"', ""{name=0, anchor=center, inner sep=0}, from=1-1, to=2-2]
      \arrow["\pP", ""{name=1, anchor=center, inner sep=0}, from=1-3, to=2-2]
      \arrow["\Rightarrow_{\connat{f}'}", shift right=10pt, draw=none, from=1, to=0]
    \end{tikzcd}
    \qquad
    \begin{tikzcd}[column sep=17, row sep=5]
      & \bB & \\
      \AS && \DS \\
      & \CS & \\
      \arrow[dashed, from=1-2, to=2-1]
      \arrow["{\Rightarrow_{\sigma'}}"{description}, shift left, draw=none, from=1-2, to=3-2]
      \arrow[dashed, from=1-2, to=2-3]
      \arrow["{\con{\bS}{g}}", from=2-3, to=3-2]
      \arrow["{\com{\bS}{\phi}}"', from=2-1, to=3-2]
    \end{tikzcd}
  \]
  \vspace{-2em}
  \[
    \begin{tikzcd}[row sep=12.5pt,column sep=17pt]
      & \bB & \\
      \AS && \DS \\
      \\
      \\
      & \bS
      \arrow[""{name=0, anchor=center, inner sep=0}, dashed, from=1-2, to=2-1]
      \arrow[dashed, from=1-2, to=2-3]
      \arrow["\pP", from=1-2, to=5-2]
      \arrow["\carA"', curve={height=12pt}, from=2-1, to=5-2]
      \arrow["\carD", curve={height=-12pt}, from=2-3, to=5-2]
      \arrow["{\Rightarrow_{\connat{f}'}}"{description}, shift left=2, draw=none, from=5-2, to=0]
    \end{tikzcd}
    \;\;\ = \;\;\;
    \begin{tikzcd}[row sep=10pt,column sep=17pt]
      & \bB & \\
      \AS && \DS \\
      & \CS \\
      \\
      & \bS
      \arrow[dashed, from=1-2, to=2-1]
      \arrow[dashed, from=1-2, to=2-3]
      \arrow["{\Rightarrow_{\sigma'}}"{description}, shift left, draw=none, from=2-1, to=2-3]
      \arrow["{\com{\bS}{\phi}}", from=2-1, to=3-2]
      \arrow["\carA"', curve={height=12pt}, from=2-1, to=5-2]
      \arrow["{\con{\bS}{g}}"', from=2-3, to=3-2]
      \arrow["\carD", curve={height=-12pt}, from=2-3, to=5-2]
      \arrow["\carC"'{pos=0.4}, from=3-2, to=5-2]
      \arrow["{\Rightarrow_{\connat{g}}}"{description}{pos=.45}, shift left, draw=none, from=5-2, to=2-3]
    \end{tikzcd}
  \]
\end{proposition}

(Moreover this proposition will still hold replacing each category
$\SA$ with any equivalent category, and replacing each category $\AS$
with the corresponding category of coalgebras. In particular it
applies using the categories $\SetpA$ in place of $\SetA$, as we have
often preferred.)

\begin{proof}[Proof sketch]
  
  %
  By \cref{prop:liftbasis}, $\clift{\cP} \cong
  \gen{\ptlift{\pP}}$. The first round of translation is by
  \cref{lem:densityformaltfae}. The translation to triangles is
  similar to \cref{prop:contfae} and \cref{prop:cmdtfae}.
\end{proof}

The generalization of this result to arbitrary subbases is actually
rather particular to $\Set$ (see \cref{cex:nonuniversal}). The plan
for our proof is as follows.

Our goal will be to show that for any $\pP \colon \bB \to \Set$ with
density comonad $\cP$ and for any function
$\one{f} \colon \xX \to \pts{\pP}$, the density comonad
$\gen{\one{f}^* \pP}$, where $\one{f}^* \pP$ is the reindexing of
$\pP$ along $\one{f}$ by pullback,\footnote{This $\one{f}^* \pP$ is
  the reindexing of $\pP$ along $\one{f}$ with respect to the
  fibration $\pts{\dash} \colon \fun{\bB}{\Set} \to \Set$, given by
  pullback --- which is defined (up to isomorphism) assuming $\bB$
  either is small or has a terminal object, so that the colimits in
  $\Set$ always exist; otherwise the domain should be replaced with
  the full subcategory of functors with a colimit. The cartesian
  arrows of the fibration are the uniformly cartesian natural
  transformations.} is the same as
$\gen{\one{f}^* \carP}$.\footnote{This $\gen{\one{f}^* \carP}$ is also
  the reindexing of $\cP$ along $\one{f}$ with respect to the
  fibration $\one{\dash} \colon \ConOSet \to \Set$ from
  \cref{prop:fibration}} With this established, we will derive the
desired universal property of $\pP$ with respect to continuous maps
from the known universal property with respect to comonad maps. (By
\cref{lem:coniscom} a continuous map $\cC \to \cD$ equivalently
amounts to a function $\one{f} \colon \pts{\cC} \to \pts{\cD}$ and an
identity-on-points comonad map $\gen{\one{f}^* \carD} \to \cC$.)

We will first show that $\Set^\xX$ is the free completion under finite
intersections of its full subcategory of globally inhabited
objects. We will then deduce that given any terminal-object-preserving
comonad on $\Set^\xX$, there is a unique
finite-intersection-preserving comonad that is defined the same way on
the globally inhabited objects: namely the sliced comonad on
$\Set^\xX$ corresponding to the induced comonad on $\Set$. Finally we
will show as desired that $\gen{\one{f}^* \pP}$ and
$\gen{\one{f}^* \carP}$ are the same by verifying that the two induced
comonads on $\Set^\xX$ agree on the globally inhabited objects.

\begin{lemma}\label{lem:surjcorefl}
  The representable functor $(\Set^\xX)\op \to \Set$ of any globally
  inhabited object sends pullback squares consisting of monomorphisms
  to pushout squares.
\end{lemma}
\begin{proof}
  We will use two facts about the category $\bS \coloneqq
  \Set^\xX$. First, unions of subobjects in $\bS$ are calculated as
  pushouts: given a cospan of monomorphisms
  $\xS_1 \hookrightarrow \xT \hookleftarrow \xS_2$ in $\bS$, the induced map
  $\xS_1 +_{(\xS_1 \times_{\xT} \xS_2)} \xS_2 \to \xT$ is a monomorphism.
  Second, the injective objects $\xA$ (meaning $\Ho{\bS}(\dash, \xA)$ sends
  monomorphisms to epimorphisms) are the globally inhabited objects.
  
  Given a pullback square of monomorphisms in $\bS$ as shown left, we
  must show that the corresponding square in $\Set$ shown right is a
  pushout.
  \[\begin{tikzcd}
      \xU & {\xS_1} \\
      {\xS_2} & \xT
      \arrow[hook, from=1-1, to=1-2]
      \arrow[hook, from=1-1, to=2-1]
      \arrow["\lrcorner"{anchor=center, pos=0.125}, draw=none, from=1-1, to=2-2]
      \arrow[hook, from=1-2, to=2-2]
      \arrow[hook, from=2-1, to=2-2]
    \end{tikzcd}
    \qquad\qquad
    \begin{tikzcd}[column sep=12.5]
      {\Ho{\bS}(\xU, \xA)} & {\Ho{\bS}(\xS_1, \xA)} \\
      {\Ho{\bS}(\xS_2, \xA)} & {\Ho{\bS}(\xT, \xA)}
      \arrow[from=1-2, to=1-1]
      \arrow[from=2-1, to=1-1]
      \arrow[from=2-2, to=1-2]
      \arrow[""{name=0, anchor=center, inner sep=0}, from=2-2, to=2-1]
      \arrow["\lrcorner"{anchor=center, pos=0.125}, draw=none, from=1-1, to=0]
    \end{tikzcd}\] That is, we must show that maps $\xU \to \xA$ are
  equivalently specified by maps $\xS_1 \to \xA$ or $\xS_2 \to \xA$ modulo
  mutual extension to $\xT \to \xA$. By injectivity of $\xA$, every map
  $\xU \to \xA$ extends to some map $\xS_1 \to \xA$.
  On the other hand, suppose $\xU \to \xA$ extends to both $\xS_1 \to \xA$ and
  $\xS_2 \to \xA$. Then it extends to the pushout $\xS_1 +_\xU \xS_2 \to \xA$ and
  hence, since $\xS_1 +_\xU \xS_2$ maps into $\xT$ monomorphically, further to
  $\xT \to \xA$.
\end{proof}

\begin{lemma}\label{lem:setcodense}
  The full subcategory of $\Set^\xX$ consisting of globally inhabited objects is codense.
\end{lemma}
\begin{proof}
  For the special case where $\xX \coloneqq 1$, note that the empty
  set is indeed the canonical limit of $\Set_{>0} \to \Set$, where
  $\Set_{>0}$ is the full subcategory of $\Set$ consisting of nonempty
  sets, so the full subcategory inclusion
  $\Set_{>0} \hookrightarrow \Set$ is codense.

  The full subcategory inclusion of globally inhabited objects in
  $\Set^\xX$ is the $\xX$-fold product power of
  $\Set_{>0} \hookrightarrow \Set$. We must show that each object
  $(\xA_i)_{i \in \xX}$ in $\Set^\xX$ is given by the canonical limit,
  namely the limit of
  $\prod_{j \in \xX}(\xA_j/\Set_{>0}) \to \Set^{\xX}$.

  For each $i \in \xX$, the product projection
  $\pi_i \colon \prod_{j \in \xX}(\xA_j/\Set_{>0}) \to
  (\xA_i/\Set_{>0})$ is cofinal (a.k.a.\ initial, the dual property of
  final as in \cref{def:final}), i.e.\ the category $(\pi_i / f)$ is
  connected for all objects $f$ in $(\xA_i/\Set_{>0})$.
  Indeed, the identity functor $\id_{(\xA_i/\Set_{>0})}$ is cofinal,
  $\bang \colon \prod_{j \in \xX \setminus \setof{i}}(\xA_j/\Set_{>0})
  \to 1$ is cofinal since
  $\prod_{j \in \xX \setminus \setof{i}}(\xA_j/\Set_{>0})$ is
  connected (it has a terminal object), and pairwise products of
  cofinal functors are cofinal.  Thus the canonical limit in
  $\Set^\xX$ is given componentwise by the canonical limits $\xA_i$ in
  $\Set$, as desired.
\end{proof}

\begin{lemma}\label{prop:freecoreflslice}
  $\Set^\xX$ is the free completion under coreflexive equalizers of its
  full subcategory consisting of globally inhabited objects.
\end{lemma}
\begin{proof}
  %
  Let $\iota \colon \sS \hookrightarrow \Set$ denote the full
  subcategory of globally inhabited objects.  The free coreflexive equalizer
  completion of $\sS$ is equivalently the opposite of the closure of
  the Yoneda embedding under reflexive coequalizers in
  $\fun{\sS}{\Set}$. It suffices to show that this full subcategory
  consists of those functors of the form
  $\Ho{\Set^\xX}(\one{f}, \iota(\dash))$, which form a full
  subcategory of $\fun{\sS}{\Set}$ equivalent to $\Set^\xX$ since
  $\iota$ is codense by \cref{lem:setcodense}. The result now follows
  by \cref{lem:surjcorefl}: reflexive coequalizers of functors
  $\Ho{\Set^\xX}(\one{f}, \iota(\dash))$ are calculated as equalizers
  in $\Set$.
\end{proof}

\begin{corollary}\label{cor:freeintslice}
  $\Set^\xX$ is the free completion under finite intersections of its
  full subcategory consisting of globally inhabited objects.
\end{corollary}
\begin{proof}
  If a category $\bS$ has regular finite intersections, then it has
  coreflexive equalizers by \cref{lem:coreflpull}.  A functor
  $\Set^\xX \to \bS$ preserves finite intersections if and only if it
  preserves coreflexive equalizers by \cref{lem:functorpreserve}. By
  \cref{prop:freecoreflslice} such a functor is equivalently specified
  by an arbitrary functor into $\bS$ from the full subcategory of
  $\Set^\xX$ consisting of globally inhabited objects.
\end{proof}



\begin{lemma}\label{lem:crepreflection}
  If $\cC$ is a comonad on $\Set^\xX$ that preserves the full
  subcategory $\sS \hookrightarrow \Set^\xX$ consisting of globally inhabited objects,
  then there is a unique crude (equivalently,
  finite-intersection-preserving) comonad $\cC'$ on $\Set$ that
  restricts to the same comonad as $\cC$ on $\sS$.
  This exhibits crude, $\sS$-preserving comonads as reflective in
  $\sS$-preserving comonads on $\Set^\xX$.
\end{lemma}
\begin{proof}
  Free coreflexive equalizer completion is 2-functorial and locally
  fully faithful. Thus given a comonad on $\Set^\xX$ that restricts to
  a comonad on $\sS$, we obtain a unique
  coreflexive-equalizer-preserving extension to the free coreflexive
  equalizer completion of $\sS$, which is $\Set^\xX$ by
  \cref{prop:freecoreflslice}. Likewise any comonad map between such
  comonads restricts to comonads on $\sS$ and we obtain a unique
  extension to the corresponding crude comonads.
  %
\end{proof}

In particular, any endofunctor on $\Set^\xX$ preserving the terminal
object $\xX$ preserves globally inhabited objects, since
an object is globally inhabited if and only if it admits a section from $\xX$.

We now have two descriptions (\cref{lem:objectivereflective} and
\cref{lem:crepreflection}) of the free crude comonad on a
terminal-object-preserving comonad on $\Set^\xX$. We hence obtain:

\begin{corollary}\label{lem:surjsufficient}
  Terminal-object preserving comonads on $\Set^\xX$ have the same
  underlying comonad on $\Set$ (the reflection from
  \cref{lem:objectivereflective}) if and only if they agree on the
  full subcategory of globally inhabited objects.\qed
\end{corollary}

\begin{lemma}\label{cor:reindexreflect}
  The reindexing functors $\reindex{\one{f}} \circ \dash \circ \Pi_\one{f}$ between
  the categories $\Com(\Set^\xX)_\xX$ commute with the reflections from
  \cref{lem:objectivereflective}.
\end{lemma}
\begin{proof}
  Observe that for any function $\one{f}\colon \xX \to \xY$ we have
  that $\Pi_\one{f}$ preserves globally inhabited objects. Thus the
  result follows from \cref{lem:surjsufficient}, since if comonads on
  $\Set^\xY$ agree on the globally inhabited objects then so do their reindexings.
\end{proof}

\begin{corollary}\label{lem:reindexbasis}
  Let $\pP \colon \bB \to \Set$ with density comonad $\cP$, and let
  $\one{f} \colon \xX \to \pts{\pP}$. Then the following comonads on $\Set$
  are isomorphic:
  \begin{enumerate}[label=(\roman*)]
  \item The reindexed comonad
    $\one{f}^*\cP \coloneqq \slicel \circ \reindex{\one{f}} \circ
    \clift{\cP} \circ \Pi_\one{f} \circ \slicer$.
  \item The density comonad $\gen{\one{f}^*\pP}$ of the reindexed
    subbasis
    $\one{f}^*\pP \coloneqq \slicel \circ \reindex{\one{f}} \circ
    \ptlift{\pP}$.\qed
  \end{enumerate}
\end{corollary}

\begin{theorem}\label{prop:subbasiscont}
  The result of \cref{prop:lccbasis} holds replacing the basis $\pP$
  with an arbitrary functor into $\Set$ admitting a density comonad
  $\cB \coloneqq \cP$.
\end{theorem}
\begin{proof}
  First we must find a correspondence between squares:
  \[
    \begin{tikzcd}[column sep=12.5]
      \ASet & \BSet \\
      {\SetpA} & {\SetpB}
      \arrow["\car{\clift{\cA}}"', from=1-1, to=2-1]
      \arrow["{\con{\Set}{f}}"', dashed, from=1-2, to=1-1]
      \arrow["\car{\clift{\cB}}", from=1-2, to=2-2]
      \arrow["{\reindex{\one{f}}}", from=2-2, to=2-1]
    \end{tikzcd}
    \qqand
    \begin{tikzcd}[column sep=12.5]
      \ASet & \bB \\
      {\SetpA} & {\SetpB}
      \arrow["\car{\clift{\cA}}"', from=1-1, to=2-1]
      \arrow[dashed, from=1-2, to=1-1]
      \arrow["{\ptlift{\pP}}", from=1-2, to=2-2]
      \arrow["{\reindex{\one{f}}}", from=2-2, to=2-1]
    \end{tikzcd}
  \]
  By \cref{prop:density} and \cref{cor:densityofleft}, the left square
  is equivalently a comonad map
  $\reindex{\one{f}} \circ \clift{\cP} \circ \Pi_{\one{f}} \to
  \clift{\cA}$ on $\SetpA$. By \cref{lem:objectivereflective} and
  \cref{cor:reindexreflect}, this is equivalently a comonad map
  $\gen{\reindex{\one{f}} \circ \ptlift{\pP}} \to \clift{\cA}$,
  corresponding to the right square by \cref{prop:density}.
  
  Next suppose given comonad maps and continuous maps as in \cref{prop:lccbasis}.
  \[
    \begin{tikzcd}[column sep=12.5]
      \ASet & \CSet \\
      {\SetpA} & {\SetpC}
      \arrow["{\com{\Set}{\phi}}", dashed, from=1-1, to=1-2]
      \arrow["\car{\clift{\cA}}"', from=1-1, to=2-1]
      \arrow["\car{\clift{\cC}}", from=1-2, to=2-2]
      \arrow["{\leindex{\one{\phi}}}"', from=2-1, to=2-2]
    \end{tikzcd}
    \quad
    \begin{tikzcd}[column sep=12.5]
      \BSet & \DSet \\
      {\SetpB} & {\SetpD}
      \arrow["{\com{\Set}{\psi}}", dashed, from=1-1, to=1-2]
      \arrow["\car{\clift{\cB}}"', from=1-1, to=2-1]
      \arrow["\car{\clift{\cD}}", from=1-2, to=2-2]
      \arrow["{\leindex{\one{\psi}}}"', from=2-1, to=2-2]
    \end{tikzcd}
    \quad
    \begin{tikzcd}[column sep=12.5]
      \ASet & \BSet \\
      {\SetpA} & {\SetpB}
      \arrow["\car{\clift{\cA}}"', from=1-1, to=2-1]
      \arrow["{\con{\Set}{f}}"', dashed, from=1-2, to=1-1]
      \arrow["\car{\clift{\cB}}", from=1-2, to=2-2]
      \arrow["{\reindex{\one{f}}}", from=2-2, to=2-1]
    \end{tikzcd}
    \quad
    \begin{tikzcd}[column sep=12.5]
      \CSet & \DSet \\
      {\SetpC} & {\SetpD}
      \arrow["\car{\clift{\cC}}"', from=1-1, to=2-1]
      \arrow["{\con{\Set}{g}}"', dashed, from=1-2, to=1-1]
      \arrow["\car{\clift{\cD}}", from=1-2, to=2-2]
      \arrow["{\reindex{\one{g}}}", from=2-2, to=2-1]
    \end{tikzcd}
  \]

  \noindent
  We would like to show a correspondence between natural transformations:
  \[
    \begin{tikzcd}[column sep=17, row sep=5]
      & \BSet & \\
      \ASet && \DSet \\
      & \CSet &
      \arrow["{\con{\Set}{f}}"'{pos=0.35}, from=1-2, to=2-1]
      \arrow["{\Rightarrow_\sigma}"{description}, shift left, draw=none, from=1-2, to=3-2]
      \arrow["{\com{\Set}{\psi}}"{pos=0.35}, from=1-2, to=2-3]
      \arrow["{\con{\Set}{g}}", from=2-3, to=3-2]
      \arrow["{\com{\Set}{\phi}}"', from=2-1, to=3-2]
    \end{tikzcd}
    \qqand
    \begin{tikzcd}[column sep=17, row sep=5]
      & \bB & \\
      \ASet && \DSet \\
      & \CSet &
      \arrow[dashed, from=1-2, to=2-1]
      \arrow["{\Rightarrow_{\sigma'}}"{description}, shift left, draw=none, from=1-2, to=3-2]
      \arrow[dashed, from=1-2, to=2-3]
      \arrow["{\con{\Set}{g}}", from=2-3, to=3-2]
      \arrow["{\com{\Set}{\phi}}"', from=2-1, to=3-2]
    \end{tikzcd}
  \]
  Note that both functors
  $\slicel \circ \leindex{\one{\phi}} \circ \reindex{\one{f}} \circ
  \car{\clift{\cB}} = \one{f}^* \carB$ and
  $\slicel \circ \leindex{\one{\phi}} \circ \reindex{\one{f}} \circ
  \ptlift{\pP} = \one{f}^* \pP$ induce the same density comonad on
  $\Set$ by \cref{lem:reindexbasis}.
  Likewise, both functors
  $\slicel \circ \reindex{\one{g}} \circ \leindex{\one{\psi}} \circ
  \car{\clift{\cB}}$ and
  $\slicel \circ \reindex{\one{g}} \circ \leindex{\one{\psi}} \circ
  \ptlift{\pP}$ induce the same density comonad on $\Set$ for the
  same reason since we may rewrite
  $\reindex{\one{g}}\circ \leindex{\one{\psi}}$ as
  $\leindex{\one{g}^*\one{\psi}}\circ
  \reindex{\one{\psi}^*\one{g}}$ by Beck-Chevalley.

  Specializations $\sigma$ as shown above are the same as
  specializations between the corresponding comonad maps between these
  two comonads on $\Set$. We then obtain the desired correspondence by
  \cref{lem:densityformaltfae}.

  It remains to verify that specializations $\sigma$
  satisfying the cube axiom
  \[
    \begin{tikzcd}[column sep=17, row sep=15]
      & \BSet & \\
      \ASet & \SetpB & \DSet \\
      \SetpA && \SetpD \\
      & \SetpC
      \arrow["{\con{\Set}{f}}"'{pos=0.4}, from=1-2, to=2-1]
      \arrow["\car{\clift{\cB}}"{pos=0.4}, from=1-2, to=2-2]
      \arrow["{\com{\Set}{\psi}}"{pos=0.4}, from=1-2, to=2-3]
      \arrow["\car{\clift{\cA}}"'{pos=0.4}, from=2-1, to=3-1]
      \arrow["{\reindex{\one{f}}}"'{pos=0.4}, from=2-2, to=3-1]
      \arrow["{\leindex{\one{\psi}}}"{pos=0.4}, from=2-2, to=3-3]
      \arrow["\Rightarrow"{description}, draw=none, from=2-2, to=4-2]
      \arrow["\car{\clift{\cD}}"{pos=0.4}, from=2-3, to=3-3]
      \arrow["{\leindex{\one{\phi}}}"', from=3-1, to=4-2]
      \arrow["{\reindex{\one{g}}}", from=3-3, to=4-2]
    \end{tikzcd}
    \;\;\; = \;\;\;
    \begin{tikzcd}[column sep=17, row sep=15]
      & \BSet & \\
      \ASet && \DSet \\
      \SetpA & \CSet & \SetpD \\
      & \SetpC
      \arrow["{\con{\Set}{f}}"'{pos=0.4}, from=1-2, to=2-1]
      \arrow["{\com{\Set}{\psi}}"{pos=0.4}, from=1-2, to=2-3]
      \arrow["{\Rightarrow_\sigma}"{description}, shift left, draw=none, from=1-2, to=3-2]
      \arrow["\car{\clift{\cA}}"'{pos=0.4}, from=2-1, to=3-1]
      \arrow["{\com{\Set}{\phi}}"'{pos=0.85}, from=2-1, to=3-2]
      \arrow["{\con{\Set}{g}}"{pos=0.85}, from=2-3, to=3-2]
      \arrow["\car{\clift{\cD}}"{pos=0.4}, from=2-3, to=3-3]
      \arrow["{\leindex{\one{\phi}}}"', from=3-1, to=4-2]
      \arrow["\car{\clift{\cC}}"{pos=0.4}, from=3-2, to=4-2]
      \arrow["{\reindex{\one{g}}}", from=3-3, to=4-2]
    \end{tikzcd}
  \]
  (where the unlabelled natural transformation $\Rightarrow$ denotes
  the mate to identity), i.e.\ identifications, correspond to natural
  transformations $\sigma'$ satisfying the cube axiom
  \[
    \begin{tikzcd}[column sep=17, row sep=15]
      & \bB & \\
      \ASet & \SetpB & \DSet \\
      \SetpA && \SetpD \\
      & \SetpC
      \arrow[dashed, from=1-2, to=2-1]
      \arrow["{\ptlift{\pP}}"{pos=0.4}, from=1-2, to=2-2]
      \arrow[dashed, from=1-2, to=2-3]
      \arrow["\car{\clift{\cA}}"'{pos=0.4}, from=2-1, to=3-1]
      \arrow["{\reindex{\one{f}}}"'{pos=0.4}, from=2-2, to=3-1]
      \arrow["{\leindex{\one{\psi}}}"{pos=0.4}, from=2-2, to=3-3]
      \arrow["\Rightarrow"{description}, draw=none, from=2-2, to=4-2]
      \arrow["\car{\clift{\cD}}"{pos=0.4}, from=2-3, to=3-3]
      \arrow["{\leindex{\one{\phi}}}"', from=3-1, to=4-2]
      \arrow["{\reindex{\one{g}}}", from=3-3, to=4-2]
    \end{tikzcd}
    \;\;\; = \;\;\;
    \begin{tikzcd}[column sep=17, row sep=15]
      & \bB & \\
      \ASet && \DSet \\
      \SetpA & \CSet & \SetpD \\
      & \SetpC
      \arrow[dashed, from=1-2, to=2-1]
      \arrow[dashed, from=1-2, to=2-3]
      \arrow["{\Rightarrow_{\sigma'}}"{description}, shift left, draw=none, from=1-2, to=3-2]
      \arrow["\car{\clift{\cA}}"'{pos=0.4}, from=2-1, to=3-1]
      \arrow["{\com{\Set}{\phi}}"'{pos=0.85}, from=2-1, to=3-2]
      \arrow["{\con{\Set}{g}}"{pos=0.85}, from=2-3, to=3-2]
      \arrow["\car{\clift{\cD}}"{pos=0.4}, from=2-3, to=3-3]
      \arrow["{\leindex{\one{\phi}}}"', from=3-1, to=4-2]
      \arrow["\car{\clift{\cC}}"{pos=0.4}, from=3-2, to=4-2]
      \arrow["{\reindex{\one{g}}}", from=3-3, to=4-2]
    \end{tikzcd}
  \]
  The domain and codomain of the natural
  transformations on either side of the first cube equation are
  $\car{\clift{\cC}} \circ \com{\Set}{\phi} \circ \con{\Set}{f}$ and
  $\car{\clift{\cC}} \circ \con{\Set}{g} \circ \com{\Set}{\psi}$. These
  themselves determine colax comonad functors, carried by
  $\slicel \circ \leindex{\one{\phi}} \circ \reindex{\one{f}}$ and
  $\slicel \circ \reindex{\one{g}} \circ \leindex{\one{\psi}}$, from the
  comonad $\clift{\bB}$ on $\SetpB$ to the comonad
  $\pts{\cC} \times \dash$ on $\Set$.
  \[
    \begin{tikzcd}[column sep=15, row sep=15]
      \SetpC & \CSet & \ASet & \BSet \\
      \Set  & \SetpC & \SetpA & \SetpB
      \arrow["{\slicel}"', from=1-1, to=2-1]
      \arrow["\car{\clift{\cB}}", from=1-4, to=2-4]
      \arrow["\car{\clift{\cC}}"', from=1-2, to=1-1]
      \arrow["{\com{\Set}{\phi}}"', from=1-3, to=1-2]
      \arrow["{\con{\Set}{f}}"', from=1-4, to=1-3]
      \arrow["{\slicel}", from=2-2, to=2-1]
      \arrow["{\leindex{\one{\phi}}}", from=2-3, to=2-2]
      \arrow["{\reindex{\one{f}}}", from=2-4, to=2-3]
    \end{tikzcd}
    \qquad
    \begin{tikzcd}[column sep=15, row sep=15]
      \SetpC & \CSet & \DSet & \BSet \\
      \Set  & \SetpC & \SetpD & \SetpB
      \arrow["{\slicel}"', from=1-1, to=2-1]
      \arrow["{\car{\clift{\cB}}}", from=1-4, to=2-4]
      \arrow["\car{\clift{\cC}}"', from=1-2, to=1-1]
      \arrow["{\con{\Set}{g}}"', from=1-3, to=1-2]
      \arrow["{\com{\Set}{\psi}}"', from=1-4, to=1-3]
      \arrow["{\slicel}", from=2-2, to=2-1]
      \arrow["{\reindex{\one{g}}}", from=2-3, to=2-2]
      \arrow["{\leindex{\one{\psi}}}", from=2-4, to=2-3]
    \end{tikzcd}
  \]
  Thus the cube axiom may itself be viewed as an equation between two
  comonad specializations between colax comonad functors.  By
  \cref{lem:densityformaltfae}, such comonad specializations are
  determined by their restriction to $\bB$, as desired.
\end{proof}

\chapter{Counterexamples}\label{app:counterexamples}

In this section we collect counterexamples, which show that various
results of the paper cannot be strengthened.

\section*{Comonads on $\Set$}

All comonads on $\Set$ preserve monomorphisms by
\cref{prop:intersections}. However, the
forgetful functor $\carC \colon \CSet \to \Set$ may not preserve
monomorphisms. (That is, not all monomorphisms in $\CSet$ need be
regular.) It is possible for a monomorphism in $\CSet$ to be sent to a
non-injection whose pullback along itself is not preserved by $\cC$,
and this limit need not be created by the forgetful functor.

\begin{counterexample}[Non-regular monomorphism]\label{cex:shear}
  Consider the comonad from \cref{ex:shear}, whose category of
  coalgebras is the category $\Shear$.
  %
  %
  %
  %
  An object in $\Shear$ consists of a pair of sets $(E_1, E_2)$, and
  an arrow $(E_1, E_2) \to (E'_1, E'_2)$ consists of a function
  $E_1 + E_2 \to E'_1 + E'_2$ sending $E_1$ entirely within
  $E'_1$. The unique map $q \colon (0,1)\to (1,0)$ in $\Shear$ is
  monic, i.e.\ $(0,1)$ is subterminal, because there is a unique map
  $(E_1,E_2)\to (0,1)$ if and only if $E_1=0$. The defining comonadic
  forgetful functor $\Shear \to \Set$ sends $(E_1, E_2)$ to the set
  $E_1 + 2E_2$, so this monomorphism is sent to the map $2 \to 1$,
  which is not injective. By \cref{lem:regularmono}, this is
  equivalent to saying that this monomorphism is not regular.
\end{counterexample}

As mentioned in the introduction and later confirmed in
\cref{sec:slices}, the category of coalgebras of a pullback-preserving
comonad on $\Set$ is a topos. This is not so for comonads that merely
preserve weak pullbacks.

\begin{counterexample}[Non-cartesian-closed coalgebra category]\label{cex:undgraphs}
  The category $\MultiGph$ from \cref{ex:undgraphs} is not cartesian
  closed. Indeed, let $V,E,L$ be the vertex, edge, and loop (the
  terminal undirected multigraph) respectively, and consider the
  coequalizer diagram $V\rightrightarrows E\to L$. The cartesian
  product with $E$ does not preserve this coequalizer. Indeed,
  $V\times E \cong V+V$ and $E\times E$ has two edges, so the
  coequalizer of $V\times E\rightrightarrows E\times E$ has two edges
  as well, while the product $E\times L$ has but a single edge.
\end{counterexample}

The following example shows that it is not possible to characterize
representable, polynomial, pullback-preserving, or
weak-pullback-preserving comonads on $\Set$ based solely on their
categories of coalgebras.

\begin{counterexample}[Non-recognizable coalgebra category I]
  Let $\bb$ be the free monoid on an idempotent element, i.e.\ the
  two-element monoid of the Booleans. Viewed as a one-object category,
  $\bb$ corresponds to a polynomial --- in fact representable ---
  comonad on $\Set$ with formula $\xX \mapsto \xX^2$ whose coalgebras
  are $\bb$-sets $\fun{\bb}{\Set}$. Let
  $U \colon \fun{\bb}{\Set} \to \Set$ denote the comonadic underlying
  set functor, also denoted $\car{\cattocmd{\bb}}$.

  We now construct another comonadic functor
  $F \colon \fun{\bb}{\Set} \to \Set$ that does not preserve
  pullbacks.  Let $U_{\mathrm{fix}} \colon \fun{\bb}{\Set} \to \Set$
  denote the functor sending a $\bb$-set to its set of fixed points,
  also denoted $\fun{\bb}{\Set}(1, \dash)$. We define $F$ to be the
  pushout of the evident inclusion
  $U_{\mathrm{fix}} \hookrightarrow U$ along itself.

  This $F$ does not preserve weak pullbacks; for example it does not
  weakly preserve the pullback of the $\bb$-set map $\bb \to 1$ along itself.
  However, as a colimit of colimit-preserving functors, $F$ preserves
  colimits and therefore is a left adjoint. (Indeed, the category
  $\fun{\bb}{\Set}$ is equivalent to the category $\fun{\bC}{\Set}$
  where $\bC$ is the idempotent splitting of $\bb$, i.e.\ the free
  category containing a split idempotent, and $U_{\mathrm{fix}}$ then
  simply picks out the value at the idempotent, thus preserving
  colimits.)  Explicitly, the right adjoint of $F$ sends a set $\xX$
  to the $\bb$-set $\xX^3$ with the action of the idempotent given by
  $(x, y, z) \mapsto (x, x, x)$.  Moreover $F$ is evidently
  conservative, and it preserves equalizers because equalizers commute
  with unions in $\Set$. Therefore it is comonadic by the crude
  comonadicity theorem (\cref{prop:crude}).
\end{counterexample}

The following is an example of a comonad on $\Set$ that preserves
finite limits as well as arbitrary intersections, but is not
polynomial.

\begin{counterexample}[Non-polynomial intersection-preserving ionad]\label{ex:dds}
  Let $\bB$ be the category whose objects are sets of the form
  $\setof{0, \ldots, n}$ for $n \in \nn$ and where an arrow
  $\setof{0, \ldots, n} \to \setof{0, \ldots, m}$ consists of a
  function of the form $x \mapsto \min\setof{x + k, m}$ that moreover sends
  $n$ to $m$, i.e.\ $n + k \geq m$. Equivalently, this is the full
  subcategory of $\fun{\nn}{\Set}$ consisting of the quotients of
  $\nn$ obtained by identifying all elements $i \geq n$ for some
  $n$. Let $\pP \colon \bB \hookrightarrow \Set$ denote the forgetful
  functor. This functor $\pP$ is readily checked to be flat; therefore
  $\cP$ preserves finite limits (and moreover $\pP$ is a basis by
  \cref{lem:extcoref}).

  We describe $\cP$ in terms of germs as in
  \cref{rem:densitygerms}. An $\xA$-valued germ $[f]_x \in \cP(\xA)$
  is given by a function $f \colon \setof{0, \ldots, n} \to \xA$ and
  an element $0 \leq x \leq n$, where we identify $[f]_x = [g]_y$
  whenever there exists a map $i$ in $\bB$ such that $\pP(i)(x) = y$
  and $f = g\circ \pP(i)$. Thus we have identified each $[g]_k$ given
  by $g \colon \setof{0, \ldots, m} \to \xA$ with any $[f]_0$ given by
  $f \colon \setof{0, \ldots, n} \to \xA$ of the form
  $f(x) = g(\min\setof{x + k, m})$ where $n + k \geq m$. Hence all germs are
  may be written as $[f]_0$, and moreover the domain of $f$ may be
  expanded arbitrarily by duplicating the value assigned its last
  element. In short, $\cP$ sends a set $\xA$ to the set of eventually
  constant sequences in $\xA$.

  A coalgebra $\hH \colon \eE \to \cP(\eE)$ assigns each element $e$
  an eventually constant sequence $(e_i)_{i \in \nn}$ valued in $\eE$,
  with $e_0 = e$ (as guaranteed by the unit law) and such that each
  $e_k$ is assigned $(e_{i + k})_{i \in \nn}$ (as guaranteed by the
  associativity law).  In short, a coalgebra is a set with an
  endofunction, a.k.a.\ discrete dynamical system or object of
  $\fun{\nn}{\Set}$, with no infinite forward trajectories. Arbitrary
  intersections of subcoalgebras are evidently subcoalgebras, so $\cP$
  preserves arbitrary intersections by \cref{prop:alextop}.
  %
\end{counterexample}

In \cref{prop:regular} we recalled that the category of coalgebras on
a weak-pullback-preserving comonad on $\Set$ is regular. The following
shows that it is not sufficient to assume tautness.

\begin{counterexample}[Non-regular taut comonad]\label{ex:monoids}
  Let $(\bM_i)_{i \in I}$ be a fixed family of monoids. There is a
  comonad on $\Set$ for which a coalgebra consists of a family of
  monoid actions $(\xX_i \colon \bM_i \to \Set)_{i \in I}$ all of whose
  fixed points (the elements $x \in \xX_i$ for which $a \act x = x$
  for all $a \in \bM_i$) are identified with the same set
  $\xX_\mathrm{fix}$.  A map of coalgebras consists of a family of
  natural transformations $(\xX_i \Rightarrow \xX'_i)_{i \in I}$ all
  of which agree on the fixed points $\xX_\mathrm{fix}$.

  The comonadic forgetful functor is the quotient of
  $\sum_{i \in I}\xX_i$ obtained by identifying the fixed points
  $\xX_\mathrm{fix}$ in each factor. In other words, the comonadic
  forgetful functor sends $(\xX_i)_{i \in I}$ to the wide pushout of
  inclusions of fixed points
  $(\xX_\mathrm{fix} \hookrightarrow \xX_i)_{i \in I}$. It has right
  adjoint sending a set $\xA$ to the family $(\xA^{\bM_i})$, with the usual
  actions, together with the canonical identifications of the constant
  sequences, which are the fixed points, with $\xA$ itself. The
  forgetful functor is evidently conservative, and as a wide pushout
  consisting of monomorphisms between representable functors it
  preserves equalizers, because equalizers commute with unions in
  $\Set$. Thus it is comonadic by the crude comonadicity theorem (\cref{prop:crude}).

  Since it preserves equalizers, the endofunctor underlying the
  comonad is taut~\cite[VII.11]{trnkova}. But when $(\bM_i)_{i \in I}$
  includes at least two nontrivial monoids $\bM_1$ and $\bM_2$, the
  category of coalgebras is not regular. Indeed, the maps
  $\bang \colon (\bM_1,\varnothing) \to 1$ and
  $\bang \colon (\varnothing,\bM_2) \to (1,1)$ are regular epimorphisms,
  each obtained by jointly coequalizing all the elements, but their
  pullback is empty.
\end{counterexample}

In \cref{prop:classifier} we recalled that if a comonad on $\Set$ is
taut, then, assuming the comonadic functor preserves monomorphisms,
the category of coalgebras has a subobject classifier. The following
example shows that the converse is not true.

\begin{counterexample}[Non-taut subobject classifier]\label{ex:unlink}
  Consider the poset-shaped diagram
  $\pP \colon \bB \to \Set$ given by
  
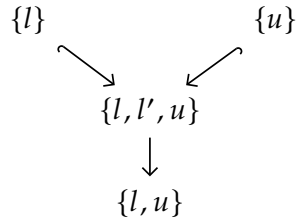
\[
    \begin{tikzcd}[row sep=15,column sep=10]
      {\setof{l}} && {\setof{u}} \\
      & {\setof{l, l', u}} \\
      & {\setof{l,u}}
      \arrow[hook', from=1-1, to=2-2]
      \arrow[hook, from=1-3, to=2-2]
      \arrow[from=2-2, to=3-2]
    \end{tikzcd}
  \]
  where the shown arrows from $\setof{l}$ and $\setof{u}$ are the inclusions
  and the shown arrow into $\setof{l,u}$ sends both $l$ and $l'$ to $l$ and
  sends $u$ to $u$.
  
  We describe the density comonad $\cP$ in terms of germs as in
  \cref{rem:densitygerms}. An $\xA$-valued germ $[f]_x \in \cP(\xA)$
  is either
  \begin{itemize}
  \item the germ of a function $f \colon \setof{l,u} \to \xA$ (an
    ordered pair) about $l$ or $u$, which we denote by $[a_l, a_u]_l$
    or $[a_l, a_u]_u$;
  \item the germ of a function $f \colon \setof{l,l',u} \to \xA$ (an
    ordered triple) about $l$, $l'$, or $u$, which we denote by
    $[a_l, a_{l'}, a_u]_l$, $[a_l, a_{l'}, a_u]_{l'}$, or
    $[a_l, a_{l'}, a_u]_u$;
  \item the germ of a function $f \colon \setof{l} \to \xA$ (an
    element) about the unique element $l$ of $\setof{l}$, which we
    denote by $[a]_l$; or
  \item the germ of a function $f \colon \setof{u} \to \xA$ (an
    element) about the unique element $u$ of $\setof{u}$, which we
    denote by $[a]_u$
  \end{itemize}
  modulo the equivalence relation generated by:
  \[
    [a]_l = [a, b, c]_l
    \quad
    [c]_u = [a, b, c]_u
    \quad
    [a, a, b]_l = [a, b]_l = [a, a, b]_{l'}
    \quad
    [a, a, b]_u = [a, b]_u
  \]
  Each equivalence class then has a canonical representative: either
  $[a]_l$, $[a]_u$, or $[a_l,a_{l'},a_u]_{l'}$ with $a_{l} \neq a_{l'}$.

  A coalgebra $\hH \colon \eE \to \cP(\eE)$ assigns each element $e$
  an $\eE$-valued germ of the form $[e]_l$, $[e]_u$, or
  $[a_l, e, a_u]_{l'}$ (as guaranteed by the unit law) and such that
  in the last case $a_l$ is always assigned $[a_l]_l$ and $a_u$ is
  always assigned $[a_u]_u$ (as guaranteed by the associativity
  law). In other words, a coalgebra is a set $\eE$ partitioned into
  three subsets $\eE_l + \eE_{l'} + \eE_u$, equipped with the
  structure of a span $\eE_l \ot \eE_{l'} \to \eE_u$. A coalgebra map
  $\eE \to \eE'$ consists of functions $\eE_l \to \eE'_l$,
  $\eE_u \to \eE'_u$, and $\eE_{l'} \to \eE'_{l'} + \eE'_l$ such that
  for each $[a_l, e, a_u]_{l'}$ sent to $[a_l', e', a_u']_{l'}$ we have that
  $[a_l]_l$ and $[a_u]_u$ are sent to $[a_l']_l$ and $[a_u']_u$, and
  for each $[a_l, e, a_u]_{l'}$ sent to $[e']_l$ we have that $[a_l]_l$ is
  also sent to $[e']_l$.
  


  Hence the objects of $\PSet$ are spans
  $\eE_l \ot \eE_{l'} \to \eE_u$, but the morphisms are allowed to
  collapse any element of $\eE_{l'}$ with its assigned element of
  $\eE_l$ thereby ``unlinking'' it from its assigned element of
  $\eE_u$. For future convenience, we denote this category by
  $\Unlink \coloneqq \PSet$.
  
  The induced comonad is not taut, since not all coalgebra maps induce
  continuous maps of the underlying topological spaces (see
  \cref{prop:tautpreimages}); for example there is a unique map from
  $\setof{l} \ot \setof{l'} \to \setof{u}$ to the terminal coalgebra
  $\setof{l} \ot \varnothing \to \setof{u}$, and moreover
  $\setof{l} \ot \varnothing \to \varnothing$ is a subcoalgebra
  (a.k.a.\ open subset) of the latter but its preimage
  $\setof{l} \ot \setof{l'} \to \varnothing$ is not a
  subcoalgebra. The category $\Unlink$ nevertheless has a subobject
  classifier, given by the inclusion of the span
  $\setof{\top} \hookleftarrow \varnothing \hookrightarrow
  \setof{\top}$ into the span
  $\setof{\bot, \top} \hookleftarrow \setof{\top} \hookrightarrow
  \setof{\bot, \top}$.

  Indeed, the classifying map $\chi_\eS$ of a subcoalgebra
  $\eS\hookrightarrow\eE$ is uniquely determined as follows. First,
  $\chi_\eS$ must send each element $e$ with germ of the form $[e]_l$
  or $[e]_u$ in $\eS$ to $\top$; likewise, $\chi_\eS$ must send each
  element $e$ with germ of the form $[e]_l$ or $[e]_u$ not in $\eS$ to
  $\bot$. Now consider an element $e$ with germ of the form
  $[a_l, e, a_u]_{l'}$. If $e$ is in $\eS$, then $\chi_\eS$ must
  collapse $e$ with $a_l$, sent to $\top$. Now suppose $e$ is not in
  $\eS$. If $a_l$ and $a_u$ are both in $\eS$, then $\chi_\eS$ must
  send $e$ to $\top$, not collapsed with $a_l$. Otherwise, $\chi_\eS$
  must collapse $e$ with $a_l$, sent to either $\top$ or $\bot$. In
  each case we have that the pullback, calculated as the largest
  subcoalgebra that is included in the pullback taken in $\Set$, is
  correct.

  The comonadic forgetful functor $\Unlink \to \Set$ also preserves
  monomorphisms. Indeed, each element in $\eE_l$, $\eE_{l'}$, or
  $\eE_u$ is picked out by a coalgebra map from the corresponding
  subbasic coalgebra $\pP(\setof{l})$, $\pP(\setof{l, l', u})$, or
  $\pP(\setof{u})$, and so if a coalgebra map sends two elements to
  the same element, then it is not monic in $\Unlink$.
  %
\end{counterexample}

The following example shows that it is not possible to characterize
taut comonads on $\Set$ based solely on their categories of
coalgebras.

\begin{counterexample}[Non-recognizable coalgebra category II]\label{ex:extensive}
  We will see that the density comonads of the two diagrams
  
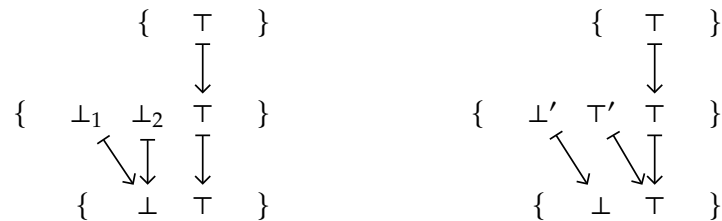
\[
    \begin{tikzcd}[column sep=0, row sep=15]
      & & {\{\;} & \top & {\;\;\}} \\
      {\{\;\;} & \bot_1 & \bot_2 & \top & {\;\;\}} \\
      & {\{\;} & \bot & \top & {\;\;\}}
      \arrow[maps to, from=1-4, to=2-4]
      \arrow[maps to, from=2-4, to=3-4]
      \arrow[maps to, from=2-3, to=3-3]
      \arrow[maps to, from=2-2, to=3-3]
    \end{tikzcd}
    \qquad\qquad\qquad
    \begin{tikzcd}[column sep=0, row sep=15]
      & & {\{\;} & \top & {\;\;\}} \\
      {\{\;\;} & {\bot'} & {\top'} & \top & {\;\;\}} \\
      & {\{\;} & \bot & \top & {\;\;\}}
      \arrow[maps to, from=1-4, to=2-4]
      \arrow[maps to, from=2-4, to=3-4]
      \arrow[maps to, from=2-3, to=3-4]
      \arrow[maps to, from=2-2, to=3-3]
    \end{tikzcd}
  \]
  are distinct, with the former taut and the latter not, but they have
  the same category of coalgebras.
  
  We describe the density comonad $\cP$ of the first diagram $\pP$ in
  terms of germs as in \cref{rem:densitygerms}. An $\xA$-valued germ
  $[f]_x \in \cP(\xA)$ is either
  \begin{itemize}
  \item the germ of a function $f \colon \setof{\bot, \top} \to \xA$ (an
    ordered pair) about $\bot$ or $\top$, which we denote by $[a_\bot, a_\top]_\bot$
    or $[a_\bot, a_\top]_\top$;
  \item the germ of a function $f \colon \setof{\bot_1,\bot_2,\top} \to \xA$ (an
    ordered triple) about $\bot_1$, $\bot_2$, or $\top$, which we denote by
    $[a_{\bot_1}, a_{\bot_2}, a_\top]_{\bot_1}$, $[a_{\bot_1}, a_{\bot_2}, a_\top]_{\bot_2}$, or
    $[a_{\bot_1}, a_{\bot_2}, a_\top]_{\top}$; or
  \item the germ of a function $f \colon \setof{\top} \to \xA$ (an
    element) about the unique element $\top$ of $\setof{\top}$, which we
    denote by $[a]_\top$
  \end{itemize}
  modulo the equivalence relation generated by:
  \[
    [x, a]_\top = [a]_\top = [x, y, a]_\top  
    \quad
    [a, a, b]_{\bot_1} = [a, b]_\bot = [a, a, b]_{\bot_2}  
  \]
  Each equivalence class then has a canonical representative: either
  $[a]_\top$, $[a_\bot, a_\top]_\bot$, or one of
  $[a_{\bot_1},a_{\bot_2},a_\top]_{\bot_1}$ or
  $[a_{\bot_1},a_{\bot_2},a_\top]_{\bot_2}$ with
  $a_{\bot_1} \neq a_{\bot_2}$.

  A coalgebra $\hH \colon \eE \to \cP(\eE)$ assigns each element $e$
  an $\eE$-valued germ of the form $[e]_\top$, $[e, a_\top]_\bot$,
  $[e,a_{\bot_2},a_\top]_{\bot_1}$ or $[a_{\bot_1},e,a_\top]_{\bot_2}$
  (as guaranteed by the unit law) and such that in the latter three
  cases $h(a_\top) = [a_\top]_\top$, in the third case
  $h(a_{\bot_2}) = [e,a_{\bot_2},a_\top]_{\bot_2}$, and in the fourth
  case $h(a_{\bot_1}) = [a_{\bot_1},e,a_\top]_{\bot_1}$ (as guaranteed
  by the associativity law). In other words, a coalgebra is a set
  $\eE$ partitioned into four subsets
  $\eE_{\bot_1} + \eE_{\bot_2} + \eE_{\bot} + \eE_{\top}$, equipped
  with functions
  $\eE_{\bot_1} \cong \eE_{\bot_2} \to \eE_\top \ot \eE_{\bot}$. A
  coalgebra map $\eE \to \eE'$ consists of functions
  $\eE_{\bot_1} \to \eE'_{\bot_1} + \eE'_\bot$,
  $\eE_{\bot_2} \to \eE'_{\bot_2} + \eE'_\bot$,
  $\eE_{\bot} \to \eE'_{\bot}$, and $\eE_{\top} \to \eE'_{\top}$ such
  that: if $[a_\bot, a_\top]_\bot$ is sent to
  $[a_\bot', a_\top']_\bot$ then $[a_\top]_\top$ is sent to
  $[a_\top]_\top$; if $[a_{\bot_1},a_{\bot_2},a_\top]_{\bot_i}$ with
  $i \in \setof{1,2}$ is sent to
  $[a_{\bot_1}',a_{\bot_2}',a_\top']_{\bot_i}$ then the other
  $[a_{\bot_1},a_{\bot_2},a_\top]_{\bot_{j}}$ is sent to
  $[a_{\bot_1}',a_{\bot_2}',a_\top']_{\bot_j}$ and $[a_\top]_\top$ is
  sent to $[a_\top']_\top$; and if
  $[a_{\bot_1},a_{\bot_2},a_\top]_{\bot_j}$ is sent to
  $[a_\bot', a_\top']_{\bot}$ then the other
  $[a_{\bot_1},a_{\bot_2},a_\top]_{\bot_j}$ is also sent to
  $[a_\bot', a_\top']_{\bot}$ and $[a_\top]_\top$ is sent to
  $[a_\top']_\top$.

  Hence the objects of $\PSet$ are cospans
  $\eE_{\bot_1} \to \eE_\top \ot \eE_\bot$, but the morphisms are
  allowed to collapse any element of $\eE_{\bot_1}$ with an element of
  $\eE_\bot$, respecting the maps into $\eE_\top$. It is readily
  checked that preimages of subcoalgebras under coalgebra maps are
  subcoalgebras, so this comonad is taut by \cref{prop:tautpreimages}.
  
  Now we describe the density comonad $\gen{\pP'}$ of the second diagram
  $\pP'$. An $\xA$-valued germ $[f]_x \in \gen{\pP'}(\xA)$ is either
  \begin{itemize}
  \item the germ of a function $f \colon \setof{\bot, \top} \to \xA$ (an
    ordered pair) about $\bot$ or $\top$, which we denote by $[a_\bot, a_\top]_\bot$
    or $[a_\bot, a_\top]_\top$;
  \item the germ of a function $f \colon \setof{\bot',\top',\top} \to \xA$ (an
    ordered triple) about $\bot$, $\top'$, or $\top$, which we denote by
    $[a_{\bot'}, a_{\top'}, a_\top]_{\bot'}$, $[a_{\bot'}, a_{\top'}, a_\top]_{\top'}$, or
    $[a_{\bot'}, a_{\top'}, a_\top]_{\top}$; or
  \item the germ of a function $f \colon \setof{\top} \to \xA$ (an
    element) about the unique element $\top$ of $\setof{\top}$, which we
    denote by $[a]_\top$
  \end{itemize}
  modulo the equivalence relation generated by:
  \[
    [x, a, a]_{\top'} = [x, a]_\top = [a]_\top = [x, y, a]_\top  
    \quad
    [x, a, a]_{\bot'} = [x, a]_\bot  
  \]
  Each equivalence class then has a canonical representative: either
  $[a]_\top$, $[a_\bot, a_\top]_\bot$, or one of
  $[a_{\bot'},a_{\top'},a_\top]_{\bot'}$ or
  $[a_{\bot'},a_{\top'},a_\top]_{\top'}$ with
  $a_{\top'} \neq a_{\top}$.

  A coalgebra $\hH \colon \eE \to \cP(\eE)$ assigns each element $e$
  an $\eE$-valued germ of the form $[e]_\top$, $[e, a_\top]_\bot$,
  $[e,a_{\top'},a_\top]_{\bot'}$ or $[a_{\bot'},e,a_\top]_{\top'}$ (as
  guaranteed by the unit law) and such that in the latter three cases
  $h(a_\top) = [a_\top]_\top$, in the third case
  $h(a_{\top'}) = [e,a_{\top'},a_\top]_{\top'}$, and in the fourth
  case $h(a_{\bot'}) = [a_{\bot'},e,a_\top]_{\bot'}$ (as guaranteed by
  the associativity law). In other words, a coalgebra is a set $\eE$
  partitioned into four subsets
  $\eE_{\bot'} + \eE_{\top'} + \eE_{\bot} + \eE_{\top}$, equipped with
  functions
  $\eE_{\bot'} \cong \eE_{\top'} \to \eE_\top \ot \eE_{\bot}$. A
  coalgebra map $\eE \to \eE'$ consists of functions
  $\eE_{\bot'} \to \eE'_{\bot'} + \eE'_\bot$,
  $\eE_{\top'} \to \eE'_{\top'} + \eE'_\top$,
  $\eE_{\bot} \to \eE'_{\bot}$, and $\eE_{\top} \to \eE'_{\top}$ such
  that: if $[a_\bot, a_\top]_\bot$ is sent to
  $[a_\bot', a_\top']_\bot$ then $[a_\top]_\top$ is sent to
  $[a_\top]_\top$; if $[a_{\bot'},a_{\top'},a_\top]_{i'}$ with
  $i \in \setof{\bot, \top}$ is sent to
  $[a_{\bot'}',a_{\top'}',a_\top']_{i'}$ then the other
  $[a_{\bot'},a_{\top'},a_\top]_{j'}$ is sent to
  $[a_{\bot'}',a_{\top'}',a_\top']_{j'}$ and $[a_\top]_\top$ is sent
  to $[a_\top']_\top$; and if $[a_{\bot'},a_{\top'},a_\top]_{i'}$ is
  sent to $[a_\bot', a_\top']_i$ then the other
  $[a_{\bot'},a_{\top'},a_\top]_{j'}$ is sent to
  $[a_\bot', a_\top']_j$ and $[a_\top]_\top$ is sent to
  $[a_\top']_\top$.

  Hence the objects of $\CoalgSet{\pP'}$ are cospans
  $\eE_{\top'} \to \eE_\top \ot \eE_\bot$, but the morphisms are
  allowed to collapse any element of $\eE_{\top'}$ with an element of
  $\eE_\bot$, respecting the maps into $\eE_\top$. This is the same
  description as $\SetP$ above. However, this comonad is not
  taut. We have a coalgebra map between sets in the diagram $\pP'$
  equipped with their tautological coalgebra structures given by
  $\setof{a_{\bot'}, a_{\top'}, a_\top} \to \setof{a_{\bot},
      a_\top}$, where $\setof{a_\top}$ is a subcoalgebra in the latter
    coalgebra, but its preimage $\setof{a_{\top'}, a_\top}$ in the
    former coalgebra is not a subcoalgebra.
\end{counterexample}

Every category of coalgebras on a comonad on $\Set$ is cototal since
its opposite is monadic over $\Set$ by \cref{lem:opmonadic}. Every
category of coalgebras on a small comonad on $\Set$ is total since it
is locally presentable by \cref{lem:presentability}. However, the
following is a large comonad on $\Set$ whose category of coalgebras is
not total.

\begin{counterexample}[Non-total comonad]
  Consider the comonad from \cref{ex:part} whose coalgebras are
  partitioned sets $\Part$. If $\Part$ were total, then any functor
  $\Part\op \to \Set$ sending colimits in $\Part$ to limits would be
  representable. Let $\xX_\kappa$ denote a part of cardinality
  $\kappa$. Consider the large product of representable functors
  $\Part(\dash, 1 + \xX_\kappa)$ over all cardinalities $\kappa$. Given
  any coalgebra (partitioned set), there is just one coalgebra map
  into $1 + \xX_\kappa$ for all but finitely many $\kappa$; hence the
  large product of functors exists. However, the large product of
  coalgebras $1 + \xX_\kappa$ does not exist.
\end{counterexample}

\section*{Density comonads}

The density comonad of any functor $\pP \colon \bB \to \Set$ with $\bB$ small
exists. For large $\bB$, the density comonad need not exist.

\begin{counterexample}[Nonexistence of density comonad]
  If the density comonad of a diagram in $\Set$ exists, then its
  colimit does as well by \cref{lem:setpointwise} and
  \cref{lem:pts}. Large diagrams in $\Set$, of course, need not have colimits;
  for instance the class of all (small) ordinals and their inclusions has 
  no union.
\end{counterexample}

In \cref{def:pointwise} we recalled the definition of pointwise
density comonad, and in \cref{def:dense} we recalled that a functor is
dense if and only if its density comonad exists, is pointwise, and is
identity. A density comonad may be identity without being pointwise.

\begin{counterexample}[Non-pointwise density comonad]\label{cex:nonpointwise}
  The density comonad of the inclusion
  $\pP \colon \bB \hookrightarrow \bS$ of $t$ in the poset $\bS$
  below is identity, although $\pP$ is not dense.
  \[\begin{tikzcd}[row sep=20pt,column sep=10pt]
      t_1 & t_2 \\
      s_1 & s_2
      \arrow[from=2-1, to=1-1]
      \arrow[from=2-2, to=1-2]
    \end{tikzcd}
  \]

  Indeed, the sources $s_1$ and $s_2$ are not initial, as required for
  density by the canonical colimit formula. Nevertheless, the identity
  $\id_{\bS}$ is the density comonad of $\pP$, i.e.\ the universal
  functor $F \colon \bS \to \bS$ equipped with a natural
  transformation
  \[
    \begin{tikzcd}[column sep=15pt, row sep=15pt]
      {\bB} &{}& {\bS} \\
      & {\bS}
      \arrow["\pP", from=1-1, to=1-3]
      \arrow["\pP"', from=1-1, to=2-2]
      \arrow["F"', dashed, from=2-2, to=1-3]
      \arrow["\Downarrow"{description, pos=0.3}, draw=none, from=1-2, to=2-2]
    \end{tikzcd}
  \]
  because any such functor $F$, i.e.\ a functor fixing $t_1$ and
  $t_2$, admits a (unique) natural transformation into it from
  identity.
\end{counterexample}

In \cref{prop:density} we recalled that the comonadic functor
associated to the density comonad of a functor $\pP$ is the free
comonadic functor extending $\pP$. But it is possible for there to
exist a free comonadic functor extending $\pP$ without the density
comonad itself existing, inducing a comonad which is not the density
comonad.

\begin{counterexample}[Non-density universal comonad]
  Every comonad on a preorder $\bS$ is an idempotent comonad, and so
  its comonadic functor is a coreflective subcategory inclusion. Thus
  in the case of preorders, \cref{prop:density} states that if the
  density comonad of a functor $\pP$ into a preorder $\bS$ exists,
  then it must be the smallest coreflective subcategory of $\bS$
  containing the image of $\pP$.
  
  Consider the empty diagram $\bang \colon \none \to \bB$ in the
  following poset $\bB$.
  \[\begin{tikzcd}[row sep=10pt,column sep=10pt]
      & t \\
      s_1 && s_2
      \arrow[from=2-1, to=1-2]
      \arrow[from=2-3, to=1-2]
    \end{tikzcd}\] Since $\bB$ itself is the only coreflective
  subcategory of $\bB$, the identity is indeed the free comonadic
  functor extending $\bang$. However, the identity is not the density
  comonad of $\bang$, i.e.\ it is not the universal functor
  $F \colon \bS \to \bS$ equipped with a natural transformation
  \[
    \begin{tikzcd}[column sep=15pt, row sep=15pt]
      {\none} &{}& {\bS} \\
      & {\bS}
      \arrow["\bang", from=1-1, to=1-3]
      \arrow["\bang"', from=1-1, to=2-2]
      \arrow["F"', dashed, from=2-2, to=1-3]
      \arrow["\Downarrow"{description, pos=0.3}, draw=none, from=1-2, to=2-2]
    \end{tikzcd}
  \]
  For instance there is no natural transformation from identity to the
  functor $F$ that swaps $s_1$ and $s_2$.
\end{counterexample}

In the proof of \cref{lem:pts} we showed that for any functor
$\pP \colon \bB \to \bS$ with density comonad $\cP$, cocones out of
$\lP \colon \bB \to \PS$ are in bijection with cocones out of the
identity on $\PS$. Even so, $\lP$ need not be final (\cref{def:final}).

\begin{counterexample}[Non-final lift]\label{cex:nonfinal}
  Let $\pP$ be the constantly $1$ functor into $\Set$ from the
  discrete category with object set $\xX$, where $\xX$ has cardinality not
  equal to $1$. The lift
  $\lP \colon \xX \to \PSet \cong \SetP \cong \Set/\xX$ is the functor
  $\ptlift{\pP}$ sending each element $x$ in $\xX$ to
  $x \colon 1 \to \xX$. But this is not final: the comma category
  $(\id_\xX, \ptlift{\pP})$ is empty, hence not connected.
\end{counterexample}

\section*{Spaces}

In \cref{thm:tfaetop} we identified topological spaces with density
comonads of diagrams of subsets of a set. However, these are not quite
the same as density comonads of preorder-shaped diagrams of
injections.

\begin{counterexample}[Topological-like non-topological comonad I]
  Consider a diagram consisting of three identities and one
  nonidentity bijection:
  
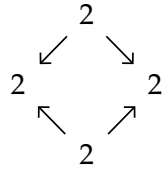
\[\begin{tikzcd}[row sep=10pt,column sep=10pt]
      & 2 \\
      2 && 2 \\
      & 2
      \arrow[from=1-2, to=2-1]
      \arrow[from=1-2, to=2-3]
      \arrow[from=3-2, to=2-1]
      \arrow[from=3-2, to=2-3]
    \end{tikzcd}\] This is a
  preordered diagram of injections, but the maps to the colimit are
  not injective, and the density comonad does not correspond to a
  topological space. Indeed, the density comonad sends a set $A$ to the colimit 
  of 
  
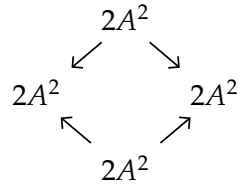
\[\begin{tikzcd}[row sep=10pt,column sep=5pt]
      & 2A^2 \\
      2A^2 && 2A^2 \\
      & 2A^2
      \arrow[from=1-2, to=2-1]
      \arrow[from=1-2, to=2-3]
      \arrow[from=3-2, to=2-1]
      \arrow[from=3-2, to=2-3]
    \end{tikzcd}\] where the only nonidentity map is a free action of
  a group of order $2$, so that the result is isomorphic to $A^2$. The
  comonad has only one point $1^2 \cong 1$, but the only topological
  space with only one point is the identity comonad. (In fact this is
  the comonad corresponding to the one-object category given by the
  group $\zz/2$.)
\end{counterexample}

Comonads arising from topological spaces cannot be characterized
solely in terms of their carrying endofunctor.

\begin{counterexample}[Topological-like non-topological comonad II]\label{cex:nogo}
  The polynomial comonads corresponding to the categories
  \[
    \begin{tikzcd}
      \smbullet & \smbullet
      \arrow[from=1-1, to=1-2, "f"]
    \end{tikzcd}
    \qqand
    \begin{tikzcd}
      \smbullet &[-10pt] \smbullet
      \arrow[from=1-1, to=1-1, loop left, looseness=12, "g=g^2"]
    \end{tikzcd}
  \]
  have the same underlying endofunctor $\xA \mapsto \xA^2+\xA$. But
  only one of these categories is a preorder, hence corresponds to a
  topological space (see \cref{rem:preorders}).
\end{counterexample}

\section*{Bases}

In \cref{def:basis} we defined a basis as a functor
$\pP \colon \bB \to \bS$ admitting a density comonad such that the
induced lift $\lP \colon \bB \to \PS$ is dense. In other words, every
coalgebra is a canonical colimit of the ``basic'' coalgebras
$\lP(\xU)$.  It is not equivalent to merely require that every
coalgebra can be expressed non-canonically as a colimit of coalgebras
$\lP(\xU)$. In particular, for a functor into $\Set$ to be a basis, it
is not sufficient for the images of subbasic coalgebras to constitute
a basis of the underlying topological space.

\begin{counterexample}[Basis-like non-basis I]\label{cex:colimitsubbasis}
  Consider the diagram $\pP \colon \bB \to \Set$, where $\bB$ is the
  walking cospan, given by
  $\setof{\top_1} \hookrightarrow \setof{\bot, \top} \hookleftarrow
  \setof{\top_2}$ where both arrows pick the element $\top$. The
  density comonad $\cP$ is the comonad corresponding to the Sierpinski
  space $\sier$. This is not a basis, but every coalgebra is a colimit
  of subbasic coalgebras.

  Indeed, let us describe the density comonad $\cP$ in terms of germs
  as in \cref{rem:densitygerms}. An $\xA$-valued germ is either the
  germ of a function $\setof{\bot, \top} \to \xA$, which is of the form
  $[a_{\bot}, a_{\top}]_\bot$ or $[a_{\bot}, a_{\top}]_\top$, or an
  element $[a_{\top_1}]_{\top_1}$, or an element
  $[a_{\top_2}]_{\top_2}$, modulo
  $[a]_{\top_1} = [b, a]_\top = [a]_{\top_2}$. Thus the germs can
  be written more simply as either $[a_{\top}]_\top$ or
  $[a_{\bot}, a_{\top}]_\top$, just as in the case of the Sierpinski
  space $\sier$.

  The only maps from the coalgebras $\setof{\top_1}$,
  $\setof{\bot, \top}$, or $\setof{\top_2}$ to the coalgebra
  $\setof{\top_1}$ are given by $\setof{\top_1} \to \setof{\top_1}$
  and $\setof{\top_2} \to \setof{\top_1}$. Thus the canonical colimit in
  $\SetP \simeq \Sh(\sier)$ indexed by such maps is the coproduct
  $\setof{\top_1} + \setof{\top_2}$, which is not
  $\setof{\top_1}$. Therefore the diagram of subbasic coalgebras
  $\lP \colon \bB \to \SetP$ is not dense, i.e.\ $\pP$ is not a basis.
\end{counterexample}

In \cref{lem:extcoref} we noted that if the Yoneda extension of a
functor $\pP$ preserves coreflexive equalizers, then $\pP$ is a
basis. The converse is not true.

\begin{counterexample}[Non-basis-like basis]\label{cex:extcoref}
  Taking $\bB$ to be the walking coreflexive equalizer candidate
  \[\begin{tikzcd}
      \xX & \xY & \xZ
      \arrow[from=1-1, to=1-2]
      \arrow[shift left=3, from=1-2, to=1-3]
      \arrow[shift right=3, from=1-2, to=1-3]
      \arrow[from=1-3, to=1-2]
    \end{tikzcd}\] the functor $\pP \colon \bB \to \Set$ defined by
  $\pP(\xX) = \varnothing$ and $\pP(\xY) = \pP(\xZ) = 1$ does not
  preserve this coreflexive equalizer and hence neither does the
  extension $\extend{\pP} \colon \Ps{\bB} \to \Set$. Nevertheless,
  $\pP$ is dense. Indeed, the canonical colimit formula says that each
  set $\xA$ is the colimit of a diagram consisting of copies of
  $\varnothing$ and $1$ featuring $\xA$ many connected components
  in the diagram consisting of copies of $1$.
\end{counterexample}

Relatedly, a functor into $\Set$ is flat, i.e.\ its category of
elements is cofiltered, if and only if its Yoneda extension preserves
finite limits. The framework of~\cite{adamek-borceux-lack-rosicky}
generalizes this and related results by replacing finite limits with
certain well-behaved classes of limits, called \emph{sound
  doctrines}. (We used this in \cref{sec:coconfluence} to show certain
functors into $\Set$ have density comonads preserving certain classes
of limits.)

As discussed in \cite[Remark 2.6]{adamek-borceux-lack-rosicky}, if a
doctrine is sound, then any functor into $\Set$ from a small category
admitting the relevant limits preserves them if and only if its
extension to the free cocompletion of the domain does. (This weaker
condition is called \emph{weakly sound} in~\cite{tendas}.)
Unfortunately, coreflexive equalizers do not form a (weakly) sound
doctrine. If they did, then \cref{lem:extcoref} would tell us that any
coreflexive-equalizer-preserving functor from a category with
coreflexive equalizers into $\Set$ is a basis.
But this is not so. For a functor into $\Set$ to be a basis, it is not
even sufficient for it to (preserve and) create intersections.

\begin{counterexample}[Basis-like non-basis II]
  Let $\arro$ be the walking arrow, and let
  $\pP \colon \arro \to \Set$ single out the constantly $\top$
  function $\setof{1, 2} \to \setof{\bot, \top}$. (Both sets are $2$
  with differently named elements for convenience.) The density
  comonad $\cP$ is somewhat similar to \cref{ex:shear}, but more
  complicated. A germ $[f]_x \in \cP(\xA)$ is either
  \begin{itemize}
  \item the germ of a function $f \colon \setof{1, 2} \to \xA$ about
    $1$ or $2$, which we denote by $[a_1, a_2]_1$ or $[a_1, a_2]_2$
  \item the germ of a function $f \colon \setof{\bot, \top} \to \xA$
    about $\bot$ or $\top$, which we denote by $[a_\bot, a_\top]_\bot$
    or $[a_\bot, a_\top]_\top$
  \end{itemize}
  modulo the equivalence relation $[a, a]_1 = [b, a]_\top = [a, a]_2$
  for all $a, b \in \xA$. Each equivalence class then has a canonical
  representative: either $[a_1, a_2]_1$ with $a_1 \neq a_2$,
  $[a_1, a_2]_2$ with $a_1 \neq a_2$, $[a_\bot, a_\top]_\bot$, or
  $[a_\top]_\top \coloneqq [a_\top, a_\top]_\top$.

  A coalgebra $\hH \colon \eE \to \cP(\eE)$ assigns each element $e$
  an $\eE$-valued germ of the form $[e, a_2]_1$, $[a_1, e]_2$,
  $[e, a_\top]_\bot$, or $[e]_\top$, (as guaranteed by the unit law)
  and such that in the first case $a_2$ is always assigned
  $[e, a_2]_2$, in the second case $a_1$ is always assigned
  $[a_1, e]_1$, and in the third case $a_\top$ is always assigned
  $[a_\top]_\top$ (as guaranteed by the associativity law). In other
  words, a coalgebra is a set $\eE$ partitioned into four subsets
  $\eE_1 + \eE_2 + \eE_\bot + \eE_\top$, equipped with an isomorphism
  $\eE_1 \cong \eE_2$ and a function $\eE_\bot \to \eE_\top$. A
  coalgebra map $\eE \to \eE'$ must restrict to a commutative square
  between the functions, and must send the elements of each
  corresponding pair $[e_1, e_2]_1 \in \eE_1$ and
  $[e_1, e_2]_2 \in \eE_2$ either to a corresponding pair
  $[e_1, e_2]_1\mapsto [e'_1, e'_2]_1 \in \eE'_1$ and
  $[e_1, e_2]_2\mapsto [e'_1, e'_2]_2\in \eE'_2$ or must send both to
  the same element $[e_1, e_2]_1\mapsto [e']_\top \in \eE'_\top$ and
  $[e_1, e_2]_2\mapsto [e']_\top\in \eE'_\top$.


  The functor $\lP \colon \arro \to \PSet$ is defined as in
  \cref{rem:lift}, with the coalgebra structure on $\setof{1, 2}$
  given by the identity germs $1 \mapsto [1, 2]_1, 2 \mapsto [1, 2]_2$
  and the coalgebra structure on $\setof{\bot, \top}$ given by the
  identity germs
  $\bot \mapsto [\bot, \top]_\bot, \top \mapsto [\top]_\top$.

  Now consider the coalgebra $1 \mapsto \cP(1)$ defined by
  $\ast \mapsto [\ast]_\top$. This coalgebra is subterminal, and there
  is only one coalgebra map from a coalgebra of the form $\lP(\xU)$
  into it, namely a map from $\setof{1, 2}$. Thus the coalgebra
  carried by $1$ is not given by the canonical colimit formula, which
  instead yields $\setof{1, 2}$.
\end{counterexample}


For ordinary subbases of topological spaces
$\topbasis{\tX} \colon \bB \to \Set$, the underlying $\bB$-presheaf of
the sheaf corresponding to a subbasic open set is its representable
presheaf. Thus given $\pP \colon \bB \to \Set$ with density comonad
$\cP$, one might expect the following diagram to commute:
\[\begin{tikzcd}[row sep=10,column sep=10]
    \bB && {\Ps{\bB}} \\
    & \PSet
    \arrow["\yo", from=1-1, to=1-3]
    \arrow["{\lP}"', from=1-1, to=2-2]
    \arrow["{\bsheaf{\pP}}"', from=2-2, to=1-3]
  \end{tikzcd}
\]
However, this is not so, even if we assume $\pP$ is a basis.

\begin{counterexample}[Non-representable basic presheaf]
  If $\bB$ is a connected category, then the terminal functor
  $1 \colon \bB \to \Set$ is dense, thus a basis of the identity
  comonad on $\Set$. The functor $\bsheaf{\pP} \circ \lP$ sends all
  objects to the terminal functor $1$, which is only representable if
  $\bB$ has an initial object.
\end{counterexample}

In \cref{ex:real_line} we defined a comonad on $\Set$ whose coalgebras
are given by ``infinitesimal $\rr$-actions''. However, such
infinitesimal $\rr$-actions can exhibit pathological behavior. Recall
that we identified the category of coalgebras with the topos whose
objects are \'etal\'e spaces over $\rr$ equipped with an action by
$\rr^\delta$ ($\rr$ with the discrete topology) such that acting by
$y$ moves the fiber over $x$ to the fiber over $x + y$.

\begin{counterexample}[Non-continuous infinitesimal $\rr$-action]\label{cex:noncontinuous}
  There is a coalgebra structure on $\setof{\bot,\top}$ whose
  infinitesimal action at $\top$ is constantly $\top$ and whose
  infinitesimal action at $\bot$ is constantly $\top$ except for the
  value $0$, which is assigned $\bot$. This corresponds to a very
  wicked sheaf on $\rr$ that reproduces the relatively familiar
  \emph{line with two origins over every point}: the \'etal\'e space
  $\tE$ over $\rr$ is $\rr\times \setof{\bot,\top}$, with the
  $\top$-line carrying the ordinary topology, but the basic open
  neighborhoods of $(x,\bot)$ are given by a punctured open interval
  around $(x,\top)$ together with $(x,\bot)$ itself. The translation
  action of $\rr^\delta$ preserves both the $\bot$ and the $\top$
  lines. This does not descend to a continuous action of $\rr$,
  because for instance the inverse image of the open neighborhood of
  $(0, \bot)$ given by
  $(-1,0)\times \{\top\} \cup \{(0,\bot)\} \cup (0,1)\times \{\top\}$
  under the action $\rr^\delta \times \tE\to \tE$ intersects the
  $\bot$ plane in a diagonal line, but every basic open in the domain
  intersecting the $\bot$ plane contains a vertical line segment.
\end{counterexample}

\section*{Subbases}

A comonad $\cC$ on $\Set$ need not be the density comonad of
itself. Thus the right adjoint functor $\rcarC \colon \Set \to \CSet$
need not be subdense, as defined in \cref{def:subbasis}.

\begin{counterexample}[Non-self-density-comonad]
  The density comonad of the functor $\cC \colon \Set \to \Set$ given
  by $F(\xA) = 2A$, underlying the discrete comonad on two points, is
  given by $\gen{\cC} \colon \xA \mapsto 2\xA^2$, the indiscrete
  comonad on two points. Indeed, $\cC$ is itself a left adjoint, with
  right adjoint $\xA \mapsto \xA^2$.
\end{counterexample}

\section*{Concrete bases}

In \cref{prop:influent} we showed that if $\El(\pP)$ has a cone to
every span (i.e.\ is co-confluent), then $\extend{\pP}$ preserves weak
limits of spans (i.e.\ weak pullbacks). The analogous statement does
not hold for arbitrary diagrams in place of spans. For example, it is
not true that if $\El(\pP)$ has a cone to every parallel pair, then
$\extend{\pP}$ preserves weak equalizers.

\begin{counterexample}[Non-weak-equalizer-preserving extension]
  Let $\bB$ be the category:
  
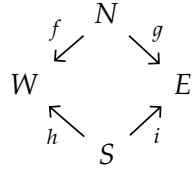
\[\begin{tikzcd}[row sep=10, column sep=10]
      & N \\
      W && E \\
      & S
      \arrow["f"',from=1-2, to=2-1]
      \arrow["g",from=1-2, to=2-3]
      \arrow["h",from=3-2, to=2-1]
      \arrow["i"',from=3-2, to=2-3]
    \end{tikzcd}\] Consider the presheaf $F \colon \bB\op \to \Set$
  sending $N$ to $2$ and all else to $1$ where $F(i)$ and $F(h)$ are
  jointly surjective. By the Yoneda lemma there are two parallel maps
  from $\Ho{\bB}(\dash, N)$ to $F$, and their equalizer is evidently
  empty.

  Now consider the terminal functor $1 \colon \bB \to \Set$. All
  diagrams of parallel pairs in $\El(1) \cong \bB$ have cones. But the
  colimit-preserving functor $\Ps{\bB} \to \Set$ extending the
  terminal functor does not preserve weak equalizers: the presheaves
  $\Ho{\bB}(\dash, N)$ and $F$ are both sent to $1$, so the equalizer
  taken in $\Set$ is $1$, but their equalizer $\varnothing$ is sent to
  $\varnothing$.
\end{counterexample}

In \cref{ex:topcat} we described ``infinitesimal categories'', which
provide a fairly general-purpose way of producing examples of
pullback-preserving comonads on $\Set$. However, not all
pullback-preserving comonads on $\Set$ arise in this way.

\begin{counterexample}[Non-infinitesimal-category ionad]
  As can be seen from the formula in \cref{ex:topcat}, the summands of
  any comonad corresponding to an infinitesimal category are
  \emph{reduced powers}~\cite{blass,pare:taut}, which are equivalently
  finite-limit-preserving quotients of representable
  functors~\cite[Theorem 1]{blass}.

  On the other hand, the finite-limit-preserving comonad from
  \cref{ex:dds} sending $\xA$ to the set of eventually constant
  sequences in $\xA$ is not given by a quotient of a representable
  functor. Indeed, there is no weakly initial object in the category
  of elements: for any particular eventually constant sequence there
  is an eventually constant sequence with greater length before the
  constant value.
\end{counterexample}

\section*{Comonads on slices}

We saw in \cref{prop:intersections} that all comonads on $\Set$ preserve
finite intersections. The same cannot be said for $\Set^\xX$ when $\xX$
has more than one element. 

\begin{counterexample}[Non-intersection-preserving comonad]
  The empty or globally inhabited objects comprise a coreflective
  subcategory of $\Set^\xX$, where the coreflection leaves the
  globally inhabited objects fixed and sends everything else to the
  initial object $\varnothing$.  If $m \colon \xA \hookrightarrow \xX$
  is a monomorphism from a non-inital non-globally-inhabited object to
  a globally inhabited object, then the induced comonad $\cC$ on
  $\Set^\xX$ does not preserve the pullback square given by the
  pushout of $m$ along itself: in this case $\xA$ is sent to the
  initial object while the other three objects in the square are
  fixed.

  Worse than this, the pushout of the comonad map
  $\bang \colon \cC \to \id_{\Set^\xX}$ along itself yields a comonad
  that does not even preserve such monomorphisms
  $\xA \hookrightarrow \xX$.
\end{counterexample}

In \cref{def:sliceable} we introduced sliceable comonads. Not all
comonads on categories with terminal objects are sliceable.

\begin{counterexample}[Non-cosliceable monad]\label{cex:sliceable}
  Consider the monad $T$ on sets of size not equal to $1$ given by the
  free meet-semilattice with both a top and bottom element. The monad
  does not lift to the coslice category under $T(0) \cong 2$ because
  there is no free algebra on any non-injective map
  $f \colon 2 \to \xX$.

  Indeed, here non-injectivity means identifying the top and bottom
  element, which entails identifying all elements, forcing the algebra
  to have one element. As there are no algebras whatsoever carried by
  any object admitting a map from $f$ in the coslice category, there
  cannot be a free algebra on $f$.
\end{counterexample}

In \cref{rem:sliceableformula} we gave a formula for the sliced
comonad $\clift{\cC}$ of a comonad $\cC$ on $\bS$ in terms of
equalizers in $\CS$. If $\bS$ itself has and both $\cC$ and
$\cC \circ \cC$ preserve these equalizers, then the comonadic functor
$\carC \colon \CS \to \bS$ creates them, allowing us to calculate
$\clift{\cC}$ in terms of equalizers in $\bS$. But this does not hold
in general.

\begin{counterexample}[Non-formulaic sliced comonad]
  We return to the comonad $\cP$ on $\Set$ from \cref{ex:unlink} whose
  category of coalgebras is $\Unlink$. The terminal coalgebra is
  described by the span $\setof{l} \ot \varnothing \to \setof{u}$,
  with underlying set $\pts{\pP} \cong \setof{l, u}$. Let
  $\xA \coloneqq \pts{\pP}$, and let $q \colon \xA \to \pts{\pP}$ be
  the function sending both elements to $l$. The cofree coalgebra on
  $q$ with respect to the lifted comonad $\clift{\cP}$ on
  $\Set^{\pts{\pP}}$, is not given by the expected equalizer formula
  in $\Set$, namely the equalizer of maps $\cP(q)$ and
  $\hH_{\pts{\pP}}\circ \bang_{\cP(\xA)}$ from $\cP(\xA)$ to
  $\cP(\pts{\pP})$.

  Indeed, as noted in \cref{rem:setsliceformula}, this equalizer
  consists of the $\xA$-valued germs lying over a $\pts{\pP}$-valued
  identity germ via $q$. As we found in \cref{ex:unlink}, the
  $\xA$-valued germs $\cP(\xA)$ are given by the disjoint union of
  the sets
  \[\setof{[a]_l \mid a \in \xA} \cong \xA
    \qquad
    \setof{[a]_u \mid a \in \xA} \cong \xA
    \qquad
    \setof{[a_l,a_{l'},a_u]_{l'} \mid a_l, a_{l'}, a_u \in \xA, a_{l} \neq a_{l'}}
  \]
  Taking $\xA \coloneqq \pts{\pP}$, we get
  \[
    \setof{[l]_l, [l]_u, [l, u, l]_{l'}, [l, u, u]_{l'}, [u]_{l}, [u]_{u}, [u, l, l]_{l'}, [u, l, u]_{l'}}
  \]
  The subset of germs lying over a $\pts{\pP}$-valued identity germ
  --- i.e.\ $[l]_l$ or $[u]_u$ --- via the constantly $l$ function $q$
  are then
  \[
    \setof{[l]_l, [l, u, l]_{l'}, [l, u, u]_{l'}, [u]_{l}, [u, l, l]_{l'}, [u, l, u]_{l'}}
  \]
  but this is not a subcoalgebra of $\cP(\xA)$. Indeed, the largest
  subcoalgebra included within it (hence the equalizer in $\Unlink$,
  which defines $\clift{\cP}$) is $\setof{[l]_l, [u]_{l}}$.
\end{counterexample}

In \cref{prop:liftbasis} we saw that if $\pP$ is a basis, then
$\gen{\ptlift{\pP}} \cong \clift{\cP}$. This does not hold for
arbitrary $\pP$. 

\begin{counterexample}[Non-slicing subbasis]\label{ex:notdense}
  Let $\bB$ be the walking parallel pair $\rightrightarrows$, and let
  $\pP' \colon \bB \to \Set$ be the functor sending the two
  nonidentity arrows to the following two maps $2 \to 3$:
  \[
    \begin{tikzcd}[column sep=0]
      {\{} & \top & \bot & {\}} & \\
      {\{} & \top & \bot_1 & \bot_2 & {\}}
      \arrow[maps to, from=1-2, to=2-2]
      \arrow[maps to, from=1-3, to=2-3]
    \end{tikzcd}
    \qquad\qquad
    \begin{tikzcd}[column sep=0]
      {\{} & \top & \bot & {\}} & \\
      {\{} & \top & \bot_1 & \bot_2 & {\}}
      \arrow[maps to, from=1-2, to=2-2]
      \arrow[maps to, from=1-3, to=2-4]
    \end{tikzcd}
  \]
  The density comonad $\cP$ is the comonad corresponding to the
  Sierpinski space $\sier$, with set of points $|\pP| = 2$. Indeed, proceeding
  as in \cref{rem:densitygerms}, an $\xA$-valued germ
  $[f]_x \in \cP(\xA)$ takes one of the forms
  \[[a_\top, a_\bot]_\top \qor [a_\top, a_\bot]_\bot \qor [a_\top, a_{\bot_1}, a_{\bot_2}]_\top \qor [a_\top, a_{\bot_1}, a_{\bot_2}]_{\bot_1} \qor [a_\top, a_{\bot_1}, a_{\bot_2}]_{\bot_2}\]
  modulo
  \[[a, b]_\top = [a, b, c]_\top = [a, c]_\top \qand [a, b,
    c]_{\bot_1} = [a, b]_\bot = [a, b, d]_{\bot_2}\] reducing to the
  forms $[a]_\top \coloneqq [a, b]_\top $ and $[a, b]_{\bot}$, which
  is as in the case of germs for the Sierpinski space.
  
  However, the density comonad of the lift
  $\ptlift{\pP} \colon \bB \to \Set/2$ is not $\clift{\cP}$. Indeed,
  $\ptlift{\pP}$ sends both objects of $\bB$ to surjections, so by the
  colimit formula for the density comonad, any non-surjection
  $\xA \to 2$ is sent to the initial object $\varnothing \to 2$. On the
  other hand, the comonad $\clift{\cP}$ on $\Set/2$ induced by the
  Sierpinski space restricts to identity on all objects $\xA \to 2$
  with image $\top$.
\end{counterexample}

The same example also shows that it is not necessarily the case that
given $\pP \colon \bB\to\Set$, the images of the maps to the colimit
$\pP(U) \to \pts{\pP}$ determine a subbasis of the underlying
topological space.

\section*{Continuous maps}

We saw in \cref{thm:topreflection} that $\Top$ is reflective in
$\ConSet$. It is not coreflective.

\begin{counterexample}[Non-preserved topological colimit]
  The pushout of the inclusion of the discrete set of two points into
  the Sierpinski space along itself in $\ConSet$ is not as in $\Top$,
  but as in $\Cat$, since $\Cat$ is coreflective in $\ConSet$.
\end{counterexample}

Neither $\ComSet$ nor $\ConSet$ is extensive.

\begin{counterexample}[Non-extensivity of comonad maps]
  For any comonad $\cC$ on $\Set$, the functor
  $\car{\clift{\cC}} \colon \CSet \to \SetpC$ corresponds a comonad
  map from $\cC$ to the comonad corresponding to the discrete
  topological space on $\pts{\cC}$, which is the coproduct of
  $\pts{\cC}$ many copies of the terminal comonad $\idSet$. But the
  domain $\cC$ need not decompose as a $\pts{\cC}$-indexed coproduct
  of one-point comonads. If it did, $\car{\clift{\cC}}$ would
  decompose a $\pts{\cC}$-indexed product of comonadic functors by
  \cref{prop:coproduct}, and therefore every coalgebra would decompose
  as a coproduct of coalgebras lying entirely over a single point in
  $\pts{\cC}$.

  For example, the functor
  $\car{\clift{\cC}} \coloneqq \slicer \colon \Set \to \Set^2$, induced by
  the indiscrete comonad on two points, is not a coproduct of two
  one-point comonads. The only coalgebra that arises as a coproduct of
  coalgebras lying entirely over a single point is the empty one.
\end{counterexample}

\begin{counterexample}[Non-extensivity of continuous maps]
  The underlying topological space of the comonad $\cP$ from
  \cref{ex:unlink}, whose category of coalgebras is $\Unlink$, is
  discrete on two points. We thus have a continuous map from $\cP$ to
  the discrete comonad on two points given by the reflection
  \cref{thm:topreflection}. But $\cP$ does not decompose as a
  coproduct of two one-point comonads. For instance, the coalgebra
  described by the span $1 \ot 1 \to 1$ is not a coproduct of
  coalgebras lying entirely over a single point.
\end{counterexample}



In \cref{prop:locallysmall} we saw that the 2-category $\ConSet$ is
locally locally small (and more generally the double category
$\CmdSet$ has a small set of 2-cells per boundaries by the same
argument). It can moreover be shown that the categories $\AccComSet$
and $\AccConSet$ are locally small, using that small comonads have
small bases. However, the larger categories $\ComSet$ and $\ConSet$
are not locally small.

\begin{counterexample}[Non-small many continuous maps]
  Consider the category $\bS$, admitting a forgetful functor to
  $\Part$ from \cref{ex:part}, described as follows. An object
  consists of a partitioned set and, for each part $\xX$, a set of
  cardinalities smaller than that of $\xX$. An arrow consists of a
  function of underlying sets such that the image of each part is a
  part where the assigned set of cardinalities is preserved, deleting
  cardinalities that are too big. The forgetful functor to $\Set$
  sending a partitioned set with such structure to the underlying set
  is comonadic.

  Indeed, the forgetful functor $\bS \to \Part$ evidently preserves
  intersections and is conservative, and it has a right adjoint
  sending a partitioned set $\eE$ to the disjoint union of all
  possible choices of a part in $\eE$ equipped with any subset of
  smaller cardinalities. Therefore the forgetful functor
  $\bS \to \Part$, as well as the composite $\bS \to \Part \to \Set$,
  is crudely comonadic by the crude comonadicity theorem
  (\cref{prop:crude}).

  Note that a set of cardinalities less than a set $\xX$ may
  equivalently be viewed as a function from the set of smaller
  cardinalities than $\xX$ to $\setof{\bot, \top}$. Choosing an
  automorphism of $\setof{\bot, \top}$ for each particular cardinality
  induces an automorphism of the category $\bS$ over $\Part$. Each
  choice determines a bijective-on-points comonad map ---
  equivalently, bijective-on-points continuous map --- from the
  comonad $\cP$ to itself, and there are a proper class of such
  choices.
\end{counterexample}

\section*{The double category of comonads}

In \cref{prop:companion} we characterized the companion comonad maps
in the double category $\CmdIdSet$ as the comonad maps carried by
cartesian natural transformations. But for general $\bS$, not all
cartesian comonad maps are companions in the double category
$\CmdId(\bS)$.

\begin{counterexample}[Non-companion cartesian comonad map]\label{cex:noncompanion}
  Let $\bS$ be the category of sets with cardinality not equal to
  six. Consider the polynomial comonad $\cC$ carried by
  $A \mapsto 2A^2$ corresponding to the indiscrete category with two
  objects $\bI$, and the polynomial comonad $\cD$ carried by
  $A \mapsto A^2$ corresponding to the one-object category
  $\zz/2$. The evident discrete opfibration $F \colon \bI \to \zz/2$
  induces a cartesian comonad map $\phi \colon \cC \to \cD$ (by
  \cref{ex:dopf} and \cref{prop:cartesianetale}).

  A comonad map is a companion in $\CmdId(\bS)$ if and only if it is
  mate to a continuous map in the sense of
  \cref{prop:companionadjoint}.
  However, there can be no right adjoint to $\com{\bS}{\phi}$ given on
  underlying sets by pullback along $\pts{F} \colon 2 \to 1$, as there
  exist $\zz/2$-sets of size $3$, and the pullback of the cospan
  $2 \to 1 \ot 3$ does not exist in $\bS$. (The pullback could not be
  anything but $6$, as seen by probing with maps from $1$.)
\end{counterexample}

In \cref{prop:cartesianetale} we saw that cartesian comonad maps and
\'etale comonad maps between pullback-preserving comonads on $\Set$
coincide. For comonad maps between more general comonads on $\Set$,
cartesian and \'etale are independent properties.

\begin{counterexample}[Non-cartesian \'etale comonad map]
  A cartesian comonad map between comonads on $\Set$ identifies each
  indecomposable summand of the domain endofunctor on $\Set$ with some
  indecomposable summand of the codomain endofunctor on $\Set$. Hence we just need
  to find a comonad admitting a coalgebra whose comonad of elements is
  not of that form. Consider the comonad $\cP$ from~\cref{ex:shear},
  whose category of coalgebras is $\Shear$. An object in $\Shear$
  consists of a pair of sets $(E_1, E_2)$, and an arrow
  $(E_1, E_2) \to (E'_1, E'_2)$ consists of a function
  $E_1 + E_2 \to E'_1 + E'_2$ sending $E_1$ entirely within $E'_1$,
  and the defining comonadic forgetful functor
  $\carP \colon \Shear \to \Set$ sends $(E_1, E_2)$ to the set
  $E_1 + 2E_2$. Thus the full subcategory of $\Shear$ consisting of
  objects of the form $(0, E_2)$, or equivalently the slice over
  $(0, 1)$, is simply $\Set$. The comonad of elements $\El(0,1)$ is
  the indiscrete comonad on $2$ points, since the associated comonadic
  composite functor $\Set \cong \Shear/(0,1) \to \Shear \to \Set$ is
  $\xA \mapsto 2\xA$. The comonad is given by the formula
  $\xA \mapsto 2\xA^2$, and the summand $\xA \mapsto \xA^2$ is not
  isomorphic as an endofunctor to $\cP$, since $\cP$ is not polynomial
  (as noted in \cref{ex:shear}, its formula is given cardinality-wise
  by $\pts{\cP(\xA)} = \pts{2 \xA^2 - \xA}$).
\end{counterexample}

\begin{counterexample}[Non-\'etale cartesian comonad map]
  Let $\cP$ be the comonad from~\cref{ex:unlink}, whose category of
  coalgebras is $\Unlink$.  The coproduct codiagonal comonad map
  $\cP + \cP \to \cP$ is cartesian (as are all structure maps for
  coproducts of comonads on $\Set$) but not \'etale. Indeed, the
  coalgebra category of $\cP + \cP$ is $\Unlink \times \Unlink$ by
  \cref{prop:coproduct}, and the induced functor
  $\Unlink \times \Unlink \to \Unlink$ takes coproducts; but this is
  not the projection from the slice of $\Unlink$ over
  $\pts{\pP} + \pts{\pP}$, as there is a map from the coalgebra $\eE$
  described by the span $\setof{l} \ot \setof{l'} \to \setof{u}$ to
  $\pts{\pP} + \pts{\pP}$ sending $l$ and $l'$ to the same element in
  one summand but $u$ to an element in the other summand, and the
  fibers do not decompose $\eE$ as a coproduct.
\end{counterexample}

For a comonad map $\phi$ and continuous map $f$ between comonads on
$\Set$ to be companions in $\CmdId$, it is not sufficient to have
$\one{f} = \one{\phi}$ and $\com{\Set}{\phi} \dashv \con{\Set}{f}$.

\begin{counterexample}[Companion-like non-companion]
  The companion of a morphism in $\CmdId$ is determined up to
  identification, i.e.\ equality. Thus we may find a counterexample by
  finding two right adjoints of $\com{\Set}{\phi}$ that are not equal,
  hence cannot both come from the companion of $\phi$, but which both
  lift $\reindex{\one{\phi}}$.

  By \cref{ex:dopf} and \cref{prop:cartesianetale}, any discrete
  opfibration $F$ of categories induces a companion comonad morphism
  $\phi$ between corresponding categories; the right adjoint of
  $\com{\Set}{\phi}$ is then the reindexing functor along $F$. Now any
  functor $F'$ naturally isomorphic to $F$ induces another naturally
  isomorphic reindexing functor, which is thus also a right adjoint of
  $\com{\Set}{\phi}$, and if $F'$ agrees with $F$ on objects then this
  reindexing functor likewise lifts $\reindex{\one{\phi}}$.

  For example, consider a nonabelian group viewed as a one-object
  category and take the unique coslice category projection onto it.
  This is naturally isomorphic to a distinct functor by conjugation.
\end{counterexample}

We saw in \cref{prop:conjoint} that the conjoints in $\Cmd(\bS)$ are
the same as those in $\CmdId(\bS)$ up to isomorphism. However, the
companions in $\CmdSet$ are more general than in $\CmdIdSet$ (even up
to isomorphism).

\begin{counterexample}[Non-companion companion]
  It is a general double-categorical fact that the companion and
  conjoint of an arrow are adjoints, and moreover if any two out of
  three of these relationships (companion, conjoint, adjoint) between
  pairs in a triple of arrows of appropriate type are satisfied then
  so is the third. In particular, in the double category of categories
  that is a sub-double-category of $\CmdSet$, any functor that is left
  adjoint to a bijective-on-objects (thus conjoint) functor is a
  companion of its right adjoint functor's conjoint retrofunctor.

  As an example, the left adjoint constantly initial functor from the
  walking isomorphism $\bI$ to the walking arrow $\arro$ (which is
  clearly not isomorphic to any discrete opfibration, i.e.\ companion
  in $\CmdIdSet$) is a companion of the retrofunctor conjoint to
  either right adjoint bijective-on-objects functor from $\arro$ to
  $\bI$.
\end{counterexample}

We saw in \cref{cor:topcom} that a comonad map $\phi$ between the
comonads on $\Set$ corresponding the topological spaces $\tX$ and
$\tY$ consists of a function
$\one{\phi} \colon \pts{\tX} \to \pts{\tY}$ and a refinement of the
domain making the function \'etale. Unlike a continuous map of
topological spaces, here the comonad map need not be uniquely
determined by the function $\one{\phi}$: there may be more than one
refinement of $\tX$ making $\one{\phi}$ \'etale.

\begin{counterexample}[Non-determined comonad map]
  Consider the function $\one{\phi}$ from
  $\pts{\tX} \coloneqq \setof{a_1, a_2, b_1, b_2}$ to
  $\pts{\tY} \coloneqq \setof{a, b}$ sending each $a_i$ to $a$ and each
  $b_i$ to $b$. Equip both sets $\pts{\tX}$ and $\pts{\tY}$ with the
  indiscrete topology. There exist two refinements of $\pts{\tX}$ making
  $\one{\phi}$ \'etale, namely the topology generated by
  $\setof{\setof{a_1, b_1}, \setof{a_2, b_2}}$ and the topology
  generated by $\setof{\setof{a_1, b_2}, \setof{a_2, b_1}}$. Thus
  there exist two comonad map structures on $\one{\phi}$.
\end{counterexample}

\section*{Limits and colimits}



In \cref{prop:comonadlimit} we saw that the category $\AccComSet$ of
small comonads on $\Set$ is complete. However, the category $\ComSet$
of arbitrary comonads on $\Set$ is not complete. By
\cref{prop:multicoalgebra}, it suffices to find a limit of comonadic
functors in $\CAT/\Set$ that is not a left adjoint.


\begin{counterexample}[Non-limit of comonads]\label{ex:nolimit}
  Consider product of the comonad from \cref{ex:part}, whose
  coalgebras are partitioned sets $\Part$, with itself. We will show
  that the functor $\carr \colon \Part \times_\Set \Part \to \Set$
  does not have a right adjoint. An object of the comma category
  $(\carr / 2)$ --- in which a terminal object is the same as a right
  adjoint to the functor $\carr$ at $2$ --- is a set equipped with two
  partitions $A$ and $B$ with each element assigned either $\bot$ or
  $\top$. We will construct a distinct object of $(\carr / 2)$ for
  each hereditary well-founded set, such that maps out of these
  objects have pairwise disjoint images. Therefore there can be no
  terminal object in $(\carr / 2)$, since it would be at least as
  large as the proper class of hereditary well-founded sets. If
  $S = \setof{S_i}_{i \in I}$ is a hereditary well-founded set, we
  define its corresponding object $F(S)$ of $(\carr / 2)$ as
  follows. The underlying object in $\Set/2$ of $F(S)$ is given by
  that of
  $\setof{\top} + I \times \setof{\bot} + \sum_{i \in I} F(S_i)$. We
  then define the $A$ and $B$ partitions so that
  ${\top} + I \times \setof{\bot}$ form an $A$-part, $\setof{\top}$
  forms a one-element $B$-part alone, each $(i, \bot)$ forms a
  two-element $B$-part with the corresponding $\top$ element of
  $F(S_i)$, and all other parts are defined as in the objects
  $F(S_i)$. We see inductively that a morphism out of such an object
  $F(S)$ in $(\carr / 2)$ can only send $\top$ to another element
  labelled $\top$ encoding the same structure. Therefore $(\carr / 2)$
  cannot have a terminal object, as it would be too large, and
  therefore $\carr$ cannot have a right adjoint.
\end{counterexample}

We saw in \cref{prop:conterminal} that for comonads on $\Set$ the
terminal object with respect to comonad maps is the same as the
terminal object with respect to continuous maps, but this is not true
in general categories besides $\Set$.

\begin{counterexample}[Non-initial initial monad]
  The identity monad on $\Set$ is initial with respect to monad maps,
  but not with respect to co-continuous maps. Indeed, there is no
  co-continuous map from identity to the constantly $1$ monad. This
  would entail a functor between algebra categories lifting the
  $\dash + 1$ functor (pushout along $0 \to 1$), but there are no
  algebras of the latter monad with more than one element.
\end{counterexample}

We saw in \cref{prop:contcocomplete} that all small colimits of
diagrams of continuous maps between small comonads exist. However, not
all colimits in $\ConOSet$ exist.

\begin{counterexample}[Continuous non-colimit]\label{ex:nocolimit}
  By \cref{prop:contcolimits} it suffices to find a connected diagram
  of comonads and bijective-on-points comonad maps whose limit does
  not exist, since this is the same as the colimit of the
  corresponding continuous maps in the opposite
  direction. \cref{ex:nolimit} does the job, since we may consider the
  product as the pullback of the maps into the terminal comonad.
\end{counterexample}

The Day convolution $\otimes$ with respect to $\times$ (see \cref{rem:dirichlet}) of
non-accessible endofunctors, or comonads, need not exist.

\begin{counterexample}[Non-product]
  Recall from \cref{ex:part} the density comonad $\cP$ of the
  subcategory inclusion $\pP \colon \Surj \hookrightarrow \Set$ of the
  category of sets and surjections in $\Set$, whose coalgebras are
  partitioned sets $\Part$. If the Day convolution $\cP \otimes \cP$
  existed, then by \cref{rem:dirichlet} it would be the density
  comonad of the external product $\Surj \times \Surj \to \Set$
  sending $(\xX_1, \xX_2)$ to $\xX_1 \times \xX_2$. However, this
  density comonad does not exist. If it did, we would have that
  $\cP(\xA)$ is the set of all germs $[f]_{(x_1,x_2)}$ of functions
  $f \colon \xX_1 \times \xX_2 \to \xA$ where $x_1 \in \xX_1$ and
  $x_2 \in \xX_2$, modulo the equivalence relation generated by
  $[f]_{(x_1,y_1)} = [f \circ (s_1 \times s_2)]_{(x'_1,y'_1)}$
  whenever $s_1 \colon \xX'_1 \to \xX_1$ is a surjection sending
  $x'_1 $ to $x_1$ and $s_2 \colon \xX'_2 \to \xX_2$ is a surjection
  sending $x'_2$ to $x_2$. But then $\cP(2)$ would be a proper class,
  since for each cardinality $\xX$ the map
  $f \colon \xX \times \xX \to 2 \coloneqq \setof{\bot, \top}$ given
  by sending the diagonal to $\top$ and all else to $\bot$ (with the
  germ centered about any particular point) is in its own equivalence
  class. Indeed, the cardinality of the set of distinct rows
  $\xX_1 \to 2$ (or columns $\xX_2 \to 2$) is invariant under the
  equivalence relation.
\end{counterexample}

Though the Day convolution $\otimes$ with respect to $\times$ agrees
with the cartesian product in $\Ionad$, it need not agree with the
cartesian product in $\AccConSet$.  Moreover, whereas the inclusion
$\Cat \hookrightarrow \Ionad$ preserves products, the inclusion
$\Cat \hookrightarrow \ConSet$ does not.

\begin{counterexample}[Non-product product]\label{cex:nonproductproduct}
  We claim that the comonad on $\Set$ from \cref{ex:monoids}, which
  we will now denote by $\cC$, is the cartesian product in $\ConSet$
  of the comonads corresponding to the individual monoids
  $(\bM_i)_{i \in I}$. Let us assume for simplicity that none of these
  monoids are terminal. Let $\bB$ denote the full subcategory of
  $\CSet$ consisting of each monoid $\bM_i$ acting on itself and the
  terminal object. Equivalently, the category $\bB$ may be obtained by
  freely adjoining a terminal object to the coproduct of one-object
  categories $\sum_{i \in I} \bM_i$. The full subcategory inclusion
  $\bB \hookrightarrow \CSet$ is dense: coalgebras embed fully
  faithfully in $\bB$-presheaves, since the underlying $\bB$-presheaf
  precisely encodes the $\bM_i$-actions and the identification between
  their fixed points. Therefore $\cC$ is the density comonad of the
  composite functor $\bB \hookrightarrow \CSet \xto{\carC} \Set$. The
  universal property of the cartesian product is then verified using
  \cref{thm:universal}: a continuous map $\cD \to \cC$ is the same as
  a functor $\CSet \to \DSet$ lifting
  $\reindex{\bang} \colon \Set \to \SetpD$, which is the same as a
  family of functors $\fun{\bM_i}{\Set} \to \DS$ lifting
  $\reindex{\bang}$, which is the same as a family of continuous maps
  $\cD \to \cattocmd{\bM_i}$.

  On the other hand, as shown in \cite[Remark 6.4]{garner:ionads}, the
  box topology formula that we used to define $\otimes$ in
  \cref{rem:dirichlet} in fact does give the cartesian product in
  $\AccIonad$. It also does so specifically in $\Cat$, since the box
  topology formula applied to the canonical bases of categories as in
  \cref{ex:catbasis} gives the canonical basis of the product
  category. However, taking $I = 2$ and $\bM_1$ and $\bM_2$ nontrivial,
  the comonad $\cC$ does not correspond to a category, or even an
  ionad, as its category of coalgebras is not regular as observed in
  \cref{ex:monoids}.
\end{counterexample}

\section*{Infinitesimal neighborhoods}

In \cref{lem:freeintslice} we saw that the canonical forgetful functor
$\LHalo{\xX} \to \Halo$ from $\xX$-labelled halos to halos becomes
faithful and conservative when restricted to those halos that arise
from comonads on $\Set$. This need not be so without this restriction.

\begin{counterexample}[Non-labelled halo isomorphism]
  Consider the formal limits of $\setof{\bot,\top}$-labelled sets given
  by the following two diagrams:
  \[
    \setof{\top}
    \qqand
    \setof{\top} \hookrightarrow \setof{\bot, \top} \hookleftarrow
    \setof{\top}
  \]
  There are two distinct maps from the former to the latter in
  $\LHalo{\setof{\top, \bot}}$, which both forget to the unique map
  $1 \to 1$ in $\Halo$. This witnesses both non-faithfulness and
  non-conservativity of the forgetful functor $\LHalo{\setof{\bot, \top}}$.
\end{counterexample}

In \cref{rem:slicefunctorcoslice} we noted that we have an
equivalence between the full subcategories of $\LHalo{X}$ and
$\Halo/X$ consisting of those halos encoded by finite-limit-preserving
functors, with the equivalence given by forgetting labels. However,
this does not constitute an equivalence $\LHalo{X} \simeq \Halo/X$.

\begin{counterexample}[Non-labelled halo]\label{cex:haloslice}
  Consider the formal limit of the one-object diagram in $\Set$
  consisting of the endomorphism on $\setof{\bot, \top}$ swapping
  $\bot$ and $\top$.  This corresponds to the functor
  $F \colon \Set \to \Set$
  \[F(\xA) \coloneqq \binom{\xA+1}{2}\] sending a set $\xA$ to the set
  of (not necessarily distinct) unordered pairs of elements in
  $\xA$. The identity on $\setof{\bot, \top}$ induces a map in $\Halo$
  from this formal limit to $\setof{\bot, \top}$ (corresponding under
  \cref{rem:halogerm} to the unordered pair
  $\setof{\bot, \top} \in F(\setof{\bot, \top})$). However, this
  underlies no map into $\setof{\bot,\top}$ in
  $\LHalo{\setof{\bot, \top}}$.
\end{counterexample}

\section*{Comonadicity}

In \cref{lem:create} we recalled that a comonadic functor creates all
limits preserved by the underlying comonad $\cC$ and its square
$\cC \circ \cC$. Dually, a monadic functor creates all colimits
preserved by the underlying monad $T$ and its square $T \circ T$. Here
preservation by just $T$ is not sufficient.

\begin{counterexample}[Preserved non-created colimit]
  Consider the algebraic theory presented by a 1-ary operation $f$ and
  a 3-ary operation $g$ with the following laws:
  \begin{itemize}
  \item $f(f(x)) = x$.
  \item $g(x,x,y) = g(x,y,x) = g(y,x,x) = x$.
  \item $g(x,f(x),y) = g(x,y,f(x)) = g(y,x,f(x)) = y$.
  \end{itemize}

  We claim that the corresponding monad on $\Set$ preserves the
  coproduct $1 + 1$, but the monadic functor does not create this
  coproduct.
  
  Indeed, for the induced monad $T$, let us construct the free algebra
  $T(\setof{a})$: this algebra receives first of all the new element
  $f(a)$.  Then $f(f(a))=a$ by the first axiom, and $g$ is determined
  by ``majority vote'' on any combination of $a$'s and $f(a)$'s by the
  second axiom.  Since of $x$ and $f(x)$ exactly one must always equal
  $y$, in this case the third axiom is implied by the second. This
  determines $T(\setof{a})=\setof{a,f(a)}$.
  
  We also calculate that $T(\setof{a,b}) = \setof{a,b,f(a),f(b)}$. In
  this case, $g$ is no longer fully determined by the second axiom:
  any of the $24$ ways to give three distinct arguments to $g$ require
  addressing via the third axiom. Since the third axiom's permutations
  of the arguments include generators of $S_3$, this axiom says in
  effect that if two of the arguments of the form $x$ and $f(x)$, then
  $g$ always returns the remaining third argument, regardless of
  ordering. Any three-element subset of $\setof{a,b,f(a),f(b)}$
  contains either $a$ and $f(a)$ or $b$ and $f(b)$, but not both. This
  shows that $g$ is uniquely determined on $\setof{a,b,f(a),f(b)}$,
  finishing the construction of the free algebra. As claimed,
  $T(1+1)=T(1)+T(1)$. However, $T(T(1)+T(1))$ is in fact infinite,
  much unlike $T^2(1)+T^2(1)$, which leaves an opening for the
  corresponding monadic functor to fail to create the coproduct $1+1$.
  
  Specifically, in the category of $T$-algebras we find that
  $1 + 1 = 1$ (where $1$ denotes the terminal $T$-algebra). Indeed,
  the $T$-algebra $1+1$ is generated by two fixed points of $f$; but
  the axioms on $g$ then say that $g$ acts via both ``minority'' and
  ``majority'' vote, which implies the two generators are identified,
  as was to be shown.
\end{counterexample}

If $T$ preserves all colimits of a certain shape $J$, then necessarily 
$T \circ T$ preserves them too, so the
monadic functor creates them. One might expect it follows that the monadic functor
also preserves colimits of shape $J$, but this is not true.

\begin{counterexample}[Preserved non-preserved colimit]
  The monad corresponding to any algebraic theory can be restricted to the 
  category of nonempty sets, where it vacuously preserves all colimits of empty 
  shape, since none exist. If the theory has a constant, then the initial algebra 
  is nonempty, so the corresponding monadic functor does \emph{not} preserve empty colimits.
\end{counterexample}

Conversely, if a monadic functor creates colimits of a certain shape,
it does not follow that the monad preserves them.

\begin{counterexample}[Created non-preserved colimit]
  Consider the poset $\bS$ sketched below. The top row
  $x_0 \to z_0 \ot y_0$ is a reflective subcategory $\bS'$, with the
  reflection projecting each column onto the top row.
  
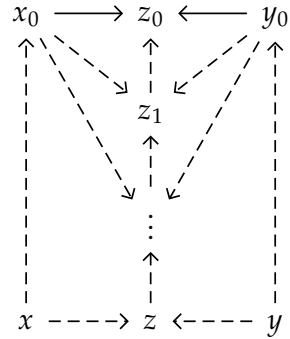
\[\begin{tikzcd}
      {x_0} & {z_0} & {y_0} \\
      & {z_1} \\
      & \vdots \\
      x & z & y
      \arrow[from=1-1, to=1-2]
      \arrow[dashed, from=1-1, to=2-2]
      \arrow[dashed, from=1-1, to=3-2]
      \arrow[from=1-3, to=1-2]
      \arrow[dashed, from=1-3, to=2-2]
      \arrow[dashed, from=1-3, to=3-2]
      \arrow[dashed, from=2-2, to=1-2]
      \arrow[dashed, from=3-2, to=2-2]
      \arrow[dashed, from=4-1, to=1-1]
      \arrow[dashed, from=4-1, to=4-2]
      \arrow[dashed, from=4-2, to=3-2]
      \arrow[dashed, from=4-3, to=1-3]
      \arrow[dashed, from=4-3, to=4-2]
    \end{tikzcd}\] The middle column consists of an infinite
  descending chain with a bottom element $z$ adjoined, and while $x_0$
  and $y_0$ are below every element of the descending chain, they are
  not below $z$, and so have no coproduct in $\bS$. The empty diagram
  in $\bS$ and all diagrams with image $\setof{x_0, y_0}$ have no
  colimit, vacuously created by the monadic inclusion $\bS' \to \bS$;
  colimits of all other diagrams in $\bS$ factoring through $\bS'$
  exist and are created by the monadic inclusion. However, the induced
  monad does not preserve the coproduct $z = x \vee y$.
  
\end{counterexample}


It is possible for a comonad to be crude, i.e.\ preserve coreflexive
equalizers, but not preserve all pullbacks of monomorphisms.

\begin{counterexample}[Non-preserved intersection]
  Consider the comonad on the poset
  
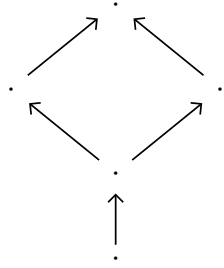
\[\begin{tikzcd}
      & \cdot & \\
      \cdot && \cdot \\
      & \cdot \\
      & \cdot
      \arrow[from=2-1, to=1-2]
      \arrow[from=2-3, to=1-2]
      \arrow[from=3-2, to=2-1]
      \arrow[from=3-2, to=2-3]
      \arrow[from=4-2, to=3-2]
    \end{tikzcd}
  \]
  acting as identity on all objects except for the second from the
  bottom, which it sends to bottom. This comonad does not preserve the
  pullback of the topmost two arrows, which are monomorphisms since
  this is a preorder, but it trivially preserves coreflexive
  equalizers.

  This counterexample can also be modified to show the same can happen
  with regular monomorphisms in place of monomorphisms. Adjoin four
  objects and parallel pairs
  
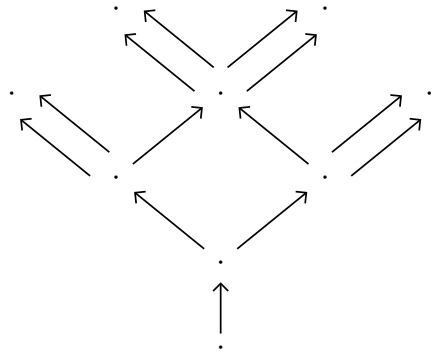
\[\begin{tikzcd}
      & \cdot && \cdot & \\
      \cdot && \cdot && \cdot \\
      & \cdot && \cdot \\
      && \cdot \\
      && \cdot
      \arrow[shift right=2, from=2-3, to=1-2]
      \arrow[shift left=2, from=2-3, to=1-2]
      \arrow[shift right=2, from=2-3, to=1-4]
      \arrow[shift left=2, from=2-3, to=1-4]
      \arrow[shift right=2, from=3-2, to=2-1]
      \arrow[shift left=2, from=3-2, to=2-1]
      \arrow[from=3-2, to=2-3]
      \arrow[from=3-4, to=2-3]
      \arrow[shift right=2, from=3-4, to=2-5]
      \arrow[shift left=2, from=3-4, to=2-5]
      \arrow[from=4-3, to=3-2]
      \arrow[from=4-3, to=3-4]
      \arrow[from=5-3, to=4-3]
    \end{tikzcd}
  \]
  and for each of these pairs impose an equation making it equalized by the
  arrow in line with it. Then the pullback square consists of regular
  monomorphisms, and again the full subcategory without the second
  object from the bottom is coreflective, with the coreflection not
  preserving the pullback.
\end{counterexample}

The dual of \cref{lem:adjointlift} states that for any a monad map
where the category of algebras of the domain monad has reflexive
coequalizers, the induced functor between categories of algebras has a
left adjoint (and hence is moreover monadic). This need not hold if
such reflexive coequalizers do not exist.

\begin{counterexample}[Non-adjoint monad map]\label{cex:noadjoint}
  \cref{cex:sliceable} or \cref{cex:noncompanion} provide examples.
\end{counterexample}

\section*{Categories without 1}

In \cref{def:slice} we defined the slice factorization of any functor
with a colimit, and in \cref{def:comprehensive} we defined the
comprehensive factorization of any functor. These agree in the case
the domain has a terminal object, but they need not agree in general.

\begin{counterexample}[Non-comprehensive slice factorization]
  Take the same example as \cref{cex:nonfinal}: the first factor in
  the slice factorization of the constantly $1$ functor into $\Set$
  from the discrete category $2$ is not final, so this is not the
  comprehensive factorization.
\end{counterexample}

In \cref{def:cartesian} we introduced uniformly cartesian natural
transformations. Not all cartesian natural transformations are
uniformly cartesian.

\begin{counterexample}[Non-uniformly cartesian]
  Consider the unique natural transformation from the diagram of shape
  $\rightrightarrows$ in $\Set$ given by the two distinct bijections
  $2 \to 2$ to the terminal diagram of the same shape. This natural
  transformation is cartesian, i.e.\ the naturality squares are
  pullbacks, because any isomorphism is a pullback of any other. The
  colimit of the first diagram is $1$, i.e. it has one indecomposable
  summand, but it is not isomorphic to the terminal diagram, so the
  natural transformation cannot be uniformly cartesian by
  \cref{prop:stronguniform}
\end{counterexample}

In \cref{def:strongcartesian} we introduced strongly uniformly
cartesian natural transformations. Not all uniformly cartesian natural
transformations between functors admitting colimits are strongly
uniformly cartesian (although we observed this is true for functors
into $\Set$).

\begin{counterexample}[Non-strongly uniformly cartesian]
  Let $\bB$ denote the discrete category with two objects
  $\setof{x, y}$, and let $\bS \coloneqq \arro$ denote the walking
  arrow category $a \to b$. Let $F \colon \bB \to \bS$ be the functor
  such that $F(x) = a$ and $F(y) = b$, and let $G \colon \bB \to \bS$
  be the functor such that $G(x) = G(y) = b$. The unique natural
  transformation $F \Rightarrow G$ is uniformly cartesian but not
  strongly uniformly cartesian.

  Indeed, the induced map between colimits $\pts{F} \to \pts{G}$ is
  simply identity on $b$; $F$ and $G$ are not isomorphic, so $F$ is
  not obtained by pulling back the diagram $G$ along an identity. Thus
  the natural transformation is not strongly uniformly
  cartesian. However it is uniformly cartesian, as is any natural
  transformation between functors from a discrete category
  $\bB$.

  Indeed, in this situation the induced maps between presheaf
  categories $\Psh{F}, \Psh{G} \colon \Set^\bB \to \Psh{\bS}$
  respectively send the $\bB$-indexed set $(X_i)_{i \in \bB}$ to
  $\sum_{i \in \bB} X_i \cdot F(i)$ and
  $\sum_{i \in \bB} X_i \cdot G(i)$, and any natural transformation
  between them is similarly defined. In particular, they respectively
  send the terminal object to $\sum_{i \in \bB} F(i)$ and
  $\sum_{i \in \bB} G(i)$, and so we see the naturality squares at
  maps into the terminal object are pullbacks, given by pulling back
  coproduct structure maps.
\end{counterexample}

As mentioned in \cref{rem:noterminal}, we may generalize definitions
and results of \cref{sec:maps} to categories without terminal objects
by replacing \emph{terminal-object-preserving} with \emph{final} and
\emph{cartesian} with \emph{uniformly cartesian}. One simply modifies
the proofs by replacing everything in sight with free (large if
necessary) cocompletions. However, it is worth noting that if $\cC$ is
a comonad on $\bS$, then the induced functor between free
cocompletions $\Coco{\carC} \colon \Coco{\CS} \to \Coco{\bS}$ is not
necessarily comonadic.

\begin{counterexample}[Non-comonadic cocompletion]\label{cex:comonadiccocompletion}
  Consider the comonadic forgetful functor
  $\fun{\zz/2}{\FinSet} \to \FinSet$, sending a finite $\zz/2$-set to
  its underlying set. We will show that the induced functor
  $\Ps{\fun{\zz/2}{\FinSet}} \to \Ps{\FinSet}$ is not conservative,
  hence not comonadic.

  There is a presheaf on $\fun{\zz/2}{\FinSet}$ sending all
  $\zz/2$-sets with no fixed points to $1$ and all $\zz/2$-sets with a
  fixed point to $\varnothing$. The map in
  $\Psh(\fun{\zz/2}{\FinSet})$ from this presheaf to the terminal
  presheaf is sent to an isomorphism in $\Psh(\FinSet)$. Indeed, this
  presheaf on $\fun{\zz/2}{\FinSet}$ is evidently obtained as a
  quotient of $\Ho{\fun{\zz/2}{\FinSet}}(\dash, \zz/2)$ by identifying
  the two product projections $\zz/2 \times \zz/2 \to \zz/2$, which is
  sent to the presheaf on $\FinSet$ obtained as a quotient of
  ${\FinSet}(\dash, 2)$ by identifying the two product projections
  $2 \times 2 \to 2$, which is terminal. Moreover, the terminal
  presheaf on $\fun{\zz/2}{\FinSet}$ is also sent to the terminal
  presheaf on $\FinSet$ since the comonadic functor
  $\fun{\zz/2}{\FinSet} \to \FinSet$ preserves the terminal object.
  %
  %
\end{counterexample}

Neither does $\Coco{\dash}$ preserve density comonads. 

\begin{counterexample}[Non-density cocompletion comonad]
  The density comonad of the unique functor from the initial category
  to the terminal category is identity. However, the induced functor
  between free cocompletions $\point \to \Set$, which picks out the
  empty set, is not dense, and so its density comonad is not identity.
\end{counterexample}

\section*{Continuous universality}

In \cref{lem:densityformaltfae} we recalled that if
$\pP \colon \bB \to \bS_1$ has density comonad $\cP$, then colax
comonad functors (\cref{def:comfun}) of the form
$F \colon \bS_1 \to \bS_2$ out of $\cP$, as well as any comonad
specializations (\cref{def:comspec}) and comonad transformations
(\cref{def:comtrans}) out of $F$, are uniquely determined by
restricting their associated commutative squares, natural
transformations, and cylinders (as in \cref{prop:formaltfae}) along
$\lP$, as long as we assume $F$ preserves the left Kan extension
$\cP$. Here the left Kan extension preservation assumption is needed.

\begin{counterexample}
  The identity comonad on $\Set$ is the density comonad $\gen{1}$ of
  the diagram $1 \colon \point \to \Set$, and we have
  $\CoalgSet{\gen{1}} \cong \Set$, with the comonadic forgetful
  functor $\car{\gen{1}}$ given by $\idSet$. The indiscrete comonad on
  two points is the density comonad $\gen{2}$ of the diagram
  $2 \colon \point \to \Set$, and we have
  $\CoalgSet{\gen{2}} \simeq \Set$, with the comonadic forgetful
  functor $\car{\gen{2}}$ given by $2 \times \dash$. We do not
  have a correspondence between squares:
  \[
    \begin{tikzcd}
      \point & \CoalgSet{\gen{2}}\\
      \Set & \Set 
      \arrow["1"', from=1-1, to=2-1]
      \arrow[dashed, from=1-1, to=1-2]
      \arrow["\car{\gen{2}}", from=1-2, to=2-2]
      \arrow["\Pow{(\dash)}"', from=2-1, to=2-2]
    \end{tikzcd}
    \qqand
    \begin{tikzcd}
      \Set & \CoalgSet{\gen{2}} \\
      \Set & \Set
      \arrow["\idSet"', from=1-1, to=2-1]
      \arrow[dashed, from=1-1, to=1-2]
      \arrow["\car{\gen{2}}", from=1-2, to=2-2]
      \arrow["\Pow{(\dash)}"', from=2-1, to=2-2]
    \end{tikzcd}
  \]
  The former exists, but the latter, i.e.\ a colax comonad functor
  from $\gen{1}$ to $\gen{2}$ carried by the covariant powerset
  functor $\Pow{(\dash)}$, does not, as there is no object of
  $\CoalgSet{\gen{2}}$ carried by the set $\Pow{\emptyset} \cong 1$.

  Now we move on to consider comonad specializations and comonad
  transformations. Note that there is a unique colax comonad functor
  from identity to identity carried by any particular functor. Thus in
  particular the terminal functor $1 \colon \Set \to \Set$ and the
  identity $\idSet \colon \Set \to \Set$ provide two examples of colax
  comonad functors from the identity comonad $\idSet$ to itself. There
  are no natural transformations $1 \Rightarrow \idSet$, hence no
  specializations between these colax comonad functors, as well as no
  comonad transformations. But both $1 \colon \Set\to \Set$ and
  $\idSet \colon \Set\to \Set$ restrict along
  $1 \colon \point \to \Set$ to the same functor
  $1 \colon \point \to \Set$, and there does exist a natural
  transformation $1 \Rightarrow 1$ between functors $\point \to
  \Set$. Hence such restriction does not provide a bijection.
\end{counterexample}

In \cref{prop:lccbasis} we saw that for comonads on a locally
cartesian closed category (with a terminal object), continuous maps
can be defined using a basis of the codomain comonad in place of its
comonadic functor. This need not hold for comonads on arbitrary
categories. For instance it does not hold for comonads on $\Set\op$,
i.e.\ monads on $\Set$.

\begin{counterexample}[Non-universal basis]


  The functor $2 \colon \point \to \Set$ is a co-basis of the double
  powerset monad, with monadic functor
  $\Pow{(\dash)} \colon \Set\op \to \Set$. Indeed, the induced lift
  $\point \to \Set\op$ can only pick out the set $1$, which is
  codense in $\Set\op$ (i.e.\ dense in $\Set$), thus giving a
  co-basis. On the other hand the functor $1 \colon \point \to \Set$
  is monadic, corresponding to the terminal monad on $\Set$. We do not
  have a correspondence between squares
  \[
    \begin{tikzcd}
      \point & \point \\
      2/\Set & 1/\Set 
      \arrow["\id_2"', from=1-1, to=2-1]
      \arrow[dashed, from=1-1, to=1-2]
      \arrow["\id_1", from=1-2, to=2-2]
      \arrow["\text{\tiny pushout}"', from=2-1, to=2-2]
    \end{tikzcd}
    \qqand
    \begin{tikzcd}
      \Set\op & \point \\
      2/\Set & 1/\Set
      \arrow[from=1-1, to=2-1]
      \arrow[dashed, from=1-1, to=1-2]
      \arrow["\id_1", from=1-2, to=2-2]
      \arrow["\text{\tiny pushout}"', from=2-1, to=2-2]
    \end{tikzcd}
  \]
  where the functor $\Set\op \to 2/\Set$ sends a set $X$ to the
  function $2 \to \Pow{X}$ picking out $\varnothing$ and $X$. The
  former exists but the latter, i.e.\ a co-continuous map of monads,
  clearly does not.
\end{counterexample}

In \cref{prop:subbasiscont} we saw that for comonads on $\Set$,
continuous maps can be defined using an arbitrary subbasis of the
codomain comonad in place of its comonadic functor. This need not hold
for comonads on locally cartesian closed categories besides $\Set$. In
fact, it need not even hold for slice categories of $\Set$.

\begin{counterexample}[Non-universal subbasis]\label{cex:nonuniversal}
  First note that the category $\Set/{\setof{a_1, a_2, b_1, b_2}}$ is
  (strictly) comonadic over $\Set/{\setof{a, b}}$ since it is the
  slice category over the object
  $\pi \colon \setof{a_1, a_2, b_1, b_2} \to \setof{a,b}$ where
  $\pi(a_i) = a$, $\pi(b_i) = b$. Next note that the empty or globally
  inhabited objects form a coreflective subcategory of
  $\Set/\setof{a, b}$. Both comonads are accessible, so by
  \cref{prop:comonadlimit} their product exists and by
  \cref{prop:multicoalgebra} its coalgebras are objects equipped with
  a coalgebra structure for both comonads.

  Therefore the full subcategory $\bS$ of
  $\Set/{\setof{a_1, a_2, b_1, b_2}}$ consisting of objects that are
  either empty or are inhabited at some $a_i$ as well as some $b_j$ is
  (strictly) comonadic over $\Set/{\setof{a, b}}$. We denote the
  comonad by $\cC$ and the comonadic functor by
  $\carC \colon \bS \to \Set/{\setof{a, b}}$ as usual. The terminal
  $\cC$-coalgebra is carried by $\pi$. The cofree $\cC$-coalgebras are
  those with the same fiber at $a_1$ and $a_2$ and with the same fiber
  at $b_1$ and $b_2$. We denote the Kleisli category of cofree
  coalgebras by $\bS'$ and the associated left adjoint by
  $L \colon \bS' \to \Set/{\setof{a, b}}$.

  The lifts
  $\ptlift{\carC} \colon \bS \to \Set/{\setof{a_1, a_2, b_1, b_2}}$
  and $\ptlift{L} \colon \bS' \to \Set/{\setof{a_1, a_2, b_1, b_2}}$
  given by slice factorization are simply the evident forgetful
  functors.  Let
  $i \colon \setof{a, b} \to \setof{a_1, a_2, b_1, b_2}$ be a section
  of $\pi$. The composite $\reindex{i} \circ \ptlift{\carC}$ does not
  factor through the coreflective subcategory of $\Set/\setof{a, b}$
  consisting of empty or globally inhabited objects; therefore there
  is no continuous map carried by $i$ from the idempotent comonad into
  $\cC$. On the other hand, the composite
  $\reindex{i} \circ \ptlift{L}$ does factor through this coreflective
  subcategory. Thus the subbasis $L$ of $\cC$ cannot be used
  to determine continuous maps.
\end{counterexample}


\bibliographystyle{alpha}
\bibliography{references}

\end{document}